%August 13, 2002
\input amstex
\documentstyle{amsppt}
\magnification=1200
\hoffset=-0.5pc
\nologo
\vsize=57.2truepc
\hsize=38.5truepc
\spaceskip=.5em plus.25em minus.20em
\define\squar{\square}
\define\rrtimes{\rtimes}
\define\lltimes{\ltimes}

\define\Bobb{\Bbb}
\define\fra{\frak}
\define\rell{r}
\define\HG{G}
\define\hg{g}
\define\GHH{H}
\define\KK{H}
\define\KKK{H}
\define\kk{h}
\define\kkk{h}
\define\ghh{h}
\define\GH{H}
\define\gh{h}
\define\armcusgo{1}
\define\armamonc{2}
\define\atibottw{3}
\define\barbsepa{4}
\define\beardboo{5}
\define\brunvett{6}
\define\burcush{7}
\define\chevaone{8}
\define\collmcgo{9}
\define\cushmone{10}
\define\djokotwo{11}
\define\elkinone{12}
\define\farakora{13}
\define\freudone{14}
\define\graueone{15}
\define\grauremm{16}
\define\guistebo{17}
\define\guistetw{18}
\define\guhujewe{19}
\define\hartsboo{20}
\define\heihuclo{21}
\define\heinloos{22}
\define\helgaboo{23}
\define\hilneors{24}
\define\hocheago{25}
 \define\howeone{26}
 \define\howetwo{27}
\define\poiscoho{28}
 \define\souriau{29}
 \define\singula{30}
\define\singulat{31}
 \define\locpois{32}
 \define\poisson{33}
   \define\modus{34}
\define\modustwo{35}
 \define\poismod{36}
 \define\smooth{37}
    \define\srni{38}
\define\claustha{39}
\define\oberwork{40}
\define\qr{41}
\define\severi{42}
\define\junigusa{43}
\define\inonwign{44}
\define\jacobexc{45}
\define\jacobstw{46}
\define\kemneone{47}
\define\kirwaboo{48}
\define\kobaswan{49}
\define\kodnispe{50}
\define\kostasix{51}
\define\kostafou{52}
\define\kraprotw{53}
\define\kralyvin{54}
\define\kronhtwo{55}
\define\kumnarra{56}
\define\lawmiboo{57}
\define\lermonsj{58}
\define\loobotwo{59}
\define\loobothr{60}
\define\marsrati{61}
\define\marswein{62}
\define\mccrione{63}
\define\mehtsesh{64}
\define\mizuntwo{65}
\define\muruschi{66}
\define\narasesh{67}
\define\naramntw{68}
\define\naramnth{69}
\define\neherone{70}
\define\neherboo{71}
 \define\nessone{72}
 \define\ohtaone{73}
\define\ortratwo{74}
\define\peteraci{75}
\define\rinehart{76}
\define\satakboo{77}
\define\satokimu{78}
\define\gwschwar{79}
\define\gwschwat{80}
\define\sekiguch{81}
\define\sjamatwo{82}
\define\sjamlerm{83}
\define\slodoboo{84}
\define\spristei{85}
\define\varoucha{86}
\define\vergnsix{87}
\define\vinberg{88}
\define\weinstwo{89}
\define\weinsone{90}
\define\weylbook{91}
\define\whitnboo{92}
\define\woodhous{93}
\topmatter
\title K\"ahler spaces, nilpotent orbits, and singular reduction
\endtitle
\author Johannes Huebschmann
\endauthor
\affil 
Universit\'e des Sciences et Technologies
de Lille
\\
U. F. R. de Math\'ematiques
\\
CNRS-UMR 8524
\\
F-59 655 VILLENEUVE D'ASCQ, France
\\
Johannes.Huebschmann\@AGAT.UNIV-LILLE1.FR
\endaffil
\date{August 12, 2002}
\enddate
\abstract{
For a stratified symplectic space, a suitable concept of stratified 
K\"ahler polarization  
encapsulates K\"ahler polarizations on the strata and the behaviour of the 
polarizations across the strata and leads to the notion of stratified 
K\"ahler space which 
establishes an intimate relationship between nilpotent 
orbits, singular reduction, invariant theory, reductive dual pairs, 
Jordan triple systems,  symmetric domains, and pre-homogeneous spaces:
The closure of a holomorphic nilpotent orbit 
or, equivalently, the closure of the 
stratum of the associated pre-homogeneous space of parabolic type 
carries a (positive) normal
K\"ahler structure. 
In the world of singular Poisson geometry, 
the closures of principal holomorphic nilpotent orbits,
positive definite hermitian JTS's, and
certain pre-homogeneous spaces
appear as different incarnations of the same structure.
The closure of the principal
holomorphic nilpotent orbit arises from a 
semisimple holomorphic orbit by contraction.
Symplectic reduction carries 
a positive K\"ahler manifold to a positive normal K\"ahler space in such a way 
that the sheaf of germs of polarized functions  coincides with the 
ordinary sheaf of germs of holomorphic functions. 
Symplectic reduction 
establishes a close relationship between singular reduced spaces and nilpotent 
orbits of the dual groups. Projectivization of 
holomorphic nilpotent orbits yields exotic (positive) stratified K\"ahler 
structures on complex projective spaces and on certain complex projective 
varieties including complex projective quadrics. 
The space of (in general twisted) representations of the fundamental group of a 
closed surface in a compact Lie group or,
equivalently, a moduli space of central Yang-Mills connections
on a principal bundle over a surface,
inherits a (positive) normal (stratified)
K\"ahler structure. 
Physical examples are provided
by certain reduced spaces arising from angular momentum zero.}
\endabstract
\address{\smallskip
\noindent
USTL, UFR de Math\'ematiques, CNRS-UMR 8524, 
\newline\noindent
F-59 655 Villeneuve d'Ascq C\'edex,
France
\newline\noindent
Johannes.Huebschmann\@agat.univ-lille1.fr}
\endaddress
\endkeywords
\endtopmatter
\document
\leftheadtext{Johannes Huebschmann}
\rightheadtext{K\"ahler spaces, nilpotent orbits, reduction}
\beginsection  Introduction

Given a quantizable system with constraints,
the issue arises whether 
reduction after quantization
coincides with quantization after reduction.
Reduction after quantization makes sense
provided the unreduced phase space is a quantizable
smooth symplectic manifold in such a way that the symmetries
can be quantized as well; it
may then be studied within the usual framework of
geometric quantization. Up to now, the available methods 
have been insufficient to attack the problem of
quantization of reduced observables, though,
once the reduced phase space is no longer a smooth manifold;
we will refer to this situation as the {\it singular case\/}.
The singular case is the rule rather than the exception.
For example, simple classical mechanical systems and
the solution spaces of classical field theories involve singularities;
see e.~g. \cite\armcusgo\ and the references there.
In the presence of singularities,
the naive restriction of the quantization problem 
to a smooth open dense part, the \lq\lq top stratum\rq\rq, 
may lead to a loss of information and in fact to
inconsistent results; see what is said below.
Trying to overcome
these difficulties on the classical level,
we were led to isolate a certain class
of \lq\lq K\"ahler spaces with singularities\rq\rq,
to which the present paper is devoted.
Examples thereof arise from symplectic reduction, applied to
K\"ahler manifolds.
Other examples arise from certain moduli spaces,
of central Yang-Mills connections over a closed surface as well as
(equivalently)
of suitable representations
of the fundamental group of a closed surface in a compact Lie group.
Still
other examples arise from taking the closure of a 
holomorphic  
nilpotent orbit
in a real semisimple 
(or reductive)
Lie algebra
of hermitian type
(see what is said below).
\smallskip
A typical example of the kind of K\"ahler space with singularities
we will study is the complex plane with its ordinary complex analytic 
structure, but with a Poisson algebra of continuous functions
which has the origin as a singularity so that geometrically
this plane should then be viewed as a half-cone rather than the ordinary plane;
the cone point,
i.~e. the singularity, cannot be detected by the complex analytic structure.
Albeit looking like a toy example, it illustrates 
the general phenomenon that, under such circumstances,
playing off against each other the complex and real (semi-analytic) 
structures yields geometric insight.
Actually, this \lq\lq exotic\rq\rq\ plane is the closure 
in $\fra{sp}(1,\Bobb R)$ of what we call a {\it holomorphic\/}
nilpotent orbit: A real semisimple Lie algebra $\fra g$
with Cartan decomposition $\fra g = \fra k \oplus \fra p$ 
together with an element $z$ in the center of $\fra k$,
referred to as an $H$-{\it element\/},
such that the restriction of $\roman{ad}(z)$ to $\fra p$ yields 
a complex structure on the latter is  said to be of {\it hermitian type\/}
\cite\satakboo; the classification of semisimple Lie algebras of hermitian 
type plainly parallels E. Cartan's classification
of semisimple hermitian symmetric spaces. We will say that an orbit $\Cal O$ 
in such a Lie algebra $\fra g$ is {\it holomorphic\/}
provided the projection from $\Cal O$ to $\fra p$ is a diffeomorphism onto its 
image in such a way that the resulting complex structure 
$J_{\Cal O}$ on $\Cal O$ combines 
with its (Kostant-Kirillov-Souriau) symplectic structure
to a positive K\"ahler structure. 
Here $\fra g$ is identified with its dual by means of an appropriate
positive multiple of the Killing form, made precise later in the paper.
The holomorphicity of an orbit does {\it not\/} depend on 
the choice of $H$-element or, equivalently, on the Cartan decomposition;
however, the induced complex structure on the orbit {\it does\/}
depend on that choice; see (3.7.5) and (3.7.6) below.
The name \lq\lq holomorphic\rq\rq\ 
is intended to hint at the fact that
the holomorphic discrete series representations arise from holomorphic 
quantization on integral semisimple holomorphic orbits $\Cal O$
(the requisite complex structure being different from 
$J_{\Cal O}$, though, see what is said below). 
Indeed, for any real positive $\lambda$, the $G$-orbit in $\fra g$
generated by $\lambda z$ is semisimple and holomorphic; 
it realizes the smooth manifold underlying the corresponding
{\it symmetric domain\/} as a symplectic manifold. In Section 3 we will prove
that a simple Lie algebra of hermitian type and split rank $r$
has exactly $r+1$ {\it holomorphic nilpotent\/} orbits
$\Cal O_0,\dots, \Cal O_r$ which are linearly ordered
in such a way that $\{0\}=\Cal O_0 \subseteq \overline{\Cal O_1}
\subseteq \ldots \subseteq \overline{\Cal O_r}$
and that the projection from $\overline{\Cal O_r}$
to the constituent $\fra p$ in the Cartan decomposition is a homeomorphism.
This entails the fact that every orbit in the closure  $\overline{\Cal O}$
of a holomorphic nilpotent orbit is itself holomorphic and that
the projection from $\overline{\Cal O}$ to $\fra p$ is a homeomorphism
onto its image which, via the complex structure on  $\fra p$
induced by 
the $H$-element $z$ (which is part of the hermitian type structure), 
turns $\overline{\Cal O}$ into a complex affine variety.
The \lq\lq exotic\rq\rq\ plane mentioned earlier is a very special case 
thereof. In general, $\overline{\Cal O_r}$ may be seen as arising from the 
corresponding semisimple holomorphic orbits by contraction;
this contraction does not change the complex analytic structure.
A comment about the complex structures is perhaps in order:
The complex structure $J_{\Cal O}$ on a semisimple holomorphic
orbit $\Cal O$ arises from the projection to $\fra p$ as does
the complex structure on $\overline{\Cal O_r}$, whence complex analytically 
these spaces are just a copy of $\fra p$. However, 
for a semisimple orbit, this is {\it not\/} the hermitian symmetric
space complex structure.
\smallskip
Our K\"ahler spaces with singularities are stratified symplectic spaces with 
additional structure. Reduced spaces stratified into smooth symplectic 
manifolds (which arise as symmetry types) occur already in
\cite\armamonc\ (where  smooth infinite dimensional pieces are not excluded).
The notion
of {\it stratified symplectic structure\/}
on a space $X$ 
which we will use is that of
\cite\sjamlerm; it is a stratification 
into smooth finite dimensional manifolds
(in a  sense stronger than that of \cite\armamonc)
together with
a  Poisson algebra $(C^{\infty}(X),\{\cdot,\cdot\})$ of continuous functions
on $X$ which, on each stratum, restricts to a symplectic Poisson algebra
of smooth functions, not necessarily consisting of all smooth functions on 
the stratum. We deliberately write $C^{\infty}(X)$ even though this algebra 
need not be an algebra of ordinary smooth  functions;
the algebra $C^{\infty}(X)$ is then referred to as a {\it smooth\/}
structure on $X$. A {\it stratified symplectic space\/}
is a stratified space $X$ together with a stratified symplectic structure
$(C^{\infty}(X),\{\cdot,\cdot\})$, referred to henceforth
as the {\it stratified symplectic Poisson algebra\/}
(of the stratified symplectic space). In \cite\sjamlerm\ it is shown that
the reduced space for a hamiltonian action of a compact Lie group on a smooth 
symplectic manifold inherits a stratified symplectic structure,
the requisite stratified symplectic Poisson algebra being that given in 
\cite\armcusgo. Stratified symplectic structures have been constructed
on various moduli spaces \cite{\guhujewe,\,\poisson--\poismod}
including moduli spaces of central Yang-Mills connections
for bundles on a surface and  spaces
of (in general twisted) 
representations of a surface group in a compact Lie group
(where \lq\lq twisted\rq\rq\ means suitable representations of
the universal central extension), 
and the stratified symplectic structure
on such a space yields a description of its symplectic singularity 
behaviour; see \cite{\srni--\oberwork} for overviews.
In the present paper we will explore  stratified symplectic spaces arising 
from holomorphic nilpotent orbits and from reduction of K\"ahler manifolds;
in physics, examples come from standard hamiltonian systems in mechanics,
in particular from angular momentum.
The relationship between angular momentum zero and certain nilpotent orbits 
in the 
real symplectic Lie algebras
has been observed already in \cite\lermonsj.
In fact, the nilpotent orbits 
described explicitly in that paper
are precisely the holomorphic ones but this was not observed there.
\smallskip
The Poisson structure of a stratified symplectic space
encapsulates the mutual positions of the symplectic structures on the strata.
The {\it new key notion\/} which we will introduce in this paper is that of 
{\it stratified polarization\/}, which likewise 
{\it encapsulates the mutual positions of polarizations on the strata\/}.
We will show that, for certain nilpotent orbits $\Cal O$ of a semisimple 
Lie algebra of hermitian type including the holomorphic ones (already 
mentioned above), the closure $\overline {\Cal O}$ carries what we will call
a {\it complex analytic stratified K\"ahler polarization\/}.
Furthermore, for a (positive) K\"ahler manifold with a hamiltonian action
of a compact Lie group $G$ which preserves the K\"ahler  polarization
and extends to an action of the complexification of $G$, the reduced space 
inherits a complex analytic stratified K\"ahler polarization.
In general, we refer to a space with a complex analytic stratified K\"ahler 
polarization as a {\it complex analytic stratified K\"ahler space\/}.
We also introduce the somewhat more general notion
of {\it stratified K\"ahler polarization\/} (not necessarily complex analytic);
on each stratum, this kind of polarization
boils down to an ordinary K\"ahler polarization, and a space equipped
with a stratified K\"ahler polarization
will be referred to as a {\it stratified K\"ahler space\/}.
While a K\"ahler polarization on an ordinary smooth symplectic manifold
is well known to be equivalent to a complex structure 
for which the symplectic structure is a K\"ahler form, we do not know whether 
a stratified K\"ahler polarization (in our sense) on a stratified symplectic 
space necessarily determines a corresponding complex analytic structure.
\smallskip
The stratification of a complex analytic stratified K\"ahler space
is in general finer than the ordinary complex analytic stratification
and the stratified symplectic structure then exhibits singularities
which are complex analytically spurious; an illustration is the complex plane 
with a cone point mentioned above. More striking examples are certain 
{\it exotic\/} complex projective spaces. By an {\it exotic\/} complex 
projective space we mean complex projective space $\Bobb C \roman P^\ell$
with its ordinary complex structure, but with a {\it non-standard\/} stratified
symplectic structure $(C^{\infty}(\Bobb C \roman P^\ell),\{\cdot,\cdot\})$
which is compatible  with the complex structure
where \lq\lq non-standard\rq\rq\ means that there are at least two strata.
For example, complex projective 3-space has an exotic structure whose
singular locus is a Kummer surface, and the stratified Poisson structure
detects the Kummer surface {\it and\/} its 16 singular points; see 4.3 below.  
Complex projective 3-space has yet another exotic structure,  having as its 
singular locus a complex quadric; see Section 10 below. Such exotic spaces 
even arise from classical constrained systems in mechanics:
The reduced space $Q_{\ell}$ (say) of $\ell$ harmonic oscillators in some
$\Bobb R^{\ell}$ with total angular momentum zero and constant energy is 
an exotic complex projective space of complex dimension 
$d=\frac {\ell (\ell+1)}2 -1$ in such a way that, for $1 \leq s \leq\ell$,
the reduced spaces $Q_s$ (say) of $\ell$ harmonic oscillators in some
$\Bobb R^s$ with total angular momentum zero and (same) constant energy
embed naturally into  $Q_{\ell}=\Bobb C \roman P^d$
and constitute an ascending chain
$Q_1 \subseteq Q_2 \subseteq \dots \subseteq Q_{\ell}$
of (compact) complex analytic stratified K\"ahler spaces which, as complex 
analytic spaces, are projective determinantal varieties, and such that the 
$Q_j \setminus Q_{j-1}$ are the strata of the decomposition of any $Q_s$
($j \leq s$).
These examples are particular cases of systematic classes of exotic stratified 
K\"ahler structures on complex projective spaces which will be given in 
the paper. The question whether there exist exotic structures on complex 
projective spaces is also prompted by the following observation:
For a compact K\"ahler manifold $N$ which is complex analytically a projective 
variety (i.~e. Hodge manifold), with reference to the standard structure
on complex projective space (i.~e. Fubini-Study metric) the Kodaira embedding 
will not in general be symplectic, but complex projective space might still 
carry an exotic structure which, via the Kodaira embedding, restricts to 
the K\"ahler structure on $N$. More generally, given a Hodge manifold $N$
together with an appropriate group of symmetries and momentum mapping,
reduction carries it to a complex analytic stratified K\"ahler space 
$N^{\roman{red}}$ which is as well a projective variety, and the question 
arises whether in such a situation the requisite complex projective space 
carries an exotic structure which, via the Kodaira embedding, restricts to the 
complex analytic stratified K\"ahler structure on $N^{\roman{red}}$.
Our approach in particular demonstrates that this situation actually occurs 
in \lq\lq mathematical nature\rq\rq. The particular case $Q_s \subseteq 
\Bobb C \roman P^d$ ($1 \leq s \leq \ell$, $d=\frac {\ell (\ell+1)}2 -1$)
arising from angular momentum is a good illustration.
\smallskip
We close the introduction with a %brief 
guide through the paper.
In Section 1 we reproduce some material from the relationship between
Poisson algebras and Lie-Rinehart algebras. In Section 2 we introduce 
stratified polarized spaces. The aim of Section 3 is to show that 
{\sl the closure of any holomorphic nilpotent orbit carries a 
normal
complex 
analytic stratified K\"ahler structure\/};
see Theorem 3.2.1 for details.
We will, furthermore, give a {\sl complete classification of all holomorphic 
nilpotent orbits, and
we will explicitly describe the
real (semi-algebraic) and the complex analytic structure
of the closure of any such holomorphic nilpotent orbit\/}.
We will also classify what we call {\it pseudoholomorphic\/} nilpotent orbits;
a pseudoholomorphic nilpotent orbit $\Cal O$ is one which under the projection 
to the constituent $\fra p$ of the Cartan decomposition is still mapped 
diffeomorphically onto its image (but the resulting K\"ahler structure
is not constrained to be positive). Along the way, we obtain several results on
holomorphic nilpotent orbits that look interesting in their own right:
For a classical Lie algebra $\fra g$, let $E$
be its standard representation,
and consider the ordinary {\it associative\/}
algebra
$\roman{End}(E)$ of endomorphisms thereof.
We will show that, in the standard cases ($\fra g = \fra{sp}(\ell,\Bobb R),
\fra{su}(p,q), \fra{so}^*(2n)$), an element $X$ of $\fra g$ generates 
a pseudoholomorphic nilpotent orbit if and only if $X^2$ is zero in 
$\roman{End}(E)$ and that, likewise, 
an element $X$ of $\fra{so}(2,q)$  generates 
a pseudoholomorphic nilpotent orbit if and only if $X^3$ is zero in 
$\roman{End}(E)$.
The Lie algebra $\fra g$ sits of course inside
$\roman{End}(E)$ 
and an element $X \in \fra g$ satisfying $X^2=0$ or $X^3=0$
is plainly nilpotent 
(i.~e. $\roman{ad}(X)$ is a nilpotent endomorphism of $\fra g$)
but the requirement that $X \in \fra g$ satisfy 
$X^2=0$ or $X^3=0$ has no direct meaning within $\fra g$ itself. 
Another interesting 
side result is the observation that,
when $\fra g$ is simple,
{\sl a pseudoholomorphic nilpotent orbit
$\Cal O$ is actually holomorphic (or antiholomorphic, which means that 
$-\Cal O$ is holomorphic) if and only if the projection from the closure 
$\overline{\Cal O}$ to $\fra p$ is a homeomorphism onto its image in\/} 
$\fra p$; see Theorem 3.2.1.
Our results are conclusive,
except that we do not 
give explicit
defining equations
for the complex analytic structures
of the closures of the rank 1 holomorphic nilpotent 
orbits in the two exceptional cases.
In \cite\vergnsix, 
given a Lie algebra of hermitian type $(\fra g,z)$,
a nilpotent orbit $\Cal O$ 
in $\fra g$ 
is defined to be holomorphic provided it corresponds, under the 
Kostant-Sekiguchi correspondence \cite\sekiguch, to an orbit in
the holomorphic constituent $\fra p^+$ of the complexification 
$\fra p^{\Bobb C} =\fra p^+ \oplus \fra p^-$, and it is then proved, via 
the Kronheimer flow \cite\kronhtwo, that the restriction of the projection 
to $\fra p$ of such an orbit is a diffeomorphism onto its image,
that is to say, that such an orbit is (pseudo)holomorphic in our sense.
The approach in the present paper is somewhat different, though: We 
{\it define\/} the property of holomorphicity of an orbit in terms of the 
projection to $\fra p$---this definition is more elementary 
than that involving the Kostant-Sekiguchi correspondence---and we classify 
{\it all\/} holomorphic nilpotent orbits (in our sense) by means of a 
geometric description which bypasses the Kostant-Sekiguchi correspondence; 
in particular, our approach does {\it not\/} involve the Kronheimer flow at all
and is therefore likely to allow for generalization over other (formally real) 
fields. It also enables us to show that 
{\sl the projection from a
holomorphic nilpotent orbit\/} $\Cal O$ {\sl to\/} 
$\fra p$ {\sl extends to a homeomorphism 
from the closure\/} $\overline{\Cal O}$ {\sl onto its image in\/} $\fra p$. 
We do not know whether
this property 
of holomorphic nilpotent orbits
can formally be deduced from properties of the Kronheimer flow since 
passing to the closure amounts to changing the topological type of the 
corresponding bundle on $S^4$. 
This 
property of holomorphic nilpotent orbits
plainly
implies the fact that,
for a holomorphic nilpotent orbit, the
Kostant-Sekiguchi correspondence is compatible with passing to the 
closure---this is actually true for any nilpotent orbit
in a real semisimple Lie algebra, cf. \cite\barbsepa---but
is considerably stronger
than 
just the 
compatibility 
of the Kostant-Sekiguchi correspondence with 
respect to the closure.
In Section 4 we show that {\sl reduction carries an 
ordinary (positive) 
K\"ahler polarization to a (positive) complex analytic stratified K\"ahler 
polarization\/}. 
As an application,
we prove that the {\sl moduli spaces mentioned earlier
inherit normal (complex analytic stratified)  K\"ahler space
structures\/}. 
Classically, such a structure was obtained only in the non-singular
case, where these spaces are smooth K\"ahler manifolds.
Our approach bypasses the usual geometric invariant theory construction
(which does not give the requisite stratified symplectic structure)
and yields an analytic construction of these moduli spaces.
A particular 
case is the space $N$ of representations of the fundamental group of
a Riemann surface of genus two in $\roman{SU}(2)$; in view of results of 
Narasimhan-Ramanan \cite{\naramntw,\,\naramnth}, this space, realized as 
the moduli space of semistable holomorphic vector bundles of rank two,
degree zero, and trivial determinant, is  complex analytically a copy of
$\Bobb C \roman P^3$, and the subspace of classes of semistable holomorphic 
vector bundles  which are not stable is a Kummer surface. We will show 
that, indeed, $N$ carries a normal (complex analytic stratified) K\"ahler 
structure having the Kummer surface as its singular locus. This is the 
exotic structure on $\Bobb C \roman P^3$ mentioned earlier, and the 
Kummer surface, being viewed as a Kodaira embedding, inherits its 
(now singular) normal K\"ahler structure from the embedding. In 
Section 5 we describe the holomorphic nilpotent orbits (for the classical cases)
as reduced spaces for suitable momentum mappings; this involves versions of 
the {\it First and Second Main Theorem of Invariant Theory\/} \cite\weylbook;
the Second Main Theorem of Invariant Theory will actually be a consequence 
of our approach to invariant theory. Our tools are the original observation 
of Kempf-Ness \cite\kemneone\ and an extension of the {\it basic construction\/}
in \cite\kraprotw; in that reference, this construction is given only for 
the case of the classical groups over the complex numbers. We extend 
this construction to the appropriate real forms of the classical groups.
The Kempf-Ness observation is then exploited to identify a holomorphic 
nilpotent orbit as a symplectic quotient with the corresponding 
complex categorical quotient (in the sense of geometric invariant theory).
\smallskip
The construction in Section 5 has a consequence for {\it singular reduction\/}
which is interesting in its own right: Let $\KK$ be a compact Lie group which 
may be written as the direct product of copies of classical groups of the kind
$\roman O(s_j,\Bobb R)$, $\roman {Sp}(s_k)$, $\roman U(s_m)$;
for each of these, let $E_{s_j}$, $E_{s_k}$, $E_{s_m}$
be the standard unitary representation, $\roman O(s_j,\Bobb R)$ 
being viewed as a subgroup of $U(s_j)$ in the usual way. Pick natural numbers
$\ell_j$, $\ell_k$, $\ell_m$, $\overline {\ell_m}$, consider the 
representations $E_{s_j}^{\ell_j}$, $E_{s_k}^{\ell_k}$, $E_{s_m}^{\ell_m}$,
$\overline E_{s_m}^{\overline{\ell_m}}$ and, finally, take the sum
$E$ of them. This is a unitary $\KK$-representation.  Here $\overline E_{s_m}$
is the conjugate representation of $E_{s_m}$. The 
representations $E_{s_j}$ 
(being real)
and $E_{s_k}$ (being quaternionic) 
are actually self-conjugate whence there is no need to take conjugate 
representations for the $\roman O(s_j,\Bobb R)$'s, nor for 
the $\roman {Sp}(s_k)$'s.
Let $\mu \colon E \to \fra \kk$ be the unique $\KK$-momentum mapping of $E$ 
having the value zero at the origin, and let $\fra g$ be the real reductive 
Lie algebra which arises, accordingly, as a sum copies of 
the $\fra {sp}(\ell_j,\Bobb R)$'s, $\fra {so}^*(2\ell_k)$'s,
$\fra u(\ell_m,\overline{\ell_m})$'s.
Corollary 5.4.4 below says that, under these circumstances, {\sl the
$\KK$-reduced space $E^{\roman{red}}$ is as a normal K\"ahler space 
isomorphic to the closure of a holomorphic nilpotent orbit in $\fra g$.\/}
Furthermore, {\sl every holomorphic nilpotent orbit in $\fra g$ arises in this 
way.\/} 
It is interesting to point out, cf. \cite\howetwo, that
the groups
$\KK$ and the direct product $G$ 
of the corresponding copies of $\roman {Sp}(\ell_j,\Bobb R)$,
$\roman {SO}^*(2\ell_k)$, $\roman U(\ell_m,\overline{\ell_m})$,
constitute a real {\it reductive dual pair\/}
in $\roman{Sp}(E)$, the symplectic structure on $E$ being determined
by its unitary one. More generally, let $H$ be a Lie group which may be 
written as the direct product of copies of classical groups of the kind
$\roman O(s'_j,s''_j)$, $\roman {Sp}(s'_k,s''_k)$, $\roman U(s'_m,s''_m)$,
and let $E$ be an $H$-representation which is a sum of standard representations.
In view of \cite\howetwo, 
the groups $G$ and $H$ then still constitute a reductive 
dual pair in $\roman{Sp}(E)$, the symplectic structure on $E$ being 
the obvious one on $E$ determined by the data (see Section 5 below for 
details). Let $\mu \colon E \to \fra h$ be the unique $H$-momentum mapping
of $E$ having the value zero at the origin. 
Corollary 5.5.4 below says that
{\sl the $H$-reduced space\/} $E^{\roman{red}}$ 
(which we define to be that of {\it closed\/} $H$-orbits in $\mu^{-1}(0)$)
{\sl is as a stratified space isomorphic
to a space $N$ which is a union of pseudoholomorphic nilpotent orbits
in a Lie algebra $\fra g$ of the above kind in such a way  that the map from
$E^{\roman{red}}$ to $N$ is a Poisson map. Furthermore,
every pseudoholomorphic nilpotent orbit of $\fra g$ arises 
in this way (as a stratum of a reduced space of the kind
$E^{\roman{red}}$). \/}
\smallskip
In Section 6 we give a construction which yields the holomorphic nilpotent 
orbits of $\fra{so}(2,q)$ by singular reduction with respect to the 
non-compact group $\roman {Sp}(1,\Bobb R)$; for arbitrary $q$, it does not 
seem possible to obtain these orbits by singular reduction with respect to 
a compact group, though. In Section 7  we relate holomorphic nilpotent 
orbits with positive definite hermitian J(ordan) T(riple) S(ystem)s,
and in Section 8 we give some relevant information for the exceptional
cases in terms of JTS's. To 
explain what we achieve at that stage, let $(\fra g,z)$ be a real reductive 
Lie algebra of hermitian type, with Cartan decomposition 
$\fra g = \fra k \oplus \fra p$ and complexification
$\fra g_{\Bobb C} = \fra p^- \oplus \fra k_{\Bobb C} \oplus \fra p^+$,
and let $G$ and $K$ denote the adjoint group and maximal compact subgroup 
with $\roman{Lie}(K) = \fra k$, respectively. 
We will establish the following facts, cf. Theorem 7.3 
for the classical cases and Theorems 8.4.1 and 8.4.2 for 
the exceptional ones:
{\sl Under the 
homeomorphism (induced by the projection from $\fra g$ to $\fra p$)
from the closure $\overline {\Cal O}$ of the principal holomorphic nilpotent 
orbit $\Cal O$ onto the complex vector space $\fra p^+$, the $G$-orbit 
stratification of $\overline {\Cal O}$ passes to the stratification of 
$\fra p^+$ by Jordan rank, with reference to the JTS-structure on $\fra p^+$,
and the latter stratification, in turn, coincides with the 
$K^{\Bobb C}$-orbit stratification of $\fra p^+$. In this fashion,
 the induced 
stratified symplectic Poisson structure on $\fra p^+$ turns the latter into a 
normal complex analytic stratified
K\"ahler space, 
and the Poisson structure
detects the stratification by Jordan rank.\/}
This reduces the study of holomorphic nilpotent orbits in $\fra g$ to that 
of the hermitian Jordan triple system $\fra p^+$: {\sl All geometric 
information of the holomorphic nilpotent orbits is encoded in that \/} JTS,
{\sl even the stratified symplectic Poisson structure since the latter
is determined
by the Lie bracket in $\fra g$, which can be reconstructed from the\/} JTS;
this parallels the well known fact that all geometric information of a 
symmetric 
domain is encoded in its JTS, cf. e.~g. \cite\mccrione. 
Our results entail in particular that the pre-homogeneous spaces
studied and classified in \cite\muruschi---these are precisely those
which arise from positive definite hermitian JTS's---carry more structure; 
indeed, such a space necessarily underlies the closure of the principal 
holomorphic nilpotent orbit (viewed as a 
complex analytic stratified K\"ahler space) of a semisimple Lie algebra of 
hermitian type. As a byproduct, we obtain the Bernstein-Sato polynomials
for the regular pre-homogeneous spaces under discussion; see Theorems 7.1 
and 8.4.1 for details. For example, for $\fra e_{7(-25)}$, 
this polynomial is Freudenthal's generalized $(3 \times 3)$-determinant 
over a complex split 
octonion algebra \cite\freudone\ or, equivalently, Jacobson's 
generic norm. In Section 9 we show how the 
{\sl closure of the principal holomorphic
nilpotent orbit arises from the corresponding 
semisimple holomorphic orbits
by contraction\/}.
In Section 10 we projectivize the holomorphic 
nilpotent orbits 
by means of {\it stratified symplectic reduction\/} with respect to 
the circle group. This gives systematic classes of examples of 
{\it exotic projective spaces\/}. 
In particular, the principal holomorphic nilpotent orbit in
$\fra{so}(2,q)$ is a copy of $\Bobb C^q$ having as singular locus in the 
sense of stratified K\"ahler spaces an affine quadric, and projectivization 
thereof yields the {\it exotic projective space\/} 
$\Bobb C \roman P^{q-1}$ having 
as {\it singular locus a projective quadric\/}. 
Thus the projective quadric inherits 
a K\"ahler, in fact Hodge structure whose symplectic constituent does 
{\it not\/} arise (via the embedding) from the Fubini-Study metric on 
projective space,
and the underlying real space arises from the non-zero non-principal
holomorphic nilpotent orbit in $\fra{so}(2,q)$ by projectivization.
 This is a particular example of a Kodaira embedding of the 
kind mentioned earlier, where the symplectic structure results
from a proper stratified symplectic Poisson algebra on projective space.
For the exceptional Lie algebra $\fra e_{7(-25)}$, we obtain an {\it exotic 
K\"ahler structure\/} on $\Bobb C \roman P^{26}$ having as 
{\it singular locus the 
projective cubic\/} defined by the generic norm over the octonions mentioned 
earlier, and this cubic hypersurface, in turn, inherits a complex analytic 
stratified K\"ahler structure with two strata
whose underlying real structure arises from the closure of the rank 2
holomorphic nilpotent orbit in
$\fra e_{7(-25)}$  by projectivization.
 Actually, this cubic has 
played a major role in the development of the exceptional Lie groups;
a version of it was studied by E. Cartan.
In the standard cases ($\fra g = \fra{sp}(\ell,\Bobb R), \fra{su}(p,q), 
\fra{so}^*(2n)$), via the principle of reduction in stages, the 
projectivized holomorphic nilpotent orbits may also be obtained by ordinary 
symplectic reduction, as explained in Theorem 10.1; 
somewhat amazingly, this does not seem to be possible for $\fra{so}(2,q)$ for 
general $q$, nor for $\fra e_{6(-14)}$ or $\fra e_{7(-25)}$; in these cases, 
we obtain the projectivized 
holomorphic nilpotent orbits only by {\it proper stratified symplectic 
reduction\/}, with respect to the circle group. 
Finally, in Section 11 we compare our notion
of stratified K\"ahler space with that of K\"ahler space with singularities
introduced by {\smc Grauert}~\cite\graueone\ and with subsequent refinements 
thereof \cite\varoucha, referred to in \cite\heihuclo\ as 
{\it stratified K\"ahlerien spaces\/}. There  is a substantial difference 
between complex analytic stratified K\"ahler spaces in our sense and 
stratified K\"ahlerien spaces, though. Our emphasis is on the Poisson structure,
and this is crucial for issues related with quantization; stratified 
K\"ahlerien spaces do not involve Poisson structures at all. 
\smallskip
Thus the notion of stratified K\"ahler space
establishes an intimate relationship between nilpotent 
orbits, singular reduction, invariant theory, reductive dual pairs, 
Jordan triple systems,  symmetric domains, and pre-homogeneous spaces;
in particular, in the world of singular Poisson geometry, 
the closures of principal holomorphic nilpotent orbits,
positive definite hermitian JTS's, and
certain pre-homogeneous spaces
appear as different incarnations of the same structure
and,
in the classical cases,
that structure arises, via invariant theory and singular reduction,
from a non-singular one 
phrased in terms of dual pairs.
This Poisson geometry is, perhaps,
already lurking behind {\smc Lie}'s theory of {\it function groups\/},
cf. the notion of dual pair of
Poisson mappings in \cite\weinstwo.
Somewhat amazingly, Jordan algebras (and generalizations 
thereof: JTS's)---which have been invented in search for new models of 
quantum mechanics (cf. e.~g. what is said in \cite\mccrione)---are 
now seen to illuminate certain classically reduced phase spaces 
(e.~g. that of finitely many harmonic oscillators with total angular 
momentum zero). 
The length of the paper
is  to some extent
explained 
by the
need for various different
descriptions of holomorphic nilpotent orbits:
To obtain a complete list thereof,
the relative root systems suffice;
to see that the list is complete,
though, we 
need the classification of nilpotent orbits
in the literature
in terms of
classical matrix realizations,
matrix realizations
are needed as well to explain the relationship with singular reduction
and,
at times, structural insight is more easily obtained from
matrix realizations than from root systems.
To describe the link with JTS's---the exceptional 
Lie algebras are actually constructed from JTS's---we 
reproduce a description of the requisite JTS's taylored to our purposes.
\smallskip
We conclude with a somewhat vague comment about normality, a notion which 
so far does not seem to have played a major role in symplectic or Poisson 
geometry.
The results of the present paper
indicate that, perhaps, normal K\"ahler spaces constitute the {\it correct\/}
category for doing K\"ahler reduction. While smoothness
is not in general preserved, in the 
circumstances
studied in our paper,
normality {\it is\/} preserved under reduction.
What remains to be done is to develop a general 
reduction procedure for (in general singular) normal K\"ahler spaces. 
\smallskip
In a follow-up paper \cite\qr,
we have extended the ordinary holomorphic quantization scheme
to stratified K\"ahler spaces.
Using this extension, we have  established the fact
that, indeed, in a suitable sense,
{\sl in the framework of K\"ahler quantization,
reduction after quantization is equivalent to
quantization after reduction, observables included\/}.
In particular, in that paper, explicit examples
are given which justify the claim made earlier
that, in
the presence of singularities, restricting quantization
to a smooth open dense part (the top stratum) leads
in general to a {\sl loss of information and in fact to
inconsistent results\/}.
In another follow-up paper \cite\severi, we have shown that
the geometry of certain {\it Severi\/} varieties may be elucidated
in terms of holomorphic nilpotent orbits.
\smallskip
I am much indebted to A.~Weinstein  and T.~Ratiu for discussions,
and for their encouragement to carry out the research program the 
present paper is part of. 
I am much indebted as well to
M.~Duflo for having 
introduced me into 
{\smc Satake}'s book \cite\satakboo\ 
and for having
suggested a possible relationship
between certain nilpotent orbits and Jordan algebras.
Thanks are also due to 
R.~Buchweitz,
R.~Cushman, 
J.~Hilgert, F.~Hirzebruch, F.~Lescure, O.~Loos, L.~Manivel, 
J.~Marsden, E.~Neher,
J.~Sekiguchi, P.~Slodowy, T.~Springer, J.~Stasheff, M.~Vergne, 
and R.~Weissauer for 
support at various stages of the project. I owe a special debt to 
{\smc Satake}'s book \cite\satakboo; without this book, the paper would 
not have been materialized in its present form. 
\bigskip
\centerline{Table of contents}
\medskip
\noindent
1. Poisson algebras and Lie-Rinehart algebras
\newline\noindent
2. Stratified polarized spaces
\newline\noindent
3. The closure of a holomorphic nilpotent orbit
\newline\noindent
4. Reduction and stratified K\"ahler spaces
\newline\noindent
5. Associated representations and singular reduction
\newline\noindent
6. Associated representations for the remaining classical case
\newline\noindent
7. Hermitian Jordan triple systems and pre-homogeneous spaces
\newline\noindent
8. The exceptional cases
\newline\noindent
9. Contraction of semisimple holomorphic orbits
\newline\noindent
10. Projectivization and exotic projective varieties
\newline\noindent
11. Comparison with other notions of K\"ahler space with singularities
\newline\noindent
References

\medskip\noindent{\bf 1. Poisson algebras and Lie-Rinehart algebras}
\smallskip\noindent
Let $R$ be a commutative ring with 1 taken as ground ring which, for the 
moment, may be arbitrary. For a commutative $R$-algebra $A$, we denote by 
$\roman{Der}(A)$ the $R$-Lie algebra of derivations of $A$, with its standard 
Lie algebra structure. Recall that a {\it Lie-Rinehart algebra\/} $(A,L)$ 
consists of a commutative $R$-algebra $A$ and  an $R$-Lie algebra $L$
together with an $A$-module structure on $L$ and an $L$-module structure 
$L \to \roman{Der}(A)$ on $A$ (that is, $L$ acts on $A$ by derivations), and 
these are required to satisfy two compatibility conditions modeled on the 
properties of the pair $(A,L)$, $A$ being the ring of smooth functions 
$C^{\infty}(N)$ on a smooth 
manifold $N$ and $L$ the 
Lie algebra  $\roman{Vect}(N)$ of smooth vector fields on $N$. Explicitly, 
these compatibility conditions read $a(\alpha(b)) = (a\alpha) b$ and
$[\alpha,a\beta] = \alpha(a) \beta + a [\alpha,\beta]$, where $a,b \in A$ 
and $\alpha,\beta \in L$. Given  an arbitrary commutative algebra $A$ over $R$,
another example of a Lie-Rinehart algebra is the pair $(A,\roman{Der}(A))$,
with the obvious action of $\roman{Der}(A)$ on $A$ and obvious $A$-module 
structure on $\roman{Der}(A)$. There is an obvious notion of morphism of
Lie-Rinehart algebras and, with this notion of morphism, Lie-Rinehart algebras 
constitute a category. 
Given a Lie-Rinehart algebra $(A,L)$ we shall occasionally refer to
$L$ as an $(R,A)$-{\it Lie algebra\/}.
More details may be found in 
{\smc Rinehart}~\cite\rinehart\ and in our papers \cite\poiscoho\ and 
\cite\souriau.
\smallskip
Let $(A,\{\cdot,\cdot\})$ be a Poisson algebra, and let $D_A$ be the $A$-module
of formal differentials of $A$. For $u,v \in A$, the assignment
$\pi (du,dv) = \{u,v\}$ yields an $A$-valued 2-form $\pi= \pi_{\{\cdot,\cdot\}}$
on $D_A$, the {\it Poisson\/} 2-{\it form\/}  for $(A,\{\cdot,\cdot\})$.
Its adjoint 
$$
\pi^{\sharp} 
\colon D_A
@>>> 
\roman{Der}(A) = \roman{Hom}_A(D_A,A)
\tag1.1$'$
$$
is a morphism of $A$-modules, and the formula
$$
[a du,b dv] = 
a\{u,b\} dv + b\{a,v\} du + ab d\{u,v\}
\tag1.2
$$
yields a Lie bracket $[\cdot,\cdot]$ on $D_A$, viewed as an $R$-module.
The $A$-module structure on $D_A$, the bracket  $[\cdot,\cdot]$,
and the morphism $\pi^{\sharp}$ of $A$-modules endow the pair $(A,D_A)$
with a Lie-Rinehart algebra structure in such a way that $\pi^{\sharp}$
is a morphism of Lie-Rinehart algebras. See \cite\poiscoho\ (3.8) for details.
We write $D_{\{\cdot,\cdot\}} = (D_A,[\cdot,\cdot],\pi^{\sharp})$ and,
for $\alpha \in D_A$, we write 
${\alpha^{\sharp} = \pi^{\sharp}(\alpha) \in \roman{Der}(A)}$.
\smallskip
We now take as ground ring that of the reals $\Bobb R$ or that of 
the complex numbers $\Bobb C$; we shall occasionally use the neutral notation
$\Bobb K$ for either of them. We shall consider spaces $N$ with an algebra
of continuous $\Bobb K$-valued functions, deliberately denoted by
$C^{\infty}(N)$, for example ordinary
smooth manifolds and ordinary smooth functions; 
such an algebra $C^{\infty}(N)$ will then be referred to as a 
{\it smooth\/} structure on $N$ and will be viewed as part 
of the structure,
and the pair $(N,C^{\infty}(N))$
will be referred to as a {\it smooth\/} space. 
A space may support {\it distinct\/} smooth 
structures, though. Given a space $N$ with a smooth structure $C^{\infty}(N)$,
we shall write $\Omega^1(N)$ for the space of formal differentials with
those differentials divided out that are {\it zero at each point\/},
cf. \cite\kralyvin.
More precisely:
Let $p \in N$,  and view $\Bobb K$ as a 
$C^{\infty}(N)$-module by means of the evaluation map which assigns to a 
function $f$ its value $f(p)$ at $p$; we shall occasionally denote this 
$C^{\infty}(N)$-module
by $\Bobb K_p$. 
A formal differential is {\it zero at the point\/}
$p$ of $N$
provided it passes to zero in
$(D_{C^{\infty}(N)})\otimes_{C^{\infty}(N)} \Bobb K_p$.
For example, 
when $N$ is an ordinary smooth  manifold and $f$ an ordinary smooth function,
in local coordinates $x_1,\dots,x_m$, a formal differential of the kind
$df - \sum \frac {\partial f} {\partial x_j} dx_j$
is non-zero but is zero at each point.
For a general smooth space $(N,C^{\infty}(N))$,
the induced $C^{\infty}(N)$-module morphism 
$$
\roman{Hom}_{C^{\infty}(N)}(\Omega^1(N),C^{\infty}(N)) 
\to \roman{Hom}_{C^{\infty}(N)}(D_{C^{\infty}(N)},C^{\infty}(N)) 
= \roman{Der}(C^{\infty}(N))
$$ 
is clearly an isomorphism. At a point $p$ of $N$, the object $\Omega^1(N)$
amounts to the ordinary space of differentials, in the following sense:
Denote by $\fra m_p$ the ideal  of functions in
$C^{\infty}(N)$  that vanish at the point $p$.
The $\Bobb K$-linear map which assigns to a function $f$ in
$\fra m_p$ its formal differential induces isomorphisms
$$
\fra m_p \big/\fra m^2_p
@>>>
(D_{C^{\infty}(N)})\otimes_{C^{\infty}(N)} \Bobb K_p
@>>>
\Omega^1(N)\otimes_{C^{\infty}(N)} \Bobb K_p
$$ 
of $\Bobb K$-vector spaces, and each of these may be taken as the space of 
{\it differentials at\/} $p$; over the reals, this space, in turn, may be 
identified with the {\it cotangent space\/} $\roman T_p^*N$. When 
$N$ is an ordinary smooth manifold,
$\Omega^1(N)$ amounts to the space of smooth sections of the cotangent bundle.
For a general smooth space $N$ over the reals, when
$A=C^{\infty}(N)= C^{\infty}(N,\Bobb R)$ is endowed with a Poisson structure,
the  formula (1.2) yields a Lie-bracket $[\cdot,\cdot]$ on the 
$C^{\infty}(N)$-module 
$\Omega^1(N)$ and the 2-form $\pi_{\{\cdot,\cdot\}}$ is still defined on 
$\Omega^1(N)$; its adjoint then yields an $C^{\infty}(N)$-linear map 
$$
\pi^{\sharp} 
\colon \Omega^1(N)
@>>> 
\roman{Hom}_{C^{\infty}(N)}(\Omega^1(N),C^{\infty}(N)) 
@>{\cong}>>\roman{Der}(C^{\infty}(N)). 
\tag1.1
$$

\proclaim{Proposition 1.3}
The $C^{\infty}(N)$-module 
structure on $\Omega^1(N)$, the bracket  $[\cdot,\cdot]$,
and the morphism  $\pi^{\sharp}$ of $C^{\infty}(N)$-modules endow the pair 
$(C^{\infty}(N),\Omega^1(N))$ 
with a Lie-Rinehart algebra structure  in such a way that
$\pi^{\sharp}$ is a morphism of Lie-Rinehart algebras.
\endproclaim

\demo{Proof}
Indeed the obvious projection map from $D_{C^{\infty}(N)}$ to $\Omega^1(N)$
is compatible with the structure. The non-triviality of the kernel of this 
projection map does not cause any problem. \qed
\enddemo

We shall write
$(\Omega^1(N),[\cdot,\cdot],\pi^{\sharp})$ as 
$\Omega^1(N)_{\{\cdot,\cdot\}}$.
When $N$ is an ordinary smooth manifold, the range 
$\roman{Der}(C^{\infty}(N))$
of the  adjoint map $\pi^{\sharp}$ from $\Omega^1(N)$ to 
$\roman{Der}(C^{\infty}(N))$
boils down to the space $\roman{Vect}(N)$ of smooth vector fields on $N$.
In this case, the Poisson structure on $N$ is symplectic, that is, arises 
from a (uniquely determined) symplectic structure on $N$, if and only if 
$\pi^{\sharp}$, which may now be written as a morphism 
of smooth vector bundles
from the cotangent bundle to the tangent bundle,
is an isomorphism. Thus the 2-form
$\pi_{\{\cdot,\cdot\}}$, which is defined for {\it every\/} Poisson algebra,
generalizes the symplectic form of a symplectic manifold;
see Section 3 of~\cite\poiscoho\  for details.

\medskip\noindent {\bf 2. Stratified polarized spaces}
\smallskip\noindent
The aim of this section is to develop a satisfactory notion of
K\"ahler polarization for a stratified symplectic space.
\smallskip
Let $N$ be a symplectic manifold, with symplectic Poisson algebra
$(C^{\infty}(N),\{\cdot,\cdot\})$, and write $\pi$ for its
Poisson 2-form.
Since $N$ is symplectic,
the complexification
$$
\pi_{\{\cdot,\cdot\}}^{\sharp}\otimes \Bobb C\colon
\Omega^1(N,\Bobb C)_{\{\cdot,\cdot\}} @>>> \roman {Vect}(N,\Bobb C)
\tag2.1
$$
of the adjoint 
$\pi_{\{\cdot,\cdot\}}^{\sharp}$
is an isomorphism of
$(\Bobb C, C^{\infty}(N,\Bobb C))$-Lie algebras.
Hence
a complex polarization 
(i.~e. integrable Lagrangian distribution)
$F \subseteq \roman T^{\Bobb C}N$ 
for $N$
corresponds to a certain 
$(\Bobb C, C^{\infty}(N,\Bobb C))$-Lie subalgebra
$P$ of
$\Omega^1(N,\Bobb C)_{\{\cdot,\cdot\}}$
which is just the pre-image under (2.1)
of the space $\Gamma F$ of sections of
$F$.
The {\it polarized\/} complex functions
determined by $F$ are
smooth complex valued functions $f$
defined on open sets of $N$ satisfying
$Zf=0$ for {\it every\/} smooth vector field
$Z$ in $\Gamma F$;
when 
$N$ is a K\"ahler manifold and
$F$ the holomorphic polarization,
these are ordinary holomorphic functions.
In general, non-zero polarized functions exist at most locally.
Following \cite\whitnboo, we shall refer to a continuous function $f$,
defined on an open subset of a space $X$, as a {\it function in\/} $X$.
For general $N$ and $F$,
given a function $f$ in $N$,
for a suitable 
compactly supported
smooth bump function $u$ on $N$
(the notion of bump function will be made precise below),
the differential $udf$ is defined everywhere on $N$.
Notice that we can extend $f$ to a smooth function 
$\widetilde f$ on $N$
in such a way that 
$udf = ud\widetilde f$, so
the differential $udf$ really lies in the $C^{\infty}(N)$-module
of differentials, for the algebra $C^{\infty}(N)$ of smooth functions 
on the whole space $N$.
Then $P$ is generated by these differentials.
A little thought reveals that,
for every such $f$ and $h$, with appropriate bump functions
$u$ and $v$
so that the differentials $udf$ 
and $vdh$
are defined everywhere on $N$,
the Lie bracket
(1.2) induces 
a bracket on $P$. This bracket is given by the expression
$$
[udf, vdh] = u\{f,v\} dh + v \{u,h\} df;
$$
notice that the third term $uvd\{f,h\}$
(which has initially to occur, cf. (1.2))
is zero since 
$f$ and $h$ are polarized.
\smallskip
Let  $X$ 
be a stratified space, and let $C^{\infty}(X)$
be  a smooth structure on $X$
which, on each stratum, restricts to an algebra of ordinary smooth
functions; 
henceforth all smooth spaces
coming into play in the paper
will be of this kind.
Given such a space $X$ and
an open subset $O$ thereof, with its induced smooth structure,
and a smooth function $f$ on $O$, we refer
to $f$ as a {\it smooth function in\/} $X$.
When $X$ is a stratified symplectic space,
we will occasionally refer to the stratification as the
{\it symplectic stratification\/}.
\smallskip
Let $X$ be a stratified symplectic space,
with complex Poisson algebra
\linebreak
$(C^{\infty}(X,\Bobb C),\{\cdot,\cdot\})$, let
$\Omega^1(X,\Bobb C)_{\{\cdot,\cdot\}}$ be the corresponding
$(\Bobb C,C^{\infty}(X,\Bobb C))$-Lie algebra, and
let
$Y\subseteq X$ be an arbitrary stratum.
Then the restriction map
from
$C^{\infty}(X,\Bobb C)$
to $C^{\infty}(Y,\Bobb C)$
is Poisson
and hence induces a morphism
$$
\left(C^{\infty}(X,\Bobb C),\Omega^1(X,\Bobb C)_{\{\cdot,\cdot\}}\right)
@>>>
\left(C^{\infty}(Y,\Bobb C),\Omega^1(Y,\Bobb C)_{\{\cdot,\cdot\}}\right)
\tag2.2
$$
of Lie-Rinehart algebras.
Given a 
$(\Bobb C, C^{\infty}(X,\Bobb C))$-Lie subalgebra
$P$ of
$\Omega^1(X,\Bobb C)_{\{\cdot,\cdot\}}$,
let $P_Y$ be the
$C^{\infty}(Y,\Bobb C)$-submodule of 
$\Omega^1(Y,\Bobb C)$
generated by the image of $P$ under (2.2);
it inherits the structure of a
$(\Bobb C, C^{\infty}(Y,\Bobb C))$-Lie subalgebra
of
$\Omega^1(Y,\Bobb C)_{\{\cdot,\cdot\}}$.
\smallskip
A  $(\Bobb C,C^{\infty}(X,\Bobb C))$-Lie subalgebra
$P$ of $\Omega^1(X,\Bobb C)_{\{\cdot,\cdot\}}$
will be said to be
a {\it stratified (complex) polarization\/}
for $X$ if, 
for every 
stratum $Y$,
under the isomorphism (2.1) of
$(\Bobb C, C^{\infty}(Y,\Bobb C))$-Lie algebras
(for the smooth symplectic manifold $Y$),
the $(\Bobb C, C^{\infty}(Y,\Bobb C))$-Lie subalgebra
$P_Y$ of
$\Omega^1(Y,\Bobb C)_{\{\cdot,\cdot\}}$
is identified with 
the space of sections of
a complex polarization in the usual sense.
Given a stratified polarization $P$ for $X$,
a smooth complex function
$f$ in
$X$, that is, a 
function which 
is defined on an open subset $O$ of $X$ and
belongs to the induced smooth structure $C^{\infty}(O)$
will be said to be {\it polarized\/}
provided it satisfies
$$
\alpha (f) = 0,\quad\text{for every}\quad \alpha \in P.
$$
Here and henceforth \lq\lq smooth\rq\rq\ refers to the structure algebras
of functions 
$C^{\infty}(X)$ on $X$ and $C^{\infty}(Y)$
on open subsets $Y$ thereof; thus a \lq\lq smooth  function\rq\rq\ 
in our sense is not necessarily
an ordinary smooth function.
In general polarized functions will exist at most locally, and we
can talk about 
{\it polarized functions in\/} $X$ or
the {\it sheaf of germs of
polarized functions\/}.
\smallskip
Infinitesimally, the notion of
stratified polarization comes down to this:
For $x \in X$,
the Poisson 2-form $\pi_{\{\cdot,\cdot\}}$
passes to a 2-form
$\Pi_x$ on the cotangent space $\roman T_x^*X
= \Omega^1(X)\otimes_{C^{\infty}(X)} \Bobb R_x$.
Given a stratified polarization,
when $x$ lies in the
stratum $Y$,
the restriction map
from
$$
{\roman T_x^*}^{\Bobb C}X =\roman T_x^*X \otimes\Bobb C
=
\Omega^1(X,\Bobb C) \otimes_{C^{\infty}(X,\Bobb C)} \Bobb C_x
$$
to
${\roman T_x^*}^{\Bobb C}Y
=\Omega^1(Y,\Bobb C)\otimes_{C^{\infty}(Y,\Bobb C)} \Bobb C_x$
induces a surjective $\Bobb C$-linear map
$$
P_x = P\otimes_{C^{\infty}(X,\Bobb C)} \Bobb C_x
@>>>
P_x(Y) = P_Y\otimes_{C^{\infty}(Y,\Bobb C)} \Bobb C_x,
$$
and $P_x(Y)$ is
a Lagrangian subspace 
of ${\roman T_x^*}^{\Bobb C}Y$
with respect to the complexified symplectic structure
$\pi^{\Bobb C}_x$ on ${\roman T_x^*}^{\Bobb C}Y$.
\smallskip
When $X$ is a smooth symplectic manifold,
viewed as a stratified symplectic space with a single stratum,
this notion of polarization manifestly boils down to the ordinary one.
For 
a general stratified symplectic space,
a stratified polarization
encapsulates the {\it mutual positions of the polarizations
on the strata\/}.
Given a stratified symplectic space
$X$,
a stratified polarization 
$P$ for $X$
will be said to be a (positive)
{\it stratified K\"ahler polarization\/}
if, for every stratum $Y$,
the image of $P_Y$ under (2.1)
is 
the space of sections of
a (positive) K\"ahler 
polarization for $Y$ 
in the usual sense.
A stratified symplectic space with
a stratified  K\"ahler polarization
will be said to be a {\it stratified K\"ahler space\/}.
Since a K\"ahler polarization on a symplectic manifold determines
a K\"ahler structure thereupon
it is obvious that
each stratum of a stratified K\"ahler space inherits a K\"ahler structure.
Recall that
a symplectic structure
$\omega$ is said to be {\it compatible\/} with a complex structure $J$
provided $J$ is a symplectic operator
so that $\omega$ defines a not necessarily positive
K\"ahler structure.

\proclaim{Proposition 2.3}
Given a real symplectic structure $\omega$ on
a complex domain $N$ of complex dimension
$n$, let $\{\cdot,\cdot\}$ be the Poisson structure determined by $\omega$.
In terms of holomorphic coordinates $w_1,\dots,w_n$,
$\omega$ is compatible with the complex structure $J$
if and only if
$\{w_j,w_k\} = 0$ for $1 \leq j,k \leq n$.
Under these circumstances, the 
K\"ahler polarization $P$ (in our sense)
is the free
$C^{\infty}(N,\Bobb C)$-submodule of
$\Omega^1(N,\Bobb C)_{\{\cdot,\cdot\}}$
generated by the differentials
$dw_1,\dots,dw_n$. For a smooth function $f$
in $N$,
the
Cauchy-Riemann equations then amount to 
$$
\{w_j,f\} = 0,\quad 1 \leq j \leq n.
$$
\endproclaim

\demo{Proof}
The differentials
$dw_1,\dots,dw_n,d\overline w_1,\dots,d\overline w_n$
freely generate $\Omega^1(N,\Bobb C)$
as a $C^{\infty}(N,\Bobb C)$-module,
and
$\roman {Vect}(N,\Bobb C)$
is, likewise, identifies the free 
$C^{\infty}(N,\Bobb C)$-module
generated by the vector fields
$\frac \partial {\partial w_1},\dots,
\frac \partial {\partial w_n},
 \frac \partial {\partial \overline w_1},\dots,
\frac \partial {\partial \overline w_n}$.
The  symplectic structure $\omega$ on $N$ is
compatible with the complex structure $J$ 
if and only if
the isomorphism (2.1)
of $(\Bobb C,C^{\infty}(N,\Bobb C))$-Lie algebras
from
$\Omega^1(N,\Bobb C)_{\{\cdot,\cdot\}}$
onto $\roman {Vect}(N,\Bobb C)$
identifies the
$C^{\infty}(N,\Bobb C)$-linear span 
of the
$dw_1,\dots,dw_n$
with the
$C^{\infty}(N,\Bobb C)$-linear span 
of the
$ \frac \partial {\partial \overline w_1},\dots,
\frac \partial {\partial \overline w_n}$
(so that this span contains the hamiltonian vector fields
$\{\cdot,w_j\}$ of the $w_j$'s ($1 \leq j \leq n$)
and therefore coincides with the span of these hamiltonian vector fields)
and, likewise, the
$C^{\infty}(N,\Bobb C)$-linear span 
of the
$d\overline w_1,\dots,d\overline w_n$
with the
$C^{\infty}(N,\Bobb C)$-linear span 
of the
$\frac \partial {\partial w_1},\dots,
\frac \partial {\partial w_n}$.
The latter property, in turn,
is equivalent
to the vanishing of $\{w_j,w_k\}$ for $1 \leq j,k \leq n$.
See for example the discussion in \cite\woodhous\ (5.4~p.~92).
Under these circumstances, the K\"ahler polarization (in our sense)
$P \subseteq \Omega^1(N,\Bobb C)_{\{\cdot,\cdot\}}$
is generated by the differentials
$dw_1,\dots,dw_n$,
and the 
{\it polarized\/}
smooth complex  functions $f$ in
$N$ 
or, equivalently, the germs of {\it polarized\/}
smooth complex  functions  
are 
(classes of)
functions $f$ in $N$ which satisfy $dw_j(f) = 0$ for $1 \leq j \leq n$,
where $dw_j(\cdot)$
refers to the action of
$dw_j \in \Omega^1(N,\Bobb C)_{\{\cdot,\cdot\}}$
on
$C^{\infty}(N,\Bobb C)$
reproduced in Section 1 above.
By construction, for every 
smooth complex function $h$ in $N$, 
given $f$, we have
$dh(f) = \{h,f\}$.
Hence the 
{\it polarized\/}
smooth complex  functions $f$ in
$N$ 
are those functions which satisfy
the equations
$\{w_j,f\} = 0$
for $1 \leq j \leq n$.
However,
for every coordinate function $w_j=q_j+ip_j$ and every 
smooth complex valued
function $f$ in $N$,
$$
\{w_j,f\} = \sum\{w_j, \overline w_k\} 
\frac{\partial f}{\partial \overline w_k}.
$$
Since $\omega$ is symplectic, the matrix
$\left[\{w_j, \overline w_k\} \right]$
is invertible
whence
the Cauchy-Riemann equations amount to
the vanishing of $\{w_j,f\}$  for $1 \leq j \leq n$. \qed
\enddemo

Under the circumstances of
Proposition 2.3,
the holomorphic functions in 
$N$
coincide with
the smooth complex valued  functions in $N$
which are polarized in our sense.
\smallskip
Let $X$ be a 
stratified K\"ahler space;
we will say that a complex analytic structure on $X$ is {\it compatible\/}
with the stratified K\"ahler structure if
(i) the stratification of $X$ 
(as a stratified symplectic space)
is a refinement of the complex analytic stratification,
if (ii) every holomorphic function in $X$ is smooth in $X$ and if,
(iii) the sheaf of germs of holomorphic functions
in $X$ is contained in its sheaf of germs of polarized functions.
Given a
stratified K\"ahler space
$X$ with a compatible complex analytic structure,
for any stratum,
endowed with 
the complex structure
coming from the induced K\"ahler polarization
on that stratum,
the restriction map from $X$
to that stratum
is plainly complex analytic. 

\proclaim{Proposition 2.4}
Given a stratified K\"ahler space $X$
together with a normal compatible complex analytic structure,
if the
stratification of $X$ 
(as a stratified symplectic space)
has an open and dense stratum,
the sheaf
of germs of holomorphic functions coincides
with that of polarized functions. 
\endproclaim

\demo{Proof}
The open and dense stratum
is necessarily contained in the top complex analytic stratum.
A continuous function
in $X$
which is holomorphic in 
the strata, in particular in the 
open and dense stratum,
is necessarily itself holomorphic,
since on a normal space
Riemann's extension theorem holds.
See e.~g. {\smc Whitney's} book \cite\whitnboo\ 
or {\smc Grauert-Remmert} \cite\grauremm\ for details.
Consequently a polarized smooth complex function in $X$
is necessarily holomorphic. \qed
\enddemo

We do not know whether for a general 
stratified K\"ahler space
having a compatible complex analytic structure
the sheaf
of germs of holomorphic functions actually coincides with that
of polarized functions,
nor do we know 
whether a stratified K\"ahler polarization
on a stratified symplectic space $X$
necessarily determines a compatible complex analytic structure
on $X$.
\smallskip
Borrowing some terminology from the language of sheaves,
we will say that a 
smooth space $X$, with smooth structure $C^{\infty}(X)$,
is {\it fine\/}
provided for an arbitrary locally finite open covering
$\Cal U=\{U_\lambda\}$ of $X$, there is a partition of unity
$\{u_\lambda\}$ subordinate to $\Cal U$
(i.~e. $\roman{supp}(u_{\lambda}) \subseteq U_{\lambda}$)
with $u_\lambda \in C^{\infty}(X)$
for every $\lambda$;
such a function $u_\lambda$ will then be referred to as
a {\it bump\/} function.
The notion of fine space will help us avoid having to talk
about the sheaf of germs of smooth functions.

\proclaim{Theorem 2.5}
Let $X$ be a fine stratified symplectic space, endowed with a
complex analytic structure.
Suppose that the symplectic stratification is a refinement of
the complex analytic
one and that germs of holomorphic functions belong to
the stratified symplectic structure, that is,
are smooth functions in $X$.
Then the
$C^{\infty}(X,\Bobb C)$-submodule $P$ of 
$\Omega^1(X,\Bobb C)_{\{\cdot,\cdot\}}$
generated by 
differentials of the kind $u df$
where $f$ is holomorphic in $X$ and $u$ a bump function is
a stratified K\"ahler polarization for $X$,
necessarily compatible with the complex analytic structure,
if and only if,
for every pair $f,h$ of holomorphic functions
in $X$,
the Poisson bracket $\{f,h\}$
vanishes.
\endproclaim

\demo{Proof}
Given differentials
$udf$ and $vdh$ where $u$ and $v$ are bump functions
and $f$ and $h$ holomorphic functions in $X$,
$$
[udf, vdh] = u\{f,v\} dh + v \{u,h\} df + uv d\{f,h\}.
\tag*
$$
Suppose that
$\{f,h\}$ vanishes for any pair of
holomorphic functions $f$ and $h$ in $X$.
The identity $(*)$ implies that 
the $C^{\infty}(X,\Bobb C)$-submodule $P$ of 
$\Omega^1(X,\Bobb C)_{\{\cdot,\cdot\}}$
is then closed under the Lie bracket.
Let $Y$ be a stratum of $X$.
We must show that  
the
$C^{\infty}(Y,\Bobb C)$-submodule 
$P_Y$
of 
$\Omega^1(Y,\Bobb C)$
generated by the image of $P$ under (2.2)
passes under (2.1) (where $Y$ plays the role of $N$
in (2.1))
to the space of 
sections of an ordinary
K\"ahler polarization on $Y$.
Since the problem is local, we may suppose that
$X$ is the zero locus of a finite set
$f_1,\dots,f_q$ of holomorphic functions on a polydisc $U$ 
in some $\Bobb C^n$, with coordinates $w_1,\dots,w_n$.
Restricted to the stratum $Y$,
the holomorphic functions
$w_1,\dots,w_n$
will not be independent
but, since $Y$ is a smooth complex
submanifold of $U$,
the
$C^{\infty}(Y,\Bobb C)$-submodule 
$P_Y$
of 
$\Omega^1(Y,\Bobb C)$
generated by
the differentials
$dw_1,\dots,dw_n$
coincides with the 
$C^{\infty}(Y,\Bobb C)$-submodule 
$P_Y$
of 
$\Omega^1(Y,\Bobb C)$
generated by the holomorphic differentials in $Y$.
More precisely, suitable holomorphic coordinate functions
$\zeta_1,\dots,\zeta_{\ell}$
on $Y$ may be written as holomorphic
functions in the
$w_1,\dots,w_n$,
whence the
$C^{\infty}(Y,\Bobb C)$-submodule 
$P_Y$
of 
$\Omega^1(Y,\Bobb C)$
generated by the 
differentials 
$dw_1,\dots,dw_n$
coincides with that
generated by the 
differentials 
$d\zeta_1,\dots,d\zeta_{\ell}$.
By virtue of Proposition 2.3,
the latter generate the
K\"ahler polarization
in our sense. 
\smallskip
Conversely,
when
$P$ is a stratified K\"ahler polarization for $X$,
given
holomorphic functions $f$ and $h$ in $X$,
in view of Proposition 2.3,
the bracket $\{f,h\}$ is zero on every stratum,
that is, $\{f,h\}$ is zero. \qed
\enddemo

Under the circumstances of Theorem 2.5,
the stratified
K\"ahler polarization $P$
is determined by the complex analytic structure;
if, on the other hand,
$P$
{\it determines\/} the complex analytic structure
in the sense that
the sheaf of germs of polarized functions actually {\it coincides\/} with that
of holomorphic functions,
$P$ will be said to be  
a {\it complex analytic stratified K\"ahler polarization\/},
and the corresponding space
will be referred to as a
{\it complex analytic stratified K\"ahler space\/}.
A complex analytic stratified K\"ahler space which,
as a complex analytic space, is normal
will be said to be a {\it normal K\"ahler space\/}.
Normal K\"ahler spaces constitute a particularly nice class
of stratified K\"ahler spaces.
All the examples of stratified K\"ahler spaces known to us are
normal K\"ahler spaces.

\smallskip
The following is now a mere observation; yet we spell it out
since it looks exceedingly attractive.

\proclaim {Proposition 2.6}
A function $f$
in a 
complex analytic stratified K\"ahler space $X$
is holomorphic if and only if
$\{w,f\} = 0$ for every holomorphic function $w$ in $X$.
Moreover, 
near any point $p$ of $X$,
the
requisite
complex analytic stratified K\"ahler polarization
$P$ has a finite set of $C^{\infty}(X,\Bobb C)$-module generators
$dw_1,\dots,d w_n$
where $w_1,\dots w_n$ are holomorphic functions
in $X$
defined near $p$, and
a function $f$ in $X$
is holomorphic (near $p$) if and only if
$$
\{w_j,f\} = 0,\quad 1 \leq j \leq n. \qed
\tag2.6.1
$$
\endproclaim

Under such circumstances,
the equations (2.6.1) may be viewed as {\it Cauchy-Riemann\/}
equations for the 
complex analytic stratified K\"ahler space $X$.
{\sl Thus the Cauchy-Riemann equations make sense
even though $X$ is not necessarily smooth\/}.
\smallskip

Henceforth every stratified symplectic space
will be assumed to be fine.
When $X$ is a smooth symplectic manifold,
viewed as a (fine) stratified symplectic space with a single stratum,
Theorem 2.5 comes down to the usual characterization
of a K\"ahler structure on $X$ in terms of a K\"ahler polarization.
However,
to attack the question whether
a stratified K\"ahler polarization
on a general stratified symplectic space $X$
necessarily determines a compatible complex analytic structure
on $X$,
one would have to 
develop a theory of obstructions to
realizing germs of stratified K\"ahler
spaces as complex analytic germs
in a compatible fashion.
In the smooth (i.~e. manifold) case,
there is no such realization problem
since manifolds are locally modeled
on affine spaces, and complex structures on
affine spaces are handled by standard linear algebra.

\beginsection 3. The closure of a holomorphic nilpotent orbit

In this section we will show that the closure of a holomorphic
nilpotent orbit inherits a normal complex analytic stratified
K\"ahler structure.
Moreover, we will give a complete classification of holomorphic
nilpotent orbits, and we will
explicitly describe the real and complex structures of the closure
of any holomorphic nilpotent orbit.
\smallskip\noindent
{\smc 3.1. Reductive Lie algebras of hermitian type.}
Following
\cite\satakboo\ (p.~54),
we define a (semisimple) Lie algebra of {\it hermitian type\/}
to be a pair $(\fra g, z)$ 
which consists of a real semisimple Lie algebra
$\fra g$ with a Cartan
decomposition $\fra g = \fra k \oplus \fra p$ 
and a central element
$z$ of $\fra k$,
referred to as an $H$-{\it element\/},
such that
$J_z = \roman{ad}(z)\big |_{\fra p}$
is a (necessarily $K$-invariant)
complex structure 
on $\fra p$.
Sometimes we will then say that $\fra g$ is of hermitian type,
without explicit reference to an $H$-element.
Slightly more generally,
a {\it reductive Lie algebra of hermitian type\/} 
is a reductive Lie algebra 
$\fra g$ together with an element $z \in \fra g$
whose constituent $z'$ (say) in the semisimple part 
$[\fra g,\fra g]$ of $\fra g$ is an $H$-element
for $[\fra g,\fra g]$
\cite\satakboo\ (p.~92).
Given two 
semisimple (or reductive) Lie algebras
$(\fra g_1, z_1)$ and
$(\fra g_2, z_2)$ 
of hermitian type,
a morphism $\phi \colon \fra g_1 \to \fra g_2$ of Lie
algebras 
is said to be an $H_1$-{\it homomorphism\/}
provided $\phi \circ \roman{ad}(z_1)= \roman{ad}(z_2)\circ \phi$. 
For a real semisimple 
Lie algebra $\fra g$,
with Cartan decomposition $\fra g = \fra k \oplus \fra p$,
we write $G$ 
for its adjoint group
and $K$ for the (compact) connected subgroup of $G$ with
$\roman{Lie}(K) = \fra k$;
the requirement that $\fra g$ be 
of hermitian type
is equivalent to
$G\big / K$ being a 
(non-compact)
hermitian
symmetric space.
Given $\fra g$, the space of $H$-elements or, equivalently,
the manifold of corresponding Cartan decompositions,
may be identified with 
a disjoint union $(G/K^+) \cup (G/K^-)$
of two copies of the
homogeneous space $G/K$ in an obvious fashion,
$G/K^+$ and  $G/K^-$ being the connected components containing
$z$ and $-z$, respectively.
For example, for $\fra g = \fra {sp}(1,\Bobb R) =\fra {sl}(2,\Bobb R)$,
$G/K^+$ and  $G/K^-$ are the two connected components of a standard
hyperboloid.
\smallskip
Given a real semisimple 
Lie algebra $(\fra g,z)$ of hermitian type,
the underlying Lie algebra $\fra g$ decomposes as 
$$
\fra g = \fra g_0 \oplus \fra g_1 \oplus \ldots  \oplus\fra g_k
\tag3.1.1
$$
where
$\fra g_0$ is the maximal compact semisimple ideal and where
$\fra g_1,\dots,  \fra g_k$ are non-compact and simple.
For a non-compact simple Lie algebra with Cartan decomposition
$\fra g = \fra k \oplus \fra p$,
the 
$\fra k$-action
on $\fra p$
coming from the adjoint representation of $\fra g$
is faithful and irreducible whence the center of $\fra k$
is then at most one-dimensional; indeed 
$\fra g$ has an $H$-element 
turning it into a Lie algebra of hermitian type
if and only if
the center of $\fra k$ has dimension one.
In view of E.~Cartan's
infinitesimal classification
of irreducible hermitian symmetric spaces,
the classical
Lie algebras 
$\fra{su}(p,q)$ ($A_n, n \geq 1$, where $n+1 = p+q$),
$\fra{so}(2,2n-1)$ ($B_n, n \geq 2 $),
$\fra{sp}(n,\Bobb R)$ ($C_n, n \geq 2$),
$\fra{so}(2,2n-2)$ ($D_{n,1}, n>2$),
$\fra{so}^*(2n)$ ($D_{n,2}, n >2$)
together with the real forms
$\fra{e}_{6(-14)}$ and $\fra{e}_{7(-25)}$
of the exceptional Lie algebras
$\fra{e}_{6}$ and $\fra{e}_{7}$,
respectively,
constitute a
complete list of
simple  Lie algebras of hermitian type.
We refer to  
$\fra{su}(p,q)$, 
$\fra{sp}(n,\Bobb R)$
and $\fra{so}^*(2n)$
as the {\it standard cases\/}.

\smallskip\noindent
{\smc 3.2. Holomorphic and pseudoholomorphic orbits.}
Let 
$(\fra g, z)$ be a semisimple Lie algebra of hermitian type.
Then,
cf.
\cite\satakboo,
$\fra k = \roman{ker}(\roman{ad}(z))$ and
$\fra p = [z,\fra g]$, that is, 
the Cartan decomposition 
$\fra g = \fra k \oplus \fra p$
is encapsulated
in the choice of $H$-element $z$. 
As usual, let $\fra p^{\Bobb C} = \fra p^+ \oplus \fra p^-$
be the decomposition of the complexification
$\fra p^{\Bobb C}$
of
$\fra p$
into $(+i)$- and $(-i)$-eigenspaces of $J_z$, respectively.
Thus $\fra p^+$ 
is the span of $w-iJ_zw$ ($=  w-i[z,w]$) and
$\fra p^-$ that of
$w+iJ_zw$ ($=  w+i[z,w]$)
where $w$ runs through $\fra p$.
The association 
$w \mapsto w-iJ_zw$
is an isomorphism of
complex vector spaces from
$\fra p$, endowed with the complex structure $J_z$,
onto  $\fra p^+$.
We will refer 
to an  adjoint orbit
$\Cal O \subseteq \fra g$
having the property that the projection map
from $\fra g$ to $\fra p$,
restricted to $\Cal O$, is a diffeomorphism onto its image,
as a {\it pseudoholomorphic\/} orbit.
In Lemma 3.3.6 below we shall show that a pseudoholomorphic nilpotent orbit
maps (diffeomorphically) onto a smooth {\it complex\/} submanifold of $\fra p$,
and in (3.4) below we shall show that a
non-nilpotent pseudoholomorphic orbit
necessarily maps diffeomorphically onto $\fra p$. Thus
a pseudoholomorphic orbit
$\Cal O$
inherits a complex structure
from the complex structure $J_z$ on
$\fra p$,
and this complex structure,
combined with the Kostant-Kirillov-Souriau
form on $\Cal O$, 
viewed as a coadjoint orbit by means of 
(a positive multiple of) the Killing form,
turns
$\Cal O$ into a (not necessarily positive)
K\"ahler manifold.
Indeed,
a choice of (complex) basis 
$\zeta_1,\dots,\zeta_k$
of 
$\fra p^+$
yields holomorphic coordinate functions
on  $\fra p^*$ which,
together with their complex conjugates
$\overline \zeta_1,\dots,\overline \zeta_k$,
may be used as coordinates
for the smooth complex functions
on 
$\fra p^*$;
since, by assumption,
the composite of the injection of $\Cal O$ into $\fra g^*=
\fra k^* \oplus \fra p^*$
(identified with $\fra g$ via a suitable multiple of the Killing form, see
below),
combined with the projection to $\fra p^*$,
is a diffeomorphism onto its image,
these functions yield coordinate functions on 
$\Cal O$ which,
with an abuse of notation,
we still write as
$\zeta_1,\dots,\zeta_k$
and
$\overline \zeta_1,\dots,\overline \zeta_k$;
we shall see below that
every smooth complex function on $\Cal O$ may be written as
a function in these variables.
Since 
the constituents $\fra p^+$ and
$\fra p^-$
are abelian subalgebras of
$\fra g^{\Bobb C}$,
cf. \cite\helgaboo\ (p.~313),
\cite\satakboo\ (p.~55),
on $\Cal O$
the Poisson-Lie brackets
$\{\zeta_j,\zeta_{j'}\}$
vanish since the Poisson bracket
on $\Cal O$
is induced from the Lie bracket on $\fra g$.
Furthermore,
since $[\fra p,\fra p] = \fra k$,
the Poisson-Lie brackets
$\{\zeta_j,\overline \zeta_{j'}\}$
are given by the restriction to
$\Cal O \subset \fra k \oplus \fra p$
of the resulting complex functions on $\fra k$
determined by the Lie bracket in $\fra k$.
In view of Proposition 2.3,
the complex structure and the symplectic structure
corresponding to the Poisson structure
thus combine to a (not necessarily positive)
K\"ahler structure on $\Cal O$;
this will be made precise in Theorem 3.7.1 below.
For nilpotent orbits, we will make
this somewhat roundabout reasoning more perspicuous
in Theorem 3.3.11
via the classification given in Theorem 3.3.3 below.
\smallskip
We now choose a
positive multiple of the Killing form.
We will say that
a pseudoholomorphic orbit
$\Cal O$
is {\it holomorphic\/} 
provided
the resulting K\"ahler structure 
on $\Cal O$
is positive.
Any positive multiple of the Killing form will do;
passing from one choice to another one
amounts to rescaling the resulting hermitian form.
The notion of
holomorphic orbit relies on the choice of $H$-element $z$
(while the notion of 
pseudoholomorphic orbit does not);
given 
a holomorphic orbit $\Cal O$,
when $z$ is replaced by $-z$,
$-\Cal O$ (and {\it not\/} 
$\Cal O$) will be holomorphic with respect to $-z$.
Sometimes we will refer to an orbit
which is holomorphic with respect to $-z$ as an
{\it antiholomorphic\/} orbit.
The name \lq\lq holomorphic\rq\rq\ 
is intended to hint at
the fact that
the holomorphic discrete series representations arise from holomorphic 
quantization on integral semisimple holomorphic orbits but,
beware, the requisite complex structure
is not the one arising from projection to the symmetric constituent
of the Cartan decomposition.
\smallskip
In particular, we can talk about {\it holomorphic nilpotent orbits\/}.
A different definition of
holomorphic nilpotent orbit
may be found in \cite\vergnsix\ 
which, as we will show in Remark 3.3.13,
is equivalent to ours.
\smallskip
Henceforth the {\it closure\/} is
always understood to be taken in the {\it ordinary\/} (not Zariski) topology.
Let $(\fra g, z)$ be a semisimple Lie algebra of hermitian type.
With reference to the decomposition (3.1.1),
any (pseudo)holomorphic nilpotent orbit $\Cal O$ in $\fra g$
may be written as
$\Cal O = \Cal O_1 \times \ldots  \times\Cal O_k$
where, for $1 \leq j \leq k$,
each $\Cal O_j$ is a 
(pseudo)holomorphic nilpotent orbit in $\fra g_j$,
with reference to the component $z_j$ of $z$
in $\fra g_j$.
Consequently, any orbit
$\Cal O'$ in the closure of
$\Cal O$ may be written as
$\Cal O' = \Cal O'_1 \times \ldots  \times\Cal O'_k$
where, for $1 \leq j \leq k$,
each $\Cal O'_j$ is  
in the closure of
$\Cal O_j$.
This reduces the study of general (pseudo)holomorphic
nilpotent orbits (and of their closures) to that of
(pseudo)holomorphic nilpotent orbits
in simple hermitian Lie algebras.
\smallskip
We now 
spell out the main result
of this section.

\proclaim{Theorem 3.2.1}
Given a
semisimple Lie algebra $(\fra g,z)$ of hermitian type,
for any holomorphic nilpotent orbit $\Cal O$,
the diffeomorphism 
from $\Cal O$ onto its image in $\fra p$
extends to a homeomorphism
from the closure
$\overline{\Cal O}$
onto its image 
in
$\fra p$, 
this homeomorphism
turns 
$\overline{\Cal O}$
into a complex affine variety
in such a way that,
$\fra g$ 
being suitably identified
with its dual so that $\overline{\Cal O}$  
is identified with the closure of 
a coadjoint orbit,
the following hold:
The induced complex analytic structure
on $\overline{\Cal O}$
combines with the resulting Poisson structure to a 
normal complex analytic
stratified K\"ahler structure,
and the closure
$\overline{\Cal O}$ is a union of finitely many
holomorphic nilpotent orbits.
Furthermore, 
when $\fra g$ is simple,
a pseudoholomorphic nilpotent orbit $\Cal O$
having the property that the diffeomorphism 
from $\Cal O$ onto its image in $\fra p$
extends to a homeomorphism
from the closure
$\overline{\Cal O}$
onto its image 
in
$\fra p$ is necessarily
holomorphic or antiholomorphic.
\endproclaim

Theorem 3.2.1 will be a consequence 
of Theorem 3.3.3 below
save that the complex analytic stratified K\"ahler structure
will be explained in Subsection 3.7 below.
In the classical cases,
the real and complex analytic structures
of the closure $\overline{\Cal O}$
of any holomorphic nilpotent orbit 
$\Cal O$ 
will be elucidated in
Theorems 3.5.4, 3.5.5, and 3.6.2 below.
The normality will be explained in Theorems 3.5.5, 3.6.3, 5.3.3, 8.4.1,
8.4.2.

\smallskip\noindent
{\smc 3.2.2. The case of $\fra{sl}(2,\Bobb R)$.}
This  Lie algebra is (well known to be) of hermitian type,
in the following fashion:
The nilpotent elements
$E = \left [ \matrix  0&1\\ 0&0 \endmatrix \right]$
and
$F = \left [ \matrix  0&0\\ 1&0 \endmatrix \right]$
together with the semisimple one
$H=\left [ \matrix  1&0\\ 0&-1 \endmatrix \right]$
span
$\fra{sl}(2,\Bobb R)$,
subject to the relations
$$
[E,F] = H,\quad [H,E] = 2E,\quad [H,F] = -2F,
\tag3.2.2.1
$$
and $Z=\frac 12 (E-F)$
is an $H$-element.
The Cartan decomposition determined by the choice of $Z$
has the form $\fra{sl}(2,\Bobb R) = \fra k \oplus \fra p$
in such a way that
$E-F$
spans $\fra k \cong \Bobb R$ and
$E+F$ and $H$ 
span $\fra p$. 
The induced complex structure $J_Z$ on $\fra p$ is given by
$$
J_Z(E+F) = H,\quad J_Z(H) = - (E+F),
$$
and $E+F$ is a basis of
$\fra p$, viewed as a complex vector space.
Thus $(\fra{sl}(2,\Bobb R),Z)$ 
is of hermitian type.
More generally,
for any $\ell\geq 1$, the product
$(\fra{sl}(2,\Bobb R)^\ell,Z^\ell)=(\fra{sl}(2,\Bobb R)^\ell,Z, \dots, Z)$ 
of $\ell$ copies of $(\fra{sl}(2,\Bobb R),Z)$ 
is of hermitian type.
\smallskip
For intelligibility we reconcile the (well known) classification of orbits
in $\fra{sl}(2,\Bobb R)$ with the notions of holomorphic 
and antiholomorphic orbit. 
The nilcone 
in $\fra{sl}(2,\Bobb R)$
consists of all 
$\alpha\in \fra{sl}(2,\Bobb R)$ with
$\roman{det}(\alpha) = 0$ and 
decomposes into the disjoint union of
the orbits 
$\Cal O^+(0)$ and $\Cal O^-(0)$ (say) of $E$ and
$F$, respectively, 
together with the orbit which consists merely of the origin.
For $\varepsilon >0$,
the two-sheeted hyperboloid
consisting of all $\alpha \in \fra{sl}(2,\Bobb R)$
with $\roman{det}(\alpha) = \varepsilon^2$
decomposes into
the two elliptic 
$\roman{SL}(2,\Bobb R)$-orbits
of $2 \varepsilon Z$ and 
$-2 \varepsilon Z$ which we denote by
$\Cal O^+(\varepsilon)$ 
and
$\Cal O^-(\varepsilon)$, respectively.
\smallskip
We will now show that
the orbits 
$\Cal O^+(\varepsilon)$ 
and $\Cal O^-(\varepsilon)$
are holomorphic 
and antiholomorphic, respectively ($\varepsilon \geq 0$).
To this end, let
$b_0=E-F$, $b_1 =E+F$, $b_2=H$. 
For $\varepsilon >0$,  $\Cal O^+(\varepsilon)$ 
and 
$\Cal O^-(\varepsilon)$ 
are the $G$-orbits of $\varepsilon b_0$ and
$-\varepsilon b_0$, respectively
and, 
when $x_0,x_1,x_2\in \fra{sl}(2,\Bobb R)^*$ 
denote the coordinate functions
in this basis,
$\roman{det}(\alpha) = x_0^2 - x_1^2 - x_2^2$. 
Thus, for $\varepsilon \geq 0$,
$$
\aligned
\Cal O^+(\varepsilon)&=
\{(x_0,x_1,x_2);
x_1^2 + x_2^2 + \varepsilon^2= x_0^2,\ x_0>0\}
\\
\Cal O^-(\varepsilon)&=
\{(x_0,x_1,x_2);
x_1^2 + x_2^2 + \varepsilon^2= x_0^2,\ x_0<0\}.
\endaligned
$$
\smallskip
The Lie algebra
$\fra{sl}(2,\Bobb R)$
may be identified
with its dual via 
{\it any positive\/} multiple of 
the Killing form $k$ 
($k({\alpha},{\beta}) = 4\, \roman{trace}({\alpha}{\beta})$);
the positivity is necessary 
in order for the isomorphism from
$\fra p$ to $\fra p^*$
to preserve the natural orientations
given by the complex structures.
To obtain 
simple formulas,
we will now identify
$\fra{sl}(2,\Bobb R)$
and $\fra p$
with their duals by means of the 
adjoint
$$
\phi \colon \fra{sl}(2,\Bobb R)
@>>>
\fra{sl}(2,\Bobb R)^*,
\quad
{\alpha}\mapsto \phi_{\alpha}\colon 
\fra{sl}(2,\Bobb R) \to \Bobb R,
$$
of the {\it trace\/} pairing
(which is one fourth of the Killing form), that is,
for 
${\alpha},{\beta} \in \fra{sl}(2,\Bobb R)$,
$\phi_{\alpha}({\beta}) = \roman{trace}({\alpha}{\beta})$. 
Then
$\phi_{b_0} = -2x_0$,
$\phi_{b_1} = 2x_1$,
$\phi_{b_2} = 2x_2$,
or, equivalently,
when $b_0,b_1,b_2$ are viewed as linear functions on
$\fra{sl}(2,\Bobb R)^*$,
$$
\frac {b_0}2 \circ \phi = - x_0,
\ 
\frac {b_1}2 \circ \phi = x_1,
\
\frac {b_2}2 \circ \phi = x_2.
$$
The Poisson structure on
$\fra{sl}(2,\Bobb R)$ induced via $\phi$ 
is determined by (3.2.2.1)
and is therefore given by
the formulas
$$
\{x_1,x_2\} = x_0,\ \{x_0,x_2\} = x_1,\ \{x_0,x_1\} = -x_2. 
\tag3.2.2.2
$$
The orthogonal projection from 
$\fra{sl}(2,\Bobb R)$ 
to the plane $\fra p$,
restricted to $\Cal O^+(\varepsilon)$ 
or to $\Cal O^-(\varepsilon)$ 
($\varepsilon \geq 0$),
is plainly a
diffeomorphism onto 
its image in $\fra p$ which is onto $\fra p$
for $\varepsilon > 0$;
for $\varepsilon = 0$, this diffeomorphism extends 
to a homeomorphism from
$\overline{\Cal O^+(0)}$
as well as from
$\overline{\Cal O^-(0)}$
onto $\fra p$.
For $\varepsilon \geq 0$, we now parametrize
the orbits $\Cal O^{\pm}(\varepsilon)$ 
by the plane $\fra p$. 
For $\varepsilon > 0$,
we
take 
$x_1$ and $x_2$
as coordinate functions for the
plane $\fra p$ while,
for $\varepsilon =0$, 
the requisite 
(non-standard) smooth structure
$C^{\infty}$
on $\fra p$
is the algebra of smooth functions
in the variables $x_0,x_1,x_2$,
subject to the relation $x_1^2+x_2^2 = x_0^2$.
The Poisson structures induced on $\fra p$
from the Poisson structures on the orbits 
via the parametrizations
are now given by
$$
\{x_1,x_2\} = x_0
= \cases 
\sqrt{\varepsilon^2 + (x_1^2+x_2^2)}&\quad\text{for}\quad \Cal O^+(\varepsilon)
\\
-\sqrt{\varepsilon^2 +(x_1^2+x_2^2)}&\quad\text{for}\quad \Cal O^-(\varepsilon)
\endcases 
$$
when $\varepsilon >0$ and
by (3.2.2.2),
subject to the relation $x_1^2+x_2^2 = x_0^2$,
when $\varepsilon = 0$.
In the latter case,
the Poisson structure extends the
symplectic Poisson structure
on $\Cal O^+(0)$
(as well as that on
$\Cal O^-(0)$)
to one on
$\overline{\Cal O^+(0)} \cong \fra p$
(or
$\overline{\Cal O^-(0)} \cong \fra p$)
(cf. \cite{\locpois,\,\srni,\,\claustha}).
The induced complex structure $J_Z$ on $\fra p$ is standard, with 
holomorphic coordinate $\zeta = x_1 +ix_2$.
Consequently,
when  
$\fra p$ is used to parametrize
$\Cal O^+(\varepsilon)$,
the resulting hermitian form 
may be written as
$\frac{d\zeta d\overline \zeta}{\sqrt{\varepsilon^2 +\zeta \overline \zeta}}$
and, when  
$\fra p$ is used to parametrize
$\Cal O^-(\varepsilon)$,
this form may accordingly be written as
$-\frac{d\zeta d\overline\zeta}{\sqrt{\varepsilon^2 +\zeta \overline \zeta}}$.
Thus, the orbits 
$\Cal O^+(\varepsilon)$ are holomorphic 
and the orbits 
$\Cal O^-(\varepsilon)$
are antiholomorphic, whatever $\varepsilon \geq 0$,
and the Poisson and complex analytic structures 
turn $\overline{\Cal O^+(0)}$ into a normal positive K\"ahler space
and, likewise, $\Cal O^-(0)$ into a normal 
\lq\lq negative\rq\rq\  K\"ahler space.
Moreover, for $\varepsilon > 0$, the K\"ahler manifold
$\Cal O^+(\varepsilon)$
{\it contracts\/} into the normal K\"ahler space
$\overline{\Cal O^+(0)}$ and, likewise,
$\Cal O^-(\varepsilon)$
contracts into $\overline{\Cal O^-(0)}$.
Hence, $\varepsilon$ being viewed as a deformation parameter,
the {\it singular\/}
normal K\"ahler structure
on $\overline{\Cal O^+(0)}$
{\it deforms into an ordinary K\"ahler structure\/};
under this deformation,
the non-standard Poisson structure
\lq\lq resolves\rq\rq\  into an ordinary smooth symplectic Poisson structure
while the complex analytic structure
is {\it not\/} deformed.
In Section 9
we will 
recall the notion of contraction
and show 
how this relationship
between 
$\overline{\Cal O^+(0)}$
and the 
$\Cal O^+(\varepsilon)$'s ($\varepsilon >0$)
extends to 
arbitrary holomorphic nilpotent orbits.
\smallskip 
For $\varepsilon >0$, 
the complex structure
$J_Z$  (say) on $\Cal O^+(\varepsilon)$
and 
$\Cal O^-(\varepsilon)$
induced from
the complex structure $J_Z$ on $\fra p$
via the projection to $\fra p$
is {\it not\/} the hermitian symmetric space structure.
In order to see this,
we note that,
for $\varepsilon > 0$,
the stereographic projection 
$$
(x_0,x_1,x_2) \longmapsto (y_1,y_2)
=\left(
\frac {x_1}{1+\frac {x_0}\varepsilon},\frac{x_2}{1+\frac{x_0}\varepsilon}
 \right)
$$
identifies each 
$\Cal O^+(\varepsilon)$ 
with the generalized Poincar\'e disc 
$D_{\varepsilon} =\{(y_1,y_2):y_1^2+y_2^2<\varepsilon\}$
in the $(y_1,y_2)$-plane,
generalized in the sense that,
in the literature, 
for $\varepsilon = 1$, this is
the familiar Poincar\'e disc;
see e.~g. (3.7.5) in \cite\beardboo.
A straightforward calculation gives
$\{y_1,y_2\} = \frac \varepsilon 4 \left(1-\frac {r^2}{\varepsilon^2}\right)^2$
where $r^2= y_1^2 + y_2^2$
whence
the induced symplectic structure $\omega$ on
$D_{\varepsilon}$ may be written as
$\omega = \frac 1\varepsilon
\frac  4 {\left(1-\frac {r^2}{\varepsilon^2}\right)^2} 
dy_1 \wedge dy_2$.
This is the K\"ahler form of the hermitian metric
$\frac 1\varepsilon \frac  4 {\left(1-\frac {r^2}{\varepsilon^2}\right)^2} 
d\sigma d\overline {\sigma}$
on 
$D_{\varepsilon}$ with holomorphic coordinate $\sigma = y_1 + i y_2$ to which
the hermitian symmetric space
hermitian form on
$\Cal O^+(\varepsilon)$
passes
under the stereographic projection.
Thus the latter is a symplectomorphism
but it is {\it not\/} compatible with
the complex structures
when $\Cal O^+(\varepsilon)$
is endowed with the complex structure $J_Z$
since then
the coordinate change
$(x_1,x_2) \mapsto (y_1,y_2)$
is not holomorphic.
Thus, 
when
endowed with the 
K\"ahler structure determined by $J_Z$,
$\Cal O^+(1)$ is only {\it symplectically\/}
a model for
the Siegel space
(the space of 
complex structures on the plane compatible with its standard symplectic 
structure).
\smallskip\noindent
{\smc Remark.}
The hermitian symmetric space K\"ahler metric 
on $\Cal O^+(\varepsilon)$ 
(equivalently: on
$D_{\varepsilon}$, with the metric
$\frac 1\varepsilon \frac  4 {\left(1-\frac {r^2}{\varepsilon^2}\right)^2} 
d\sigma d\overline {\sigma}$)
has constant negative curvature
equal to
$-\frac 1{\varepsilon^2}$ while
the hermitian form
$\frac {d\zeta d \overline\zeta}{\sqrt {\varepsilon^2+ \zeta\overline\zeta}}$
has strictly negative curvature which is
constant {\it only\/} along $K$-orbits, $K=\roman U(1)$ 
being the maximal compact subgroup of
$\roman{SL}(2,\Bobb R)$ determined by the choice of $J_Z$.
Indeed, a calculation shows that, for $\varepsilon \geq 0$,
$(\Cal O^+(\varepsilon),J_Z)$ has curvature equal to
$- \frac {\varepsilon^2} 
{\left(\varepsilon^2 + \zeta\overline\zeta\right)^{\frac 52}}$;
notice that $\zeta \ne 0$ when $\varepsilon = 0$.
This shows again that,
for $\varepsilon >0$, 
$J_Z$
does not yield the hermitian symmetric space structure
on $\Cal O^+(\varepsilon)$.
Moreover, this formula entails
that, in particular, 
$\Cal O^+(0)$ 
(as well as $\Cal O^-(0)$)
is flat (has zero curvature).

\smallskip\noindent
{\smc 3.3. Classification of pseudoholomorphic nilpotent orbits.}
Let $\fra g$ be a real Lie algebra.
As usual, we will refer to a triple
$(e,f,h)$ of elements of $\fra g$
as an
$\fra{sl}(2)$-{\it triple\/} in $\fra g$ provided
$[h,e] = 2e,
[h,f] = -2f,
[e,f] = h$.
An $\fra{sl}(2)$-triple $(e,f,h)$ plainly determines and is determined by
a homomorphism
$\kappa \colon \fra{sl}(2,\Bobb R)\to\fra g$ of real Lie algebras
via $\kappa(E) = e,\, \kappa(F) = f,\,\kappa(H) = h$.
Suppose in addition that $\fra g$ is semisimple,
with Cartan involution $\vartheta$ and Cartan decomposition
$\fra g = \fra k \oplus \fra p$.
We will say that an $\fra{sl}(2)$-triple
$(e,f,h)$ 
in $\fra g$
is {\it invariant\/}
provided 
$e-f \in \fra k$ and $e+f, h \in \fra p$.
An $\fra{sl}(2)$-triple
$(e,f,h)$ 
in $\fra g$
is plainly invariant if and only if
the corresponding
Lie algebra homomorphism
$\kappa$ is compatible with the Cartan involutions.
\smallskip
Let $(\fra g,z)$ be a simple
Lie algebra of Hermitian type,
with Cartan decomposition $\fra g = \fra k \oplus \fra p$.
We will say that an 
invariant $\fra{sl}(2)$-triple 
is {\it holomorphic\/}
or an $H_1$-triple
provided the Lie algebra homomorphism
from
$\fra{sl}(2,\Bobb R)$ to $\fra g$ 
it determines is an $H_1$-homomorphism.
An invariant
$\fra{sl}(2)$-triple 
$(e,f,h)$ 
is plainly an $H_1$-triple
if and only if
$$
J_z(e+f) (= [z, e+f]) = h,
\quad J_z(h) (=[z,h]) = - (e+f). 
\tag3.3.1
$$

We will deduce the classification of holomorphic nilpotent orbits
in $\fra g$ from 
the relative root system.
We will systematically denote the split rank,
also called real rank (the dimension
of a maximal abelian subalgebra of $\fra p$) by $r$.
The following result is a version of 
an observation spelled out as
Remark 3 on p.~62 of \cite\satakboo\ 
and we therefore label it as a Proposition.
This observation, in turn, is a consequence of 
the existence of $r$ strongly orthogonal (absolute) roots
relative to the complexification of a compact Cartan subalgebra
(referred to as well as {\it compactly embedded\/} 
Cartan subalgebra in the literature),
cf. Lemma II.4.3 on p.~60 of \cite\satakboo.

\proclaim{Proposition 3.3.2}
Given a simple
Lie algebra $(\fra g,z)$  
of Hermitian type and split rank $r$,
with Cartan decomposition
$\fra g = \fra k \oplus \fra p$,
a choice of maximal 
abelian subalgebra $\fra a_{\fra p}$ of $\fra p$
determines
a system
$(e_1,f_1,h_1)$,\dots,$(e_r,f_r,h_r)$
of $H_1$-triples which
combine to
an $H_1$-embedding
of $(\fra{sl}(2,\Bobb R)^r,Z^r)$
into  $(\fra g,z)$
in such a way that the following hold.
\newline\noindent
{\rm 1)} The elements
$h_1,\dots, h_r$ constitute a basis of 
$\fra a_{\fra p}$. 
\newline\noindent
{\rm 2)} 
With the notation $(\xi_1,\dots,\xi_r)$ for the basis of
$\fra a^*_{\fra p}$
dual to 
$h_1,\dots, h_r$,
the relative root system
$S = S(\fra g, \fra a_{\fra p}) \subseteq \fra a_{\fra p}^*$
determined by 
$\fra a_{\fra p}$ is given by
$$
S= 
\cases \{ \pm \xi_j \pm \xi_k\,(j \ne k),\ \pm 2 \xi_j \} \quad
&\text{type}\ (\roman{C}_r),\quad\text{or}\\
\{ \pm \xi_j \pm \xi_k\,(j \ne k),\ \pm 2 \xi_j, \pm  \xi_j \} \quad
&\text{type}\ (\roman{BC}_r) .
\endcases
$$
{\rm 3)} 
The $H$-element $z$ of $\fra g$ equals
$\frac 12 (e_1-f_1 + \dots + e_r - f_r)$
if and only if $S$ is of type $(\roman C_r)$.
\newline\noindent
{\rm 4)} 
The $H_1$-embedding
of $(\fra{sl}(2,\Bobb R)^r,Z^r)$
into  $(\fra g,z)$
extends to
an $H_1$-embedding
of $\fra{sp}(r,\Bobb R)$
into $\fra g$, 
$\fra{sp}(r,\Bobb R)$ being endowed with the appropriate
$H$-element, and this embedding
identifies the two $H$-elements
if and only if 
$S$ is of type $(\roman C_r)$.
\endproclaim

We will say that
a simple Lie algebra
of hermitian type 
is {\it regular\/}
provided its relative root system is
of type $(\roman C_r)$,
and we will say that it is
{\it non-regular\/} otherwise.
This notion of regularity 
fits with that
of regularity of pre-homogeneous spaces
to be examined in Section 7 below. 
For the corresponding symmetric domain,
regularity amounts to the property of being a tube domain.
 
\demo{Proof}
To adjust the notation on p.~110 of \cite\satakboo\ 
to ours,
recall that in
\cite\satakboo\  (p.~92 and p.~109)
certain root vectors
$X_k,Y_k,H_k \in \fra g_{\Bobb C}$
($1 \leq k \leq r$)
are given
satisfying the equations
$$
[H_k,X_k] = -2 Y_k,
\quad
[H_k,Y_k] = 2 X_k,
\quad
[X_k,Y_k] = 2 H_k.
$$
Then the triples $(e_k,f_k,h_k)$ given by
$$
e_k = \frac 12 (Y_k+H_k),
\quad
f_k = \frac 12 (Y_k-H_k),
\quad
h_k = X_k
$$
(cf. III.1.6. on p.~91 and III.2.9 on p.~97 of \cite\satakboo)
will have the  properties
asserted in 1) and 2).
\smallskip
The statement 3) is a consequence of Corollary
III.1.6 on p.~94, combined
with Remark 1 in III.8 on p.~150 of \cite\satakboo. 
To justify
statement 4) we observe that,
in the regular case, 
appropriate root vectors for
the relative root system
$S$
span a copy of 
$\fra{sp}(r,\Bobb R)$ in $\fra g$.
In the non-regular case, the relative root system
${S=\{ \pm \xi_j \pm \xi_k\,(j \ne k),\ \pm 2 \xi_j, \pm  \xi_j \}}$
contains the subsystem
$S'= \{ \pm \xi_j \pm \xi_k\,(j \ne k),\ \pm 2 \xi_j \}$
and,
likewise,
appropriate root vectors for the latter span a copy of 
$\fra{sp}(r,\Bobb R)$
in $\fra g$. \qed
\enddemo

\smallskip\noindent
{\smc Remark.} For intelligibility,
we list the various types of the relative root systems
that occur here,
cf. \cite\satakboo\ p.~115: 
\newline\noindent
$\fra{su}(q,q)$: $(\roman C_q)$; $\fra{su}(p,q)$, $p>q$: $(\roman {BC}_q)$;
$\fra{so}(2,q)$: $(\roman C_2)$ $(q \geq 3)$;
$\fra{sp}(n,\Bobb R)$:  $(\roman C_n)$;
\newline\noindent
$\fra{so}^*(2n)$: $\cases 
&(\roman C_{\frac n2}),\quad n \ \text{even} 
\\ 
&(\roman {BC}_{[\frac n2]}),\quad n>2 \ \text{odd} 
\endcases$;
$\fra{e}_{6(-14)}$: $(\roman {BC}_2)$;
$\fra{e}_{7(-25)}$: $(\roman C_3)$.
\smallskip
Given a simple
Lie algebra $(\fra g,z)$,  
with adjoint group $G$,
for $0 \leq t,u \leq r,\, t+u \leq r$, let
$e_{t,u} = e_1 +\dots+e_t -e_{t+1} - \dots -e_{t+u} \in \fra g$
and let
$\Cal O_{t,u} = Ge_{t,u}\subseteq \fra g$;
this is plainly a nilpotent
orbit 
in $\fra g$.
We will refer to $(t,u)$ as the {\it type\/}
of
$\Cal O_{t,u}$,
to $t+u$ as its {\it rank\/},
and to
$t-u$ as its {\it signature\/}.
For later reference, we will write
$f_{t,u} = f_1 +\dots+f_t -f_{t+1} - \dots -f_{t+u} \in \fra g$
($0 \leq t,u \leq r,\, t+u \leq r$).

\proclaim{Theorem 3.3.3}
For a simple  Lie algebra
$(\fra g,z)$ 
of hermitian type
and split rank $r$,
the $\frac {(r+1)(r+2)}2$ nilpotent 
orbits
$\Cal O_{t,u}$
$(0 \leq t,u \leq r,\, t+u \leq r)$
constitute a complete list
of the pseudoholomorphic
nilpotent orbits.
Given such an orbit 
$\Cal O_{t,u}$,
an orbit $\Cal O'$
distinct from
$\Cal O_{t,u}$
is in the closure of
$\Cal O_{t,u}$
if and only if
$\Cal O' =\Cal O_{t',u'}$
for some $t' \leq t$ and $u' \leq u$
with $t'+u'<t+u$.
In particular, the orbits $\Cal O_{t,0}$ 
$(t \leq r)$
and $\Cal O_{0,u}$ 
$(u \leq r)$ are precisely
the holomorphic and antiholomorphic ones, respectively.
\endproclaim

\proclaim{Addendum} 
For a nilpotent orbit $\Cal O$ of a simple  Lie algebra
$(\fra g,z)$ of hermitian type,
the property of being holomorphic, antiholomorphic, or pseudoholomorphic
does not depend on the choice of $H$-element $z$.
More precisely: Once an $H$-element $z$ has been chosen,
when $\Cal O$ is holomorphic or antiholomorphic with respect to
$z$, it will be so with respect to any $\widetilde z$ in the space
$G/K$ of $H$-elements;
when $\Cal O$ is pseudoholomorphic with respect to
$z$, it will be so with respect to any $\widetilde z$ in the space
$(G/K^+)\cup (G/K^-)$.
\endproclaim

\smallskip
The proof requires some preparation:
Following the procedure in III.4 of \cite\satakboo\ 
(p.~110), we choose the basis 
$\Delta = \{\gamma_1,\dots,\gamma_r\}$
of $S$
where $\gamma_j = \xi_j-\xi_{j+1}$ for
$1 \leq j \leq r-1$
and where 
$\gamma_r = \xi_r$ 
in the regular case
and
$\gamma_r = 2\xi_r$
in the non-regular case.
Let
$S= S^+ \cup S^-$
be the resulting
decomposition
into positive and negative relative roots,  so that
$$
S^+= 
\cases \{ \xi_j \pm \xi_k\,(j < k),\ 2 \xi_j \} \quad
&\text{regular case}\\
       \{ \xi_j \pm \xi_k\,(j < k),\ 2 \xi_j, \xi_j \} \quad
&\text{non-regular case} .
\endcases
$$
For 
$\lambda \in S$, 
let $\fra g^{\lambda}$ be
its root space and, for
any $S' \subseteq S$, write
$\fra g^{S'}=\sum_{\lambda \in S'} \fra g^{\lambda}$.
Thus $\fra g$ decomposes as
$\fra g = \fra l \oplus \fra g^S$
where $\fra l = \fra m + \fra a_{\fra p}$
is the centralizer of $\fra a_{\fra p}$
in $\fra g$,
$\fra m$ being the centralizer of
$\fra a_{\fra p}$
in $\fra k$. With the notation
$\fra n^+ = \fra g^{S^+}$, 
the resulting Iwasawa decomposition
has the form
$\fra g = \fra k + \fra a_{\fra p} + \fra n^+$; 
accordingly, we write
$G=KA_{\fra p}N$ for 
the corresponding global Iwasawa decomposition---in particular, this
fixes the nilpotent group $N$---and,
for later reference,
we write $\fra n^- = \fra g^{S^-}$.
\smallskip
For $k = 1,\dots,r$, let
$$
e^k = e_1 + \dots + e_k,\ 
f^k = f_1 + \dots + f_k,\ 
h^k = h_1 + \dots + h_k;
$$ 
this yields the
$H_1$-triple $(e^k,f^k,h^k)$ and
we write
$\kappa_k \colon \fra{sl}(2,\Bobb R) \to \fra g$
for the resulting $H_1$-homomorphism.
Here and below the ordering of the $H_1$-triples 
$(e_k,f_k,h_k)$ ($1 \leq k \leq r$)
does not matter.
For $k = 1,\dots,r$, 
the parabolic subalgebra
$\fra q_k$
of $\fra g$
(written as $\fra b_{\kappa}$
in \cite\satakboo\ III.2 (p.~95)) which is
determined by
$(e^k,f^k,h^k)$ or, 
what amounts to the same, by
$\kappa_k$,
has Levi decomposition 
$\fra q_k = \fra l_k \oplus \fra n_k^+$
where 
$\fra l_k$ is the (reductive) centralizer of
$h^k$ in $\fra g$,
and 
$$\fra g = \fra n_k^-\oplus \fra l_k \oplus \fra n_k^+.
\tag3.3.4(k)
$$
Let
${\fra n^+(0) = \{x \in \fra n^+; [h^r,x] = 0\}}$,
${\fra n^-(0) = \{x \in \fra n^-; [h^r,x] = 0\}}$,
for $j=1,2$, let
$$
{\fra n(j) = \{x \in \fra n^+; [h^r,x] = jx\}},
\quad 
{\fra n(-j) = \{x \in \fra n^-; [h^r,x] = -jx\}}
$$ 
and,
more generally,
given $1 \leq k \leq r$,
let
$$
\fra n_k(j) = \{x \in \fra n^+; [h^k,x] = jx\},
\
\fra n_k(-j) = \{x \in \fra n^-; [h^k,x] = -jx\}.
$$
We note that
$\fra n_k(1)$ and $\fra n_k(2)$
correspond to $V_{\kappa}$ 
and $U_{\kappa}$ 
in III.2  of \cite\satakboo\ (p.~95 ff.), respectively.
The following is well known (and easy to prove).

\proclaim{Lemma 3.3.5}
As a vector space, the nilpotent Lie algebra
$\fra n^+$ decomposes as
$$
\fra n^+ = \fra n^+(0) \oplus \fra n(1) \oplus \fra n(2),
$$
and the decomposition yields a graded Lie algebra
(\lq\lq graded\rq\rq\  being understood in the 
obvious naive sense)
whence, in particular,
$\fra n(2)$ is abelian.
Furthermore,
%$\fra n(1) \oplus \fra n(2) =\fra n^+_r$ and, 
in the regular case,
$\fra n(1)$ is zero
while, in the non-regular case,
$\fra n^+_r$ which, as a vector space,
is the direct sum $\fra n(1) \oplus \fra n(2)$,
may be written as a 
central Lie algebra extension
$0\to \fra n(2)\to \fra n^+_r \to \fra n(1)\to 0$ 
with abelian
quotient,
i.~e. as a Heisenberg algebra.
More generally,
for $1 \leq k \leq r$, as a vector space, $\fra n^+_k$ decomposes as
$\fra n^+_k=\fra n_k(1) \oplus \fra n_k(2)$. \qed
\endproclaim

To spell out the next lemma we recall that the linear $K$-action
on 
$\fra p$, viewed as a complex vector space via the complex structure
$J_z$, extends canonically to a linear $K^{\Bobb C}$-action
on this vector space.

\proclaim{Lemma 3.3.6}
Each nilpotent orbit
$\Cal O_{t,u} (= Ge_{t,u})$ (where $0 \leq t,u \leq r,\, t+u \leq r$)
is pseudoholomorphic, that is,
the projection from
$\Cal O_{t,u}$
to $\fra p$
is a diffeomorphism
onto its image in $\fra p$, and the image of
$\Cal O_{t,u}$ in 
$\fra p$ is the $K^{\Bobb C}$-orbit
of $\frac 12(e_{t,u} +f_{t,u})$.
\endproclaim

\demo{Proof}
Write
$\fra n^+_{\fra p} \subseteq \fra n(2)$ and
$\fra n^-_{\fra p}  \subseteq \fra n(-2)$
for the span of $e_1, e_2,\dots,e_r$ and $f_1, f_2,\dots,f_r$,
respectively.
We claim that
$G\fra n^+_{\fra p}=K\fra n^+_{\fra p}$.
Indeed,
$G\fra n^+_{\fra p}$
is the union of the orbits $\Cal O_{t,u}$
($0 \leq t,u \leq r,\ t+u \leq r$).
By Lemma 3.3.5,
$[\fra n^+,\fra n(2)] \subseteq \fra n(2)$;
hence 
$N\fra n^+_{\fra p}= \fra n^+_{\fra p}$
since $N$ acts trivially on $\fra n^+(2)$.
Furthermore,
$\fra n^+_{\fra p}$
is closed under the $A_{\fra p}$-action.
Thus,
by virtue of
the Iwasawa decomposition $G=KA_{\fra p}N$
of $G$,
$G\fra n^+_{\fra p}=K\fra n^+_{\fra p}$ as asserted. 
Consequently
the orbit $\Cal O_{t,u}=G(e_1 + \dots + e_t -e_{t+1} - \dots -e_{t+u})$ 
decomposes
into the $K$-orbits of the elements
$$
e_{t,u,\lambda} = \lambda_1 e_1 + \dots + \lambda_te_t -\lambda_{t+1}e_{t+1} - 
\dots -\lambda_{t+u}e_{t+u},
$$
where $\lambda=(\lambda_1,\dots,\lambda_{t+u})\in \Bobb R^{t+u}$ 
with $\lambda_j > 0$ ($1 \leq j \leq t+u$).
In view of the action of the relative Weyl group,
which contains a copy of the symmetric group on $r$ letters,
we normalize 
$\lambda$
by requiring
$\lambda_1\leq \lambda_2 \leq \dots \leq \lambda_t$
and $\lambda_{t+1} \leq \dots \leq \lambda_{t+u}$.
Likewise
the orbit $G(f_1 + \dots + f_t -f_{t+1} - \dots -f_{t+u})$ 
decomposes
into the $K$-orbits of the elements
$$
f_{t,u,\mu} = \mu_1 f_1 + \dots + \mu_tf_t -\mu_{t+1}f_{t+1} - 
\dots -\mu_{t+u}f_{t+u}
$$
where
$\mu =(\mu_1,\dots,\mu_t,\mu_{t+1},\dots,\mu_{t+u})\in \Bobb R^{t+u}$
with $\mu_j > 0$ ($1 \leq j \leq t+u$).
\smallskip
The image in $\fra p$
of $e_{t,u,\lambda}$
under the orthogonal projection
to $\fra p$
is the element
$\frac 12 (e_{t,u,\lambda}+f_{t,u,\lambda})\in \fra p $.
Using the (known) structure
of the stabilizer of any nilpotent element
in a semisimple Lie algebra,
cf. e.~g. Theorem 1.7 in \cite\hilneors,
we conclude that the stabilizer 
$Z_K(\frac 12 (e_{t,u,\lambda}+f_{t,u,\lambda}))$
in $K$ of
$\frac 12 (e_{t,u,\lambda}+f_{t,u,\lambda})\in \fra p $
coincides with the
stabilizer 
$Z_K(e_{t,u,\lambda})$ in $K$ of
$e_{t,u,\lambda}$.
More precisely:
The nilpotent element $e = e_{t,u,\lambda}$
belongs to the invariant 
$\fra{sl}(2)$-triple 
$(e,f,h)$ 
where
$$ 
f= \frac 1 {\lambda_1} f_1 + \dots + \frac 1 {\lambda_t}f_t 
-\frac 1 {\lambda_{t+1}}f_{t+1} - 
\dots -\frac 1 {\lambda_{t+u}}f_{t+u}
$$
and $h = h_1 + \dots+ h_t+h_{t+1} +\dots + h_{t+u}$;
we note that 
$(e,f,h)$ is an $H_1$-triple
if and only if $u=0$
and if $\lambda_1 = \dots = \lambda_{t+u} = 1$.
The stabilizer $Z_G(e)$ of $e$ in the adjoint group $G$
of $\fra g$ is contained in the Jacobson-Morozow parabolic subgroup
$Q=UL$ (say) associated with $(e,f,h)$
(where $U$ is the unipotent radical and $L$ the Levi subgroup),
and the reductive constituent 
of $Z_G(e)$ is contained in the (reductive) stabilizer
$Z_G(\fra s)$ of the entire triple
(where $\fra s$ refers to the copy of $\fra{sl}(2,\Bobb R)$ in $\fra g$
generated by $e,f,h$) whence
$Z_K(e)=Z_G(e) \cap K$ is contained in $Z_G(\fra s)$,
in fact, equals $Z_G(\fra s)\cap K$;
the notation in \cite\hilneors\ 
is $(X,Y,H)$ for our $(e,f,h)$.
Noting that,
with reference to the Cartan involution $\vartheta$ on $\fra g$,
$f_{t,u,\lambda} = -\vartheta (e_{t,u,\lambda}) =-\vartheta (e)$,
we must show that
the inclusion
$Z_K(e) \subseteq Z_K(\frac 12 (e-\vartheta e))$
is the identity.
In order to justify this assertion, we note first that,
when $u=0$ and $\lambda_1 = \dots = \lambda_t = 1$,
we have $e = e_1 + \dots +e_t= e^t$,
$f=f_1 + \dots +f_t=f^t$,
and
$J_z(e+f) = h = h_1+\dots+ h_t=h^t$;
this implies that
$Z_K(e+f) =Z_K(h)= K \cap L_t$
(since $J_z$ commutes with the adjoint action of $K$ on $\fra g$)
where $L_t$ refers to the Levi subgroup
of the parabolic subgroup $Q_t = U_t L_t$ 
associated with the $H_1$-triple $(e^t,f^t,h^t)$.
Thus 
$Z_K(e+f)$ now coincides with the stabilizer 
in $K$ of $h$. The latter, in turn, equals
the stabilizer 
$Z_K(e)$ of $e$ in $K$;
since the factor $\frac 12$ is irrelevant,
we conclude that
$Z_K(\frac 12(e+f))$ is equal to $Z_K(e)$
in this special case.
Next we consider the case where
$e = e_1 + \dots +e_t-e_{t+1} - \dots - e_{t+u}$
and
$f=f_1 + \dots +f_t -f_{t+1} - \dots - f_{t+u}$.
Noting that
$e$ and $f$
are conjugate to
$e_1 + \dots +e_t +e_{t+1} + \dots + e_{t+u}$
and
$f_1 + \dots +f_t +f_{t+1} + \dots + f_{t+u}$,
respectively,
in the complexification $\fra g_{\Bobb C}$,
taking suitable conjugates of stabilizers
in the complexification of $G$, we conclude
that
$Z_K(\frac 12(e+f))$ coincides with $Z_K(e)$.
Finally, we may reduce the general case to the one just established
by noting that
$e_{t,u,\lambda}$ and
$f_{t,u,\lambda}$ 
arise from
$e_1 + \dots +e_t -e_{t+1} - \dots - e_{t+u}$
and
$f_1 + \dots +f_t -f_{t+1} - \dots - f_{t+u}$
respectively
as the result of the adjoint action with a suitable element
from the constituent $A_{\fra p}$
in the Iwasawa-decomposition $G=KA_{\fra p} N$.
Consequently
$Z_K(\frac 12(e+f))$ coincides with $Z_K(e)$
in the general case.
\smallskip
We conclude that
the projection from
$Ke_{t,u,\lambda}$
onto its image in
$\fra p$
is a diffeomorphism.
Furthermore,
with the above normalization,
the orbits
$K(\frac 12 (e_{t,u,\lambda}+f_{t,u,\lambda})) \subseteq \fra p$
($0<\lambda_1\leq \lambda_2 \leq \dots \leq \lambda_t$
and $0<\lambda_{t+1} \leq \dots \leq \lambda_{t+u}$)
are mutually disjoint,
and their union
coincides with the $K^{\Bobb C}$-orbit
of $\frac 12(e_{t,u} +f_{t,u})$ in $\fra p$
as may be see from writing out the $K^{\Bobb C}$-action on
$\fra p$ in terms of
the polar decomposition $K^{\Bobb C} = K\roman{exp}(i\fra k)$
of $K^{\Bobb C}$.
In particular,
given
$0<\lambda_1\leq \lambda_2 \leq \dots \leq \lambda_t$
and $0<\lambda_{t+1} \leq \dots \leq \lambda_{t+u}$,
for some $X \in \fra k$, we have
$$
e_{t,u,\lambda} +f_{t,u,\lambda}=
\roman{exp}(i X) (e_{t,u} +f_{t,u}) .
$$
Hence, 
for $0 \leq t,u \leq r$ with $t+u \leq r$,
the projection from
$\Cal O_{t,u}=Ge_{t,u}$
onto its image in $\fra p$
is a diffeomorphism
onto the $K^{\Bobb C}$-orbit
of $\frac 12(e_{t,u} +f_{t,u})$
in $\fra p$.
In particular,
$\Cal O_{t,u}$ is pseudoholomorphic as asserted. \qed
\enddemo

\noindent
{\smc Remark.} 
Let $L_r$ be the (reductive) centralizer of
$h^r$ in the adjoint group $G$ and let
$K_r= K \cap L_r$.
The union of the orbits $\Cal O_{t,u}$
($0 \leq t,u \leq r,\ t+u \leq r$)
equals $K \fra n^+_r \subseteq \fra g$, and the
projection to $\fra p$
is induced by the canonical surjection from
$K \times_{K_r} \fra n^+_r$ to $\fra p$
given by the assignment to
$(x,\alpha) \in K \times \fra n^+_r$ of
$x\overline\alpha x^{-1} \in \fra p$
where 
$\overline\alpha \in \fra p$
refers to the image of
$\alpha \in \fra n^+_r$ 
under the projection to $\fra p$.
From this observation,
one may as well deduce that
the projection from
$\Cal O_{t,u}=Ge_{t,u}$
onto its image in $\fra p$
is a diffeomorphism
onto the $K^{\Bobb C}$-orbit
of $\frac 12(e_{t,u} +f_{t,u})$
in $\fra p$.
For example, 
for $\fra g = \fra{sp}(\ell,\Bobb R)$,
the real rank $r$ equals $\ell$,
$\fra n^+_r$ amounts to the space $\roman S^2[\Bobb R^{\ell}]$
of real symmetric $\ell \times \ell$-matrices,
$\fra p$ to that 
of complex symmetric $\ell \times \ell$-matrices
$\roman S^2[\Bobb C^{\ell}]$,
associating 
$\overline\alpha \in \fra p$
to $\alpha \in \fra n^+_r$ 
amounts to associating to a
real symmetric $\ell \times \ell$-matrix
the same matrix, 
viewed as a complex symmetric $\ell \times \ell$-matrix,
and the canonical surjection from
$K \times_{K_r} \fra n^+_r$ to $\fra p$
amounts to the familiar map from
$\roman U(\ell) \times_{\roman O(\ell,\Bobb R)} \roman S^2[\Bobb R^{\ell}]$ 
to $\roman S^2[\Bobb C^{\ell}]$.
The restriction of this map to 
the $\roman U(\ell)$-span of any
$\roman O(\ell,\Bobb R)$-orbit of non-degenerate real symmetric 
$\ell \times \ell$-matrices
is a diffeomorphism onto
the space
of non-degenerate complex symmetric
$\ell \times \ell$-matrices.

\demo{Proof of Theorem 3.3.3}
Lemma 3.3.6 says that
the orbits
$\Cal O_{t,u}$ are pseudoholomorphic.
It remains to justify the claim as to the closures 
of the orbits $\Cal  O_{t,u}$
and to explain
how the holomorphic and
antiholomorphic orbits are singled out.
To this end we observe that
every nilpotent orbit of
the semisimple Lie algebra of hermitian type
$((\fra{sl}(2,\Bobb R)^r,(z_1,\dots,z_r))$
is pseudoholomorphic;
these orbits are generated by
the $(\varepsilon_1 e_1,\dots,\varepsilon_r e_r)$'s where
$\varepsilon_j = -1,0,+1$ ($1 \leq j \leq r$).
The classification of holomorphic orbits
of $\fra {sl}(2,\Bobb R)$ reproduced in (3.2.2) above entails at once
that
the orbits 
involving only
$\varepsilon_j = 0, +1$
are the holomorphic ones
and that, likewise, the orbits 
involving only
$\varepsilon_j = -1,0$
are the antiholomorphic ones.
Furthermore,
when $\Cal O$ 
is the
$(\roman{SL}(2,\Bobb R)^r)$-orbit of
$(\varepsilon_1 e_1,\dots,\varepsilon_r e_r)$,
the orbits in the closure
of $\Cal O$
are precisely the orbits
generated by
those $r$-tuples which arise from
$(\varepsilon_1 e_1,\dots,\varepsilon_r e_r)$
by replacing some $\varepsilon_j \ne 0$
with zero.
However,
within $\fra g$,
the relative Weyl group $W_{\fra p}\cong (\Bobb Z\big /2)^r \rrtimes S_r$
permutes the $r$ copies of
$\fra{sl}(2,\Bobb R)$
and hence in particular
the entries of the $(\varepsilon_1 e_1,\dots,\varepsilon_r e_r)$
whence,
within $\fra g$,
what remains invariant under 
the relative Weyl group,  
are the rank $t+u$ and the signature $t-u$.
Thus an orbit $\Cal O'$
distinct from
$\Cal O_{t,u}$
is in the closure of
$\Cal O_{t,u}$
if and only if
$\Cal O' =\Cal O_{t',u'}$
for some $t' \leq t$ and $u' \leq u$
with $t'+u'<t+u$.
Furthermore, since the complex structure $J_z$ on $\fra p$
is $K$-invariant,
the orbits $\Cal O_{t,0}$ ($t \leq r$) 
and $\Cal O_{0,u}$ ($u \leq r$)  are precisely
the holomorphic and antiholomorphic ones, respectively.
\smallskip
Finally we establish the fact that 
the orbits $\Cal O_{t,u}$
$(0 \leq t,u \leq r,\, t+u \leq r)$
exhaust the pseudoholomorphic nilpotent orbits.
Inspection of the possible nilpotent orbits in $\fra g$
shows that, given a nilpotent orbit
which is different from one of the kind
$\Cal O_{t,u}$,
the projection to 
the constituent $\fra p$ of the Cartan decomposition of $\fra g$
is no longer injective.
For the classical cases
we  postpone the details to the proofs of Theorems 3.5.3 and 3.6.2
below. 
We now settle
the two exceptional cases 
$\fra{e}_{6(-14)}$ and $\fra{e}_{7(-25)}$
by inspection.
As before, $G$ denotes the adjoint group of the corresponding 
simple Lie algebra $\fra g$ under discussion.
\smallskip
The Lie algebra $\fra g =\fra e_{6(-14)}$ 
has 12 non-zero nilpotent orbits.
The orbits numbered (1) to (5) in Table X of \cite\djokotwo\ (p.~512) are
precisely the
pseudoholomorphic ones;
with the appropriate choice of $H$-element,
the orbits (1) and (3)
coincide 
with the orbits $Ge^1$ and $Ge^2$,
respectively,
the orbits (2) and (4)
are antiholomorphic
and amount to the orbits
$G(-e^1)$ and $G(-e^2)$,
while the orbit (5)
is our pseudoholomorphic nilpotent orbit 
$G(e^1-e^2)$
which  is neither holomorphic nor antiholomorphic.
There are no other pseudoholomorphic nilpotent orbits
in $\fra e_{6(-14)}$.
Indeed, in this case, $\dim_{\Bobb C} \fra p = 16$,
and column 4) of Table X of \cite\djokotwo\
gives the complex dimensions of the stabilizers of the
corresponding $K^{\Bobb C}$-orbits
in $\fra p$ where $K\cong \roman{Spin}_{\roman c}(10)$ 
is the maximal compact subgroup
of the adjoint group $G$.
Since $\dim K^{\Bobb C} = 46$, 
inspection of the dimensions of the stabilizers
shows that, under the projection from $\fra g$ to $\fra p$,
none of the orbits (6)--(12)
can map diffeomorphically onto its image in $\fra p$.
\smallskip
The Lie algebra $\fra g =\fra e_{7(-25)}$
has 22 non-zero nilpotent orbits.
The orbits numbered (1) to (9) in Table XIII of \cite\djokotwo\ (p.~516) are
precisely the
pseudoholomorphic ones;
with the appropriate choice of $H$-element,
the orbits (1),  (3), and (6)
coincide 
with,
respectively, the orbits $Ge^1$, $Ge^2$, and $Ge^3$,
the orbits (2), (4), and (7)
are antiholomorphic
and amount to, respectively, the orbits
$G(-e^1)$, $G(-e^2)$, and $G(-e^3)$,
while the remaining orbits (5), (8), (9)
amount to our pseudoholomorphic nilpotent orbits 
respectively $G(e^1-e^2)$,
$G(e^1+e^2 - e^3)$,
$G(e^1-e^2 - e^3)$,
which  are neither holomorphic nor antiholomorphic.
There are no other pseudoholomorphic nilpotent orbits
in $\fra e_{7(-25)}$.
Indeed, in this case, $\dim_{\Bobb C} \fra p = 27$,
and
the complex dimensions of the stabilizers of the
corresponding $K^{\Bobb C}$-orbits
in $\fra p$
may be found in column 4) of Table X of \cite\djokotwo\
where $K$
is the maximal compact subgroup
of the adjoint group $G$;
the group 
$K$ is locally isomorphic to $E_{6(-78)} \times \roman{SO}(2,\Bobb R)$
($E_{6(-78)}$ being a compact form of $E_6$) 
and $\dim_{\Bobb R}K = 79$.
Inspection of the dimensions of the stabilizers
shows that, under the projection from $\fra g$ to $\fra p$,
none of the orbits (10)--(22)
can map diffeomorphically onto its image in $\fra p$. \qed 
\enddemo

\demo{Proof of the Addendum}
Let $z$ and $\widetilde z$ be two $H$-elements, and use the notation
$\widetilde e_j$, $\widetilde e_{t,u}$, etc.
with reference to $\widetilde z$.
Then, for 
$0 \leq t,u \leq r,\, t+u \leq r$, 
$\Cal O_{t,u} = Ge_{t,u} = G\widetilde e_{t,u}=\widetilde{\Cal O_{t,u}}$
provided 
$z$ and $\widetilde z$
are in the same component of the space of $H$-elements, and
$\Cal O_{t,u} = Ge_{t,u} = G\widetilde e_{u,t}=\widetilde{\Cal O_{u,t}}$
otherwise. \qed
\enddemo

\noindent
{\smc Remark 3.3.7.} In the standard cases, 
the holomorphicity of the
orbits $\Cal O_{t,0}$ ($1 \leq t \leq r$),
in particular, the assertion in the above proof
of Lemma 3.3.6 involving the stabilizers,
will also be a consequence of Theorem 5.3.3 below.

\proclaim{Corollary 3.3.8} 
The projection
from 
the closure 
$\overline {\Cal O}$
of a pseudoholomorphic nilpotent orbit
$\Cal O$
to the constituent $\fra p$ of the Cartan decomposition is injective
if and only if 
$\Cal O$
is holomorphic or
antiholomorphic.
\endproclaim

\demo{Proof}
Indeed, 
a pseudoholomorphic nilpotent orbit
$\Cal O$
is one of the kind $\Cal O_{t,u}$. 
If $u=0$ or $t=0$,
for each 
$0 \leq s< \roman{rank}(\Cal O_{t,u})$,
there is only a single 
rank $s$ nilpotent
orbit 
in the closure
$\overline {\Cal O_{t,u}}$.
If $u\ne 0$ and $t\ne 0$,
the projection
from 
the closure $\overline {\Cal O_{t,u}}$
to $\fra p$
is no longer injective since
the closure contains at least two orbits
of the same rank but having different signature;
the projection 
to $\fra p$
identifies the two. \qed
\enddemo

\proclaim{Corollary 3.3.9}
An 
orbit which belongs to the closure
of 
a (pseudo)holomorphic 
nilpotent orbit
is necessarily 
(pseudo)holomorphic,
and the closure of a 
(pseudo-)\linebreak
holomorphic 
nilpotent orbit
is a union of finitely many
(pseudo)holomorphic nilpotent orbits.
\endproclaim

\demo{Proof} This is an almost immediate consequence of 
Theorem 3.3.3. \qed\enddemo

\demo{Proof of Theorem 3.2.1}
This follows readily from Theorem 3.3.3. We leave the details
to the reader. \qed \enddemo

\smallskip
Let
$(\fra g,z)$ be a simple  Lie algebra
of hermitian type and real rank $r$;
for $ 1 \leq t \leq r$,
we will henceforth denote the holomorphic
nilpotent orbit $\Cal O_{t,0}$ by $\Cal O_t$.
Theorem 3.3.3 entails that
the holomorphic nilpotent orbits
$\Cal O_0,\dots, \Cal O_r$ are linearly ordered
in such a way that
$$
\{0\}=\Cal O_0 \subseteq \overline{\Cal O_1} 
\subseteq \ldots \subseteq \overline{\Cal O_r}.
\tag3.3.10
$$
Borrowing and extending terminology from \cite\kostasix,
where the regular (complex) nilpotent orbit is said to be
principal,
we will refer to the top orbit
$\Cal O_r$ as the 
{\it principal holomorphic\/} nilpotent orbit.
The principal holomorphic nilpotent (in our sense)
may be viewed as the {\it regular holomorphic nilpotent orbit\/};
it is in particular unique, and its closure contains every
holomorphic nilpotent orbit.

\proclaim{Theorem 3.3.11}
For
the principal holomorphic nilpotent orbit
$\Cal O_r$,
the composite of the projection from
the closure $\overline {\Cal O_r}$ to $\fra p$
with the complex linear isomorphism
from
$\fra p$
to $\fra p^+$ is a
homeomorphism $\overline {\Cal O_r} \to \fra p^+$.
Under this homeomorphism,
the $G$-orbit stratification
of 
$\overline {\Cal O_r}$
passes to the 
$K^{\Bobb C}$-orbit stratification
of $\fra p^+$.
Thus, for $1 \leq s \leq r$,
restricted to $\Cal O_s = Ge^s$,
this homeomorphism
is a $K$-equivariant diffeomorphism
from $\Cal O_s$ onto the
$K^{\Bobb C}$-orbit 
in $\fra p^+$ 
of $x_s = \frac 12 (e^s + f^s - i h^s) \in \fra p^+$. 
\endproclaim

\demo{Proof}
This is a consequence of Lemma 3.3.6,
except perhaps the fact
that the indicated map is onto $\fra p^+$;
the latter, in turn, is implied by the observation that
$\fra p^+$ is the disjoint union of its $K^{\Bobb C}$-orbits. \qed
\enddemo

\noindent
{\smc Remark 3.3.12.}
In \cite\hilneors\ (p.~185)
a nilpotent orbit $\Cal O$ of $\fra g$
is said to be of convex type provided
it is contained in a proper generating invariant cone,
and a complete classification of orbits of convex type is given.
The orbits of convex type in this sense are precisely
the holomorphic and antiholomorphic ones.

\smallskip\noindent
{\smc Remark 3.3.13.}
It is known 
(Theorem 1.9 in \cite\sekiguch) that,
for a general real semisimple Lie algebra $\fra g$,
the
assignment to
$(e,f,h)$ of  $\frac 12 (e+f-ih)$
where $(e,f,h)$ runs through 
invariant $\fra {sl}(2)$-triples in $\fra g$
yields a bijective correspondence
$\Cal O \mapsto c(\Cal O)$,
referred to in the literature
as {\it Kostant-Sekiguchi correspondence\/}, 
between
nilpotent ($G$-)orbits in $\fra g$ and
nilpotent $K_{\Bobb C}$-orbits in $\fra p_{\Bobb C}$
(i.~e. $K_{\Bobb C}$-orbits 
in $\fra p_{\Bobb C}$ which arise as intersection of
a nilpotent orbit in 
$\fra g_{\Bobb C}$
with $\fra p_{\Bobb C}$).
In particular, 
for a semisimple Lie algebra $(\fra g,z)$ of hermitian type,
a little thought reveals that
the holomorphic nilpotent orbits $\Cal O$ in $\fra g$
are precisely those
which 
have the property that the orbit $c(\Cal O)$ 
lies in $\fra p^+$.
In \cite\vergnsix,
a nilpotent orbit $\Cal O$ of $\fra g$
is defined to be holomorphic provided
$c(\Cal O)$ lies in $\fra p^+$;
it is then observed that, 
for a 
nilpotent orbit $\Cal O$ 
having the property that
$c(\Cal O)$ lies in $\fra p^+$,
the projection from
$\Cal O$ 
to $\fra p$ is a diffeomorphism
onto its image.
The argument in \cite\vergnsix\  relies
on the 
fact, established in \cite\vergnsix\ 
that,
for an arbitrary nilpotent orbit $\Cal O$,
the Kostant-Sekiguchi correspondence 
may be realized by a diffeomorphism
from $\Cal O$ to $c(\Cal O)$;
this diffeomorphism involves a certain flow  
constructed by
Kronheimer in \cite\kronhtwo\ 
and is
referred to in the literature as
{\it Kronheimer-Vergne diffeomorphism\/}.
Our approach avoids the detour
via 
this diffeomorphism.

\smallskip\noindent{Remark 3.3.14.}
The classification of nilpotent orbits in the
real exceptional Lie algebras given in \cite\djokotwo\ 
(coming into play in the proof of Theorem 3.3.3)
relies
on the Kostant-Sekiguchi 
correspondence. 
We now indicate briefly how, in the two exceptional cases,
the 
holomorphic nilpotent orbits 
can be classified
without reference to the Kostant-Sekiguchi correspondence:
Over the complex numbers,
the classification
of nilpotent orbits
is in terms of 
\lq\lq regular\rq\rq\ subalgebras 
of $\fra g$ of minimal rank containing a representative
of the nilpotent orbit being classified,
the requisite subalgebra being written as $L_1$ in \cite\elkinone\ 
and $L(x)$ in \cite\mizuntwo.
These subalgebras 
actually play a role similar to that of
the indecomposable types in \cite\burcush\ 
for the classical cases which we exploited for the proofs of Theorems 3.5.3 
and 3.6.2. As before, $G$ refers to the adjoinnt group of the Lie algebra
$\fra g$ under discussion.
\smallskip\noindent
$\fra g =\fra e_{6(-14)}$:
The orbits $Ge^1$ and $Ge^2$ complexify to the orbits
in 
$\fra g_{\Bobb C}$
given in terms of the subalgebras
$L_1$ and/or $L(x)$
written
as $A_1$ and $2A_1$, respectively,
on p.~152 of
\cite\elkinone\ 
and in Table 1 on p.~447 of
\cite\mizuntwo.
Inspection of these tables shows that
no other nilpotent orbit 
in $\fra g$
can be holomorphic:
given any such orbit,
it complexifies to an orbit
involving
a regular subalgebra 
of minimal rank
which 
does {\it not\/}
arise by complexification from a
reductive Lie algebra of hermitian type.
Hence 
$\fra g =\fra e_{6(-14)}$ has no 
holomorphic nilpotent orbits
other than $\{0\},\Cal O_1=Ge^1, \Cal O_2=Ge^2$.
A slight extension of this reasoning shows that
the  orbits $\Cal O_{t,u}$ where $0 \leq t,u \leq 2$ and
$t+u \leq 2$ constitute a complete list of pseudoholomorphic
nilpotent orbits.
\smallskip\noindent
$\fra g =\fra e_{7(-25)}$:
The orbits $Ge^1$, $Ge^2$, and $Ge^3$ complexify to the orbits
in 
$\fra g_{\Bobb C}$
given in terms of the subalgebras
$L_1$ and/or $L(x)$
written, respectively,
as $A_1$, $2A_1$,
$(3A_1)''$,
on p.~153 of
\cite\elkinone\ 
and in Table 2 on p.~447 of
\cite\mizuntwo.
Inspection of these tables shows that
no other nilpotent orbit 
in $\fra g$
can be holomorphic
except possibly one which complexifies to an orbit
in 
$\fra g_{\Bobb C}$
given in terms of the subalgebra $L_1$
and/or $L(x)$
written as $(3A_1)'$.
However, inspection of Table 9 on p.~454 of \cite\mizuntwo\ 
reveals that the stabilizer of this orbit 
has dimension 69 whence the dimension of this orbit is strictly
larger than that of $\fra p$ and hence this orbit cannot possibly arise
from a (real) pseudoholomorphic nilpotent orbit.
Actually, this orbit does not arise from any real orbit at all.
Hence 
$\fra g =\fra e_{7(-25)}$ has no other
holomorphic nilpotent orbits
than $\{0\},\Cal O_1=Ge^1, \Cal O_2=Ge^2, \Cal O_3=Ge^3$.
A slight extension of this reasoning shows that
the  orbits $\Cal O_{t,u}$ where $0 \leq t,u \leq 3$ and
$t+u \leq 3$ constitute a complete list of pseudoholomorphic
nilpotent orbits.

\smallskip\noindent
{\smc 3.4. Non-nilpotent pseudoholomorphic orbits.}
Let $(\fra g,z)$ be a simple Lie algebra of hermitian type
and split rank $r$, with Cartan decomposition
$\fra g = \fra k \oplus \fra p$. Let $\varepsilon >0$
and let
$\Cal O^+(\varepsilon)$ 
and $\Cal O^-(\varepsilon)$
be the 
(semisimple elliptic) orbits of
$2 \varepsilon z$ and 
$-2 \varepsilon z$, respectively;
these generalize the non-nilpotent pseudoholomorphic orbits 
in $\fra {sl}(2,\Bobb R)$, cf. (3.2.2) above.

\proclaim{Theorem 3.4.1}
An orbit of the kind
$\Cal O^+(\varepsilon)$ 
is holomorphic
and one of the kind
$\Cal O^-(\varepsilon)$
is antiholomorphic.
\endproclaim

\demo{Proof} Indeed, for simplicity, let
$\varepsilon = \frac 12$, so that
$\Cal O^+(\varepsilon)= Gz$.
Since
$G=KA_{\fra p}K$,
the orbit $\Cal O^+(\varepsilon)$ may be written as
$\Cal O^+(\varepsilon)= Gz=KA_{\fra p}z$.
We now assert that
the projection from
$Gz$ to $\fra p$ is a diffeomorphism.
In order to see this,
exploiting
Proposition 3.3.2, pick a maximal 
abelian subalgebra $\fra a_{\fra p}$ of $\fra p$,
let
$(e_1,f_1,h_1)$,\dots,$(e_r,f_r,h_r)$
be the corresponding  system
of $H_1$-triples
so that, in particular, 
the elements $h_1,\dots,h_r$ constitute a basis of
$\fra a_{\fra p}$, and consider the resulting $H_1$-embedding
of $(\fra{sl}(2,\Bobb R)^r,Z^r)$
into  $(\fra g,z)$.
Let $\widetilde z =\frac 12 (e_1-f_1 + \dots + e_r - f_r)$.
Notice that, in view of Proposition 3.3.2(3), in the regular case,
$\widetilde z = z \in \fra g$.
Whether or not $(\fra g, z)$ is regular,
since the embedding 
of $(\fra{sl}(2,\Bobb R)^r,Z^r)$
into  $(\fra g,z)$
is an $H_1$-embedding,
we have $[z,h_j] =[\widetilde z,h_j]$ for every $j$ ($1 \leq j \leq r$). 
Consequently
$A_{\fra p}$ leaves the point $z-\widetilde z$ of $\fra g$
invariant.
Hence
$A_{\fra p}z = A_{\fra p}\widetilde z + z-\widetilde z$.
Furthermore,
since $\fra a_p$ is the span of
$h_1,\dots,h_r$,
the orbit $A_{\fra p}\widetilde z$
consists of 
the elements
$\frac 12 (c_1 + \dots+  c_r)\in \fra g$ 
where 
$c_j =\roman{ch}(\lambda_j)(e_j-f_j) + \roman{sh}(\lambda_j)(e_j+f_j)$,
for $1 \leq j \leq r$,
with $\lambda_1,\dots,\lambda_r \in \Bobb R$. Thus
the orbit $A_{\fra p} z$
consists of 
the points
$\frac 12 (c_1 + \dots+  c_r)+ z-\widetilde z\in \fra g$. 
Since 
$e_j-f_j \in \fra k$
($1 \leq j \leq r$) and since
$z-\widetilde z\in \fra k$,
the projection from
the orbit $A_{\fra p}z$ to $\fra p$
sends each
$\frac 12 (c_1 + \dots+  c_r) + z-\widetilde z$
to
$\frac 12 (\roman{sh}(\lambda_1)(e_1+f_1)+\dots+ 
\roman{sh}(\lambda_r)(e_r+f_r))$,
and this projection
is a diffeomorphism onto its image.
Since
$$
Gz=KA_{\fra p}z = K(A_{\fra p}\widetilde z + z-\widetilde z)
=
K(A_{\fra p}\widetilde z -\widetilde z) + z,
$$
the elements $\frac 12 (c_1 + \dots+  c_r)+z-\widetilde z$
constitute a complete set of representatives of the
$K$-orbits of
$Gz$.
Furthermore, it is well known that
$\fra p = \cup_{y \in K} \roman{Ad}(y) \fra a_p$
whence the elements
$\frac 12 (\roman{sh}(\lambda_1)(h_1)+\dots+ 
\roman{sh}(\lambda_r)(h_r))$
constitute a complete set of representatives of the
$K$-orbits of $\fra p$; this fact entails that
the elements
$\frac 12 (\roman{sh}(\lambda_1)(e_1+f_1)+\dots+ 
\roman{sh}(\lambda_r)(e_r+f_r))$
constitute a complete set of representatives of the
$K$-orbits of $\fra p$ as well
since,
for $1 \leq j \leq r$, $J_z(e_j+f_j) = h_j$, and since the operator
$J_z$ may be realized by the adjoint action with an element of $K$.
Consequently
the projection from
the orbit $Gz$ to $\fra p$ is a bijection onto its image,
in fact, a diffeomorphism onto 
$\fra p$.
Moreover, since the orbit
$\Cal O^+(\varepsilon)(1) \times \dots \times \Cal O^+(\varepsilon)(r)$
in
$\fra {sl}(2,\Bobb R)^r$ is holomorphic where
$\Cal O^+(\varepsilon)(j)$ refers to the corresponding orbit
in the $j'$th copy of
$\fra {sl}(2,\Bobb R)$ in $\fra {sl}(2,\Bobb R)^r$,
cf. Proposition 3.3.2,
we conclude that
the orbit
$Gz$
is holomorphic in $\fra g$.
Likewise,
the orbit
$G(-z)$
is antiholomorphic in $\fra g$. \qed
\enddemo

\proclaim{Theorem 3.4.2}
The only non-nilpotent pseudoholomorphic orbits
are the 
semisimple 
holomorphic and antiholomorphic
ones, i.~e. those of the kind 
$\Cal O^+(\varepsilon)$ and
$\Cal O^-(\varepsilon)$.
\endproclaim

Thus there are no other holomorphic orbits than the 
semisimple
and the nilpotent ones, 
and the latter ones are
classified in Theorem 3.3.3 above.

\demo{Proof}
Given a (co)adjoint orbit $\Cal O$
of a general real or complex semisimple 
Lie algebra $\fra g$,
the closure
$\overline{\Cal O}$ contains a single semisimple
orbit $\Cal O_{\roman{ss}}$, and any other orbit in 
$\overline {\Cal O}$ fibers over 
$\Cal O_{\roman{ss}}$
with typical fiber  a nilpotent orbit
for some reductive subalgebra $\fra g'$ of $\fra g$
which is a proper subalgebra when
$\Cal O_{\roman{ss}}$ is non-trivial.
More precisely,
the assignment to
$x \in \overline {\Cal O}$ 
of the semisimple constituent
$s$ in its  Jordan decomposition
$x = s + u$ is a fiber bundle map
$\pi \colon \overline {\Cal O} \to 
\Cal O_{\roman{ss}}$;
for $s\in \Cal O_{\roman{ss}}$,
the  stabilizer
$Z_G(s)$ of $s$ in the adjoint group $G$
is reductive, and
the fiber $\pi^{-1}(s)$
is (isomorphic to) the 
closure
$\overline {(Z_G(s)) u} \subseteq \fra z_{\fra g}(s)$
of the orbit 
$(Z_G(s)) u$
in $\fra z_{\fra g}(s) = \roman{Lie}(Z_G(s))$
of any $u$ with 
$x = s + u\in \Cal O$.
Thus,
given a (real) simple Lie algebra of hermitian type
$(\fra g,z)$,
let
$\Cal O$
be a non-nilpotent pseudoholomorphic orbit,
and let  
$\Cal O_{\roman{ss}}$
be the unique semisimple
orbit in $\overline {\Cal O}$.
Since $\fra k_{\Bobb C}$ has the same rank as
$\fra g_{\Bobb C}$ and since
the span of $e_1-f_1,\dots, e_r -f_r$ is a Cartan subalgebra
of
$\fra k_{\Bobb C}$
and hence of 
$\fra g_{\Bobb C}$,
there are complex numbers $\lambda_1,\dots,\lambda_s$
($s \leq r$)
such that 
the complexification $\Cal O_{\roman{ss}}^{\Bobb C}$
is the $G_{\Bobb C}$-orbit of
$w=\lambda_1(e_1-f_1) +\dots +\lambda_s (e_s -f_s)$.
Since $\Cal O$ is pseudoholomorphic, the projection
from $\Cal O$ to $\fra p$ is a diffeomorphism onto its image,
and the induced tangent map
$\roman T_{w+u}(\Cal O^{\Bobb C}) \to \fra p_{\Bobb C}$
is injective for any $u$ with $w+u \in \Cal O^{\Bobb C}$,
the complexification of $\Cal O$.
Moreover, 
the tangent space
$\roman T_{w+u}(\Cal O^{\Bobb C})$
contains 
the tangent space
$\roman T_w(\Cal O_{\roman{ss}}^{\Bobb C})$.
Hence $[y,w]= 0$ for every $y \in \fra k_{\Bobb C}$.
But 
$w \in \fra k_{\Bobb C}$ whence
$w$ lies in the center of
$\fra k_{\Bobb C}$, that is,
$w = \lambda z$ for some non-zero complex number $\lambda$.
Consequently
$\Cal O_{\roman{ss}}$ is the $G$-orbit of $\varepsilon z$
for some non-zero real number $\varepsilon$,
that is to say, 
$\Cal O_{\roman{ss}}$ is a 
semisimple holomorphic or antiholomorphic orbit,
and for dimensional reasons,
$\Cal O_{\roman{ss}}$
coincides with $\Cal O$. \qed
\enddemo

\proclaim{Addendum}
For a semisimple orbit of a simple  Lie algebra
$(\fra g,z)$ of hermitian type,
the property of being holomorphic or antiholomorphic
does not depend on the choice of $H$-element $z$.
\endproclaim

This is established in much the same way as the
Addendum to Theorem 3.3.3.

\smallskip\noindent{Remark 3.4.3.}
The holomorphic discrete series representations
arise from holomorphic quantization on 
integral holomorphic semisimple
orbits (but with reference to a complex structure different
from that coming from the projection to $\fra p$).
Thus our terminology \lq\lq holomorphic semisimple
orbit\rq\rq\ is consistent
with established terminology in the literature.

\smallskip\noindent
{\smc 3.5. The standard cases} $\fra{sp}(n,\Bobb R)$, $\fra{su}(p,q)$, 
$\fra{so}^*(2n)$.
\smallskip\noindent
{\smc 3.5.1. Explicit descriptions of these Lie algebras.}
We shall need them for ease of exposition and
later reference. Thus,
let $\Bobb K = \Bobb R, \Bobb C, \Bobb H$,
and let $V$ be an $n$-dimensional
{\it right\/} $\Bobb K$-vector space,
endowed with 
a positive definite hermitian form
$(\cdot,\cdot)$, that is 
$$ 
(x,y) =\overline{(y,x)},\quad 
(x,y\lambda) =(x,y)\lambda,
\quad
(x\lambda,y) =\overline \lambda(x,y),
\quad x,y \in V,\ \lambda \in \Bobb K,
$$
where the bar indicates conjugation as usual.
For right vector spaces
over the quaternions, 
requiring  linearity in the second variable
is more appropriate than
the more usual 
requirement
of linearity in the first variable.
After a suitable choice of basis,
$V$ gets identified with
$\Bobb K^n$, the algebra $\roman {End}_{\Bobb K}(V)$
may be identified with
the (left $\Bobb K$-vector) space
of $(n \times n)$-matrices
with entries from $\Bobb K$
which act on column vectors from the left in the usual fashion,
and 
$(\cdot,\cdot)$
is given by $(x,y) = \sum \overline {x_j}  y_j$ ($x,y \in V$);
in particular,
composition of endomorphisms now corresponds to
multiplication of matrices.
The
Lie group $\roman U(V,(\cdot,\cdot))$
of $(\cdot,\cdot)$-isometries
and its Lie algebra
$\fra u(V,(\cdot,\cdot))$
amount 
to
$\roman O(n,\Bobb R)$
and
$\fra {so}(n,\Bobb R)$
for $\Bobb K = \Bobb R$,
to
$\roman U(n)$
and
$\fra u(n)$
for $\Bobb K = \Bobb C$,
and to
$\roman {Sp}(n)$
and
$\fra {sp}(n)$
for $\Bobb K = \Bobb H$.
\smallskip
The requisite additional ingredient is
a 
($\Bobb K$-linear)
complex structure
$J_V$ on $V$ 
which is compatible with 
$(\cdot,\cdot)$,
that is,
$J_V \in \fra u(V,(\cdot,\cdot))$
or, equivalently,
for every $u,v \in V$, $(J_Vu,v) + (u,J_Vv) = 0$.
Define the non-degenerate skew-hermitian form $\Cal B$ on $V$ by
${
\Cal B(u,v) = (J_V u,v)}$, where $u,v \in V$.
We write
$\roman U(V,\Cal B)$
for the Lie group 
of $\Cal B$-isometries,
$\roman U(V,\Cal B)^0$
for its connected component of the identity,
and
$\fra u(V,\Cal B)=\roman {Lie}(\roman U(V,\Cal B))$
for its Lie algebra
(which consists of
all $X\in \roman {End}_{\Bobb K}(V)$ 
satisfying 
$\Cal B(X v, w) + \Cal B( v,X w) = 0$ 
for $ v, w \in V$).
Thus  
$X\in \roman {End}_{\Bobb K}(V)$
lies in
$\fra u(V,\Cal B)$
if and only if
the assignment
$$
\Cal B_X(v,w) = -\Cal B(X v, w)
\quad (=-(J_V X  v, w)),\quad  v, w \in V,
\tag3.5.1.1
$$
yields a hermitian form $\Cal B_X$ on $V$.
Then
$J_V$ lies in
$\fra u(V,\Cal B)$, and
$\Cal B_{J_V} = (\cdot,\cdot)$  
(the latter identity
justifies, in particular, the minus sign in (3.5.1.1)).
This sign also entails that,
cf. what is said after Theorem 3.5.4,
$\Cal B_X$
being non-negative
(in a sense to be explained)
is equivalent to
the holomorphicity of the orbit generated by $X$. 
The maximal compact subgroup
of $\roman U(V,\Cal B)$
and, accordingly, the Lie algebra thereof,
are the intersections
$\roman U(V,\Cal B)\cap \roman U(V,(\cdot,\cdot))$
and $\fra u(V,\Cal B)\cap \fra u(V,(\cdot,\cdot))$,
respectively.
Further, with the standard notation $X^*$ for the adjoint
of
$X\in \roman {End}_{\Bobb K}(V)$
so that
$(X u, v) =(u,X^* v)$ for every $u, v \in V$,
we have $J_V^* = -J_V$, and
$X\in \roman {End}_{\Bobb K}(V)$
lies in
$\fra u(V,\Cal B)$
if and only if $J_VX + XJ_V = 0 \in \roman {End}_{\Bobb K}(V)$.
Now $z_{\Cal B} = \frac 12 J_V$ is an $H$-element,
so that
$(\fra g,z_{\Cal B}) =(\fra u(V,\Cal B), z_{\Cal B})$
is a 
simple Lie algebra 
of hermitian type
in case $\Bobb K = \Bobb R$
or 
$\Bobb K = \Bobb H$,
and
a reductive Lie algebra
of hermitian type
in case $\Bobb K = \Bobb C$
whose semisimple part is simple
and whose center is a copy of $\Bobb R$.
We will write
$\fra {su}(V,\Cal B)$
for the corresponding simple Lie algebra
(for $\Bobb K = \Bobb R$
and $\Bobb K = \Bobb H$
there is no difference between
$\fra u(V,\Cal B)$ and
$\fra {su}(V,\Cal B)$).
With the notation $r$ for the corresponding
split rank
and $1,\Cal I, \Cal J, \Cal K$
for the unit quaternions, i.~e. standard real basis
of $\Bobb H$,
we now recall the
standard descriptions
of the three classes of {\it standard \/}
simple Lie algebras 
$(\fra g,z_{\Cal B})$
of hermitian type,
where $V$ is the {\it standard representation\/},
and where
the Cartan decomposition is written as $\fra k \oplus \fra p$.
Here and henceforth we will often write 
a matrix
$\left[\matrix 
\bold v_1 
\\
\ldots \\
\bold v_t
\endmatrix
\right]$ 
having row vectors
$\bold v_1, \dots,\bold v_t$
in the form
$[\bold v_1,\dots,\bold v_t]^{\tau}$ where the superscript $\tau$ 
stands for \lq\lq transpose\rq\rq.
\smallskip\noindent
(3.5.1.2) $\Bobb K=\Bobb R,\ n = 2 \ell,\ V = \Bobb R^n,\ r=\ell, \ 
J_V=\left [\matrix 0 & -I_{\ell}\\ I_{\ell}&0 \endmatrix \right]$,
\newline 
$\Cal B(\bold x,\bold y) 
= \sum_{j=1}^{\ell}( x_j y_{j+\ell}-x_{j+\ell} y_j),\
\bold x,\bold y \in V,\,
\roman U(V,\Cal B)=\roman U(V,\Cal B)^0= \roman {Sp}(\ell,\Bobb R)$
$$
\align
\fra g &=\fra u(V,\Cal B)=\fra{sp}(\ell,\Bobb R) =\left\{\left [\matrix
A&B\\
C&-A^{\roman t} 
\endmatrix
\right];
B^{\roman t} = B, \, C^{\roman t} = C,
\ A,B,C \in \roman M_{\ell,\ell}(\Bobb R)
\right\}
\\
\fra k &= \fra u(\ell) =
\left\{\left [\matrix
A&B\\
-B&A 
\endmatrix
\right];
A^{\roman t} = -A, \, B^{\roman t} = B \right\}
\\
\fra p &=  
\left\{\left [\matrix
A&B\\
B&-A 
\endmatrix
\right];
A^{\roman t} = A, \, B^{\roman t} = B \right\}
= \roman S_{\Bobb R}^2[\Bobb R^\ell] \oplus \roman S_{\Bobb R}^2[\Bobb R^\ell]
\\
\fra g 
&
= \fra k \oplus \fra p:
\quad
\left [\matrix
A&B\\
C&-A^{\roman t} 
\endmatrix
\right]
=
\left [\matrix
\frac {A-A^{\roman t}}2&\frac {B-C}2\\
\frac {C-B}2 & \frac {A-A^{\roman t}}2
\endmatrix
\right]
+
\left [\matrix
\frac {A+A^{\roman t}}2&\frac {B+C}2\\
\frac {C+B}2 & -\frac {A+A^{\roman t}}2
\endmatrix
\right] .
\endalign
$$

\smallskip\noindent
{\smc Remark.}
Our present sign convention is guided by the standard description
of angular momentum,
that is,
with the notation
$\bold x = [\bold q_1, \bold p_1]^{\tau}$ and
$\bold y = [\bold q_2, \bold p_2]^{\tau}$ 
where $\bold q_1, \bold q_2,\bold p_1, \bold p_2 \in \Bobb R^{\ell}$,
the given expression for $\Cal B(\bold x, \bold y)$ 
equals the expression for
$\omega([\bold q_1, \bold p_1]^{\tau},[\bold q_2, \bold p_2]^{\tau})$
for
the standard symplectic structure $\omega$ on $\Bobb R^{2\ell}$ 
evaluated at these vectors,
that is,
the skew form $\Cal B$ comes down to the standard symplectic form $\omega$ 
on $\Bobb R^{2\ell}$; for example, when $\ell =1$,
$\Cal B([q^1,p^1]^{\tau},[q^2,p^2]^{\tau}) = q^1 p^2 - q^2 p^1$.
See also what is said in Section 5 below.
When $(V, J_V)$ is viewed as a flat K\"ahler manifold,
the present sign for $J_V$ etc. is consistent with the standard
convention in K\"ahler geometry, to the effect that, in particular,
when $V$  is the ordinary plane, with coordinates $x,y$, 
the standard complex structure
is given by $J\frac {\partial}{\partial x} = \frac {\partial}{\partial y} $
and
$J\frac {\partial}{\partial y} = -\frac {\partial}{\partial x}$.
The sign convention for $J_V$ and $\Cal B$
in (3.2.2) above with 
$V= \Bobb R^2$ and
$J_V = E-F = \left[\matrix 0 & 1 \\ -1 & 0 \endmatrix \right]$,
as well as
in Satake's book \cite\satakboo\  and elsewhere in the literature,
is opposite to the present one.
Changing the sign convention 
amounts to interchanging holomorphic and antiholomorphic orbits. 
Under the circumstances of (3.2.2),
the present sign convention
may be realized via the $\fra {sl}(2)$-triple $(-E,-F,H)$ 
in $\fra {sl}(2,\Bobb R)$,
so that $\frac 12(F-E)$ is the corresponding $H$-element.
\smallskip\noindent
(3.5.1.3) $\Bobb K=\Bobb C,\ V = \Bobb C^n,\ 
n = p +q,\ p \geq q=r,\ 
J_V=i\left [\matrix I_p &0\\ 0& -I_q \endmatrix \right]$,  
\newline
$\Cal B(\bold x,\bold y) 
= i(\sum_{j=p+1}^{p+q}x_j \overline y_j-\sum_{j=1}^px_j \overline y_j),\ 
\bold x,\bold y \in V,\,
\roman U(V,\Cal B) =\roman U(V,\Cal B)^0 
=\roman U(p,q)$ 
$$
\align
\fra g=&\fra {su}(V,\Cal B)=\fra{su}(p,q) =
\left\{\left [\matrix
A&B\\
B^*&D 
\endmatrix
\right];
A^* = - A, \, D^* = -D, \, \roman{tr}(A) +\roman{tr}(D) =0\right\}
\\
\text{where\ } 
A &\in \roman M_{p,p}(\Bobb C),
D \in \roman M_{q,q}(\Bobb C),
B \in \roman M_{q,p}(\Bobb C)
\\
\fra  k=&(\fra u(p) \oplus \fra u(q))\cap \fra{su}(p,q),
\ 
\fra p =  
\left\{\left [\matrix
0&B\\
B^*&0 
\endmatrix
\right] \right\}\cong \roman M_{q,p}(\Bobb C) .
\endalign
$$
When $p>q$, $J_V$ lies in 
$\fra{u}(p,q)$ but not in $\fra{su}(p,q)$.
However, with reference to the decomposition 
$\fra{u}(p,q) = \fra z \oplus \fra{su}(p,q)$
where $\fra z \cong \Bobb R$ is the center of
$\fra{u}(p,q)$, 
the element $J_V$ decomposes as
$$
J_V =  i \frac {p-q}{p+q} I_{p+q}+ J'_V,\quad \text{where} \quad 
J'_V = i
\left [\matrix 
\frac {2q}{p+q} I_p &0\\ 
                   0& -\frac {2p}{p+q}I_q 
\endmatrix \right],
$$
and $\frac 12 J_V' \in \fra{su}(p,q)$ is an $H$-element;
since, in case $p>q$, the constituent
$i \frac {p-q}{p+q} I_{p+q}$ in the decomposition of $J_V$
lies in the center of $\fra{u}(p,q)$, it acts trivially
on $\fra p$ and the complex structure on $\fra p$
given by $\frac 12 \roman{ad}_{J_V'}$
within $\fra{su}(p,q)$
is the same as that 
given by $\frac 12 \roman{ad}_{J_V}$
within $\fra{u}(p,q)$.
See also Ex.~III.2.1 on p.~98 of \cite\satakboo.
\smallskip\noindent
(3.5.1.4) $\Bobb K=\Bobb H,\ V = \Bobb H^n,\ r = [\frac n2],\ 
J_V=\Cal J I_n,\  
\Cal B(\bold x,\bold y) = \sum_{j=1}^n \overline {x_j} \overline{\Cal J} y_j,\,
\bold x,\bold y \in V$,
\newline\noindent
$\roman U(V,\Cal B)= \roman O^*(2 n)$
(the non-compact dual of $\roman O(2 n,\Bobb R)$), 
$\roman U(V,\Cal B)^0= \roman {SO}^*(2 n)$
\linebreak
($=\roman {SU}^-(n,\Bobb H)$)
$$
\align
\fra g&=\fra u(V,\Cal B)=\fra{so}^*(2n)
=
\left\{\left [\matrix
A&-B\\
\overline B & \overline A 
\endmatrix
\right];
A^{\roman t} = -A,
\ 
B^* = B,\  A,B \in \roman M_{n,n}(\Bobb C)
\right\}
\\
\fra k &=\fra u(n) = 
\left\{\left [\matrix
A&-B\\
B&A 
\endmatrix
\right] \in \fra{so}^*(2n);
A= \overline A =-A^{\roman t} , \, B=\overline B= B^{\roman t}\right\}
\\
\fra p &=
\left\{i\left [\matrix
V&W\\
W&-V
\endmatrix
\right];
V^{\roman t} = -V, \, W^{\roman t} = -W,
\, V,W \in  \roman M_{n,n}(\Bobb R) \right\}
= i\Lambda^2_{\Bobb R}[\Bobb R^n] \oplus i\Lambda^2_{\Bobb R}[\Bobb R^n]
\\
\fra g 
&= \fra k \oplus \fra p:
\quad
\left [\matrix
A&-B\\
\overline B & \overline A 
\endmatrix
\right]
=
\left [\matrix
\frac {A+\overline A}2&- \frac {B+\overline B}2\\
\frac {B+\overline B}2 & \frac {A+\overline A}2
\endmatrix
\right]
+
\left [\matrix
\frac {A-\overline A}2&- \frac {\overline B-B}2\\
\frac {\overline B-B}2 & \frac {\overline A -A}2
\endmatrix
\right] .
\endalign
$$

\noindent
{\smc Remarks.}
(i) In each of these cases, the complex structure
$J_{z_{\Cal B}}$ on $\fra p$ is simply given by matrix multiplication
by $J_V$ from the left.
\newline\noindent
(ii) In case (3.5.1.3), the form $\Cal B$, written out as
$\Cal B = \Cal B_1 +i\Cal B_2$ with
real forms
$\Cal B_1$, $\Cal B_2$,
has 
$\fra u(V,\Cal B_1)=\fra {sp}(p+q,\Bobb R)$ and
$\fra u(V,\Cal B_2)=\fra {so}(2p,2q)$ whence
$\fra u(p,q)= \fra {so}(2p,2q) \cap\fra {sp}(p+q,\Bobb R)$.
The resulting inclusion
$$
\fra {su}(p,q) @>>> \fra {sp}(p+q,\Bobb R)
\tag3.5.1.5
$$
is an $H_1$-embedding.
Alternatively, $\fra u(p,q)$ 
arises as the centralizer of 
$$
\left[\matrix   0  &  0  & I_p & 0  \\
                0  &  0  & 0   & -I_q\\
              -I_p &  0  & 0   &  0  \\
                0  & I_q & 0   &  0
\endmatrix
\right] \in \fra {sp}(p+q,\Bobb R).
$$
\newline\noindent
(iii)
Still in case (3.5.1.3),
with reference to the Cartan decomposition
\linebreak
$\fra {sp}(q,\Bobb R) = \fra u(q) \oplus \fra p$
of $\fra {sp}(q,\Bobb R)$,
cf. (3.5.1.2),
the assignment to
$$
\left [\matrix
A&B\\
C&-A^{\roman t} 
\endmatrix
\right] = 
\left [\matrix
A'&B'\\
-B'&A' 
\endmatrix
\right] 
+
\left [\matrix
A''&B''\\
B''&-A'' 
\endmatrix
\right]
\in \fra {sp}(\ell,\Bobb R)
$$
of
$
\left[\matrix A'+iB' & 0 & iA''+B''
\\            0 & 0 & 0
\\            -iA''+B'' & 0 & A'-iB'
\endmatrix\right] 
\in \fra {su}(p,q)
$
when $p>q$ and of
$$
\left[\matrix    A'+iB' &  iA''+B''
\\            -iA''+B'' &  A'-iB'
\endmatrix\right] 
\in \fra {su}(p,q)
$$
when $p=q$
where
$\left [\matrix
A'&B'\\
-B'&A' 
\endmatrix
\right] 
 \in \fra u(q)$ and 
$\left [\matrix
A''&B''\\
B''&-A'' 
\endmatrix
\right] \in \fra p$
yields
a rank preserving $H_1$-injection
$$
\fra {sp}(q,\Bobb R) @>>> \fra {su}(p,q).
\tag3.5.1.6
$$
This is an embedding of the kind asserted in Proposition 3.3.2(4);
it identifies the two $H$-elements in the regular case $(p=q)$.
\newline\noindent
(iv)
In case (3.5.1.4), the form $\Cal B$, written out as
$\Cal B = \Cal B_1  +i\Cal B_2 + \Cal B_3 j$ with
real forms
$\Cal B_1$, $\Cal B_2$
and complex form 
$\Cal B_3$
has 
$\fra u(V,\Cal B_1)=\fra {sp}(2n,\Bobb R)$,
$\fra u(V,\Cal B_2)=\fra {so}(2n,2n)$,
and
$\fra u(V,\Cal B_3)=\fra {so}(2n,\Bobb C)$.
Hence
$\fra {so}^*(2 n)$ 
is the intersection 
in $\fra {gl}(4n,\Bobb R)$
of the Lie subalgebras
$\fra {sp}(2n,\Bobb R)$,
$\fra {so}(2n,2n)$,
and 
$\fra {so}(2n,\Bobb C)$,
and the resulting inclusion
$$
\fra {so}^*(2 n) @>>> \fra {sp}(2n,\Bobb R)
\tag3.5.1.7
$$
is  an $H_1$-injection.
\newline\noindent
(v)
Still in case (3.5.1.4), 
with reference to the Cartan decomposition
\linebreak
$\fra {sp}(\ell,\Bobb R) = \fra u(\ell) \oplus \fra p$
of $\fra {sp}(\ell,\Bobb R)$,
cf. (3.5.1.2),
the assignment to
$$
\left [\matrix
A&B\\
C&-A^{\roman t} 
\endmatrix
\right] = 
\left [\matrix
A'&B'\\
-B'&A' 
\endmatrix
\right] 
+
\left [\matrix
A''&B''\\
B''&-A'' 
\endmatrix
\right]
\in \fra {sp}(\ell,\Bobb R)
$$
of
$$
\left[\matrix U & -V  
\\            \overline V & \overline U 
\endmatrix\right] 
\in \fra {so}^*(4 \ell),
\quad
U= \left[\matrix A' & i B''
\\             -iB''&   A'
\endmatrix
\right],
\ 
V= \left[\matrix -B' & i A''
\\              -iA''&  -B'
\endmatrix
\right]
$$
where
$\left [\matrix
A'&B'\\
-B'&A' 
\endmatrix
\right] 
 \in \fra u(q)$ and 
$\left [\matrix
A''&B''\\
B''&-A'' 
\endmatrix
\right] \in \fra p$
yields
a rank preserving $H_1$-injection
of $\fra {sp}(\ell,\Bobb R)$
into $\fra {so}^*(4 \ell)$;
combination with the obvious injection
of 
$\fra {so}^*(4 \ell)$
into
$\fra {so}^*(4 \ell+2)$,
gives
a rank preserving $H_1$-injection
of $\fra {sp}(\ell,\Bobb R)$
into $\fra {so}^*(4 \ell+2)$ as well.
The resulting rank preserving $H_1$-morphism
$$
\fra {sp}\left(\left [\frac n2\right],\Bobb R\right)
@>>>
\fra {so}^*(2n),
\tag3.5.1.8
$$ 
for any $n\geq 2$,
is an embedding of the kind asserted in Proposition 3.3.2(4);
in the regular case ($n$ even), this embedding
identifies the $H$-element of the source with that of the target.
\newline\noindent
{\smc 3.5.2. Explicit $H_1$-triples.}
With reference to the notation introduced in (3.5.1.2),
for $1 \leq k \leq \ell$,
consider the  embedding
$\kappa_k \colon \fra{sp}(1,\Bobb R)@>>>\fra{sp}(\ell,\Bobb R)$
given by the assignment to
$\left [\matrix
a&b\\
c&-a 
\endmatrix
\right]\in
\fra{sp}(1,\Bobb R)$
of
$\left [\matrix
A_k&B_k\\
C_k&-A_k 
\endmatrix
\right]$
where
$A_k = \roman {diag}(0,\dots,0,a,0,\dots,0)$, 
$B_k = \roman {diag}(0,\dots,0,b,0,\dots,0)$, 
$C_k = \roman {diag}(0,\dots,0,c,0,\dots,0)$, 
($a,b,c$ in the $k'$th position),
and let $e_k = \kappa_k(E)\in \fra{sp}(\ell,\Bobb R)$,
$f_k = \kappa_k(F)\in \fra{sp}(\ell,\Bobb R)$,
$h_k = \kappa_k(H)\in \fra{sp}(\ell,\Bobb R)$.
Each $(e_k,f_k,h_k)$ 
($1 \leq k \leq \ell)$
is an $H_1$-triple
in $\fra{sp}(\ell,\Bobb R)$,
cf. Ex.~III.3.2 on p.~105/106 of \cite\satakboo,
and these combine to an $H_1$-embedding
of
$\fra{sl}(2,\Bobb R)^{\ell}$
into $\fra{sp}(\ell,\Bobb R)$
of the kind spelled out in Proposition 3.3.2.
A corresponding $H_1$-embedding
of
$\fra{sl}(2,\Bobb R)^q$
into $\fra{su}(p,q)$ ($p \geq q$)
results
from the corresponding
$H_1$-embedding
of
$\fra{sl}(2,\Bobb R)^q$
into $\fra{sp}(q,\Bobb R)$, combined
with the
$H_1$-embedding
(3.5.1.6)
of
$\fra{sp}(q,\Bobb R)$
into
$\fra{su}(p,q)$.
Alternatively, noting that
$\fra{su}(1,1)$ is isomorphic to
$\fra{sl}(2,\Bobb R)$,
we obtain an obvious embedding of
$\fra{sl}(2,\Bobb R)^q\cong \fra{su}(1,1)^q$
into $\fra{su}(q,q)$; when $p>q$,
this embedding, combined with
the obvious embedding of
$\fra{su}(q,q)$
into $\fra{su}(p,q)$,
yields the requisite embedding of
$\fra{sl}(2,\Bobb R)^q$ into
$\fra{su}(p,q)$.
Likewise,
the composite
of the
$H_1$-embedding
of
$\fra{sl}(2,\Bobb R)^{[\frac n2]}$
into $\fra{sp}([\frac n2],\Bobb R)$
with the
$H_1$-embedding
(3.5.1.8)
of
$\fra{sp}([\frac n2],\Bobb R)$
into $\fra{so}^*(2n)$
yields
a corresponding $H_1$-embedding
of
$\fra{sl}(2,\Bobb R)^{[\frac n2]}$
into
$\fra{so}^*(2n)$.

 \smallskip
We will say that a matrix $X$ in 
$\fra g$
has {\it square zero\/}
provided $X^2 = 0$,
the square being taken in the algebra 
$\roman {End}_{\Bobb K}(V)$
of endomorphisms
of 
the standard representation $V$ of $\fra g$.
A square zero matrix $X$ in $\fra g$
is a nilpotent element of
$\fra g$ 
(i.~e. $\roman{ad}(X)$ is a nilpotent endomorphism of $\fra g$) 
since $X^2=0$ entails $(\roman{ad}(X))^3 =0$.

\proclaim{Theorem 3.5.3}
A nilpotent orbit
of a standard simple  Lie algebra
of hermitian type
is pseudoholomorphic
if and only if it is generated by
a square zero matrix.
More precisely, 
given a standard simple  Lie algebra
$(\fra g,z)$ 
of hermitian type
and split rank $r$,
for $0 \leq t,u \leq r, t+u \leq r$,
the pseudoholomorphic nilpotent orbit $\Cal O_{t,u}$
is the smooth connected manifold in 
$\fra g$
which consists of 
square zero
matrices 
in $\fra g$
which have $\Bobb K$-rank $t+u$
and $\Bobb K$-signature $t-u$,
and 
the $\frac {(r+1)(r+2)}2$ pseudoholomorphic nilpotent orbits
$\Cal O_{t,u}$ constitute a complete list of all pseudoholomorphic
nilpotent orbits. 
\endproclaim

\demo{Proof} 
We will explain in detail the case 
$\fra g = \fra{sp}(\ell,\Bobb R)$
and will  give some hints how the other cases
may be reduced to this one.
Thus,
let $X \in \fra{sp}(\ell,\Bobb R)$ and suppose that $X^2=0$.
With reference to the Cartan decomposition
of $\fra{sp}(\ell,\Bobb R)$, cf. (3.5.1.2),
write
$$
X = \left[\matrix A & B \\
-B & A \endmatrix
\right]
+
\left[\matrix S & T \\
T & -S \endmatrix
\right],\quad
A,B,S,T \in \roman M_{\ell,\ell}(\Bobb R)
\tag3.5.3.1
$$
where $A$ is skew-symmetric and $B,S,T$ are symmetric.
Since every element in $\fra u(\ell)$ lies in the Lie algebra
of a maximal torus, we may assume that
$A=0$ and that $B$ is a diagonal matrix.
Then a straightforward calculation shows that
$X^2=0$
implies that $S$ and $T$ are diagonal matrices as well,
and that $B^2 = S^2 + T^2$.
Hence
$X \in \fra g=\fra{sp}(\ell,\Bobb R)$ is determined by its projection
$\left[\matrix S & T \\
T & -S \endmatrix
\right]$ to 
$\fra p=S_{\Bobb R}^2[\Bobb R^{\ell}] \oplus S_{\Bobb R}^2[\Bobb R^{\ell}]$,
up to the signs in the entries of $B$,
and a choice of signs
corresponds precisely to
a choice of signature for the hermitian form
$\Cal B_X$, once the rank is fixed.
\smallskip
Consider
the subset
$\widetilde {\Cal O}_{t,u}$ of $\fra g$
which consists of 
square zero
matrices 
in $\fra g$
which have rank $t+u$
and signature $t-u$.
Inspection shows that
$X=e_1+ \dots+ e_t - e_{t+1} - \dots -e_{t+u}$,
cf. Proposition 3.3.2 and Subsection 3.5.2,
is a square zero matrix in $\fra g$.
Hence
the pseudoholomorphic nilpotent orbit $\Cal O_{t,u}$
is contained in $\widetilde {\Cal O}_{t,u}$.
Moreover, for $X \in \widetilde {\Cal O}_{t,u}$,
written out as in (3.5.3.1),
after diagonalization so that $A=0$ and so that $B$ is a diagonal matrix,
the signature being constant
(and equal to $t-u$)
amounts
to making a possible choice of signs for $B$.
Consequently the restriction
to 
$\widetilde {\Cal O}_{t,u}$
of the projection
from $\fra g = \fra k \oplus \fra p$
to $\fra p$
is a bijection
onto its image in $\fra p$.
We will show shortly that
$\widetilde {\Cal O}_{t,u}$
coincides with $\Cal O_{t,u}$.
\smallskip
In view of the classification of
real adjoint orbits by {\it types\/} \cite\burcush,
nilpotent orbits are classified by 
sums of {\it indecomposable (nilpotent) types\/},
where a(n indecomposable)
\lq\lq type\rq\rq\ is a generalized Jordan normal form (Jordan block);
cf. also (9.3) in \cite\collmcgo.
Under our circumstances, this means the following:
Given a nilpotent orbit $\Cal O$ of $\fra g$ and $X \in \Cal O$,
with reference to the operator $X$
on $V= \Bobb R^{2\ell}$,
the space $V$ decomposes into a sum
$V= V_1 \oplus \dots \oplus V_{\tau}$
of $X$-modules $V_j$
which are mutually isotropic with respect to the symplectic form
$\omega$ (say) on $V$, that is,
$\omega(V_j, V_{j'}) = 0$
for $j \neq j'$,
and which are either $X$-{\it indecomposable\/}
or are decomposed into two indecomposable $X$-stable 
subspaces of the same dimension,
and each $V_j$ 
which is $X$-indecomposable,
i.~e. not decomposed
into two $X$-stable subspaces of the same dimension,
carries as an additional ingredient
a sign $\pm$
which encapsulates the behaviour of the hermitian form
$\Cal B_X$.
A space of the kind $\widetilde {\Cal O}_{t,u}$
corresponds to a sum of indecomposable types of the kind
$$
V = \Delta_1^+(0)\oplus \dots\oplus\Delta_1^+(0)
\oplus
\Delta_1^-(0)
\oplus \dots\oplus\Delta_1^-(0)
\oplus
\Delta_0(0,0)
\oplus \dots\oplus\Delta_0(0,0)
$$
involving $t$ copies
of $\Delta_1^+(0)$,
$u$ copies
of $\Delta_1^-(0)$,
and
$r-t-u$
copies of $\Delta_0(0,0)$;
here the notation is that in \cite\burcush,
each
summand is a real symplectic vector space
of dimension two,
for $X\in \widetilde {\Cal O}_{t,u}$
the form $\Cal B_X$
is positive on each copy of $\Delta_1^+(0)$,
negative on each copy of $\Delta_1^+(0)$,
and zero 
on each copy of 
$\Delta_0(0,0)$.
This shows in particular that each
$\widetilde {\Cal O}_{t,u}$
is a nilpotent orbit
which necessarily coincides with
$\Cal O_{t,u}$
and that, when $X \in \fra g$
has square zero,
the only possible indecomposable types
in the decomposition of $V$
are
$\Delta_1^+(0)$,
$\Delta_1^-(0)$,
and $\Delta_0(0,0)$.
Indeed, when a higher type occurs,
the operator $X^2$ is non-zero on this type.
\smallskip
The proof is then completed by the observation that,
if a nilpotent orbit $\Cal O$ involves 
a type $\Delta$ of a kind
different from those spelled out above,
the 
projection from $\Cal O$ to $\fra p$
is no longer injective.
Indeed, it suffices to consider the special case
where $V$ 
is the underlying space of such a higher type $\Delta$,
that is, where $\fra g = \fra{sp}(V)=\fra{sp}(m,\Bobb R)$ 
for some $1<m \leq \ell$.
There are only two kinds of types:
\smallskip
\noindent
1) $X^{2m -1} \ne 0$.
In this case, $\Cal O$ is 
one of the two {\it real\/} regular nilpotent orbits
in 
$\fra{sp}(m,\Bobb R)$.
However, the dimension of each of the two regular orbits
is plainly strictly larger than that of
the constituent $\fra p$ in the Cartan decomposition
of
$\fra{sp}(m,\Bobb R)$ 
unless $m=1$ but
the case
$m=1$ is excluded here
since we suppose
that
the type $\Delta$ under consideration is
different from those spelled out above.
Hence 
the projection from $\Cal O$ to $\fra p$
cannot then be injective. 
Indeed,
$\dim (\fra{sp}(m,\Bobb R)) = 2 m^2 + m$
and
$\dim(\Cal O) = 2 m^2$
whereas
$\dim(\fra p) = m^2 + m$. 
In the literature, this type is written
as $\Delta^{\pm}_k(0)$, with $k=2m-1$.
\newline\noindent
2) As an $X$-module, $V$ decomposes into a sum $V=V' \oplus V''$ of
two Lagrangian subspaces
$V' \cong \Bobb R^m$ and
$V'' \cong \Bobb R^m$, each of which is $X$-indecomposable,
in such a way that the orbit $\Cal O$
arises from
the 
(nonnegative or from the nonpositive)
regular nilpotent orbit of $\fra {gl}(m,\Bobb R)$
in the following fashion:
Consider the canonical injection
of $\fra {gl}(m,\Bobb R)$ into
of $\fra {sp}(m,\Bobb R)$
which sends a matrix
$A \in \fra {gl}(m,\Bobb R)$
to the matrix
$\widehat A =
\left[\matrix A & 0 \\ 0 & -A^t\endmatrix \right] \in \fra {sp}(m,\Bobb R)$.
Up to a choice of basis in $V'$ and $V''$,
the endomorphism $X$ of $V$
is now of the kind
$X=\widehat A$,
where $A$ is a regular nilpotent element of 
$\fra {gl}(m,\Bobb R)$,
and it does not matter whether
$A$ is taken from the nonnegative or nonpositive  
regular nilpotent orbit of
$\fra {gl}(m,\Bobb R)$.
For example, $A$ could be the ordinary Jordan normal form
of a nilpotent rank $m-1$ matrix in
$\fra {gl}(m,\Bobb R)$.
In the literature, this type is written
as $\Delta_k(0,0)$, with $k=m-1$.
(It yields new types only
for $k$ even.)
To verify that
the 
projection from $\Cal O$ to the constituent $\fra p$
of the Cartan decomposition
$\fra{sp}(m,\Bobb R) = \fra u(m) \oplus \fra p$ 
is no longer injective,
it suffices to 
prove that
the 
projection from one
of the two regular nilpotent orbits
$\Cal O'$ (say) 
of
$\fra {gl}(m,\Bobb R)$
to the constituent $\fra p' \cong \roman S_{\Bobb R}^2[\Bobb R^m]$
of the Cartan decomposition
$\fra{gl}(m,\Bobb R) = \fra {so}(m) \oplus \fra p'$ 
is not injective.
However,
$\fra{gl}(m,\Bobb R)$ is reductive and not simple, and
$\Cal O'$ 
is actually an orbit
in  
$\fra{sl}(m,\Bobb R)$ whence
$\Cal O'$  projects to
the constituent $\fra p'' \subseteq \fra p'$
of the Cartan decomposition 
$\fra{sl}(m,\Bobb R) = \fra {so}(m) \oplus \fra p''$;
notice that
$\fra p''$ consists of the  symmetric trace zero
$(m\times m)$-matrices.
Thus, when $m\geq 3$,
$$
\dim(\Cal O') = m^2-m >\frac 12 m(m+1)-1=\dim (\fra p'')
$$
whence the 
projection from 
$\Cal O'$  
to $\fra p''$ cannot be injective.
The case $m=2$
(i.~e. the type $\Delta_1(0,0)$)
is already covered by the type
$\Delta_1^+(0)\oplus\Delta_1^-(0)$
examined earlier. 
\smallskip
A similar reasoning establishes the claim in the other standard cases.
Alternatively, the  embeddings 
(3.5.1.4)
of $\fra {su}(p,q)$ 
into
$\fra {sp}(p+q,\Bobb R)$
and 
(3.5.1.6)
of
$\fra {so}^*(2n)$
into $\fra {sp}(2n,\Bobb R)$
are compatible with the Cartan decompositions
and $H_1$-embeddings.
By means of these embeddings,
we may reduce the case
of $\fra {su}(p,q)$ 
and
that of
$\fra {so}^*(2n)$
to that of
$\fra {sp}(\ell,\Bobb R)$
(for suitable $\ell$).
More precisely: For
$\fra {su}(p,q)$, we have the chain of $H_1$-injections
$$
\fra {sp}(q,\Bobb R) \subseteq \fra {su}(p,q) \subseteq \fra {sp}(p+q,\Bobb R),
$$ 
and things may be arranged in such a way that,
cf. Proposition 3.3.2 above,
a complete system 
$(e_1,f_1,h_1)$,\dots,$(e_q,f_q,h_q)$
of $H_1$-triples
for
$\fra {sp}(q,\Bobb R)$
constitutes a complete system of
$H_1$-triples for
$\fra {su}(p,q)$ and that,
together with suitable additional
$H_1$-triples
$(e_{q+1},f_{q+1},h_{q+1})$,\dots,$(e_{p+q},f_{p+q},h_{p+q})$,
the $H_1$-triples
$(e_1,f_1,h_1)$,\dots,$(e_{p+q},f_{p+q},h_{p+q})$
constitute a complete system
of $H_1$-triples for 
$\fra {sp}(p+q,\Bobb R)$.
Then, for $0 \leq t,u \leq q, t+u \leq q$,
the 
intersection
$\Cal O_{t,u} \cap \fra {su}(p,q)$ of the
corresponding
pseudoholomorphic nilpotent orbit $\Cal O_{t,u}$
for
$\fra {sp}(p+q,\Bobb R)$
is the 
corresponding
pseudoholomorphic nilpotent orbit 
for $\fra {su}(p,q)$.
The same kind of reasoning applies to $\fra {so}^*(2n)$.
We leave the details to the reader.
Furthermore, the statement
referring to the closures
of the orbits
$\Cal O_{t,u}$
is a consequence of the observation that
passing to the closure
amounts to dropping
summands in the decomposition of types. \qed
\enddemo

As noted by Springer-Steinberg \cite\spristei,
nilpotent orbits
in $\fra{sp}(\ell,\Bobb C)$
arise as the intersections of
nilpotent orbits in $\fra{sl}(2\ell,\Bobb C)$
with 
$\fra{sp}(\ell,\Bobb C)$.
Theorem
3.5.3 entails that,
given a nilpotent orbit $\widetilde {\Cal O}$ in
$\fra{sp}(\ell,\Bobb C)$
arising as the complexification of a pseudoholomorphic
nilpotent orbit in
$\fra{sp}(\ell,\Bobb R)$,
the intersection $\widetilde {\Cal O}\cap \fra{sp}(\ell,\Bobb R)$
decomposes into connected components
which are distinguished
by the signature
of the form $\Cal B_X$ for 
$X \in \widetilde {\Cal O}\cap \fra{sp}(\ell,\Bobb R)$.
The same kind of remark applies to 
$\fra{su}(p,q)$ and
$\fra{so}^*(2n)$.
In particular, 
in the standard cases,
we can thus single out
the holomorphic nilpotent orbits
among the pseudoholomorphic ones.
We will say that $X \in \fra u(V,\Cal B)$
is {\it non-negative\/}
({\it non-positive\/})
provided the hermitian form
$\Cal B_X$ is non-negative
(non-positive);
plainly, the property
of $X$ being non-negative
(non-positive)
depends only on the
$\roman U(V,\Cal B)^0$-orbit
of $X$,
and we can talk about
{\it non-negative\/}
({\it non-positive\/})
$\roman U(V,\Cal B)^0$-orbits.
The non-negative elements $X$ in 
$\fra u(V,\Cal B)$
constitute a (real) cone in
$\fra u(V,\Cal B)$
which is invariant under
$\roman U(V,\Cal B)^0$
(and the same is true, of course, of 
the non-positive elements).

\proclaim{Theorem 3.5.4}
Given a standard simple Lie algebra $(\fra g,z)$
of hermitian type and real rank $r$,
for $1 \leq s \leq r$,
the holomorphic nilpotent orbit $\Cal O_s$
{\rm (= $\Cal O_{s,0}$)}
is the (connected) smooth manifold 
of non-negative rank $s$ nilpotent
matrices 
in $\fra g$
(real rank $s$ for $\fra g = \fra{sp}(r,\Bobb R)$;
complex rank $s$ for $\fra g = \fra {su}(p,q)$, where $p \geq q = r$;
quaternionic rank  $s$
for $\fra g = \fra{so}^*(2n)$, where $r= [\frac n2]$).
Furthermore, as a real semi-algebraic set in $\fra g$, 
the closure $\overline {\Cal O_s}$ 
is the space of
non-negative nilpotent 
matrices 
in $\fra g$
which have rank at most $s$.
\endproclaim

\demo{Proof}
By Theorem 3.5.3, the 
holomorphic nilpotent orbit $\Cal O_s$
($1 \leq s \leq r$)
is the smooth connected manifold 
of non-negative rank $s$ nilpotent
square zero
matrices 
in $\fra g$.
Theorem 5.4.1 below implies that a non-negative nilpotent matrix
generates a holomorphic nilpotent orbit and is therefore of square zero,
and the \lq\lq Furthermore\rq\rq\ statement is implied by
that theorem as well. \qed
\enddemo

Thus, in the definition of the 2-form $\Cal B_X$, 
cf. (3.5.1.1),
the
sign has been adjusted in such a way that
$\Cal B_X$ being non-negative corresponds to 
the holomorphicity of the orbit generated by $X$.
Explicit equations and inequalities
describing
$\overline {\Cal O_s}$
as a real semi-algebraic set in $\fra g$ 
will be given in Remark 5.4.2.
\smallskip
In the standard cases, the complex analytic structures
of the holomorphic nilpotent orbits 
and their closures are elucidated by the following
result which, apart from the last statement referring to the normality
of the resulting complex varieties,
makes explicit 
the statement of Theorem 3.3.11 above
for the standard cases and in particular spells out 
the $K^{\Bobb C}$-stratification.
\smallskip
Recall that, with reference to the Cartan decomposition
$\fra g = \fra k \oplus \fra p$,
$\fra p^+=\roman S_{\Bobb C}^2[\Bobb C^{r}]$ 
for $\fra g=\fra {sp}(r,\Bobb R)$,
$\fra p^+=\roman M_{q,p}(\Bobb C)$
for $\fra g=\fra {su}(p,q)$ $(p \geq q = r)$, 
$\fra p^+=\Lambda^2_{\Bobb C}[\Bobb C^n]$
for $\fra g=\fra {so}^*(2n,\Bobb R)$ $(r= [\frac n2])$
where $r$ denotes the real rank.
For $1 \leq s < r$, let
$V_s$
be the smooth complex affine 
variety 
in 
$\fra p^+$
which consists of complex 
matrices 
in $\fra p^+$
having rank $s$
for $\fra g = \fra{sp}(r,\Bobb R)$
and
$\fra g = \fra {su}(p,q)$
and
rank  $2s$
for $\fra g = \fra {so}^*(2n)$;
likewise,
let
$V_{\leq s}$
be the complex affine 
determinantal variety 
in 
$\fra p^+$
which consists of complex 
matrices 
in $\fra p^+$
having rank at most $s$
for $\fra g = \fra{sp}(r,\Bobb R)$
and
$\fra g = \fra {su}(p,q)$
and
rank at most $2s$
for $\fra g = \fra {so}^*(2n)$.

\proclaim{Theorem 3.5.5}
For 
a standard simple rank $r$ Lie algebra
of hermitian type
$(\fra g,z)$,
under the 
projection from $\fra g = \fra k \oplus \fra p$ to
$\fra p$, followed by the complex linear isomorphism from
$\fra p$ to $\fra p^+$,
the closure
$\overline{\Cal O_r}$ 
of the principal holomorphic nilpotent orbit
is identified with $\fra p^+$ and,
for $1 \leq s < r$,
$\Cal O_s$ 
is mapped diffeomorphically onto
the smooth affine complex subvariety $V_s$ of $\fra p^+$.
Likewise,
for $1 \leq s < r$,
as a complex analytic space, 
$\overline {\Cal O_s}$ 
is identified with
the complex affine 
determinantal variety 
$V_{\leq s}$ in $\fra p^+$.
An explicit system of equations
defining
$\overline {\Cal O_s}$ as an affine complex variety is thus given by 
the determinantal equations which 
in the first two cases 
say that
all $((s+1) \times (s+1))$-minors
are zero
and which in the third case say
that
all 
$(2(s+1) \times 2(s+1))$-minors are zero.
As a complex variety, each $\overline {\Cal O_s}$ is normal. \qed
\endproclaim

Determinantal varieties have been studied in the literature;
see e.~g. \cite\brunvett, \cite\hocheago;
they are known to be normal
\cite\brunvett\ (Theorem 2.11), \cite\hocheago.
Thus,
at the present stage, the proof of Theorem 3.5.5 is complete.
A result 
more general than Theorem 3.5.5
will be given as Theorem 5.3.3 below
which will, in particular, entail the normality of the closures of the strata
and hence the normality of complex determinantal
varieties.
The reasoning 
aimed at unravelling the complex analytic structures
of the categorical quotients 
$W_J\big /\big / {\KKK}^{\Bobb C}$
in the proof of Theorem 5.3.3 below
establishes the claim of Theorem 3.5.5 as well.

\smallskip\noindent
{\smc Remark 3.5.6.}
The proof of Theorem 3.5.3 is independent of those of Theorems 3.3.3 and
3.3.11 and thus yields another proof of the statement of 
Theorem 3.3.3 and
of that of Theorem
3.3.11 
for the special case where
$(\fra g,z)$ is a standard simple rank $r$ Lie algebra
of hermitian type.

\smallskip\noindent
{\smc 3.6. The Lie algebra $\fra g=\fra{so}(2,q),\ q \geq 2$.}
This is the remaining classical case. 
We note that $\fra{so}(2,1)\cong \fra{sp}(1,\Bobb R)$,
$\fra{so}(2,3)\cong \fra{sp}(2,\Bobb R)$,
$\fra{so}(2,4)\cong \fra{su}(2,2)$,
$\fra{so}(2,6)\cong \fra{so}^*(8)$,
and that
$\fra{so}(2,10)$ arises from
the real Cayley division algebra in a similar fashion; see Section 8 below.
Recall that
$$
\fra g=\fra{so}(2,q)=
\left\{
\left [\matrix
A&B\\
B^{\roman t}&D 
\endmatrix\right];
A \in \fra {so}(2,\Bobb R),
D \in \fra {so}(q,\Bobb R),
B \in \roman M_{q,2}(\Bobb R)
\right\}.
$$
Its split rank $r$ equals 2, and 
$z = \left [\matrix
J&0\\
0&0 
\endmatrix
\right]$ where $J=\left [\matrix 0&1\\ -1&0 \endmatrix\right]$
is an $H$-element, so that
$(\fra g,z)$ is a simple Lie algebra of hermitian type
for $q \geq 3$.
The constituents $\fra k$ and $\fra p$ of its Cartan decomposition
have the form
$\fra k = \fra {so}(2,\Bobb R) \oplus \fra {so}(q,\Bobb R)$
and
$\fra p = 
\left\{
\left [\matrix
0&B\\
B^{\roman t}&0 
\endmatrix\right] 
\right\}\cong \roman M_{q,2}(\Bobb R)$
and,
for $x =
\left [\matrix
0&B\\
B^{\roman t}&0 
\endmatrix\right] 
\in \fra p$,
the 
value $J_z(x)$
under the complex operator
$J_z$ on 
$\fra p$
is easily seen to be given by matrix multiplication
by $J$ on $B$ from the left.
\smallskip
The Lie algebra $\fra{so}(2,2)$
is generated by
$$
\align
X&=\left [\matrix
J&0\\
0&J 
\endmatrix
\right],
\quad
A_1=\left [\matrix
0&0&1&0
\\
0&0&0&-1
\\
1&0&0&0
\\
0&-1&0&0
\endmatrix
\right],
\quad
A_2=\left [\matrix
0&0&0&1
\\
0&0&1&0
\\
0&1&0&0
\\
1&0&0&0
\endmatrix
\right]
\\
Y&=\left [\matrix
J&0\\
0&-J 
\endmatrix
\right],
\quad
B_1=\left [\matrix
0&0&0&1
\\
0&0&-1&0
\\
0&-1&0&0
\\
1&0&0&0
\endmatrix
\right],
\quad
B_2=\left [\matrix
0&0&1&0
\\
0&0&0&1
\\
1&0&0&0
\\
0&1&0&0
\endmatrix
\right],
\endalign
$$
these matrices satisfy the relations
$$
\align
[A_1,A_2] &= 2X,
\quad
[X,A_1] = -2A_2,
\quad
[X,A_2] = 2A_1,
\\
[B_1,B_2] &= 2Y,
\quad
[Y,B_1] = -2B_2,
\quad
[Y,B_2] = 2B_1,
\endalign
$$
and 
$(X,A_1,A_2)$
and
$(Y,B_1,B_2)$
each generate a copy of
$\fra{sl}(2,\Bobb R) \cong \fra{so}(2,1)$ in
$\fra{so}(2,2)$
in such a way that
$\fra{so}(2,2)= \fra{so}(2,1) \oplus \fra{so}(2,1)$.
Let 
$$
\gather
e_1 = \frac 12 (A_2+X),
\quad
f_1 = \frac 12 (A_2-X),
\quad
h_1 = A_1,
\\
e_2 = \frac 12 (B_2+Y),
\quad
f_2 = \frac 12 (B_2-Y),
\quad
h_2 = B_1 .
\endgather
$$
These lie in
$\fra{so}(2,2)$ and,
via the standard embedding of
$\fra{so}(2,2)$
into $\fra{so}(2,q)$ for $q \geq 3$, 
we view
$e_1,f_1,h_1,e_2,f_2,h_2$
as elements of 
$\fra{so}(2,q)$; 
inspection 
shows that
$(e_1,f_1,h_1)$
and $(e_2,f_2,h_2)$ 
are holomorphic
$\fra{sl}(2)$-triples
in $\fra{so}(2,q)$, 
and these combine
to an $H_1$-embedding
of $\fra{sl}(2,\Bobb R)^2$  into $\fra{so}(2,q)$
of the kind spelled out in Proposition 3.3.2.
By Theorem 3.3.3, with the notation
$G=\roman {SO}(2,q)^0$,
for $q \geq 3$,
the orbits
$$
\gathered
\Cal O_{0,0}= \{0\},
\ \Cal O_{1,0} = Ge_1,
\ \Cal O_{0,1} = G(-e_1),
\\
\ \Cal O_{2,0} = G(e_1+e_2),
\ \Cal O_{1,1} = G(e_1-e_2),
\ \Cal O_{0,2} = G(e_1+e_2),
\endgathered
\tag3.6.1
$$
are pseudoholomorphic,
$\Cal O_{1,0},\, \,\Cal O_{2,0}$
are holomorphic,
$\Cal O_{0,1},\, \,\Cal O_{0,2}$
are antiholomorphic,
and $\Cal O_{1,1}$
is neither holomorphic nor antiholomorphic.
Recall that
$\roman {SO}(2,q)$
and
$\roman O(2,q)$
have  two and four connected components, respectively,
and that 
$\roman O(2,q)^0=\roman {SO}(2,q)^0$.
Notice that a matrix $X$ in $\fra g$
with $X^3 = 0$
satisfies $(\roman{ad}(X))^4 =0$
and is therefore a nilpotent element of
$\fra g$.

\proclaim{Theorem 3.6.2}
For $q \geq 3$,
{\rm (3.6.1)}
constitutes a complete list
of the pseudoholomorphic nilpotent $\roman {SO}(2,q)^0$-orbits
in $\fra{so}(2,q)$.
Furthermore:
\newline\noindent
{\rm (i)} The orbits $\Cal O_{1,0}$ and $\Cal O_{0,1}$
are the two connected components of
the single 
$\roman {SO}(2,q)$-orbit 
(and hence $\roman O(2,q)$-orbit)
in 
$\fra{so}(2,q)$
which consists of
square zero matrices
in 
$\fra{so}(2,q)$.
\newline\noindent
{\rm (ii)} The
three $\roman {SO}(2,q)^0$-orbits 
$\Cal O_{2,0}$, $\Cal O_{1,1}$, $\Cal O_{0,2}$
are the connected components of
the smooth manifold
of matrices $X$
in 
$\fra{so}(2,q)$ satisfying $X^3 = 0$
in such a way that,
as $\roman {SO}(2,q)$-orbits
as well as
as $\roman O(2,q)$-orbits,
$\Cal O_{2,0}$  and $\Cal O_{0,2}$
are the two connected components of
one 
orbit and
$\Cal O_{1,1}$
is the other 
orbit.
\newline\noindent
{\rm (iii)} 
For 
$q \geq 2$, the orbits
$\Cal O_{1,0}$
and
$\Cal O_{0,1}$
are the two connected components
of the smooth manifold 
of rank two  
nilpotent
matrices 
$\left [\matrix
A&B\\
B^{\roman t}&D 
\endmatrix
\right]$
in $\fra{so}(2,q)$ 
which,
when $B$ is written as $B=\left[\matrix x_1 & x_2 & \cdots & x_q\\
y_1 & y_2 & \cdots & y_q\endmatrix \right]$,
satisfy the equations
$$
\align
x_1^2 - y_1^2 + x_2^2 - y_2^2 + \dots + x_q^2 - y_q^2 &= 0
\\
x_1 y_1 + x_2 y_2 + \dots + x_q y_q &= 0 .
\endalign
$$
\endproclaim

\demo{Proof}
Statements (i) and (ii) are established by a reasoning 
involving the decomposition of
the type of a nilpotent orbit into
a sum of indecomposable 
$\roman O(2,q)$-types
\cite \burcush,
similar to that given in the proof of Theorem 3.5.3,
together with a closer look at
the 
decomposition 
of the resulting $\roman O(2,q)$-orbits
into
$\roman {SO}(2,q)$- and 
$\roman {SO}(2,q)^0$-orbits;
we leave the details to the reader.
Statement (iii)
is a consequence of Theorem 3.6.3 below. 
It may as well
be deduced from 
the discussion in Section 6 below. \qed
\enddemo

The next result 
makes explicit 
the statement of Theorem 3.3.11 above
for $\fra g = \fra{so}(2,q)$.

\proclaim{Theorem 3.6.3}
Under the 
projection from $\fra g = \fra k \oplus \fra p$ to
$\fra p$, followed by the canonical isomorphism from
$\fra p$ to $\fra p^+$,
the closure
$\overline{\Cal O_{2,0}}$ 
of the principal holomorphic nilpotent orbit
$\Cal O_{2,0}$ 
is identified with $\fra p^+\cong \Bobb C^q$, and
the closure $\overline{\Cal O_{1,0}}$ 
of 
$\Cal O_{1,0}$ 
amounts to the smooth complex affine 
quadric
in $\Bobb C^q$ 
given by the equation
$w_1^2 + w_2^2 + \dots +w_q^2 = 0$. 
Furthermore, as
a complex analytic space,
$\overline{\Cal O_{1,0}}$ is normal, and so is, of course,
$\overline{\Cal O_{2,0}}$ since it is an affine space.
\endproclaim

\demo{Proof} The first statement is established by inspection.
See also Remark 6.3 below.
Being an affine quadric, the normality of 
$\overline{\Cal O_{1,0}}$
is immediate. \qed
\enddemo

\smallskip\noindent
{\smc Remark 3.6.4.}
The proof of Theorem 3.6.2 is independent of those of Theorems 3.3.3 and
3.3.11 and thus yields proofs of the statements of 
these Theorems for the present special case as well.

\smallskip
Since $\fra{so}(2,3)$
is isomorphic to
$\fra{sp}(2,\Bobb R)$, for $q \geq 3$, the standard inclusion
of $\fra{so}(2,3)$
into 
$\fra{so}(2,q)$
amounts to an injection
$$
\fra{sp}(2,\Bobb R) @>>> \fra{so}(2,q).
\tag3.6.5
$$
This is an $H_1$-embedding of the kind asserted in Proposition 3.3.2(4);
it identifies the $H$-element of the source with that of the target.

\smallskip\noindent
{\smc 3.7. The complex analytic stratified K\"ahler structure.}
We will now show that
the closure of any 
holomorphic 
nilpotent orbit 
$\Cal O$
in a
semisimple Lie algebra 
of hermitian type
carries the structure of a 
complex analytic stratified K\"ahler space
which, by virtue of Theorems 3.5.5, 3.6.3, 5.3.3, 8.4.1, 8.4.2,
is actually normal.
Thus, let $(\fra g,z)$ be a 
semisimple Lie algebra 
of hermitian type,
and let $\Cal O$ be a holomorphic nilpotent orbit.
Its closure $N=\overline{\Cal O}$
in $\fra g$ is a union of finitely many 
holomorphic
nilpotent orbits and, 
as we have seen above,
the 
restriction to $N$ of the 
projection from $\fra g$ to
the constituent
$\fra p$
of the Cartan decomposition of $\fra g$
is
an injective map, 
even though the projection
from $\fra g$ to $\fra p$ is not injective;
this injection
turns $N$ into a complex analytic space
in such a way that each $G$-orbit in $N$ is a
smooth complex manifold.
Independently of the structure results obtained 
above,
the decomposition of $N$ into $G$-orbits 
may be shown to be a stratification,
since $\fra g$ is semisimple;
see \cite\lermonsj\ (Theorem 2.23).
Furthermore, the
{\it Tarski-Seidenberg\/} Theorem entails that,
as a subset of $\fra g$, $N$ is real semi-algebraic.
\smallskip
Upon identifying $\fra g$ with its dual $\fra g^*$
by means of 
an appropriate positive multiple of the
Killing form
 we view $N$ as the closure
of a 
(nilpotent)  
{\it coadjoint\/} orbit
(i.~e. orbit in $\fra g^*$).
At this stage any positive multiple 
of the
Killing form
will do; in Section 5 below
a more 
specific choice will be made; see also (3.2.2) above.
Restricting the 
algebra 
$C^{\infty}(\fra g^*)$
of
ordinary smooth 
functions 
on $\fra g^*$
to
$N$ yields 
a smooth structure $C^{\infty}(N)$ on $N$; a
smooth structure
of this kind is
usually referred to as an algebra
of {\it Whitney smooth functions\/}.
The existence of smooth partitions of unity on $\fra g^*$
entails that
this smooth structure
is fine.
Furthermore,
since $N$ is the closure
of a 
nilpotent orbit,
the Lie-Poisson structure  $\{\cdot,\cdot\}$
on 
$C^{\infty}(\fra g^*)$
passes
to a Poisson bracket $\{\cdot,\cdot\}$
on $C^{\infty}(N)$ which,
on each stratum of $N$
(i.~e. coadjoint orbit contained in $N$), 
restricts to the ordinary
symplectic Poisson algebra.
Thus
$(C^{\infty}(N),\{\cdot,\cdot\})$
{\it turns $N$ into a stratified symplectic space\/}.
\smallskip
To spell out the requisite 
complex analytic stratified
K\"ahler 
polarization,
we note that, as a module over the algebra
$C^{\infty}(\fra g^*,\Bobb C)$,
the module $\Omega^1(\fra g^*,\Bobb C)$ of 1-forms may be written
in the form $C^{\infty}(\fra g^*,\Bobb C) \otimes_{\Bobb C}\fra g^{\Bobb C}$.
Hence
the decomposition
$\fra g^{\Bobb C} = \fra k^{\Bobb C} \oplus \fra p^+ \oplus \fra p^-$
entails a direct sum decomposition
$$
\Omega^1(\fra g^*,\Bobb C)
=
C^{\infty}(\fra g^*,\Bobb C) \otimes_{\Bobb C}\fra k^{\Bobb C} 
\oplus 
C^{\infty}(\fra g^*,\Bobb C) \otimes_{\Bobb C}\fra p^+ 
\oplus 
C^{\infty}(\fra g^*,\Bobb C) \otimes_{\Bobb C}\fra p^-
$$
as a module over the algebra
$C^{\infty}(\fra g^*,\Bobb C)$.
The Lie-Rinehart structure 
on
$\Omega^1(\fra g^*,\Bobb C)$
coming from the Poisson structure on
$C^{\infty}(\fra g^*,\Bobb C)$
(explained in Section 1 above)
is actually a \lq\lq crossed product
structure\rq\rq, cf. \cite\poiscoho.
Since $\fra k^{\Bobb C}$ is a (complex) Lie algebra,
the constituent
$C^{\infty}(\fra g^*,\Bobb C) \otimes_{\Bobb C}\fra k^{\Bobb C}$
is closed under the Lie bracket
in 
$\Omega^1(\fra g^*,\Bobb C)$
and hence inherits a Lie-Rinehart structure,
a crossed product structure as well.
We have already observed that the constituents $\fra p^+$ and
$\fra p^-$
are known to be abelian subalgebras of
$\fra g^{\Bobb C}$,
cf. \cite\helgaboo\ (p.~313),
\cite\satakboo\ (p.~55). 
Hence the
constituents
$C^{\infty}(\fra g^*,\Bobb C) \otimes_{\Bobb C}\fra p^+$ 
and
$C^{\infty}(\fra g^*,\Bobb C) \otimes_{\Bobb C}\fra p^-$
are
closed under the Lie bracket
in 
$\Omega^1(\fra g^*,\Bobb C)$
and inherit Lie-Rinehart structures;
these are crossed product structures again.
Even though 
$\fra p^+$ 
and
$\fra p^-$ 
are abelian,
the brackets on
$C^{\infty}(\fra g^*,\Bobb C) \otimes_{\Bobb C}\fra p^+$ 
and
$C^{\infty}(\fra g^*,\Bobb C) \otimes_{\Bobb C}\fra p^-$
are non-trivial since they are induced from brackets
of the kind $[u,v]$ where $u \in \fra p^+$ 
and $v \in \fra g^{\Bobb C}$.
\smallskip
Write $A(N)=C^{\infty}(N,\Bobb C)$, and let
$I$ be the ideal of smooth functions on $\fra g^*$
which vanish on $N$.
By the general theory of K\"ahler differentials,
the projection map
from
$C^{\infty}(\fra g^*,\Bobb C)$ to $A(N)$
gives rise to the exact sequence
$$
I \big / I^2
@>{\delta}>>
A(N)\otimes_{C^{\infty}(\fra g^*,\Bobb C)}
\Omega^1(\fra g^*,\Bobb C)
@>>>
\Omega^1(N,\Bobb C)
@>>>
0
$$
of $A(N)$-modules
where $\delta$ is induced by the association
$a \mapsto 1 \otimes da\ (a \in I)$.
Furthermore, 
as an
$A(N)$-module,
$A(N)\otimes_{C^{\infty}(\fra g^*,\Bobb C)}
\Omega^1(\fra g^*,\Bobb C)$
may plainly be written as
$$
A(N)\otimes_{C^{\infty}(\fra g^*,\Bobb C)}
\Omega^1(\fra g^*,\Bobb C)
\cong
A(N)
\otimes_{\Bobb C}\fra k^{\Bobb C} 
\oplus 
A(N) \otimes_{\Bobb C}\fra p^+ 
\oplus 
A(N) \otimes_{\Bobb C}\fra p^- .
$$
Let
$P\subseteq \Omega^1(N,\Bobb C)$
be the image of
$A(N) \otimes_{\Bobb C}\fra p^+$
in
$\Omega^1(N,\Bobb C)$.
This is a $(\Bobb C,A(N))$-Lie subalgebra of 
$\Omega^1(N,\Bobb C)_{\{\cdot,\cdot\}}$
which,
as an $A(N)$-module,
is generated by the differentials $dw$ where $w \in \fra p^+$;
here the elements of $\fra g^{\Bobb C}$ are viewed as
functions on $(\fra g^{\Bobb C})^*$.

\proclaim{Theorem 3.7.1}
The
$(\Bobb C,A(N))$-Lie subalgebra $P$ of 
$\Omega^1(N,\Bobb C)_{\{\cdot,\cdot\}}$
is a 
complex analytic stratified K\"ahler polarization for
the (fine) stratified symplectic space
$(N,C^{\infty}(N),\{\cdot,\cdot\})$;
the underlying complex analytic structure
is that induced from the embedding of $N$ into
$\fra p (\cong \fra p^+)$.
\endproclaim

\demo{Proof}
Since each stratum of $N$ inherits a smooth complex  structure
from the complex analytic structure
of $N$,
the symplectic stratification is a refinement of
the complex analytic one.
Let $w_1,\dots,w_{\ell}$ be a (complex) basis of 
$\fra p^+$ and view the $w_j$'s as holomorphic coordinate functions
on $(\fra p^+)^*$; restricting them to $\Cal O$
yields holomorphic functions on $\Cal O$,
and every holomorphic function in $\Cal O$ may be written
as a holomorphic function in these coordinates.
This implies that, 
given
two holomorphic functions $f$ and $h$ in $\Cal O$,
the Poisson bracket $\{f,h\}$ vanishes,
since the Lie bracket on $\fra p^+$ (which yields the brackets of the kind
$\{w_j,w_k\}$) is zero.
Consequently
the sheaf
of germs of holomorphic functions
on $N$ is contained in the sheaf of germs of
polarized functions.
In view of Theorem 2.5,
$P$ is a stratified K\"ahler polarization which is
compatible with the complex analytic
structure; this may also be seen directly.
\smallskip
To see that the 
sheaf of germs of holomorphic functions coincides with that
of polarized ones, we argue as follows:
The top stratum
$\Cal O_r$ (i.~e. the principal holomorphic nilpotent orbit)
is Zariski open in $\fra p^+$ 
and has 
$\fra p^+$ as its closure.
The complement 
$\overline{\Cal O_{r-1}}$
of $\Cal O_r$
in $\fra p^+$ 
is a thin set in $\fra p^+$.
Hence a continuous function 
in $\fra p^+$
which is holomorphic in
$\Cal O_r$ is holomorphic
in $\fra p^+$
\cite\grauremm.
This implies
the assertion, since the complex analytic structure
on the closure of any smaller holomorphic 
nilpotent  orbit
is induced from the embedding thereof into
$\fra p^+$. \qed
\enddemo

\smallskip\noindent
{\smc Remark 3.7.2.}
While $\fra p$ 
carries a 
symplectic structure as well which
combines with its complex structure  $J_z$
to a K\"ahler structure,
the symplectic structures 
on the strata of $N$
{\it cannot\/}  be induced from
a symplectic structure on $\fra p$, though.

\smallskip\noindent
{\smc Remark 3.7.3}
The image of
$A(N) \otimes_{\Bobb C}\fra p^-$
in
$\Omega^1(N,\Bobb C)$
is a stratified complex polarization
$\overline P \subseteq \Omega^1(N,\Bobb C)$
as well
which 
is generated by the differentials $dx$ where $x \in \fra p^-$;
it is 
just the complex 
conjugate of $P$.
The 
$C^{\infty}(N,\Bobb C)$-submodule
of $\Omega^1(N,\Bobb C)$
generated by
$P$ and $\overline P$
decomposes into a
direct sum
$P \oplus\overline P$.
However, 
$\Omega^1(N,\Bobb C)$
does  {\it not\/}
exhaust this direct sum. In fact,
the inclusion of
$P\oplus\overline P$
into
$\Omega^1(N^{\roman{red}},\Bobb C)$
gives rise to
an exact 
sequence
$$
0
@>>>
P \oplus\overline P
@>>>
\Omega^1(N,\Bobb C)
@>>>
C
@>>> 0
\tag3.7.4
$$
of $C^{\infty}(N,\Bobb C)$-modules,
the cokernel
$C$ being 
generated by (classes of) the 
differentials $dw$
where $w \in \fra k^{\Bobb C}$.
The precise information about holomorphic nilpotent orbits
spelled out in Theorems 3.5.4, 3.5.5, and 3.6.2 above
yields
at once explicit examples where $C$ is non-zero.
Thus a complex polarization can {\it no longer
be thought of as being given by the $(0,1)$-vectors
of a complex structure\/}.
On any stratum,
the 
constituent $C$ is zero, though, and the
module of 
differentials
does of course 
decompose into the 
direct sum of the 
holomorphic and antiholomorphic 
polarization:
Restricted to a stratum
$S$,
the exact sequence 
(3.7.4) comes down to an isomorphism
$P_S \oplus\overline P_S
@>>>
\Omega^1(S,\Bobb C)$.

\smallskip\noindent
{\smc  3.7.5. Dependence on the choice of $H$-element}.
Given a  nilpotent orbit
$\Cal O$, we have observed above
(Addendum to Theorem 3.3.3)
that the property of being holomorphic
does not depend on the choice of $H$-element.
The induced complex structure on the closure
$\overline {\Cal O}$  depends on 
that choice, though.
Letting the $H$-element
$z$ vary, we obtain a moduli space of
complex structures compatible with the other structure,
i.~e. a moduli space of complex analytic stratified K\"ahler
structures on 
$\overline {\Cal O}$,
the stratified symplectic structure being fixed.
This moduli space amounts to the homogeneous space
$G/K$.

\smallskip\noindent
{\smc  3.7.6. The  holomorphic semisimple orbits.}
In the same fashion, 
for any semisimple holomorphic orbit $\Cal O$,
a choice of $H$-element $z$ induces a 
smooth $K$-invariant complex structure
on $\Cal O\cong G/K$ turning the latter into a K\"ahler manifold.
This is {\it not\/} the standard hermitian symmetric space structure
on 
$G/K$, though, since
only $K$ acts by isometries,
and the curvature tensor is not
parallel. 
A special case has been illustrated in (3.2.2) above.

\medskip\noindent{\bf 4. Reduction and stratified K\"ahler spaces}
\smallskip\noindent
In this section we will show that reduction is another source
for stratified K\"ahler spaces. For later reference, it will be convenient
to extend the usual reduction scheme in symplectic geometry.
\smallskip
Let 
$(N,C^{\infty}(N))$ be a 
stratified space, with smooth structure
$C^{\infty}(N)$,
and let  $G$ be a Lie group, with Lie algebra $\fra g$.
We will say that an action 
of $G$ on $N$
is {\it smooth\/} provided it preserves strata and carries 
smooth functions on $N$
to smooth functions on $N$,
that is, 
(i) given 
$h \in C^{\infty}(N)$ 
and $y \in G$, the function $f \circ y$ 
(where we do not distinguish in notation between $y \in G$
and the homeomorphism it induces on $N$)
is in
$C^{\infty}(N)$ and (ii),
given 
$h \in C^{\infty}(N)$ 
and $Y \in \fra g$, the function $Y_Nh$ 
is in $C^{\infty}(N)$; 
here 
$Y_N$ refers to the derivation 
which, on each stratum, is given by the vector field
$Y_N$ induced from $Y$ on that stratum.
As a smooth  function,
$Y_Nh$ initially exists only stratumwise,
whence the requirement (ii).
We refer to a derivation of the kind
$Y_N$ as the {\it stratified vector field\/}
on $N$ induced from $Y$.
\smallskip
Let $\{\ ,\ \}$ be
a Poisson structure on
$C^{\infty}(N)$
which, on each stratum of $N$, restricts to a
smooth 
(not necessarily symplectic) Poisson structure.
We will refer to
$(N,C^{\infty}(N),\{\ ,\ \})$
as a {\it stratified Poisson space\/}.
We will say that the $G$-action on $N$
is {\it hamiltonian\/},
with  {\it stratified Poisson momentum mapping\/}
$\mu \colon N \to \fra g^*$,
provided (i) $\mu$ is a smooth $G$-equivariant map,
and (ii) for every $Y \in \fra g$ and every smooth function
$f$ on $N$,
$$
\{f, \mu^Y\} = Y_N f,
$$
where
$\mu^Y\colon N \to \Bobb R$
is the composite of $\mu$ with $Y$, viewed as a linear map on $\fra g^*$,
and where $Y_N$
refers to the
stratified vector field
on $N$
induced by the $G$-action.
Here \lq\lq smooth\rq\rq\  means
that, for every ordinary smooth function $h$ on $\fra g^*$,
the composite
$h \circ \mu$ is in $C^{\infty}(N)$.
The requirement
(ii) is equivalent
to
the adjoint $\delta \colon \fra g \to C^{\infty}(N)$,
which is given by $\delta(Y) = \mu^Y$ ($Y\in \fra g$),
having the property that
(ii$'$) 
for every $Y\in \fra g$, the hamiltonian vector field $\{\cdot,\mu^Y\}$
of $\delta(Y)$ coincide with the vector field $Y_N$;
the mapping $\delta$ is then occasionally referred to as {\it comomentum\/}.
Furthermore,
when (ii) holds,
the $G$-equivariance of $\mu$ implies that
(i$'$) $\delta[Y,Z] = \{\delta(Y),\delta(Z)\}$,
for every $Y,Z\in \fra g$ and, for connected $G$,
this is actually equivalent to the $G$-invariance.
Notice that (i$'$) says that $\mu$ is a Poisson map where $\fra g^*$
is endowed with its Lie-Poisson structure.
When $N$ is a smooth symplectic manifold
or,
more generally, 
a smooth Poisson manifold,
this notion of Poisson
momentum mapping amounts to the ordinary one in symplectic geometry
or, more generally, 
to that of Poisson momentum mapping
\cite\marsrati.
When
$(N,C^{\infty}(N),\{\ ,\ \})$ is a stratified symplectic space
acted upon by a Lie group $G$
in a hamiltonian fashion with
stratified Poisson momentum mapping
$\mu \colon N \to \fra g^*$,
we will refer to $\mu$
as a {\it stratified symplectic space momentum mapping\/}.
Examples which are not ordinary momentum mappings
will be given in Section 8 below.
\smallskip
Let $(N,C^{\infty}(N),\{\ ,\ \})$ be a 
stratified Poisson space,
and let $G$ be a Lie group which acts
smoothly on $N$ in a hamiltonian fashion
with  stratified Poisson momentum mapping
$\mu \colon N \to \fra g^*$.
Generalizing the standard construction, define the
{\it reduced space\/} $N^{\roman{red}}$ by
$N^{\roman{red}} = \mu^{-1}(0) \big / G$.
The construction of the
{\smc Arms-Cushman-Gotay}-algebra \cite\armcusgo\ 
extends as well: Define
$$
C^{\infty}(N^{\roman{red}})
\quad\text{as}\quad
\left(C^{\infty}(N)^G\right) \big / I^G,
$$
the algebra $C^{\infty}(N)^G$ of smooth $G$-invariant functions
on $N$, modulo the ideal $I^G$
of functions in $C^{\infty}(N)^G$ that vanish on the zero locus
$\mu^{-1}(0)$, so that
$C^{\infty}(N^{\roman{red}})$
yields 
a {\it smooth structure\/}
$C^{\infty}(N^{\roman{red}})$, i.~e.
algebra of continuous functions
on $N^{\roman{red}}$ in an obvious fashion.

\proclaim{Proposition 4.1}
The Poisson structure on 
$C^{\infty}(N)$
induces a Poisson structure
$\{\cdot,\cdot\}^{\roman{red}}$ on
$C^{\infty}(N^{\roman{red}})$.
\endproclaim

\demo{Proof} 
Let $f$ be a $G$-invariant function in $C^{\infty}(N)$.
Noether's theorem is still true and says that
$\mu$ is constant along the trajectories
of $X_f = \{\cdot,f\}$; 
indeed, any trajectory lies in a stratum.
Hence, given a function $h$ in
$C^{\infty}(N)$
which vanishes on $\mu^{-1}(0)$, 
for any $q \in\mu^{-1}(0)$, since $\mu(\roman{exp}(tX_f)q) = \mu(q)$,
we have
$$
\{h,f\}(q) = X_f(h)|_q = \frac d{dt}(h(\roman{exp}(tX_f)q))|_{t=0} = 0,
$$
that is,
$\{h,f\}$ vanishes on $\mu^{-1}(0)$ as well.
Consequently $I^G$ is a Poisson ideal in
$C^{\infty}(N)^G$. \qed
\enddemo

We now suppose that $G$ is compact.
Given a symplectic manifold $N$ acted upon by 
$G$ in a hamiltonian fashion,
each connected component of
the reduced space $N^{\roman{red}}$ is a stratified
symplectic space \cite\sjamlerm;
to simplify the exposition somewhat we shall not distinguish
in notation between 
the reduced space and any one of its connected components.
The existence of $G$-invariant partitions of unity on $N$
with respect to arbitrary locally finite coverings
by  $G$-invariant open subsets entails that
$C^{\infty}(N^{\roman{red}})$ is fine.

\proclaim{Proposition 4.2}
Given a positive K\"ahler manifold $N$
with a holomorphic $G^{\Bobb C}$-action 
whose restriction to (a compact real form)
$G$
preserves the K\"ahler structure
and
is  hamiltonian, 
the K\"ahler polarization $F$
induces a 
(positive) normal K\"ahler
structure
on the reduced space $N^{\roman{red}}$.
\endproclaim

\demo{Proof} This follows readily from the results
in \cite\heinloos\ 
and
\cite\sjamatwo,
combined with Theorem 2.5 above.
Indeed, the zero locus $\mu^{-1}(0)$
is a {\it Kempf-Ness\/} set
of $N$,
and the inclusion
of $\mu^{-1}(0)$ into 
$N$  
induces a homeomorphism
$$
N^{\roman{red}}
= \mu^{-1}(0) \big / G
@>>>
N \big /\big / G^{\Bobb C} 
$$ 
from the reduced space onto
the  G(eometric) I(nvariant) T(heory) quotient
$N \big /\big / G^{\Bobb C}$.  
The latter inherits a complex analytic structure
which, locally, amounts to that of a Stein
quotient in such a way that the symplectic stratification
of $N^{\roman{red}}$ is a refinement of the complex analytic one
and that each stratum,
with its induced complex analytic structure,
inherits a K\"ahler metric;
see (2.3) of \cite\sjamatwo\ 
(under the additional technical hypothesis
that the momentum mapping be \lq\lq admissible\rq\rq) 
and \cite\heinloos.
For every pair $f$ and $h$ of holomorphic functions in $N$,
the Poisson bracket $\{f,h\}$ manifestly vanishes, and this property
descends to $N^{\roman{red}}$.
By Theorem 2.5, the $C^{\infty}(N^{\roman{red}},\Bobb C)$-submodule
$P$ of 
$\Omega^1(N^{\roman{red}},\Bobb C)_{\{\cdot,\cdot\}}$
described there (where $X= N^{\roman{red}}$)
is a complex analytic stratified K\"ahler polarization.
Normality is a consequence of the fact that, locally,
the GIT-quotient may be written as a quotient of a smooth Stein space
and, a smooth complex manifold being necessarily normal,
passing to the quotient preserves normality. 
Since, complex analytically, the space
$N^{\roman{red}}$ is normal,
we might as well  first observe that
$P$ is a stratified K\"ahler polarization on
$N^{\roman{red}}$
such that the complex analytic structure is compatible with it
and then deduce,
using Proposition 2.4, that $P$ is actually a normal K\"ahler structure.
The positivity of $P$ is obvious. \qed
\enddemo

\smallskip\noindent
{\smc Example 4.3.}
Let 
$\Sigma$
be a closed Riemann surface, 
let $G$
be a compact Lie group, 
and let $\xi \colon P \to \Sigma$ be a principal $G$-bundle.
Pick a Riemannian metric on $\Sigma$
and an invariant positive definite
symmetric bilinear form on the Lie algebra $\fra g$ of $G$, and
let $N$ be the moduli space of central Yang-Mills connections
on $\xi$ determined by these data \cite\atibottw.
To recall its well known representation theory interpretation, 
write $Z$ for the center of $G$, let $\fra z = \roman{Lie}(Z)$,
and pick $X \in \fra z$.
Let $\pi$ be the fundamental group of $\Sigma$,
let $0 \to \Bobb Z \to \Gamma \to \pi \to 1$
be its universal central extension,
let $\roman{Hom}_X(\Gamma, G)$ 
be the space of homomorphisms
from
$\Gamma$ to $G$
which send a preferred generator of the kernel $\Bobb Z$ to 
$\roman{exp}X \in G$, 
and consider the space 
$\roman{Rep}_X(\Gamma, G) =\roman{Hom}_X(\Gamma, G)\big / G$
of $G$-orbits (under conjugation);
depending on the choice of $X$, the space
$\roman{Hom}_X(\Gamma, G)$ may be empty.
The topological type of $\xi$
determines an element $X \in \fra z$ such that
taking holonomies identifies the space $N$ with
a connected component of $\roman{Rep}_X(\Gamma, G)$
\cite\atibottw,
\cite{\locpois--\modustwo}.
We will refer to
any of the connected components
of $\roman{Rep}_X(\Gamma, G)$
as a {\it space of twisted representations\/}
of $\pi$ in $G$. As final ingredient,
we pick a complex structure on $\Sigma$.

\proclaim{Theorem 4.3.1} 
The data 
induce on $N$ a 
normal (complex analytic positive) K\"ahler structure
turning $N$ into an (in general) exotic projective variety.
\endproclaim

Explicit constructions of the stratified symplectic  structure 
(which does not depend on the complex structure on $\Sigma$)
may be found in
\cite{\locpois--\modustwo}.
For $G= \roman U(n)$, 
the requisite complex analytic structure on $N$
was constructed by 
Narasimhan and Seshadri
\cite\narasesh\ 
by means of geometric invariant theory;
the space $N$ may then be identified with the moduli space
of  semistable holomorphic rank $n$ vector bundles
on $\Sigma$
of a certain degree $k$
which corresponds to the 
topological type of $\xi$ or, equivalently,
to the twisting 
of the representations of $\pi$;
for $k=0$ the representations of $\pi$
in $G$ are ordinary ones, i.~e. are untwisted.
For general $G$, a construction of the complex analytic structure
may be found e.~g. in \cite\kumnarra.
\smallskip
It is not obvious that the two structures
combine to a complex analytic K\"ahler structure, though, except in the 
(complex analytically or algebraically) non-singular case
where the space $N$ is known to be a K\"ahler manifold
\cite\narasesh.
We now 
explain the requisite analytical details.
\smallskip
Let  $A$ be an arbitrary connection
(not necessarily Yang-Mills)
on $\xi$, and let $d_A$ be its operator of covariant derivative.
It is well known that $A$ 
determines (and is determined by) a  complex structure
on 
the complexification $\xi^{\Bobb C} \colon P^{\Bobb C} \to \Sigma$ 
of $\xi$
which 
turns $\xi^{\Bobb C}$ 
into a holomorphic principal 
$G^{\Bobb C}$-bundle;
we write the latter as $\xi_A^{\Bobb C}$.
The 
chosen invariant positive definite
symmetric bilinear form on $\fra g$
induces a hermitian structure on $\roman{ad}(\xi^{\Bobb C})$
which, in turn, determines a Hodge decomposition
$$
\Omega^{0,0}(\Sigma,\roman{ad}(\xi^{\Bobb C}))
= \Cal H_A^{0,0} \oplus L^{0,0},\quad
\Omega^{0,1}(\Sigma,\roman{ad}(\xi^{\Bobb C}))
= \Cal H_A^{0,1} \oplus d_A''(L^{0,0}),
$$
the spaces of harmonic forms being denoted by
$\Cal H_A^{0,*}$ and the Cauchy-Riemann operator by
$d_A''$.
A standard construction, cf. Lemma 2.1 (ii) in \cite\narasesh\ 
for $G=\roman U(n)$, yields a complete
analytic family of holomorphic
principal $G^{\Bobb C}$-bundles
$\{\xi_y\}_{y \in Y_A}$
parametrized by a based ball
$(Y_A,o)$ in the complex vector space $\Cal H_A^{0,1}$
so that $\xi_o$ \lq\lq coincides\rq\rq\  with 
$\xi_A^{\Bobb C}$,
that is,
a holomorphic principal  $G^{\Bobb C}$-bundle 
over $Y_A \times \Sigma$
whose Kodaira-Spencer map
$\roman T_o Y_A @>>> \Cal H_A^{0,1}$
is an isomorphism,
together with an isomorphism
of holomorphic principal bundles between
$\xi_o$ and
$\xi_A^{\Bobb C}$.
\smallskip
Suppose that $A$ is a central Yang-Mills connection, 
and consider the real differential graded Lie algebra
$(\Omega^*(\Sigma,\roman{ad}(\xi)),d_A)$.
It controls the infinitesimal 
central Yang-Mills connections variations of $(\xi,A)$.
The star operator 
$* \colon \Omega^*(\Sigma,\roman{ad}(\xi)) 
\to \Omega^{2-*}(\Sigma,\roman{ad}(\xi))$
determined by the data
endows
$\Omega^1(\Sigma,\roman{ad}(\xi))$
with a complex structure, and associating
$\eta + i*\eta$ to $\eta$ 
identifies
$\Omega^1(\Sigma,\roman{ad}(\xi))$
with 
$\Omega^{0,1}(\Sigma,\roman{ad}(\xi^{\Bobb C})) $
as a complex vector space.
Under this isomorphism,
the space $\Cal H_A^1$
of real harmonic $1$-forms
gets identified with $\Cal H_A^{0,1}$.
The space of connections $\Cal A(\xi)$
being endowed with the standard symplectic structure $\omega$ given by 
$\omega(\alpha,\beta) = \int_{\Sigma}\langle \alpha \wedge \beta \rangle$
where $\alpha,\beta \in \Omega^1(\Sigma,\roman{ad}(\xi))$ and
where 
$\langle\, , \rangle$ is the pairing induced from the given
invariant symmetric bilinear form on $\fra g$,
the assignment to a connection of its curvature
is a momentum mapping for the action of the group of gauge transformations.
Consider the symplectic slice
$\Cal M_A\subseteq \Cal A(\xi)$ given in (2.16) of \cite\singula.
The symplectic structure $\omega$
induces a symplectic structure $\omega_A$ on
$\Cal M_A$.
Write $Z_A$ for the stabilizer of $A$
in the group of gauge transformations of $\xi$,
and let $\fra z_A = \roman{Lie}(Z_A)$.
The momentum mapping for the action of the group of gauge transformations
restricts to a momentum mapping
$\mu_A \colon \Cal M_A \to \Cal H_A^2$ 
where
$\Cal H_A^2 \cong \roman H^2(\Sigma,\roman{ad}(\xi))$
is identified with the dual of
$\fra z_A$
in the standard fashion.
After having $\Cal M_A$ suitably cut to size (if need be), 
when we associate to
a connection $A+\eta$ in $\Cal M_A$
its induced Cauchy-Riemann operator $d''_{A+\eta}$,
we obtain a map from
$\Cal M_A$
to
the GIT  quotient
$Y_A\big / \big / Z^{\Bobb C}_A$ 
which induces a homeomorphism from the reduced space
$\mu_A^{-1}(F_{\xi})\big /Z_A$
onto $Y_A\big / \big / Z^{\Bobb C}_A$;
here 
$F_{\xi}$ denotes the corresponding constant central curvature.
On the other hand,
an explicit map from
$\mu_A^{-1}(F_{\xi})\big /Z_A$
to
$N$,
realized as a connected
component of a space of the kind
$\roman{Rep}_X(\Gamma, G)$,
is given by the assignment to a point
$A+\eta$ of its holonomies with respect to a system of curves
representing the generators of the fundamental group.
See \cite\smooth\ for details.
Letting $A$ vary
(among central Yang-Mills connections), we obtain these structures on suitable
\lq\lq coordinate patches\rq\rq\ of the space $N$.
The Poisson structure comes from
the reduced spaces of the kind $\mu_A^{-1}(F_{\xi})\big /Z_A$,
while the complex analytic one
is inherited from the
complex GIT  quotient's
$Y_A\big / \big / Z^{\Bobb C}_A$. 
The global stratified symplectic Poisson structure 
on $\roman{Rep}_X(\Gamma, G)$ is that constructed in
\cite\modus.
Thus we obtain a normal K\"ahler structure
on $N$, the normality being an immediate consequence
of the local 
complex GIT  quotient constructions.
More details will be given at another occasion.
\smallskip
Apart from the Hodge decomposition 
of $\Omega(\Sigma,\roman{ad}(\xi^{\Bobb C}))$
on which the existence of the complete families of 
holomorphic principal $G^{\Bobb C}$-bundles relies
and the Banach space techniques which yield the slice $\Cal M_A$,
{\sl this construction of the moduli space
$N$ as a normal K\"ahler space is finite dimensional and
entirely analytic, that is, it avoids the usual
detour via geometric invariant theory.\/}
Normality is a consequence of the fact that, locally, the space is an affine
geometric invariant theory quotient.
The {\sl complex analytic stratified K\"ahler structure\/}
isolated here is somewhat finer
than the complex analytic structure alone, though, and
cannot be obtained from geometric invariant theory.
{\sl This stratified K\"ahler structure gives a precise description of the 
K\"ahler singularity behaviour.\/}
The issue of singularities was raised in \cite\atibottw. 
When $G=\roman{SU}(2)$,
as already pointed out in the introduction, the space
$N$ is an exotic $\Bobb C \roman P^3$,
and the complement of the top stratum is
a Kummer surface which, in turn,
is the singular locus, 
in the sense of 
complex analytic stratified K\"ahler spaces,
of this exotic $\Bobb C \roman P^3$.
For various special cases,
the local structure of the space $N$ near any of its points
has been examined in \cite\locpois;
see also \cite{\srni--\oberwork}.
\smallskip
More generally,
certain moduli spaces including  
moduli spaces
of
semistable holomorphic vector bundles
with parabolic structure
\cite\mehtsesh\ 
inherit normal K\"ahler structures,
the requisite stratified
symplectic structure
being that constructed in
\cite\guhujewe.
We will explain the details elsewhere.

\medskip\noindent
{\bf 5. Associated representations and singular reduction}
\smallskip\noindent
In this section we will describe the holomorphic nilpotent orbits
for the classical cases 
as reduced spaces for suitable momentum mappings.
\smallskip\noindent
{\smc 5.1. Invariants and momentum mapping.}
Consider a real finite-dimensional symplectic vector space $W$, with 
symplectic form $\omega$. 
The  algebra $\Bobb R[W^*]$
of real polynomial functions 
on $W$---after a choice of basis
$b_1,\dots,b_n$: the algebra of polynomials in the
coordinate functions $x_1,\dots,x_n$---inherits 
a Poisson structure $\{\cdot,\cdot\}$
from the symplectic structure
which is given by $\{x_j,x_k\} = \omega(b_j,b_k)$, $1 \leq j,k\leq n$;
furthermore, the homogeneous quadratic part
$\Bobb R_2[W^*]$ 
of $\Bobb R[W^*]$
is closed under the (Poisson) bracket and, in fact,
as a Lie algebra, isomorphic to
$\fra{sp}(W,\omega)$.
An explicit isomorphism $\delta \colon
\fra{sp}(W,\omega) \to \Bobb R_2[W^*]$
is given by the assignment to 
$X \in \fra{sp}(W,\omega)$
of the quadratic polynomial function $f_X$ on $W$ where
$f_X(v) =\frac 12\omega(X v,v)$ ($v \in W$).
The inverse map sends a quadratic polynomial
$f$ to its Hamiltonian vector field $X_f = \{\cdot,f\}$.
More precisely, the assignment to $X \in \fra{sp}(W,\omega)$
of the Hamiltonian vector field of $f_X$ is the linear vector field
on $W$ induced by $X$, the corresponding map
from $\fra{sp}(W,\omega)$ to $\roman{Vect}(W)$
being 
an {\it anti\/} Lie isomorphism onto the Lie algebra of linear
vector fields.
For $W=\Bobb R^2$, with coordinates $q,p$ and $\omega$ the standard form
so that $\{q,p\} = 1$, with the notation $E,F,H$ in (3.2.2),
$$
\delta(E) = \frac {p^2}2,\ \delta(F) = -\frac {q^2} 2,\ 
\delta(H) = pq. 
$$
\smallskip
For general $W$ and $\omega$,
the unique momentum mapping 
$\mu \colon W \to (\Bobb R_2[W^*])^*$
(having the value zero at the origin)
for
the $\roman{Sp}(W,\omega)$-action on
$W$ 
is 
just the adjoint of $\delta$ and therefore
simply given by the formula
$$
(f \circ \mu) (v) = f(v),
\quad v \in W,\ f \in \Bobb R_2[W^*],
$$
which, in terms of
$\fra{sp}(W,\omega)$,
comes down to the well known formula
$$
(X \circ \mu) (v) = \frac 12 \omega(X v,v),
\quad v \in W,\ X \in \fra{sp}(W,\omega).
$$
Under the canonical isomorphism between
$(\Bobb R_2[W^*])^*$ 
and   the symmetric square 
$\roman S_{\Bobb R}^2[W]$ of $W$,
this momentum mapping amounts to the canonical squaring map
for $W$. 
Given $v \in W$, define
$v^\dagger \colon W \to \Bobb R$ 
by $\omega(u,v) = v^\dagger u$; the composite
$v v^\dagger$, viewed as a linear endomorphism
of $W$, then lies in
$\fra{sp}(W,\omega)$.
Plainly, for $v \in W$,
$$
(X \circ \mu) (v) =
\frac 12 v^\dagger X v = \frac 12 \roman{trace}(X \circ (v v^\dagger)),
\tag5.1.1
$$
that is to say,
when we identify $\fra{sp}(W,\omega)$ with its dual by means of the 
$\roman{Sp}(W,\omega)$-invariant
pairing 
$(A,B) \mapsto \frac 12 \roman{trace}(AB)$,
referred to henceforth as the {\it half-trace pairing\/}---it 
is a positive multiple
of the Killing form---the momentum mapping
$\mu$ amounts to the assignment to $v \in W$ of
$v v^\dagger \in \fra{sp}(W,\omega)$.

\smallskip\noindent
{\smc Example.} Let $W = \Bobb R^4$, with 
standard symplectic structure which,
in the coordinates
$q^1,q^2,p^1,p^2$,
amounts to the Poisson brackets $\{q^j, p^j\} = 1$ ($j=1,2$) etc.
With the notation
$v=[q^1,q^2,p^1,p^2]^{\tau}$,
we obtain
$$
\mu (v) =
v v^\dagger
=\left[
\matrix
q^1 p^1 &q^1 p^2 & -q^1 q^1 &-q^1 q^2 \\
q^2 p^1 &q^2 p^2 & -q^2 q^1 &-q^2 q^2 \\
p^1 p^1 &p^1 p^2 & -p^1 q^1 &-p^1 q^2 \\
p^2 p^1 &p^2 p^2 & -p^2 q^1 &-p^2 q^2 
\endmatrix
\right].
$$
We may identify $W$ with the total space $\roman T^*(\Bobb R^2)$ 
of the cotangent 
bundle of the ordinary 
$(q^1,q^2)$-plane $\Bobb R^2$. 
Ordinary angular momentum in the plane
has the infinitesimal generator
$X= \left[\matrix 
0 & -1 & 0 &  0 \\
1 &  0 & 0 &  0 \\
0 &  0 & 0 & -1 \\
0 &  0 & 1 &  0 
\endmatrix
\right]$, and
$$
(X \circ \mu) (v) =
\frac 12 \roman{trace}(X \circ (v v^\dagger)) =q^1 p^2 - q^2 p^1.
$$
This is the angular momentum mapping
for a single particle in the plane.
\smallskip
For general $W$ and $\omega$,
given   a 
subgroup $H$ of $\roman{Sp}(W,\omega)$,
the induced $H$-action  on $W$ 
induces an $H$-action on 
$\Bobb R_2[W^*]$
preserving the Poisson structure,
and the invariants
$\fra g=\Bobb R_2[W^*]^H$
constitute a Lie subalgebra
of $\Bobb R_2[W^*]$.
Indeed,
in $\fra{sp}(W)$,
$\fra g$
amounts to the $H$-invariant subspace
of $\fra{sp}(W)$, with reference to the adjoint representation.
The $H$-invariant map
$W \to \fra g^*$
whose adjoint is the injection $\fra g \to \Bobb R[W^*]$
is plainly a momentum mapping for the induced action of
the Lie group $G \subseteq \roman{Sp}(W,\omega)$ 
with $\roman {Lie}(G) = \fra g$.
For example, when $H$ is the obvious copy
of $\Bobb Z/ 2$ in 
$\roman{Sp}(W,\omega)$ 
whose generator acts on $W$ via multiplication by $-1$,
$\Bobb R_2[W^*]^H =\Bobb R_2[W^*]$, the homogeneous degree 2 part
$\Bobb R_2[W^*]$ generates the algebra of $\Bobb Z/2$-invariants,
and the 
momentum mapping
$\mu$ amounts to the Hilbert map of invariant theory
for the $\Bobb Z/2$-action on $W$.
\smallskip\noindent
{\smc 5.2. Associated representations.}
Given a unitary representation $E$ of a compact Lie group ${\KKK}$,
with $\fra {\kkk} = \roman{Lie}({\KKK})$,
the construction in (5.1) applies to
$W=E$ endowed with the symplectic structure
coming from the hermitian form on $E$
and yields  
the unique momentum mapping $\mu_{\KKK} \colon E \to \fra {\kkk}^*$
having the value zero at the origin of $E$.
Given a 
nilpotent orbit
$\Cal O$
in a real semisimple Lie algebra $\fra g$,
we will say that a unitary representation $E$
of a compact Lie group ${\KKK}$ is a {\it compact associated representation\/}
for $\Cal O$
provided
$\Bobb R_2[E^*]^{\KKK} \cong \fra g$
in such a way that the ${\KKK}$-invariant map
$\mu_G \colon E \to \fra g^*$
whose adjoint is the injection $\fra g \to \Bobb R[E^*]$
induces an isomorphism 
of stratified symplectic spaces
from the ${\KKK}$-reduced space
$E^{\roman{red}}=\mu_{\KKK}^{-1}(0) \big / {\KKK}$ 
onto the closure $\overline{\Cal O}$ of $\Cal O$
in $\fra g^*$.
\smallskip
More generally, let $\GHH$
be an arbitrary Lie group, 
not necessarily compact,
with Lie algebra
$\fra \ghh$,
which acts 
linearly and symplectically on a symplectic vector space
$E$, with
momentum mapping $\mu_{\GHH} \colon E \to \fra \ghh^*$
having the value zero at the origin of $E$.
We then refer to $E$ as a {\it symplectic representation\/} of $\GHH$
and {\it define\/} the reduced space
$E^{\roman{red}}$ 
to be the space 
of {\it closed\/} orbits in
the zero locus $\mu_{\GHH}^{-1}(0)$.
Given   a 
nilpotent orbit
$\Cal O$
in a real semisimple Lie algebra $\fra g$,
we will say that a symplectic representation $E$
of a  Lie group $H$ is an {\it associated representation\/}
for $\Cal O$
provided
$\Bobb R_2[E^*]^H \cong \fra g$
in such a way that the $H$-invariant map
$\mu_G\colon E \to \fra g^*$
whose adjoint is the injection $\fra g \to \Bobb R[E^*]$
induces a Poisson
map
$\overline\mu_G\colon E^{\roman{red}}  \to \fra g^*$
which is an
isomorphism 
of stratified spaces
from
$E^{\roman{red}}$ 
onto the closure $\overline{\Cal O}$ of $\Cal O$
in $\fra g^*$.

\smallskip\noindent
{\smc 5.3. Compact associated representations\/}.
We will now establish their existence 
for {\it holomorphic\/} nilpotent orbits 
in the {\it standard cases\/}.
Our tool will be an extension of the
{\it basic construction\/}
in \cite\kraprotw;
in that reference,
this construction is given only for the case
of the classical groups over
the complex numbers.
We now extend this construction to 
the appropriate real forms of the classical groups.
A key step will be an extension of the \lq\lq First Main Theorem
of Invariant Theory\rq\rq\ to the real forms of the classical groups,
given as Theorem 5.3.1 below.
The significance of
the First Main Theorem
of Invariant Theory
for orbits
of complex classical groups has been noticed
in \cite\kraprotw.
\smallskip
Let $\Bobb K = \Bobb R, \Bobb C, \Bobb H$,
let $V^s =\Bobb K^s$,
the standard 
(right) $\Bobb K$-vector space of dimension $s$,
endowed with 
a non-degenerate hermitian form
$(\cdot,\cdot)$ 
of type $(s',s'')$, that is, of rank $s=s'+s''$
and signature $s'-s''$;
further, let $V = \Bobb K^n$, 
endowed with  
a 
complex structure $J_V$ and
skew form 
$\Cal B$ of the kind reproduced in (3.5),
let $\GHH =\GHH(s',s'')= 
\roman U(V^s,(\cdot,\cdot))$,
$\fra h = \roman{Lie}(\GHH)$, $\HG = \roman U(V,\Cal B)$,
and
$\fra g = \roman{Lie}(\HG)$
(cf. (3.5)),
and denote the split rank of 
$\HG = \roman U(V,\Cal B)$ by $r$.
When $s=s'$ and $s''=0$
we will write 
$\GHH(s)$ instead of $\GHH(s',s'')$.
More explicitly:
\newline\noindent
(1) $\Bobb K=\Bobb R,\ V^s = \Bobb R^s,\ V = \Bobb R^n,\
n = 2 \ell,\ 
\GHH= \roman O(s',s''),\ 
\HG= \roman {Sp}(\ell,\Bobb R),\ r = \ell$; 
\newline\noindent
(2) $\Bobb K=\Bobb C,\ V^s = \Bobb C^s,\ V = \Bobb C^n,\ 
n = p +q,\ 
\GHH= \roman U(s',s''),\ 
\HG= \roman U(p,q)$, $p \geq q$, $r = q$;
\newline\noindent
(3) $\Bobb K=\Bobb H,\ V^s = \Bobb H^s,\ V = \Bobb H^n,
\GHH= \roman U(s',s'',\Bobb H) = \roman {Sp}(s',s''),
\HG= \roman O^*(2 n)$, $r = [\frac n2]$.
\newline
When $s=s'$ and $s''=0$, $\KK=\KK(s)$ is the group
$\roman O(s,\Bobb R)$,
$\roman U(s)$, or
$\roman {Sp}(s)$ 
as appropriate.
\smallskip
Let $W=W(s)=\roman{Hom}_{\Bobb K}(V^s,V)$.
The groups $\KK$ and $\HG$ act on $W$ in an obvious fashion: Given
$x \in \KK, \ \alpha \in \roman{Hom}_{\Bobb K}(V^s,V),\ 
y \in \HG$,
the action is given by the assignment to
$(x,y,\alpha)$ of $y \alpha x^{-1}$.
Given a $\Bobb K$-linear map $\alpha \colon V^s \to V$, define
the $\Bobb K$-linear map $\alpha^\dagger \colon V \to V^s$ by
$$
(\alpha^\dagger \bold u,\bold v) = \Cal B(\bold u,\alpha \bold v),
\ \bold u \in V,\, \bold v \in V^s;
$$
given
$x \in \KK, y \in \HG$, plainly
$(y \alpha x^{-1})^\dagger =x \alpha^\dagger y^{-1}$.
In case (1), for $s=1$, this construction comes down to the 
assignment to $v \in W$ of $v^\dagger \colon W \to \Bobb R$
given in (5.1).
For intelligibility we give explicit formulas for 
$\alpha^{\dagger}$
for the special cases
where $s'=s$ 
(arbitrary value)
and $s''=0$.
\newline\noindent
Case (1): 
$\alpha = \left[\bold q_1, \ldots, \bold q_{\ell},
\bold p_1 ,\ldots ,\bold p_{\ell}\right]^{\tau}$:
$\alpha^{\dagger} = [\bold p_1,\dots,\bold p_{\ell}, 
-\bold q_1, \dots, -\bold q_{\ell}]
\colon \Bobb R^n \to \Bobb R^s$,
\linebreak 
$\bold q_1, \ldots, \bold q_{\ell},
\bold p_1 ,\ldots ,\bold p_{\ell} \in \Bobb R^s$
being the column vectors of $\alpha^{\dagger}$.
\newline\noindent
Case (2): 
$\alpha = \left[\bold w_1,\ldots, \bold w_{p+q}\right]^{\tau}$:
$\alpha^{\dagger} = i[\overline{\bold w_1},\dots,\overline{\bold w_p}, 
-\overline{\bold w_{p+1}},\dots,-\overline{\bold w_{p+q}}]
\colon \Bobb C^{p+q} \to \Bobb C^s$,
\linebreak 
$\bold w_1,\ldots, \bold w_{p+q} \in \Bobb C^s$
being the column vectors of $\alpha^{\dagger}$.
\newline\noindent
Case (3): 
$\alpha = \left[\bold w_1,\ldots, \bold w_n\right]^{\tau}$:
$\alpha^{\dagger} 
= [\overline{\bold w_1}\Cal J,\dots,\overline{\bold w_n}\Cal J] 
\colon \Bobb H^n \to \Bobb H^s$,
$\bold w_1,\ldots, \bold w_n \in \Bobb H^s$
being the column vectors of $\alpha^{\dagger}$.
\smallskip
The 2-forms $(\cdot,\cdot)$ and 
$\Cal B$
induce a symplectic structure
$\omega_W$ on $W$ 
by means of the assignment
$$
\omega_W(\alpha,\beta) = \roman{trace}_{\roman r}(\beta^{\dagger}\alpha),
\quad \alpha,\beta \in W,
$$
where
$\roman{trace}_{\roman r}(\beta^{\dagger}\alpha) \in \Bobb R$
is the {\it real\/} constituent of the value
of the trace of $\beta^{\dagger}\alpha$.
In case (1) where $\Bobb K= \Bobb R$ the function 
$\roman{trace}_{\roman r}$ is the ordinary trace function.
In case (2) and case (3),
the hermitian form
$(\cdot,\cdot)$
on $V^s = \Bobb K^s$
determines a 
non-degenerate real symmetric bilinear form
$(\cdot,\cdot)_r$ on $\Bobb R^s$
of the same signature as
$(\cdot,\cdot)$
such that
the {\it real\/} constituent 
of $(\cdot,\cdot)$
may be written as
$(\Bobb R^s,(\cdot,\cdot)_r) \otimes_{\Bobb R} \Bobb K$,
and the real 
constituent
$\Cal B_1$ 
of $\Cal B$
is a real
symplectic structure
on $V$ as observed earlier.
Writing
$W$ in the form $W=\roman{Hom}_{\Bobb R}(\Bobb R^s, V)$,
given $\alpha \colon \Bobb R^s \to V$,
we may define $\alpha^{\dagger\dagger} \colon V \to \Bobb R^s$ by
${
(\alpha^{\dagger\dagger} \bold u,\bold v)_r = 
\Cal B_1(\bold u,\alpha \bold v)
}$ 
where
$\bold u \in V$ and $\bold v \in \Bobb R^s$;
then
$$
\omega_W(\alpha,\beta) =
\roman{trace}(\beta^{\dagger\dagger}\alpha),
\quad \alpha,\beta \in W,
$$
where \lq\lq trace\rq\rq\ now refers to the ordinary real trace on 
$\roman{End}(\Bobb R^s)$,
and where we do not distinguish in notation between
an $\Bobb R$-linear map $\alpha$ from $\Bobb R^s$ to $V$
and its
$\Bobb K$-linear extension to $V^s$.
\smallskip
The actions of $\KK$ and $\HG$ on $W$ 
preserve the symplectic structure $\omega_W$ and
are hamiltonian.
Up to the requisite identifications
of $\fra h$ and $\fra g$ with their duals
by means of the half-trace pairing,
cf. (5.1) above,
the momentum mappings $\mu_{\KK}$ and
$\mu_{\HG}$
for the $\KK$- and $\HG$-actions are given by
$$
\aligned
\mu_{\KK} &\colon W @>>> \fra {\kk},
\quad
\mu_{\KK}(\alpha)= -\alpha^\dagger \alpha \colon V^s \to V^s,
\\
\mu_{\HG} &\colon W @>>> \fra {\hg},
\quad
\mu_{\HG}(\alpha) =
 \alpha \alpha^\dagger \colon V \to V,
\endaligned
$$
respectively. 
Indeed, since the $\KK$-action on $W$
is given by
the assignment
to
$(x,\alpha)$ of $\alpha x^{-1}$
($x \in \KK,\, \alpha \in W$)
the $\fra{\kk}$-action
on $W$ is given by
the formula
$$
X(\alpha) =-\alpha X,\quad X \in \fra{\kk},\, \alpha \in W.
$$
Hence, cf. (5.1.1),
$$
(X \circ \mu_{\KK}) (\alpha) =
\frac 12 \omega_W(X(\alpha),\alpha) =
-\frac 12 \omega_W(\alpha X,\alpha) =
-\frac 12 \roman{trace}_{\roman r}((\alpha^{\dagger}\alpha)\circ X).
$$
This shows that
$\mu_{\KK}$ is indeed the $\KK$-momentum mapping.
A similar reasoning justifies the claim
that
$\mu_{\HG}$ is the $\HG$-momentum mapping.
\smallskip
We could have identified 
$\fra h$ and $\fra g$ with their duals
by means of the {\it trace\/} pairing,
but we prefer {\it not\/} to do so
since,
up to sign,
with our conventions,
the momentum mappings will appear as the
Hilbert maps of invariant theory.
\smallskip
\noindent
{\smc Illustration.} Consider
the case (1) with
$s'=s$ and $s''=0$ 
so that $\KK = \roman O(s,\Bobb R)$.
Given
$\bold x = (x^1,\dots, x^s)\in \Bobb R^s$ and
$\bold y = (y^1,\dots, y^s) \in \Bobb R^s$,
with the standard notation
$$
\bold x \wedge \bold y\in \Lambda_{\Bobb R}^2[\Bobb R^s] =\fra {so}(s,\Bobb R)
$$
for the skew-symmetric matrix 
having $x^j y^k - x^k y^j$ as its $(j,k)$'th entry,
when
the linear map
$\alpha \colon V^s \to V$ 
is written as
${\alpha = \left[\bold q_1, \ldots, \bold q_{\ell},
\bold p_1 ,\ldots ,\bold p_{\ell}\right]^{\tau}}$
with row vectors $\bold q_1,\dots,\bold q_{\ell},
\bold p_1,\dots,\bold p_{\ell}  \in \Bobb R^s$,
the values
$\mu_{\KK}(\alpha) =-\alpha^\dagger \alpha 
\in \fra {so}(s,\Bobb R)$
and
$\mu_{\HG}(\alpha) = \alpha \alpha^\dagger\in \fra {sp}(\ell,\Bobb R)$
admit the following descriptions:
$$
\align
\mu_{\KK}(\alpha) =-\alpha^\dagger \alpha 
&=\bold q_1\wedge \bold p_1 +\dots +\bold q_{\ell}\wedge\bold p_{\ell}
\in \Lambda_{\Bobb R}^2[\Bobb R^s] =\fra {so}(s,\Bobb R)
\\
\mu_{\HG}(\alpha) = \alpha \alpha^\dagger
&=
\left[\matrix 
\left[\bold q_j \bold p_k\right] & -\left[\bold q_j \bold q_k\right]\\
\left[\bold p_j \bold p_k\right] & -\left[\bold p_j \bold q_k\right]
\endmatrix
\right]
\in \fra {sp}(\ell,\Bobb R).
\endalign
$$
Here 
$\fra {so}(s,\Bobb R)$
(as well as $\fra {sp}(\ell,\Bobb R)$)
is identified with its dual
via the
half-trace pairing 
$(A,B) \mapsto \frac 12 \roman{trace}(AB)$.
The map $\mu_{\KK}$  comes down to the ordinary angular momentum mapping for
a constrained system of $\ell$ particles in $\Bobb R^s$ whence the notation
$\bold q_1,\dots,\bold q_{\ell},
\bold p_1,\dots,\bold p_{\ell}$.
For example, when $s=2$ and $\ell=1$, 
$W = \roman{End}_{\Bobb R}(\Bobb R^2)$ and,
taking 
$\alpha = \left[\matrix \bold q\\ \bold p\endmatrix\right]
= \left[\matrix q^1&q^2\\ p^1&p^2\endmatrix\right]$,
we have
$\alpha^\dagger = \left[\matrix p^1&-q^1\\ p^2&-q^2\endmatrix\right]$
and
$
-\alpha^\dagger \alpha =
\left[\matrix 0&q^1p^2- p^1q^2\\ q^2p^1-p^2q^1&0\endmatrix\right].
$
When we identify
$W$ with the total space $\roman T^*(\Bobb R^2)$ 
of the cotangent 
bundle of the ordinary 
$(q^1,q^2)$-plane $\Bobb R^2$, 
taking the infinitesimal generator
$X= \left[\matrix 
0 & -1  \\
1 &  0 
\endmatrix
\right]$
of ordinary angular momentum, we have
$$
(X \circ \mu_{\KK}) (\alpha) =
-\frac 12 \roman{trace}((\alpha^\dagger\alpha)\circ X) =q^1 p^2 - q^2 p^1 .
$$
This is again the angular momentum mapping
for a single particle in the plane,
the present description being a variant of
the example in (5.1) above. 
\smallskip
Let $E$ be a (real or complex) symplectic vector space.
Recall that two subgroups $\Gamma$ and $\Gamma'$
of $\roman {Sp}(E)$
are
said to constitute a {\it reductive dual pair\/} provided
$\Gamma$ and $\Gamma'$ are each other's full centralizer in
$\roman {Sp}(E)$ \cite\howeone.
In each of the cases (1),(2),(3) above,
the two groups $G$ and $H$ constitute a real reductive dual pair in
$\roman {Sp}(W)$ \cite\howetwo;
these are precisely the pairs
(5.2)(i),(iii),(iv) in that reference.
Below $(G,H)$ refers to any of these three pairs.

\proclaim{Theorem 5.3.1}
The ring of all $\KK$-invariant polynomials
in $\Bobb R[W^*]$ is generated by its homogeneous quadratic 
part, and the canonical map from $\Bobb R_2[W^*]^{\KK}$
to $\fra {\hg}$ is an isomorphism
in such a way that
the adjoint of the injection 
of $\fra {\hg}$
into 
$\Bobb R[W^*]$
amounts to the
$\HG$-momentum mapping $\mu_{\HG}\colon W \to  \fra {\hg}^*$.
Likewise,
the ring of all $\HG$-invariant polynomials
in $\Bobb R[W^*]$ is generated by its homogeneous quadratic 
part,
and the canonical map from $\Bobb R_2[W^*]^{\HG}$
to $\fra {\kk}$ is an isomorphism
such that the
$\KK$-momentum mapping $\mu_{\KK}\colon W \to  \fra {\kk}^*$
is the adjoint of the injection 
of $\fra {\kk}$
into 
$\Bobb R[W^*]$.
\endproclaim

\proclaim{Corollary 5.3.2}
The components of the
$\mu_{\KK}$-momentum mapping
generate the algebra of $\HG$-invariants
in $\Bobb R[W^*]$
and the components of the
$\mu_{\HG}$-momentum mapping
generate the algebra of $\KK$-invariants
in $\Bobb R[W^*]$.
\endproclaim

This Corollary extends a result in \cite\lermonsj.
Its statement
is a version of the \lq\lq First Main Theorem
of Invariant Theory\rq\rq\ 
\cite\weylbook, cf. also \cite\howeone.
However this theorem does {\it not\/}
say that the 
$\KK$-invariants combine to a
momentum mapping for $\HG$
nor that
the $\HG$-invariants combine to a
momentum mapping for $\KK$.
Moreover, this theorem is usually phrased only for the
classical groups
$\roman O(n,\Bobb R),\roman O(n,\Bobb C),
\roman{Sp}(n,\Bobb R),\roman{Sp}(n,\Bobb C),
\roman{GL}(n,\Bobb R), \roman{GL}(n,\Bobb C)$.
There is as well a close relationship between
the situation of Theorem 5.3.1 and the notion of dual pair of
Poisson mappings \cite\weinstwo.

\demo{Proof of Theorem 5.3.1}
Since $\KK$ and $\HG$
constitute a reductive dual pair
in $\roman {Sp}(W)$,
the homogeneous degree 2 part
$\Bobb R_2[W^*]^{\KK}$
coincides with (a copy of)
$\fra {\hg}$ and, likewise, 
the homogeneous degree 2 part
$\Bobb R_2[W^*]^{\HG}$
coincides with (a copy of)
$\fra {\kk}$.
We note that,
by construction,
$\fra {\hg}$ and 
$\fra {\kk}$
are thus 
realized within 
$\Bobb R_2[W^*] =\fra{sp}(W)$,
which equals 
$\fra{sp}(s\ell,\Bobb R)$ for $\Bobb K = \Bobb R$,
$\fra{sp}(sn, \Bobb R)$   for $\Bobb K = \Bobb C$,
$\fra{sp}(2sn,\Bobb R)$   for $\Bobb K = \Bobb H$.
\smallskip
To see that
the subalgebra $\Bobb R[W^*]^{\KK}$ 
of real invariant polynomials
is generated by its homogeneous degree 2 part
$\Bobb R_2[W^*]^{\KK}$,
we complexify the $\KK$-action on $W$
and examine the three cases separately.
\smallskip\noindent
Case (1): $W=\roman{Hom}_{\Bobb R}(\Bobb R^s, \Bobb R^{2\ell}),\ 
\roman{Sp}(W)=\roman{Sp}(s\ell, \Bobb R),\ 
\KK=O(s',s'',\Bobb R),\, \HG = \roman{Sp}(\ell,\Bobb R)$.
By construction,
$W_{\Bobb C} =W \otimes_{\Bobb R} \Bobb C = \Bobb C^{2s\ell} 
= \Bobb C^s\otimes_{\Bobb C} \Bobb C^{2\ell}$,
that is,
$W_{\Bobb C}$
may be written 
as
the product 
of $2\ell$ copies of
$\Bobb C^s$,
and the ${\KK}$-action on $W$ complexifies to
the diagonal action of $\roman{O}(s,\Bobb C)=\KK^{\Bobb C}$
on 
the product 
of $2\ell$ copies of
$\Bobb C^s$.
\smallskip\noindent
Case (2): $W= \roman{Hom}_{\Bobb C}(\Bobb C^s,\Bobb C^n)\cong\Bobb C^{sn},\ 
\roman{Sp}(W)=\roman{Sp}(sn, \Bobb R),\ \KK=\roman U(s',s''),
\ \HG=\roman U(p,q)$, $n=p+q$.
Now
$
W_{\Bobb C} =W \otimes_{\Bobb R} \Bobb C \cong \Bobb C^{2sn} 
\cong \Bobb C^{2s}\otimes_{\Bobb C} \Bobb C^{n}, 
$
that is,
$W_{\Bobb C}$
may be written 
as
the product 
of $n$ copies of
$\Bobb C^{2s}$,
and the ${\KK}$-action on $W$ complexifies to
the diagonal action of $\roman{GL}(s,\Bobb C)=\KK^{\Bobb C}$
on 
the product 
of $n$ copies of
$\Bobb C^s \times (\Bobb C^s)^*$
where $(\Bobb C^s)^*$ is identified with
$\overline {\Bobb C^s}$
by means of the hermitian form.
\smallskip
\noindent
Case (3): $W= \roman{Hom}_{\Bobb H}(\Bobb H^s, \Bobb H^n) \cong \Bobb H^{sn},
\roman{Sp}(W) =\roman{Sp}(2sn, \Bobb R),\KK=\roman{Sp}(s',s''),
\HG= \roman O^*(2 n)$.
Now
$
W_{\Bobb C} =W \otimes_{\Bobb R} \Bobb C \cong \Bobb C^{2sn} 
\cong \Bobb C^{2s}\otimes_{\Bobb C} \Bobb C^n, 
$
that is,
$W_{\Bobb C}$
may be written 
as
the product 
of $n$ copies of
$\Bobb C^{2s}$,
and the $\KK$-action on $W$ complexifies to
the diagonal action of $\roman{Sp}(s,\Bobb C)={\KK}^{\Bobb C}$
on 
the product 
of $n$ copies of
$\Bobb C^{2s}$.
\smallskip
By the first main theorem of invariant theory, 
the algebra $\Bobb C[W_{\Bobb C}^*]^{\KK^{\Bobb C}}$
of complex
${\KK^{\Bobb C}}$-invariants
in
$\Bobb C[W_{\Bobb C}^*]$
is generated by its homogeneous degree 2 part
$\Bobb C_2[W_{\Bobb C}^*]^{\KK^{\Bobb C}}$.
Since $\KK$ and $\HG$
constitute a reductive dual pair
in $\roman {Sp}(W)$,
the homogeneous degree 2 part
$\Bobb C_2[W_{\Bobb C}^*]^{\KK^{\Bobb C}}$ is the complexification of
the homogeneous 
degree 2 part
$\Bobb R_2[W^*]^{\KK}$.
This implies that the subalgebra $\Bobb R[W^*]^{\KK}$ 
of real invariant polynomials
is generated by its homogeneous degree 2 part
$\Bobb R_2[W^*]^{\KK}$.
Interchanging the roles of
$\KK$ and $\HG$ we see that
the subalgebra $\Bobb R[W^*]^{\HG}$ of real invariant polynomials
is generated by its homogeneous degree 2 part
$\Bobb R_2[W^*]^{\HG}$. \qed
\enddemo

\noindent
{\smc Remark.}
Our description
gives the
\lq\lq Second Main Theorem
of Invariant Theory\rq\rq\ 
\cite\weylbook\ 
for free; indeed, this theorem provides
defining relations among those generators
of the invariants which occur above as
the components of the corresponding momentum mapping,
and these defining relations 
just fix the rank of the matrices
whose entries are the generators.

\noindent
{\smc Remark.}
In case (2), letting $s=1$ so that $H=\roman U(1)$, the resulting injection
of $\fra g=\fra u(p,q)= \Bobb R_2[W^*]^{\KK}$ into 
$\Bobb R_2[W^*] = \fra{sp}(p+q,\Bobb R)$ is one of the kind (3.5.1.5),
up to the central factor in $\fra u(p,q)$.
Likewise,
in case (3), letting $s=1$ so that $H=\roman {Sp}(1)$, the resulting injection
of $\fra g=\fra {so}^*(2n)= \Bobb R_2[W^*]^{\KK}$ into 
$\Bobb R_2[W^*] = \fra{sp}(2n,\Bobb R)$ is one of the kind (3.5.1.7).
\smallskip
We will now exploit Theorem 5.3.1 to relate certain
reduced spaces with
nilpotent orbits.
The basic observation is that,
under the circumstances of that theorem,
given $\alpha \colon V^s \to V$,
the vanishing of $\mu_{\KK}(\alpha)$ entails
$\mu_{\HG}(\alpha)\mu_{\HG}(\alpha) = 0$, that is,
$\mu_{\HG}(\alpha)$ is then a 
square zero
nilpotent matrix in $\fra {\hg}$.
Since $\alpha$ has rank at most $s$,
the matrix  $\alpha \alpha^\dagger \in \fra {\hg}$,
viewed as a $\Bobb K$-linear map from $\Bobb K^n$ to
$\Bobb K^n$, has $\Bobb K$-rank at most $s$.
As a side remark we note that, likewise,
$\mu_{\HG}(\alpha)=0$ entails
$\mu_{\KK}(\alpha)\mu_{\KK}(\alpha) = 0$;
in particular,
when $\KK$ is compact, 
$\mu_{\KK}(\alpha)$ is then necessarily zero.
As before, we take $W^{\roman{red}}$ to be
the space of {\it closed\/} $\KK$-orbits in the zero locus
of $\mu_{\KK}$. 
Whether or not
$\KK$ is compact,
by $\KK$-invariance, the momentum mapping
$\mu_{\HG}$ induces a map
from the space 
of closed $\KK$-orbits in $W$
to $\fra {\hg}$ and hence a map
$\overline\mu_{\HG}\colon W^{\roman{red}} \to \fra {\hg}$
the image of which, in turn, lies in
the subspace
of square zero nilpotent matrices in $\fra {\hg}$
which have rank at most $s$
since for any $\alpha \colon V^s \to V$
the image $\mu_{\HG}(\alpha) =\alpha\alpha^\dagger \colon V \to V$
has rank at most $s$.
We will show that 
$\overline\mu_{\HG}$ identifies
the $\KK$-reduced space $W^{\roman{red}}$
with the closure
of 
finitely many pseudoholomorphic
nilpotent orbits
in $\fra {\hg}$
of rank 
$s_r= \roman{min}(\rell,s)$. 
In particular,
when
$\KK$ is compact,
i.~e. $s=s'$ and $s'' = 0$,
$\overline\mu_{\HG}$ identifies
$W^{\roman{red}}$
with the closure
of 
the holomorphic
nilpotent orbit
$\Cal O_{s_r}$
which consists of
non-negative 
nilpotent matrices
in $\fra {\hg}$
of rank 
$s_r= \roman{min}(\rell,s)$. 

\smallskip
Recall that
$V$ 
carries the complex structure
$J=J_V$ 
which is induced by the element $J_V \in \fra {\hg}$.
(We remind the reader that $\frac 12 J_V$ is the requisite
$H$-element in 
the standard cases, cf. (3.5)).
Define $J_W \colon W \to W$ by
$J_W (\alpha) = J_V \circ \alpha$, where $\alpha \colon V^s \to V$.
This yields a complex structure $J_W$ on $W$
which,
together with the symplectic structure
$\omega_W$,
turns $W$ into a flat
(not necessarily positive)
K\"ahler manifold.
When we wish to emphasize in notation that $W$ is endowed
with this complex structure we write $W_J$.
\smallskip
Let $V^{s-1}= \Bobb K^{s-1}\subseteq V^s$,
and endow it with the induced hermitian form
which we still write as $(\cdot,\cdot)$.
Let $W(s-1)=\roman{Hom}_{\Bobb K}(V^{s-1},V)$,
$\KK(s-1)= \roman U(s-1,\Bobb K)$.
The inclusion
$V^{s-1}\subseteq V^s$
induces a projection
$V^s \to V^{s-1}$
of hermitian $\Bobb K$-vector spaces which,
in turn, induces an injection
$W(s-1) \to W=W(s)$ of $\Bobb K$-vector spaces
compatible with all the structure:
It is $\KK(s-1)$-linear where $\KK(s-1)$ is viewed as a subgroup of
$\KK(s)$ in the obvious fashion; it is $\HG$-linear;
it is compatible with the complex structure $J_W$;
and it is compatible with the momentum mappings in the sense that the diagrams
$$
\CD
W(s-1) @>>> W(s)   @.\phantom {longlonglong} @.    W(s-1) @>>> W(s)
\\
@V{\mu_{\KK(s-1)}}VV       
@VV{\mu_{\KK(s)}}V        
@.
@V{\mu_{\HG}}VV
@VV{\mu_{\HG}}V
\\
\fra {\kk}(s-1) @>>> \fra {\kk}(s) @. 
\phantom {longlonglong} @. \fra {\hg} @>{\roman{Id}}>> \fra {\hg}
\endCD
$$
are commutative.

\smallskip
We now suppose that $\KK$ is compact.
Thus ${\KKK}= \roman O(s,\Bobb R)$,
${\KKK}= \roman U(s)$,
or
${\KKK}= \roman {Sp}(s)$, according to the case under consideration.
The complex structure $J=J_W$ on $W$
now combines with the symplectic structure
$\omega_W$
to a {\it positive\/} 
flat K\"ahler structure on $W$---maintaining notation
established earlier, we write $W_J$ for $W$,
endowed with this K\"ahler structure---, and the ${\KKK}$-representation
on $W_J$ is unitary;
in particular, it extends to a holomorphic action of
${\KKK}^{\Bobb C}$ on $W_J$.
Thus, 
in view of Proposition 4.2,
via the Kempf-Ness homeomorphism
from $W^{\roman{red}}=\mu_{\KKK}^{-1}\big/{\KKK}$ onto the categorical quotient
$W_J\big /\big / {\KKK}^{\Bobb C}$, 
the ${\KKK}$-reduced space 
$W^{\roman{red}}$
inherits a (positive) normal K\"ahler structure.
\smallskip
The next result 
establishes 
the existence of a compact associated representation
for any holomorphic nilpotent orbit
in a standard simple Lie algebra of hermitian type.
As a by-product, 
independently of
the classification of pseudoholomorphic
nilpotent orbits given in (3.3) above,
we then deduce
the 
holomorphicity
of the nilpotent orbits spelled out there.
Thus 
the next result
does {\it not\/} rely
on the fact,
obtained in (3.3) via the classification
of pseudoholomorphic nilpotent orbits, that
non-negative nilpotent elements
in $\fra \hg$
(i.~e. matrices $X$ in $\fra g$ with $\Cal B_X$ non-negative, see
(3.5) above)
generate holomorphic nilpotent orbits.

\proclaim{Theorem 5.3.3}
The induced map $\overline \mu_{\HG}$
from the ${\KKK}$-reduced space
$W(s)^{\roman{red}} = \mu_{\KKK}^{-1}(0)\big / {\KKK}$
to
$\fra {\hg}$  
is a proper embedding 
of
$W(s)^{\roman{red}}$ into $\fra {\hg}$
and
induces an isomorphism of 
stratified symplectic spaces
from $W(s)^{\roman{red}}$ onto the closure $\overline {\Cal O_{s_r}}$
of the unique
nilpotent orbit $\Cal O_{s_r}$
which consists of non-negative nilpotent matrices in $\fra {\hg}$ of rank 
$s_r=\roman{min}(r,s)$,
and $\overline {\Cal O_{s_r}}$
consists of non-negative nilpotent matrices in $\fra {\hg}$ of rank 
at most $\roman{min}(r,s)$.
The orbit 
$\Cal O_{s_r}$
is necessarily holomorphic,
and the isomorphism
$\overline \mu_{\HG}$
from 
$W(s)^{\roman{red}}$
onto 
$\overline {\Cal O_{s_r}}$
is one of
complex analytic stratified K\"ahler spaces
whence $\overline {\Cal O_{s_r}}$ is a normal K\"ahler space.
Furthermore, when $s > r$,
the inclusion of $W(s-1)$ into $W(s)$ induces an isomorphism
of $W(s-1)^{\roman{red}}$ onto $W(s)^{\roman{red}}$
while when
$s \leq r$,
so that $s_r = s$, 
the inclusion of $W(s-1)$ into $W(s)$ induces an 
injection 
of $W(s-1)^{\roman{red}}$ into $W(s)^{\roman{red}}$
which, under the identifications with the closures of 
holomorphic
nilpotent orbits,
amounts to the inclusion
of  $\overline {\Cal O_{s-1}}$
into $\overline {\Cal O_s}$
in such a way that  
$\overline {\Cal O_{s-1}}$
is the complement of the top stratum
$\Cal O_s$ of
$\overline {\Cal O_s}$.
Finally,
under the
projection from $\fra g =\fra k \oplus \fra p$
to $\fra p$, followed by the identification of the latter
with $\fra p^+$,
the image of the reduced space
$W^{\roman{red}}$ in $\fra g$
is identified
with the affine complex variety consisting of
matrices in $\fra p^+$ 
of rank at most $s$
in such a way that the complex structure 
on $\overline {\Cal O_s}$
which arises from the complex algebraic quotient
construction gets identified with that resulting
from the projection of 
$\overline {\Cal O_s}$
to $\fra p$;
in particular, when
$s =r$,
the image of the reduced space
$W^{\roman{red}}$ in $\fra g$
is identified with 
$\fra p^+$,
that is,
the categorical quotient
$W_J\big /\big / {\KKK}^{\Bobb C}$
is identified with 
$\fra p^+$.
\endproclaim

Thus we see once more that, for $s \leq r$, 
$\{0\},\Cal O_1,\dots,\Cal O_s$
are the strata of
$\overline {\Cal O_s}$,
and that their closures constitute a chain of inclusions
$\{0\}\subseteq
\overline {\Cal O_1}
\subseteq
\dots
\subseteq
\overline {\Cal O_{s-1}}
\subset
\overline {\Cal O_s}$.
A somewhat more explicit description of the complex analytic structures
of the strata was given in Theorem 3.5.5 above.
For the special case
where ${\KKK}=\roman O(s,\Bobb R)$ and $\HG=\roman {Sp}(\ell,\Bobb R)$
(case (1) above),
the observation that
$\overline \mu_{\HG}$
identifies the reduced space with the closure of a certain nilpotent
orbit may be found in \cite\lermonsj.

\demo{Proof}
We have already observed that
the momentum mapping $\mu_{\HG}$
induces a map 
$\overline \mu_{\HG}$
from $W^{\roman{red}}$ to 
the space of nilpotent matrices
in $\fra{\hg}$
which have rank at most $s$.
Since ${\KKK}$ and $\HG$ centralize each other
in $\roman{Sp}(W) =\roman{Sp}(\frac {dsn}2, \Bobb R)$,
$W^{\roman{red}}$
inherits a 
$\HG$-action,
and $\overline \mu_{\HG}$
is $\HG$-equivariant.
By Theorem 5.3.1,
the algebra 
$\Bobb R[W^*]^{\KKK}$ of ${\KKK}$-invariants 
is generated by the components
of 
$\mu_{\HG} \colon W \to \fra{\hg}$.
Hence the map
from the real categorical quotient
$W\big/\big/{\KKK}$ to $\fra {\hg}$ induced by $\mu_{\HG}$
is the ordinary (real) Hilbert map;
since for $X = \alpha \alpha^{\dagger}$
and $v,w \in V$,
$$
\Cal B_X(v,w) = - \Cal B(\alpha \alpha^{\dagger}v,w)
=
\overline{\Cal B(w,\alpha \alpha^{\dagger}v)}
=
\overline{(\alpha^{\dagger}w,\alpha^{\dagger}v)}
=
(\alpha^{\dagger}v,\alpha^{\dagger}w),
$$
this Hilbert map realizes
$W\big/{\KKK}$
(in a $\HG$-equivariant fashion)
as the space of non-negative matrices in 
$\fra {\hg}$ of rank at most $s$.
Consequently the induced map
$\overline \mu_{\HG}$ from the ${\KKK}$-reduced space
$W^{\roman{red}} = \mu_{\KKK}^{-1}(0)\big / {\KKK}$ to $\fra {\hg}$
is injective as well.
Furthermore, the square of the distance to the origin in
$W$ is a ${\KKK}$-invariant polynomial function.
Hence the Hilbert map is proper, and
$\overline \mu_{\HG}$ is a 
$\HG$-equivariant
homeomorphism 
from $W^{\roman{red}}$
onto its image.
\smallskip
Let $\Cal O_{s_r}$ be the space
of non-negative nilpotent matrices in $\fra {\hg}$
of rank 
$s_r=\roman{min}(r,s)$.
We claim that
$\Cal O_{s_r}$
is a single $\HG$-orbit, that
$\overline \mu_{\HG}$ 
maps $W^{\roman{red}}$ onto $\overline{\Cal O_{s_r}}$,
and that
$\overline{\Cal O_{s_r}}$
is the space of non-negative nilpotent matrices in $\fra {\hg}$
of rank 
at most $s_r=\roman{min}(r,s)$.
To justify this claim
we observe first that, for any $\alpha \colon V^s \to V$, 
the rank of
$\alpha \alpha^{\dagger}$ is at most $r$, whatever $s$. 
Indeed, the ${\KKK}$-reduced space $W^{\roman{red}}$
consists of orbits of $\alpha \colon V^s \to V$
such that
$\alpha^{\dagger} \alpha \colon V^s \to V^s$ is zero;
for such an $\alpha$,
necessarily
$\alpha \alpha^{\dagger} \alpha =0$,
that is,
$\alpha (V^s)$ lies in
the kernel of
$\alpha \alpha^{\dagger}$.
Thus, if $s>r$,  
if $\alpha^{\dagger} \alpha =0$,
and if $\roman{rank}(\alpha)>r$,
the kernel of
$\alpha \alpha^{\dagger}$
contains a space of dimension $>r$
whence
$\roman{rank}(\alpha \alpha^{\dagger})<r$.
Consequently
$\overline \mu_{\HG}$ 
maps $W^{\roman{red}}$ into the  space
of nilpotent orbits in
$\fra {\hg}$
of rank 
at most $s_r=\roman{min}(r,s)$.
Furthermore, 
this observation entails that,
when $s > r$,
the inclusion of $W(s-1)$ into $W(s)$ induces an isomorphism
of $W(s-1)^{\roman{red}}$ onto $W(s)^{\roman{red}}$.
\smallskip
Next we show that,
for some $\alpha_0$
with $\mu_{\KKK}(\alpha_0)=0$,
the image $\mu_{\HG}(\alpha_0) \in \fra {\hg}$
has rank
$s_r=\roman{min}(r,s)$.
To simplify the exposition,
we will henceforth suppose that $s \leq \rell$.
\smallskip\noindent
Case (1): $W=\roman{Hom}_{\Bobb R}(\Bobb R^s, \Bobb R^{2\ell}),\ 
\roman{Sp}(W)=\roman{Sp}(s\ell, \Bobb R),\ 
{\KKK}=\roman O(s,\Bobb R),\, \HG = \roman{Sp}(\ell,\Bobb R)$,
$r=\ell$.
Pick
vectors $\bold q_1,\dots, \bold q_\ell$
in $\Bobb R^s$ and consider the linear map 
$\alpha \colon \Bobb R^s \to \Bobb R^{2\ell}$
whose matrix has the
$\bold q_1,\dots, \bold q_\ell$
as first $\ell$ row vectors and whose
other row vectors are zero.
Then 
$\alpha^{\dagger} \colon \Bobb R^{2\ell} \to \Bobb R^s$
is given by the matrix which has the first $\ell$ column vectors
zero and which has,
for $1 \leq j \leq \ell$, the vector $\bold q_j$
as $(\ell + j)$'th column vector.
Then 
$\mu_{\KKK}(\alpha) =\alpha^{\dagger} \alpha = 0$ and
${
\mu_{\HG}(\alpha) =\alpha \alpha^{\dagger} 
= \left[\matrix 0 & Q
\\
0 & 0
\endmatrix\right]
}$
where $Q= [\bold q_j \bold q_k]_{1 \leq j,k\leq \ell}$.
Thus, when we choose 
$\bold q_1,\dots, \bold q_s$
orthonormal and 
$\bold q_{s+1},\dots, \bold q_\ell$
zero
and write $\alpha_0$ for the resulting linear map
from
$\Bobb R^s$  to $\Bobb R^{2\ell}$, the image
$\mu_{\HG}(\alpha_0)$ is a non-negative
nilpotent matrix in $\fra{sp}(\ell,\Bobb R)$
of rank $s$.
\smallskip\noindent
Case (2): $W= \roman{Hom}_{\Bobb C}(\Bobb C^s,\Bobb C^n)\cong\Bobb C^{sn},\ 
\roman{Sp}(W)=\roman{Sp}(sn, \Bobb R),\ {\KKK}=\roman U(s),\ \HG=\roman U(p,q)$,
$p+q=n, p \geq q =r$.
Pick vectors $\bold w_1,\dots, \bold w_q$
in $\Bobb C^s$ and consider the linear map 
$\alpha \colon \Bobb C^s \to \Bobb C^n$
given by
the transpose of the matrix
$[
\bold w_1,
\ldots, 
\bold w_q,
0
\ldots 
0,
\bold w_{p+1}, 
\ldots ,
\bold w_{n}]$
where, for $1 \leq j \leq q$,
$\bold w_{p+j} = i\bold w_j$.
Then 
$\alpha^{\dagger} = i[
\overline{\bold w_1},\dots, \overline{\bold w_q},
0,\dots,0,
-i\overline{\bold w_1},\dots, -i\overline{\bold w_q}]$,
$\mu_{\KKK}(\alpha) =\alpha^{\dagger} \alpha = 0$, and
$$
\mu_{\HG}(\alpha) =\alpha \alpha^{\dagger} 
= \left[\matrix 
iA & 0& -A
\\
0 & 0& 0
\\
-A & 0& -iA
\endmatrix\right]
$$
where 
$A= [\bold w_j \overline{\bold w_k}]_{1 \leq j,k\leq q}$.
Thus, when we choose 
$\bold w_1,\dots, \bold w_s$
orthonormal and 
$\bold w_{s+1},\dots, \bold w_q$
zero
and write $\alpha_0$ for the resulting linear map
from
$\Bobb C^s$  to $\Bobb C^{2\ell}$, the image
$\mu_{\HG}(\alpha_0)$
in
$\fra u(p,q)$
is a non-negative nilpotent matrix in
$\fra u(p,q)$ of rank $s$.
\smallskip
\noindent
Case (3): $W= \roman{Hom}_{\Bobb H}(\Bobb H^s, \Bobb H^n) \cong \Bobb H^{sn},
\roman{Sp}(W) =\roman{Sp}(2sn, \Bobb R),{\KKK}=\roman{Sp}(s),
\HG= \roman O^*(2 n)$.
The reasoning is essentially the same as in the two other cases.
We leave the details to the reader.
\smallskip
Write $\Cal O'_s$
for the $\HG$-orbit
in $\fra \hg$ of $\mu_{\HG}(\alpha_0)$;
by construction, $\Cal O'_s$ 
is contained in the space
$\Cal O_s$ 
of non-negative nilpotent matrices in $\fra {\hg}$
of rank 
$s$ ($\leq r$),
and the 
reasoning given so far implies that the
closure
$\overline{\Cal O'_s}$
lies in the image
$\overline \mu_{\HG}(W^{\roman{red}})$
of
$W^{\roman{red}}$
in $\fra \hg$. 
Thus 
$\overline{\Cal O'_s}\subseteq
\overline \mu_{\HG}(W^{\roman{red}})
\subseteq
\overline{\Cal O_s}$.
We claim that
$\Cal O'_s= \Cal O_s$.
This follows for example from the classification of
holomorphic nilpotent orbits established in (3.3)  above.
However, it may be verified directly,
by means of an argument in the proof of Theorem 5.3.1
of \cite\lermonsj:
Suppose for the moment that $\fra g = \fra{sp}(\ell,\Bobb R)$
($r=\ell$) and let $X$ be a non-negative nilpotent matrix in $\fra {\hg}$
of rank $s$ ($\leq r$).
Then $X$ is of square zero and hence,
when $X$ is viewed as an endomorphism of $V$
($= \Bobb R^{2\ell}$),
the image $\roman{im}(X)$
of $X$ is $\Cal B$-isotropic,
since
$\Cal B(Xu,Xv) = \Cal B(u,X^2v)$.
Furthermore, $\roman{im}(X) \subseteq\roman{ker}(X)$.
Write $V= L_1 \oplus L_2$ as a direct sum of transverse
Lagrangian subspaces $L_1$ and $L_2$,
with 
$\roman{im}(X) \subseteq L_1$, and
pick $A \in \roman{GL}(L_1)$ such that
$B=\left[\matrix A & 0 \\ 0 & (A^t)^{-1}\endmatrix\right] \in 
\roman{Sp}(V)
$
carries
the decomposition
$V= L_1 \oplus L_2$
into the standard decomposition 
$\Bobb R^{2\ell} =\Bobb R^{\ell} \oplus\Bobb R^{\ell}$.
Then $BXB^{-1}$ has the form
${
Y=\left[\matrix 0 & S \\ 0 & 0\endmatrix\right] 
}$
where $S$ is a symmetric $(\ell \times \ell)$-matrix
of rank $s$.
Thus we may assume that $X$ is of the kind $Y$.
Furthermore,
when we choose $B$ suitably,
the matrix $S$ will be a diagonal matrix.
Hence the orbit of $X$ is determined
by the signature of $S$.
The positivity assumption means that
the signature of $S$ is zero.
The matrix $Y$ plainly lies in
$\Cal O'_s$
whence
$\overline{\Cal O_s} \subseteq
\overline{\Cal O'_s}$.
This kind of reasoning may be applied to the other cases as well.
Alternatively,
we may reduce these cases to that
of
$\fra {sp}(\ell,\Bobb R)$
(for suitable $\ell$)
by means of the embeddings
(3.5.1.5) and (3.5.1.7).
A closer look shows that the argument comes down
to Gram-Schmidt orthogonalization of
the complex subspace of $V$ generated by the image of $X$,
with respect to the complex structure $J_V$.
\smallskip
Finally we examine the complex analytic structures.
\smallskip
\noindent
Case (1):
The 
Cartan decomposition 
$\fra {sp}(\ell,\Bobb R) = \fra u(\ell) \oplus \fra p$
of 
$\fra {sp}(\ell,\Bobb R)$
complexifies to
$\fra {sp}(\ell,\Bobb C) = \fra {gl}(\ell,\Bobb C) 
\oplus \fra p^+ \oplus \fra p^-$, where
$\fra p^+=\roman S_{\Bobb C}^2[\Bobb C^{\ell}]$
and $\fra p^-=\overline{\roman S_{\Bobb C}^2[\Bobb C^{\ell}]}$.
Being endowed
with  the complex structure $J=J_V$ on $W\cong\Bobb R^{2s\ell}$,
this space  is just a copy of $\Bobb C^{s\ell}$,
${\KKK}^{\Bobb C} = \roman{O}(s,\Bobb C)$, and
the
${\KKK}^{\Bobb C}$-action on $W_J=\Bobb C^{s\ell}$
amounts to the standard diagonal action of
$\roman{O}(s,\Bobb C)$
on
a product 
$\Bobb C^s\times \dots \times \Bobb C^s$
of $\ell$ copies of $\Bobb C^s$.
By the first main theorem of invariant theory
for the group $\roman{O}(s,\Bobb C)$,
for $s \geq 1$, 
the invariants 
in $\Bobb C[W_J^*]^{{\KKK}^{\Bobb C}}$
are generated by
the homogeneous degree 2 part
$\Bobb C_2[W_J^*]^{{\KKK}^{\Bobb C}}$ thereof,
and the complex algebraic map---the (complex) Hilbert
map---from $W_J=(\Bobb C^s)^\ell$
to the complex vector space
$\roman S_{\Bobb C}^2[\Bobb C^{\ell}]$
of complex symmetric
$(\ell \times \ell)$-matrices
which assigns 
the symmetric matrix
$\left[\bold w_j \bold w_k\right]_{1 \leq j, k \leq \ell}$
to
$\bold w = (\bold w_1,\dots,\bold w_\ell) \in W_J$
identifies
the 
complex algebraic quotient
$W_J\big /\big / \roman O(s,\Bobb C)$  
and hence the reduced space $W^{\roman{red}}$
as a complex analytic space
with the affine complex variety consisting of
matrices in $\roman S_{\Bobb C}^2[\Bobb C^{\ell}]$
of rank at most $s$.
\smallskip
\noindent
Case (2): The Cartan decomposition
$\fra u(p,q) = \fra u(p) \oplus \fra u(q) \oplus \fra p$
of
$\fra u(p,q)$ 
has $\fra p = \roman M_{q,p}(\Bobb C)$
and complexifies to
$
(\fra u(p,q))^{\Bobb C} =\fra {gl}(n,\Bobb C)
= \fra {gl}(p,\Bobb C) \oplus 
\fra {gl}(q,\Bobb C) \oplus
\fra p^+ \oplus \fra p^-$,
where
\linebreak
$\fra p^+=\roman M_{q,p}(\Bobb C)$ and 
$\fra p^-=\overline{\roman M_{q,p}(\Bobb C)}\cong \roman M_{p,q}(\Bobb C)$.
Being endowed
with  the complex structure $J=J_V$ on $W\cong\Bobb C^{sn}$,
this space is  a copy of $\Bobb C^{sn}$,
${\KKK}^{\Bobb C} = \roman{GL}(s,\Bobb C)$, and
the
${\KKK}^{\Bobb C}$-action on $W_J=\Bobb C^{sn}$
amounts to the standard diagonal action of
$\roman{GL}(s,\Bobb C)$
on
a product 
$\Bobb C^{s}\times \dots \times \Bobb C^{s}
\times (\Bobb C^s)^* \times \dots \times (\Bobb C^s)^*$ 
of $p$ copies of $\Bobb C^{s}$
and $q$ copies of
$(\Bobb C^s)^*$.
By invariant theory,
just as before, the invariants 
in $\Bobb C[W_J^*]^{{\KKK}^{\Bobb C}}$
are generated by
the homogeneous degree 2 part
$\Bobb C_2[W_J^*]^{{\KKK}^{\Bobb C}}$ thereof,
and the Hilbert map
from $W_J$ to
$\roman M_{q,p}(\Bobb C)$ 
assigns 
the matrix
$\left[(\xi_j ,\bold w_k)\right]_{1 \leq j \leq q, 1 \leq k \leq p}$
to
$\bold w = (\bold w_1,\dots,\bold w_p,\xi_1,\dots,\xi_q) \in W_J$
and thereby
identifies
the 
complex algebraic quotient
$W_J\big /\big / \roman {GL}(s,\Bobb C)$  
and hence the reduced space $W^{\roman{red}}$
as a complex analytic space
with 
the
complex affine variety in
$\roman M_{q,p}(\Bobb C)$
which consists of matrices
of rank at most
$\roman{min}(s,q)$.
In particular, when $s \geq q$,
there are no relations among the invariants;
indeed, any complex
$(p \times q)$-matrix
may be written as
$[(\xi_j,\bold w_k)]_{1\leq j\leq q, 1 \leq k\leq p}$ 
for suitable vectors
$\bold w_1,\dots, \bold w_p \in \Bobb C^{s}$ 
and 
$\xi_1,\dots, \xi_q \in (\Bobb C^{s})^*$ 
provided 
$s \geq q$.
Consequently,
when $s \geq q$,
the algebraic quotient 
$W_J\big /\big / {\KKK}^{\Bobb C}$
actually coincides with
$\fra p^+=\roman M_{q,p}(\Bobb C)$.

\smallskip\noindent
Case (3): The Cartan decomposition
$\fra {so}^*(2 n) = \fra u(n) \oplus \fra p$
of $\fra {so}^*(2 n)$
has
\linebreak
$\fra p = i\Lambda_{\Bobb R}^2[\Bobb R^n]\oplus i\Lambda_{\Bobb R}^2[\Bobb R^n]$
and complexifies to
$$(\fra {so}^*(2 n))^{\Bobb C} =\fra {so}(2 n,\Bobb C)
= \fra {gl}(n,\Bobb C) \oplus \fra p^+ \oplus \fra p^-$$
where
$\fra p^+=\Lambda_{\Bobb C}^2[\Bobb C^n]$
and
$\fra p^- =\overline{\Lambda_{\Bobb C}^2[\Bobb C^n]}$.
With respect to the complex structure $J=J_V$ on $W\cong \Bobb H^{sn}$,
the latter is a copy of $\Bobb C^{2sn}$,
${\KKK}^{\Bobb C} = \roman{Sp}(s,\Bobb C)$, and
the
${\KKK}^{\Bobb C}$-action on $W_J=\Bobb C^{2sn}$
amounts to the standard diagonal action of
$\roman{Sp}(s,\Bobb C)$
on
a product 
$\Bobb C^{2s}\times \dots \times \Bobb C^{2s}$
of $n$-copies
of
$\Bobb C^{2s}$.
By invariant theory,
just as before, the invariants 
in $\Bobb C[W_J^*]^{{\KKK}^{\Bobb C}}$
are generated by
the homogeneous degree 2 part
$\Bobb C_2[W_J^*]^{{\KKK}^{\Bobb C}}$ thereof,
and from the theory of reductive dual pairs we know that
$\Bobb C_2[W_J^*]^{{\KKK}^{\Bobb C}}$
is a copy of 
$\fra {so}(n,\Bobb C) \cong \Lambda_{\Bobb C}^2[\Bobb C^n]=\fra p^+$.
With reference to the notation $\omega$ for the standard complex
symplectic form on
$\Bobb C^{2s}$,
the Hilbert map
from $W_J$ to
$\Lambda_{\Bobb C}^2[\Bobb C^n]$ 
assigns 
the 
skew-symmetric complex
$(n \times n)$-matrix
$[\omega(\bold w_j,\bold w_k)]_{1\leq j,k\leq n}$
to
$(\bold w_1,\dots,\bold w_n) \in W_J$
and thereby
identifies
the complex algebraic quotient
$W_J\big /\big / {\KKK}^{\Bobb C}$,
and hence the reduced space
$W^{\roman{red}}$ as a complex analytic space,
with a complex affine variety in
$\Lambda_{\Bobb C}^2[\Bobb C^n] =\fra{so}(n,\Bobb C)$.
When $s \geq \rell (= [\frac n2])$,
there are no relations among the invariants;
indeed, any skew-symmetric complex
$(n \times n)$-matrix
may be written as
$[\omega(\bold w_j,\bold w_k)]_{1\leq j,k\leq n}$
for suitable vectors
$\bold w_1,\dots, \bold w_n \in \Bobb C^{2s}$ provided 
$2s \geq n$.
Consequently,
when $s \geq \rell$,
the algebraic quotient 
$W_J\big /\big / {\KKK}^{\Bobb C}$
coincides with $\fra p^+=\Lambda_{\Bobb C}^2[\Bobb C^n]$.
When $2s < n$,
the Hilbert map
from $W_J$ to
$\Lambda_{\Bobb C}^2[\Bobb C^n]$ 
identifies
the complex algebraic quotient
$W_J\big /\big / {\KKK}^{\Bobb C}$ with
the complex affine variety 
of matrices in
$\Lambda_{\Bobb C}^2[\Bobb C^n] =\fra{so}(n,\Bobb C)$
having rank at most $2s$. 
\smallskip
Thus, in each of the three standard cases,
when $1 \leq s \leq r$,
under the
projection from $\fra g = \fra k \oplus \fra p$
to $\fra p$, followed by the identification of the latter
with $\fra p^+$,
the (image of the) reduced space
$W^{\roman{red}}$ in $\fra g$
is identified
with the affine complex variety consisting of
matrices in $\fra p^+$ 
of rank at most $s$
for $\fra g = \fra{sp}(\ell,\Bobb R)$
and
$\fra g = \fra {su}(p,q)$
and of rank at most $2s$
for $\fra g = \fra{so}^*(2n)$;
in particular, when
$s =r$,
the (image of the) reduced space
$W^{\roman{red}}$ in $\fra g$,
that is, (that of) the categorical quotient
$W_J\big /\big /{\KKK}^{\Bobb C}$ 
is identified with 
$\fra p^+$.  
Furthermore,
the complex structure 
on $\overline {\Cal O_s}$
which arises from the algebraic quotient
construction gets identified with that resulting
from the projection of 
$\overline {\Cal O_s}$
to $\fra p$.
Consequently,
for $1 \leq s \leq r$,
the restriction 
to the closure
$\overline{\Cal O'_s}$
of the nilpotent ($\HG$-)orbit $\Cal O'_s$
of
the projection from $\fra g = \fra k \oplus \fra p$
to $\fra p$
is injective, whence
$\Cal O'_s$ is a holomorphic nilpotent orbit. \qed
\enddemo

We mention in passing that
the reasoning 
aimed at unravelling the complex analytic structures
of the categorical quotients written
as $W_J\big /\big / {\KKK}^{\Bobb C}$,
given
in the proof of Theorem 5.3.3,
establishes the claim of Theorem 3.5.5 above as well.

\noindent
{\smc 5.4. The real semi-algebraic structure.}
The proof of Theorem 5.3.3 actually yields more explicit information,
which we now spell out.
To this end,
we write $\rho$ for the (ordinary) rank of $\fra {\hg}$.
Recall that
$\rho(\fra{sp}(\ell,\Bobb R)) = \ell$,
$\rho(\fra u(p,q)) = p+q=n$,
$\rho(\fra {so}^*(2n)) =n$.
Let $\beta_1,\dots,\beta_{\rho}$
be a set of generators for the $\HG$-invariants
on $\fra {\hg}$ and view them as polynomial functions
on $\fra {\hg}$ as usual
by means of the 
half-trace pairing (the appropriate multiple of the Killing form).
Let
$\beta =(\beta_1,\dots,\beta_{\rho})\colon \fra {\hg} \to \Bobb R^{\rho}$.

\proclaim{Theorem 5.4.1}
In the standard cases,
for $1 \leq s \leq r$, 
the real algebraic quotient $W \big / \big / {\KKK}$
is the space
of matrices 
in $\fra {\hg}$
which have rank at most $s$
(real rank $s$ for $\fra {\hg} = \fra{sp}(\ell,\Bobb R)$, $\ell = r$;
complex rank $s$ for $\fra {\hg} = \fra {su}(p,q)$, $p \geq q = r$;
quaternionic rank $s$ for $\fra {\hg} = \fra{so}^*(2n)$, $r= [\frac n2]$),
that is, 
$W \big / \big / {\KKK}$
is a real determinantal variety in $\fra {\hg}$.
Furthermore, the image
$p(W)$ of $W$ under the projection
from $W$ to
$W \big / \big / {\KKK}$
sits inside 
$W \big / \big / {\KKK}$
as the space of non-negative matrices 
in $W \big / \big / {\KKK}$,
and
the reduced space
$W^{\roman{red}}$
sits inside
$p(W)$
as the zero locus  of 
$\beta$, that is, of
the $\HG$-invariants on $\fra \hg$,
viewed as functions 
on $W \big / \big / {\KKK}$
(and hence on $p(W)$).
Finally, for $s \geq 1$, the Zariski tangent space
$\roman T_0(\overline {\Cal O_s})$
at the origin
(for the smooth structure
$C^{\infty} (\overline {\Cal O_s})$)
coincides with the entire ambient space
$\fra {\hg}$
whence $\overline {\Cal O_s}$ cannot be realized in any proper linear subspace
of $\fra {\hg}$.
In particular,
the stratified symplectic Poisson structure on
each $\overline {\Cal O_s}$ necessarily involves {\rm all}
the elements of a basis of
$\fra {\hg}$,
i.~e. cannot be described by fewer generators than the
dimension of $\fra {\hg}$.
\endproclaim

\demo{Proof}
Since the momentum mapping 
$\mu_{\HG} \colon W \to \fra {\hg}$
amounts to the
Hilbert map for the ${\KKK}$-action on $W$, 
$\mu_{\HG}$ induces an
injection 
of the real algebraic quotient $W \big / \big / {\KKK}$ into
$\fra {\hg}$
which identifies 
$W \big / \big / {\KKK}$
with the 
subspace of 
$\fra {\hg}$
of matrices of rank at most $s$
(real rank $s$ for $\fra {\hg} = \fra{sp}(\ell,\Bobb R)$;
complex rank $s$ for $\fra {\hg} = \fra {su}(p,q)$;
quaternionic rank $s$ for $\fra {\hg} = \fra{so}^*(2n)$).
\smallskip
The image $p(W)$ of $W$ in
$W \big / \big / {\KKK}$ 
under the projection
$
p \colon W 
@>>>
W \big / \big / {\KKK}
\subseteq
\fra{\hg}
$
is the
space
$W \big / {\KKK}$
of ${\KKK}$-orbits
in 
$W \big / \big / {\KKK}$;
it 
is the real semi-algebraic set
of 
non-negative 
matrices
in $\fra {\hg}$
of rank at most $s$ ($\leq \ell$).
Recall that a matrix $X$
in $\fra {\hg}$ is non-negative if the hermitian form
$-\Cal B(X \cdot,\cdot) = -(J_VX\cdot,\cdot)$
is non-negative or, equivalently,
if the hermitian matrix
$-J_VX$ is non-negative.
The property of 
$-J_VX$ being non-negative
may be expressed in terms of suitable inequalities
involving the appropriate minors.
Thus, as a real semi-algebraic subset of
$\fra {\hg}$,
the image $p(W)$ is given by
the equations
which say that all
$(s+1) \times (s+1)$-minors vanish
and the 
inequalities
saying that
a matrix $X$ in $\fra {\hg}$ 
of rank at most $s$ is non-negative.
\smallskip
Since the nilcone in $\fra {\hg}$ is the zero locus of the
$\HG$-invariants $\beta_1,\dots,\beta_{\rho}$ in $\fra {\hg}$,
viewed as functions on $\fra {\hg}$,
the zero locus of these invariants,
viewed as functions on
$W \big / \big / {\KKK}$,
is the space of nilpotent matrices in $\fra {\hg}$
of rank $\leq s$.
Hence
the zero locus of the invariants,
viewed as functions on
$p(W) \subseteq W \big / \big / {\KKK}$,
is the space  of non-negative nilpotent matrices in $\fra {\hg}$
of rank $\leq s$.
Thus, in view of Theorem 5.3.3, the reduced space
$W^{\roman{red}}$
or, equivalently,
the 
closure of
the holomorphic nilpotent orbit
$\Cal O_s$
which consists
of 
non-negative nilpotent matrices in $\fra {\hg}$
of rank $\leq s$,
sits inside
$p(W)$
as the zero locus  of 
$\beta$, that is, of
the $\HG$-invariants on $\fra \hg$,
viewed as functions 
on $W \big / \big / {\KKK}$
(and hence on $p(W)$). 
\smallskip
To justify the claim
involving the Zariski tangent spaces at the origin,
it will suffice to prove that
the Zariski tangent space
$\roman T_0(\overline {\Cal O_1})$
at the origin
(for the smooth structure
$C^{\infty} (\overline {\Cal O_1})$)
coincides with the entire ambient space
$\fra {\hg}$.
In case (1),
for $s=1$,
$H=\roman O(1,\Bobb R)$,
the $H$-momentum mapping is trivial, and
the $G$-momentum mapping
$\mu_G \colon (\Bobb R^2)^{\ell} \to \fra{sp}(\ell,\Bobb R)$
sends the vector
$[q_1,p_1,\ldots,q_{\ell},p_{\ell}]^{\tau}$
(where $q_j, p_j \in \Bobb R$)
to the matrix
$$
\left[\matrix 
\left [q_j  p_k\right ] &  -\left [q_j  q_k\right ]\\
\left [p_j  p_k\right ] &  -\left [p_j  q_k\right ]
\endmatrix
\right]
\in \fra {sp}(\ell,\Bobb R);
$$
it is thus straightforward to find suitable
vectors
$[q_1,p_1,\ldots,q_{\ell},p_{\ell}]^{\tau}$
whose images in 
$\fra {sp}(\ell,\Bobb R)$
constitute a basis thereof.
Likewise,
in case (2),
for $s=1$,
$H=\roman U(1)$,
the $H$-momentum mapping 
$\mu_H \colon \Bobb C^{p+q} \to \fra u(1)$
sends the vector
$\left[w_1,\ldots, w_{p+q}\right]^{\tau}$
(where $w_1,\ldots, w_{p+q} \in \Bobb C$)
to
$i(w_1 \overline {w_1} + \dots +w_p \overline {w_p}
-w_{p+1} \overline {w_{p+1}} - \dots -w_{p+q} \overline {w_{p+q}})
\in \fra u(1)$,
and the $G$-momentum mapping
$\mu_G \colon \Bobb C^{p+q} \to \fra u(p,q)$
sends 
$\left[w_1,\ldots, w_{p+q}\right]^{\tau}$
to
$$
\left[\matrix A   & B\\
              B^* & D
\endmatrix
\right]
\in \fra u(p,q)
$$
where
$A=i\left[w_j\overline w_k\right]_{1\leq j,k\leq p}$,
$B=-i\left[w_j\overline w_{p+k}\right]_{1\leq j \leq p,1\leq k \leq q}$,
$D=-i\left[w_{j+p}\overline w_{k+p}\right]_{1\leq j,k\leq q}$.
In particular,
when $\mu_H(\left[w_1,\ldots, w_{p+q}\right]^{\tau}) = 0$,
the value $\mu_G(\left[w_1,\ldots, w_{p+q}\right]^{\tau})$
lies in $\fra {su}(p,q)$.
A little thought reveals that,
for suitable choices of
$w_1,\ldots, w_{p+q}$ with 
$\mu_H(\left[w_1,\ldots, w_{p+q}\right]^{\tau}) = 0$,
the images
$\mu_G(\left[w_1,\ldots, w_{p+q}\right]^{\tau}) \in  \fra {su}(p,q)$
constitute a basis.
For example,
for $\fra {su}(2,1)$,
the vectors
$[1,0,1]^{\tau}$, $[0,1,1]^{\tau}$,
$[1,0,i]^{\tau}$, $[0,1,i]^{\tau}$,
$[i,0,1]^{\tau}$, $[0,i,1]^{\tau}$,
$[1,1,\sqrt 2]^{\tau}$, $[i,1,\sqrt 2]^{\tau}$ 
will have the asserted properties.
The same kind of reasoning establishes the claim in case (3).
We leave the details to the reader. \qed 
\enddemo

\smallskip
\noindent
{\smc Remark 5.4.2.}
Theorem 5.4.1 yields, for $1 \leq s \leq r$,
explicit equations and inequalities
for the closure $\overline {\Cal O_s}$,
realized as a semi-algebraic set in $\fra {\hg}$.
Indeed, $\overline {\Cal O_s}$ consists of all matrices
$\alpha \in \fra g$
(taken in its matrix realization spelled out in (3.5))
which satisfy
\roster
\item
the real determinantal equations
which say that
all 
$((s+1) \times (s+1))$-minors 
of $\alpha$
vanish;
\item
the equations
$\beta_1 (\alpha)=0,\dots, \beta_\rho(\alpha)=0$
(which say that $\alpha$ lies in the real nilvariety of $\fra g$),
and  
\item the inequalities which say that
the matrix $\alpha$ is non-negative.
\endroster

\smallskip
\noindent
{\smc Remark 5.4.3.}
Theorems 5.3.3 and 5.4.1 show once more that,
for a standard simple Lie algebra of hermitian type
$(\fra g,z)$ of real rank $r$,
even though 
as a complex analytic space, the closure $\overline {\Cal O_r}$
of the top stratum
is affine and hence non-singular,
as a stratified symplectic space,
it has {\it singularities\/}, the {\it singular locus\/}
(in the sense of stratified symplectic spaces)
being the complement $\overline {\Cal O_{r-1}}$
of the top stratum
$\Cal O_r$.
Furthermore,
Theorem 5.3.3
entails at once
the fact, implied as well by Corollary 3.3.8 above
that, for a holomorphic nilpotent orbit
$\Cal O$,
the diffeomorphism 
from $\Cal O$
onto its image 
in $\fra p$
extends to a homeomorphism
from the closure
$\overline{\Cal O}$
onto its image 
in
$\fra p$.

\proclaim {Corollary 5.4.4}
Let ${\KKK}$ be a compact Lie group which decomposes into a direct product
of copies of $\roman O(s_j,\Bobb R)$,
$\roman {Sp}(s_k)$, 
$\roman U(s_m)$, 
and let $E$ be a representation of ${\KKK}$ which decomposes
into a product of standard representations
and of the conjugate representation of the standard
one for the case of $\roman U(s_m)$. 
Then the reduced space $E^{\roman{red}}$
is as a normal K\"ahler space
isomorphic to a holomorphic nilpotent orbit
in a real reductive Lie algebra $(\fra g,z)$
of hermitian type
which decomposes into a sum of copies
of $\fra{sp}(\ell_j,\Bobb R)$,
$\fra{so}^*(2\ell_k)$,
$\fra u(p_m, q_m)$,
for suitable $\ell_j, \ell_k, p_m, q_m$.
Furthermore,
given a real reductive Lie algebra 
of hermitian type
which decomposes into a sum of Lie algebras
of the standard kind,
every holomorphic nilpotent
orbit 
in such a Lie algebra
arises in this fashion.
\endproclaim

\demo{Proof} This is an immediate consequence of Theorem 5.4.1. \qed
\enddemo

\noindent
{\smc 5.5. Existence of associated representations
for {\it arbitrary\/} pseudoholomorphic nilpotent orbits.}
Return to the
real symplectic vector space $W=W(s)=\roman{Hom}_{\Bobb K}(V^s,V)$,
endowed with
the symplectic actions
of the groups $\KK$ and $\HG$ 
and corresponding momentum mappings
where 
$\KK$ is no longer assumed to be compact.
As before, write $(s',s'')$
for the type of the hermitian form
$(\cdot,\cdot)$
on $W$, so that $s'+s''$ is its rank and
$s'-s''$ its signature.

\proclaim{Theorem 5.5.1}
The induced map $\overline \mu_{\HG}$
from the $\KK$-reduced space
$W(s)^{\roman{red}}$ 
to
$\fra {\hg}$  
is an 
embedding 
of
$W(s)^{\roman{red}}$ into $\fra {\hg}$
and
induces an isomorphism of 
stratified  spaces
from $W(s)^{\roman{red}}$ onto the 
union $N$ of all
pseudoholomorphic 
nilpotent orbits 
$\Cal O_{(t,u)}$
in $\fra {\hg}$
with $t \leq \min (r,s')$,
$u \leq \min (r,s'')$
(and, of course, $t+u \leq r$).
Furthermore,
the map from $W(s)^{\roman{red}}$ to $N$ is a Poisson map.
\endproclaim

\demo{Proof}
This is proved in much the same way as
Theorem 5.3.3. We leave the details to the reader. \qed
\enddemo

Thus when $s=2r$ and
the hermitian form $(\cdot,\cdot)$
on $V^s$ has type $(s',s'') = (r,r)$,
the reduced space
$W(s)^{\roman{red}}$
recovers {\it all\/} pseudoholomorphic nilpotent orbits
in $\fra {\hg}$.

\smallskip\noindent
{\smc Remark 5.5.2.}
By symmetry,
the above construction identifies
the $G$-reduced space
with a union of square zero nilpotent orbits
in $\fra{\kk}$, and every 
square zero nilpotent orbit
orbit 
in $\fra{\kk}$ arises in this way.
This yields the square zero nilpotent orbits
for
the classical Lie algebras
$\fra {so}(s',s'')$,
$\fra u(s',s'')$,
and $\fra {sp}(s',s'')$.
The case
$\fra u(s',s'')$
does not provide anything new, though.
The case
of $\fra {so}(2,q)$
will be examined in Section 6 below.

\smallskip\noindent
{\smc Remark 5.5.3.}
Earlier 
we made
the 
observation
that,
for a pseudoholomorphic nilpotent orbit $\Cal O$
which is neither holomorphic nor antiholomorphic,
the diffeomorphism from
$\Cal O$ onto  its image in $\fra p$ under the projection
from $\fra g$ to $\fra p$
does not extend to a homeomorphism.
This observation translates
to the fact that
the obvious map from
the reduced space
$W(s)^{\roman{red}}$
to the corresponding
categorical quotient
is {\it not\/}
a homeomorphism.
{\sl Thus the Kempf-Ness observation 
(relating the symplectic quotient for a hamiltonain action of a 
compact group with the categorical quotient for the complexified group)
does
not in general extend to the case
where a 
complex group acts
holomorphically
with hamiltonian
action 
of
the restriction to a
non-compact real form.\/}

\proclaim {Corollary 5.5.4}
Let $\KK$ be a Lie group which decomposes into a direct product
of copies of $\roman O(s'_j,s''_j)$,
$\roman {Sp}(s'_k,s''_k)$, 
$\roman U(s'_m,s''_m)$, 
and let $E$ be a representation of $\KK$ which decomposes
into a product of standard representations
and of the conjugate representation of the standard
one for the case of $\roman U(s'_m,s''_m)$. 
Then the reduced space $E^{\roman{red}}$
is as a stratified space
isomorphic to a union $N$ of pseudoholomorphic nilpotent orbits
in a real reductive Lie algebra $(\fra g,z)$
of hermitian type
which decomposes into a sum of copies
of $\fra{sp}(\ell_j,\Bobb R)$,
$\fra{so}^*(2\ell_k)$,
$\fra u(p_m, q_m)$,
for suitable $\ell_j, \ell_k, p_m, q_m$.
Furthermore,
the map from 
$E^{\roman{red}}$
to $N$ is a Poisson map,
and every pseudoholomorphic nilpotent
orbit 
in a real reductive Lie algebra 
of hermitian type
which decomposes into a sum of Lie algebras
of the standard kind
arises in this fashion
as a constituent of a reduced space. 
\endproclaim

\demo{Proof} This is again an immediate consequence of Theorem 5.4.1. \qed
\enddemo

\smallskip\noindent
{\smc Example 5.5.5.} Let 
$V^1$ be the standard representation
of $\KK = \roman O(1,1)$ 
on $\Bobb R^2$
and
$V$ the standard representation of
$G=\roman{Sp}(1,\Bobb R)$
on $\Bobb R^2$.
Then $W= \roman{Hom}_{\Bobb R}(V^1,V)$
has dimension four,
as an $\KK$-representation, $W$ amounts to the sum
$V^1 \oplus V^1$ of two copies of the standard representation,
the $\KK$-momentum 
mapping sends $\alpha= (\bold q,\bold p) \in \Bobb R^2 \times \Bobb R^2$
to $\bold q \wedge \bold p\,\left[\matrix 0 & 1 \\ 1 & 0 \endmatrix \right]$,
and
the $G$-momentum mapping
$\mu_G \colon W \to \fra {sp}(1,\Bobb R)$
sends $\alpha= (\bold q,\bold p) \in \Bobb R^2 \times \Bobb R^2$
to 
$$
\alpha \alpha^{\dagger}=
\left[\matrix (\bold q,\bold p) & -(\bold q,\bold q) 
\\ (\bold p,\bold p) & -(\bold q,\bold p) \endmatrix \right]
\in \fra{sp}(1,\Bobb R).
$$
The zero locus $\mu_{\KK}^{-1}(0)$
consists of the points $(\bold q,\lambda \bold q)$ 
and
$(\lambda\bold q,\bold q)$
($\bold q =(q_1,q_2)\in \Bobb R^2$,\linebreak $\lambda \in \Bobb R$),
and the non-closed ${\KK}$-orbits in
$\mu_{\KK}^{-1}(0)$
are those which correspond to 
non-zero 
isotropic elements $\bold q$, that is, to non-zero
$\bold q$'s with 
$(\bold q,\bold q) = q_1^2 - q_2^2= 0$.
The $\KK$-reduced space
is the {\it whole\/} nilcone in
$\fra{sp}(1,\Bobb R)$,
that is, the closure of the holomorphic {\it and\/} antiholomorphic
2-dimensional nilpotent orbits
in 
$\fra{sp}(1,\Bobb R)$.

\smallskip\noindent
{\smc 5.6. The semisimple holomorphic orbits.} 
Maintaining the notation established earlier,
let $s= 2 \ell, p+q, n$,
and  $s'=s, s'' = 0$,
so that, according to the case considered: 
\newline\noindent
(1) $W= \roman{End}_{\Bobb R}(\Bobb R^{2\ell})$,
$(\HG,\GHH,K) = \left(\roman{Sp}(\ell,\Bobb R), \roman{O}(2\ell,\Bobb R),
\roman{U}(\ell)\right)$
\newline\noindent
(2) $W= \roman{End}_{\Bobb C}(\Bobb C^{p,q})$,
$(\HG,\GHH,K) = \left(\roman{U}(p,q), \roman{U}(p+q),
\roman{U}(p)\times \roman{U}(q)\right)$
\newline\noindent
(3) 
$W= \roman{End}_{\Bobb H}(\Bobb H^{n})$,
$(\HG,\GHH,K) = \left(\roman{O}^*(2n), \roman{Sp}(n),
\roman{U}(n)\right)$.
\newline\noindent
In particular,
the group $\GHH$ is compact.
Then 
$J_V$ lies in $\fra {\kkk}$,
the pre-image
$\mu_{\KKK}^{-1}(J_V)$ may be identified with the group 
$G$---in fact, this is the very definition of $G$ as the group of
$\Cal B$-isometries of $V$, cf. (3.5.1) above---, the 
stabilizer of $J_V$ in ${\KKK}$
may be identified with the maximal compact subgroup $K$ of $G$
determined by the choice of $H$-element,
and the reduced space
$\mu_{\KKK}^{-1}(J_V)\big / K$
amounts to the homogeneous space $G/K$.
Thus the
$\HG$-momentum mapping $\mu_{\HG}\colon W \to  \fra {\hg}$
($\fra {\hg}$ being identified with its dual as explained earlier)
identifies
the reduced space
$\mu_{\KKK}^{-1}(J_V)\big / K\cong G/K$
with the semisimple holomorphic orbit $G(2z)$ of $2z$ in $\fra {\hg}$.
More generally,
for $\varepsilon >0$,
the 
$\HG$-momentum mapping $\mu_{\HG}\colon W \to  \fra {\hg}$
identifies
the reduced space
$\mu_{\KKK}^{-1}(\varepsilon J_V)\big / K$
with the semisimple holomorphic orbit $G(2\varepsilon z)$
of $2 \varepsilon z$
in $\fra {\hg}$,
and this orbit is still canonically diffeomorphic to $G/K$.

\beginsection 6. Associated representations for the remaining classical  case

While the construction given in Section 5
does not yield the holomorphic nilpotent orbits in $\fra{so}(2,q)$
by symplectic reduction with respect to a {\it compact group\/},
these orbits may be obtained by symplectic reduction with respect to  
the non-compact group $\roman{Sp}(1,\Bobb R)$,
as already hinted at in Remark 5.9 above.
\smallskip\noindent
{\smc 6.1. The real structure.}
Let $V= \Bobb R^n$ ($n \geq 2$), endowed
with the standard Lorentz form 
$(\cdot,\cdot)$
of signature $(p,q)$ ($p+q \geq 2$)---we indicate this with the notation 
$\Bobb R^{p,q}$, so that
$\roman U(V,(\cdot,\cdot)) = \roman O(p,q)$---and, 
for $s \geq 1$, consider
$\Bobb R^{2s}$, endowed with
its standard symplectic structure
which,
for consistency of exposition,
we write $\Cal B(\cdot,\cdot)$,
so that 
$\roman U(\Bobb R^{2s},\Cal B) = \roman {Sp}(s,\Bobb R)$;
the notation $\roman U(\cdot,\cdot)$ was introduced in (3.5) above.
Let $W=W(s)=\roman{Hom}_{\Bobb R}(\Bobb R^{2s},V)$.
As before,
the groups $\GH=\roman {Sp}(s,\Bobb R)$
and 
$\HG=\roman O(p,q)$
act on $W$ 
by means of the assignment to
$(x,y,\alpha)$ of $y \alpha x^{-1}$
where
$x \in \roman {Sp}(s,\Bobb R), \ \alpha \in \roman{Hom}_{\Bobb K}(V^s,V),\ 
y \in \roman O(p,q)$.
Given an $\Bobb R$-linear map
$\alpha \colon \Bobb R^{2s} \to \Bobb R^{p,q}$, 
define
$\alpha^\dagger \colon \Bobb R^{p,q} \to \Bobb R^{2s}$ by
$$
(\alpha \bold u,\bold v) = \Cal B(\bold u,\alpha^\dagger \bold v),
\ \bold u \in \Bobb R^{2s},
\ \bold v \in \Bobb R^{p,q};
$$
we note that here $\Cal B$ and $(\cdot,\cdot)$
are defined on the domain and range,
respectively, of $\alpha$
whereas in Section 5 above
$\Cal B$ and $(\cdot,\cdot)$
are defined on the  range and domain, 
respectively, of the $\alpha$
coming into play there.
The forms $\Cal B$ 
and $(\cdot,\cdot)$ induce a symplectic structure
on $W\cong \Bobb R^{2sn}$, 
and the actions of $\GH$ and $\HG$ on $W$ are hamiltonian,
having momentum mappings
$$
\align
\mu_{\GH} \colon W &@>>> \fra {\gh},
\quad
\mu_{\GH}(\alpha)= -\alpha^\dagger \alpha \colon \Bobb R^{2s} \to \Bobb R^{2s},
\\
\mu_{\HG} \colon W &@>>> \fra {\hg},
\quad
\mu_{\HG}(\alpha) =
 \alpha \alpha^\dagger \colon \Bobb R^{p,q} \to \Bobb R^{p,q}
\endalign
$$
respectively where, as before,
$\fra {\gh}$ and $\fra {\hg}$
are identified with their duals by means of the half-trace pairings.
The groups $\roman {Sp}(s,\Bobb R)$
and $\roman O(p,q)$ constitute a dual pair
in 
$\roman {Sp}(W,\Bobb R)\cong\roman {Sp}(sn,\Bobb R)$
\cite\howetwo\ (5.2(i)) and,
up to sign,
the $\GH$-momentum mapping
$\mu_{\GH}$ 
amounts to the Hilbert map for the $\GH$-action on $W$,
that is to say,
it passes to an embedding
of the real algebraic quotient
$W \big / \big / \GH$ 
into 
$\fra {so}(p,q)$.
As before, let $W^{\roman{red}}$
be the space of {\it closed\/} $H$-orbits
in $\mu_H^{-1}(0)$,
viewed as a subspace of
$W \big / \big / \GH$.
\smallskip
We will now unravel the real orbit structure which
is relevant for our purposes:
Let $s=1,\  p=2,\ q > 2$;
since 
$\fra {so}(2,2) \cong \fra {so}(2,1)\oplus \fra {so}(2,1)$,
$\fra {so}(2,3) \cong \fra {sp}(2,\Bobb R)$,
$\fra {so}(2,4) \cong \fra {su}(2,2)$,
the first truly interesting case actually occurs for $q=5$.
Given a linear map 
$$
\alpha=[\bold w_1,\dots,\bold w_n]^{\tau}
\colon \Bobb R^2 \to \Bobb R^{2,q},
$$
with the notation
$
\bold w_1=(w_1^1,w_1^2),\bold w_2=(w_2^1,w_2^2),
\bold w_3=(w_3^1,w_3^2),
\dots,
\bold w_n=(w_n^1,w_n^2)
$
(where $(\cdot,\cdot)$ refers to a 2-component vector and {\it not\/}
to the Lorentz form),
the adjoint $\alpha^{\dagger} \colon
 \Bobb R^{2,q} \to \Bobb R^2$, or rather its negative,
is given by the matrix
$$
-\alpha^{\dagger} =
\left[
\matrix
w_1^2
&
w_2^2
&
-w_3^2
&
\ldots
&
-w_n^2
\\
-w_1^1
&
-w_2^1
&
w_3^1
&
\ldots
&
w_n^1
\endmatrix
\right].
$$
Hence, with the notation $\bold w^1 = (w_1^1,\dots,w_n^1)$ etc.,
$$
\mu_{\GH}(\alpha) 
=-\alpha^{\dagger} \alpha
= \left[\matrix (\bold w^1,\bold w^2) & (\bold w^2,\bold w^2)\\
-(\bold w^1,\bold w^1) & -(\bold w^1,\bold w^2)\endmatrix \right] 
\in \fra{sp}(1,\Bobb R).
$$
Under these circumstances,
$\mu_H$ induces a continuous 
$\roman{SO}(2,q)$-equivariant
bijection
from
$W^{\roman{red}}$
onto the closure of the
$\roman{SO}(2,q)$-orbit $\Cal O_{1,0}\cup\Cal O_{0,1}$
which coincides with the subspace of square zero matrices
in $\fra {so}(2,q)$.
As
$\roman{SO}(2,q)^0$-orbits,
$\Cal O_{1,0}$ is holomorphic and $\Cal O_{0,1}$
antiholomorphic.
For completeness, we note that
$e_1=\mu_H(\alpha_1)$ and $e_2=\mu_H(\alpha_2)$ where
$\alpha_{1,2} \colon \Bobb R^2 @>>> \Bobb R^n$
are given by
$$
\alpha_1 
= 
\left [\matrix
 1& 0& 1& 0& \ldots&0
\\
 0& 1& 0& -1& \ldots&0
\endmatrix \right]^{\tau},
\quad
\alpha_2 
=
\left [\matrix
 0&-1& 1& 0& \ldots&0
\\
 1& 0& 0& 1& \ldots&0
\endmatrix \right]^{\tau}.
$$

The remaining pseudoholomorphic nilpotent orbits 
$\Cal O_{2,0}$, $\Cal O_{0,2}$, and $\Cal O_{1,1}$
of
$\fra{so}(2,q)$ 
are the connected components
of the
$\roman{Sp}(1,\Bobb R)$-reduced space
$W_C$, the space of closed
$\roman{Sp}(1,\Bobb R)$-orbits
in the pre-image $\mu_{\GH}^{-1}(C)$
of
the nilcone $C$ in 
$\fra{sp}(1,\Bobb R)$.
Indeed,
the nilcone $C$ 
is the disjoint union
of the three nilpotent orbits
$\Cal O_0=\{0\},\Cal O_0^+,\Cal O_0^-$,
the notation being that introduced in (3.2.2).
Inspection shows that any matrix in $C$
has square zero;
this is also a very special case of Theorem 3.5.3
since
$\Cal O_0^+$ and $\Cal O_0^-$
are pseudoholomorphic.
Thus, given $\alpha \in \mu_{\GH}^{-1}(C)$,
$$
(\mu_G(\alpha))^3 =
(\alpha\alpha^{\dagger})^3 =
\alpha(\alpha^{\dagger}\alpha)^2\alpha^{\dagger}
=\alpha (\mu_H(\alpha))^2 \alpha^{\dagger} = 0
$$
whence
the momentum mapping $\mu_G\colon W \to \fra{so}(2,q)$
maps the pre-image
$\mu_{\GH}^{-1}(C)$
into the space of
matrices $X$ in 
$\fra{so}(2,q)$
with $X^3=0$.
The classification of nilpotent orbits
in 
$\fra{so}(2,q)$
\cite\burcush\ 
implies
that
this map is onto
the space of
matrices $X$ with $X^3=0$,
in fact,
induces
a continuous bijection from
$W_C$ onto 
$\Cal O_{2,0}\cup \Cal O_{1,1}\cup \Cal O_{0,2}$.
A closer look reveals that
$$
\mu_G(\mu_{\GH}^{-1}(\Cal O_0^+)) = \Cal O_{2,0}\cup \Cal O_{0,2},
\quad
\mu_G(\mu_{\GH}^{-1}(\Cal O_0^-)) = \Cal O_{1,1},
$$
where the notation is that in (3.6).
From this observation, it is not hard to deduce that
the orbits
$\Cal O_{2,0}$, $\Cal O_{0,2}$, and $\Cal O_{1,1}$
are pseudoholomorphic.
For example,
since
we may change $\alpha$ within its
$\roman {Sp}(1,\Bobb R)$-orbit without
changing
$\alpha \alpha^{\dagger} \in \fra{so}(2,q)$,
to study
$\Cal O_{2,0}\cup \Cal O_{0,2}$,
we may  
consider
the smooth submanifold
$\widetilde {\Cal O_0^+}$
of $W= \Bobb R^{2n}$
which consists of those 
$\alpha =[\bold w_1,\bold w_2,\dots,\bold w_n]^{\tau}
\colon \Bobb R^2 \to \Bobb R^n$ 
which satisfy the equation
$$
\mu_{\GH}(\alpha) = 
-\alpha^{\dagger} \alpha
= \left[\matrix (\bold w^1,\bold w^2) & (\bold w^2,\bold w^2)\\
-(\bold w^1,\bold w^1) & -(\bold w^1,\bold w^2)\endmatrix \right] 
=
E=\left[\matrix 0& 1\\ 0 & 0 \endmatrix \right].
$$
The stabilizer $Z_E$ 
of $E \in \fra{sp}(1,\Bobb R) =\fra{sl}(2,\Bobb R)$ 
is a one-dimensional unipotent subgroup
in 
$\roman {Sp}(1,\Bobb R)$;
it acts freely on
$\widetilde {\Cal O_0^+}$
in such a way that
the space of
$Z_E$-orbits
coincides with $\Cal O_{2,0}\cup \Cal O_{0,2}$.
We claim that
the $\roman{SO}(q,\Bobb R)$-equivariant map
$$
\pi \colon W= \Bobb R^{2n} @>{\mu_G}>> 
\fra{so}(2,q)=\fra g = \fra k \oplus \fra p 
@>{\roman {proj}}>> \fra p =\Bobb R^{2q},
$$
restricted 
to $\widetilde {\Cal O_0^+}$,
is a fibration onto its image in
$\fra p =\Bobb R^{2q}$
with fiber the real line.
Indeed,
by construction,
$$
\pi [\bold w_1,  \bold w_2,\dots, \bold w_n]^{\tau}
=\left[\matrix 
\bold w_3 \wedge \bold w_1 & \ldots \ \bold w_n \wedge \bold w_1\\
\bold w_3 \wedge \bold w_2 & \ldots \ \bold w_n \wedge \bold w_2\endmatrix
\right]
$$
Plainly,
$\roman{SO}(q,\Bobb R)$-equivariance
reduces the study of $\pi$ to that of
the special case where $q = 3$ since,
for any 
$\alpha =[\bold w_1,  \bold w_2,\dots, \bold w_n]^{\tau}$,
there is some $y \in \roman{SO}(q,\Bobb R)$
such that
$$
y[\bold w_1, \bold w_2, \bold w'_3,\bold w'_4,\bold w'_5,0,\dots, 0]^{\tau}
=
[\bold w_1,  \bold w_2,\dots, \bold w_n]^{\tau}
$$
for some 
$\bold w'_3,\bold w'_4,\bold w'_5 \in \Bobb R^2$.
However,
$\fra {so}(2,3) \cong \fra {sp}(2,\Bobb R)$ and,
from the description  of 
holomorphic nilpotent orbits in
$\fra {sp}(2,\Bobb R)$
given in Section 3 above
we know already that, for $q=3$,
the restriction of the map $\pi$
to $\widetilde {\Cal O_0^+}$
is a fibration of
$\widetilde {\Cal O_0^+}$
onto its image in
$\fra p =\Bobb R^6$
with fiber the real line.
By $\roman{SO}(q,\Bobb R)$-equivariance,
this implies that, for general $q \geq 3$,
the restriction of $\pi$ 
to $\widetilde {\Cal O_0^+}$
is a fibration of
$\widetilde {\Cal O_0^+}$
onto its image in
$\fra p =\Bobb R^{2q}$
with fiber the real line.
Consequently
the restriction 
to $\Cal O_{2,0}$
of
the projection
from $\fra g = \fra k \oplus \fra p$
to $\fra p$
is a diffeomorphism onto its image,
and the same is true of
$\Cal O_{0,2}$.
The same kind of reasoning
applies
to $\Cal O_{1,1}$. 
Thus the orbits $\Cal O_{2,0}$,
$\Cal O_{1,1}$,
and $\Cal O_{0,2}$
are pseudoholomorphic.
Henceforth we write
$\Cal O_1 =\Cal O_{1,0}$
and
$\Cal O_2 =\Cal O_{2,0}$.
\smallskip\noindent
{\smc 6.2. The complex structure.}
The constituent $\fra p$ in the Cartan decomposition
$\fra {so}(2,q) = \fra k \oplus \fra p$
consists of matrices
in $\fra {so}(2,q)$
of the kind
$\left[\matrix 0 &B \\ B^t & 0 \endmatrix \right]$
where 
$B$ is an arbitrary real $(q \times 2)$-matrix
and thus
amounts to the space
$\roman M_{q,2}(\Bobb R)$
of real $(q \times 2)$-matrices.
The requisite complex structure
on
$\roman M_{q,2}(\Bobb R)$ is given by
matrix multiplication
$$
J \left[\matrix x_1 & \ldots &x_q\\ y_1 & \ldots &y_q \endmatrix \right]
= 
\left[\matrix  y_1 & \ldots &y_q \\ -x_1 & \ldots &-x_q\endmatrix \right]
$$
whence complex analytically,
$\overline{\Cal O_2} \cong \fra p^+$
amounts to a copy of $\Bobb C^q$.
In this description,
$$
e_1=\mu_{\GH}(\alpha_1) =(i,1,0,\dots, 0) \in \Bobb C^q,\ 
e_2=\mu_{\GH}(\alpha_2) =(1,i,0,\dots,0) \in \Bobb C^q,
$$ 
and the action of the complexification
$K^{\Bobb C} = \roman{SO}(2,\Bobb C) \times \roman{SO}(q,\Bobb C)$
of the maximal compact subgroup 
$K = \roman{SO}(2,\Bobb R) \times \roman{SO}(q,\Bobb R)$
of $\roman{SO}(2,q)$
is given by the ordinary 
$\roman{SO}(q,\Bobb C)$-action
on 
$\Bobb C^q$
together with the
multiplication by non-zero complex numbers,
when 
$\roman{SO}(2,\Bobb C)$ is identified with
$\Bobb C^*$.
After introduction of the complex coordinates $w_1,\dots, w_q$
where $w_j = x_j + i y_j$ ($1 \leq j \leq q$),
under the homeomorphism from
$\overline{\Cal O_2}$ 
onto 
$\fra p^+$,
complex analytically,
the closure $\overline{\Cal O_1}$ 
is identified with
the complex quadric 
in $\Bobb C^q$
given by the equation
$f(w_1,\dots, w_q) = 0$ where
$$
f(w_1,\dots, w_q) = w_1^2 + \dots + w_q^2 .
$$
In fact,
$\Cal O_1$ is generated by $e_1$,
and the $K^{\Bobb C}$-orbit 
of
$(i,1,0,\dots, 0) \in \Bobb C^q$
consists of the non-zero points
of the quadric 
$w_1^2 + \dots + w_q^2= 0$.
The significance of the quadratic form $f$ will be elucidated further
in Theorem 7.1 below.
\smallskip\noindent
{\smc Remark 6.3.} The above observations combine to  a proof of
Theorem 3.6.3 and hence
of statement (iii) of Theorem 3.6.2.  

\smallskip \noindent
{\smc Remark 6.4.} Just as in the standard cases, even though 
as a complex analytic space, $\overline {\Cal O_2}$
is affine and hence non-singular,
as a stratified symplectic space
and hence, as a normal K\"ahler space,
it has {\it singularities\/}, the {\it singular locus\/}
(in the sense of stratified K\"ahler spaces)
being the complex quadric $\overline {\Cal O_1}$.

\beginsection 7. Hermitian Jordan triple systems and 
pre-homogeneous spaces

In this section we will relate holomorphic nilpotent orbits
with positive definite hermitian Jordan triple systems.

Recall that a pair $(\fra g,\vartheta)$
consisting of a semisimple  Lie algebra $\fra g$
together with a 3-grading
$\fra g = \fra g_{-1} +\fra g_0+\fra g_1$
and an involution $\vartheta$ of $\fra g$ such that
$\vartheta \fra g_{j} =\fra g_{-j}$
($j=-1,0,1$)
is called a ({\it non-degenerate\/}) {\it symmetric Lie algebra\/}
(\cite\satakboo\ I.7 p.~28)
provided 
$\fra g_0 = [\fra g_{-1},\fra g_1]$.
Given
a symmetric Lie algebra
$(\fra g,\vartheta)$,
the abelian subalgebra
$V=\fra g_{-1}$
(as well as 
the abelian subalgebra
$\fra g_1$),
endowed with the trilinear operation $\{\cdot,\cdot,\cdot\}$ defined by
$$
\{x,y,z\} =-\frac 12 [[x,\vartheta y],z]
$$
is what is called a {\it non-degenerate J(ordan) T(riple) 
S(ystem)\/}---the defining axioms are not reproduced 
here, cf. \cite\satakboo\  p.~22---and
every non-degenerate JTS arises in this fashion,
see \cite\satakboo\  p.~28 and the references there.
For a simple Lie algebra $(\fra g,z)$ of hermitian type,
the decomposition (3.3.4($r$))
(where $r$ is the real rank of $\fra g$)
 yields a
real symmetric Lie algebra if and only if
$(\fra g,z)$ is regular,
the requisite involution being the Cartan involution.
\smallskip
A (non-degenerate) real JTS
$(V,\{\cdot,\cdot,\cdot\})$, together with a complex structure (on $V$)
such that 
$\{\cdot,\cdot,\cdot\}$
is $\Bobb C$-linear in the first and third variable and
$\Bobb C$-antilinear in the second one
is said to be a {\it hermitian\/} JTS,
and it is said to be {\it positive definite\/}
provided 
the hermitian form $\tau'$ on $V$ ($\Bobb C$-linear in the second variable)
given by the trace $\tau'(x,y)$ (taken over $\Bobb C)$
of the $\Bobb C$-linear transformation
$y \squar x$ of $V$ given by $(y \squar x)z=\{y,x,z\}$
is positive definite ($x,y,z \in V$).
\smallskip
Let $(\fra g,z)$ be a
simple Lie algebra of hermitian type,
with Cartan decomposition
$\fra g = \fra k \oplus \fra p$.
Its complexification
$\fra g^{\Bobb C} =  \fra p^+\oplus \fra k^{\Bobb C} \oplus \fra p^-$
inherits a real symmetric Lie algebra structure,
the requisite involution being complex conjugation
(but, beware, this is {\it not\/} the 
symmetric Lie algebra structure mentioned earlier).
Thus 
$(\fra g,z)$ determines  
the (non-degenerate) JTS $(V,\{\cdot,\cdot,\cdot\})$
where $V=\fra p^+$.
It is known that
this JTS is hermitian and positive definite
and that every
positive definite hermitian JTS arises in this fashion,
cf. \cite\satakboo\ (Proposition II.3.3 on p.~56);
indeed, with the appropriate notions of morphism,
associating with
a (real) semisimple Lie algebra of hermitian type
$(\fra g,z)$
the JTS $(\fra p^+,\{\cdot,\cdot,\cdot\})$
yields an equivalence of categories between
(real) semisimple Lie algebras of hermitian type
and positive definite hermitian JTS's
\cite\satakboo\  (Proposition II.8.2 on p.~85).
Under this correspondence,
the regular semisimple Lie algebras of hermitian type
correspond to positive definite hermitian JTS's
which \lq\lq arise\rq\rq\ from a Jordan algebra (with unit element),
in the following sense:
A finite dimensional (non-associative) algebra
$U$ satisfying
$xy = yx$ and $x^2(xy)= x(x^2y)$ 
($x,y \in U$)
is called a {\it Jordan algebra\/}.
Any Jordan algebra $U$ becomes a JTS when the triple product
$\{\cdot,\cdot,\cdot\}$
is defined by
$$
\{x,y,z\} = (xy)z + x(yz) -y(xz),
\quad x,y,z \in U,
$$
and we then say that the JTS {\it arises\/} from
the Jordan algebra $U$.
When $U$ has a unit element $e$, the Jordan product is determined by
the triple product by means of 
$$
\{x,y,e\} (= \{x,e,y\}=\{e,x,y\} )=xy,\quad x,y \in U.
$$
\smallskip
A representation $E$ of a complex reductive Lie group is
said to be a
{\it pre-homogeneous space\/}
\cite\vinberg\ 
provided it has a Zariski-open orbit, 
referred to henceforth as the {\it top\/} orbit,
and it
is said to be {\it regular\/}
\cite\satokimu\ 
provided the complement of the
top orbit is a hypersurface or, equivalently,
provided the isotropy subgroups of the points in the
top orbit are reductive (cf. p. 96 of \cite\muruschi).
For an irreducible
regular pre-homogeneous space,
where \lq\lq irreducible\rq\rq\ 
means that the representation is irreducible,
the hypersurface
arising as the complement of the top
orbit is known to be irreducible;
the polynomial which it defines
is a {\it relative invariant\/}
(i.~e. it is a \lq\lq twisted\rq\rq\ invariant, with reference
to a suitable character)
of the representation
(p. 96 of \cite\muruschi)
which is unique up to a constant
and referred to as the {\it fundamental relative
invariant\/};
this 
fundamental relative
invariant then generates 
the
relative invariants.
\smallskip
Given
a simple Lie algebra $(\fra g,z)$ of hermitian type,
the space
$\fra p^+$ is plainly a pre-homogeneous
space for the complexification 
$K^{\Bobb C}$
of the compact constituent
$K$ in the Cartan decomposition of 
the adjoint group
$\HG$ of $\fra g$.
These pre-homogeneous spaces
are precisely those which arise from
positive definite hermitian JTS's.
In \cite\muruschi, they 
have been classified 
in terms of root systems
(without reference to the notion of Lie algebra of hermitian type).
For intelligibility, we reproduce explicit descriptions
of the
$K^{\Bobb C}$-actions on
$\fra p^+$ in the classical cases and recall which ones are regular:
\smallskip\noindent
$\fra g= \fra{sp}(\ell,\Bobb R)$: 
$\fra p^+=\roman S_{\Bobb C}^2[\Bobb C^{\ell}]$ is regular;
$K^{\Bobb C} = \roman{GL}(\ell,\Bobb C)$, the action on
$\roman S_{\Bobb C}^2[\Bobb C^{\ell}]$ being given by
$$
x\cdot S = x S x^t,\quad x \in \roman{GL}(\ell,\Bobb C),
\ S \in \roman S_{\Bobb C}^2[\Bobb C^{\ell}].
$$
This is the second symmetric power 
of the standard $\roman{GL}(\ell,\Bobb C)$-representation.
\newline\noindent
$\fra g=\fra {u}(p,q)$: $\fra p^+=\roman M_{q,p}(\Bobb C)$
is regular if and only if $p=q$;
$K^{\Bobb C} = \roman{GL}(p,\Bobb C) \times \roman{GL}(q,\Bobb C)$, 
the action on
$\roman M_{q,p}(\Bobb C)$ being given by
$$
(x,y)\cdot M = x M y^t,\quad x \in \roman{GL}(p,\Bobb C),
\ y \in \roman{GL}(q,\Bobb C),
\ M \in \roman M_{q,p}(\Bobb C);
$$
\newline\noindent
$\fra g=\fra {so}^*(2n)$: $\fra p^+=\Lambda^2[\Bobb C^n]$
is regular if and only if $n$ is even;
$K^{\Bobb C} = \roman{GL}(n,\Bobb C)$, the action on
$\Lambda^2[\Bobb C^n]$ being given by
$$
x\cdot A = x A x^t,\quad x \in \roman{GL}(n,\Bobb C),
\ A \in \Lambda^2[\Bobb C^n].
$$
This is the second exterior power 
of the standard $\roman{GL}(n,\Bobb C)$-representation.
\newline\noindent
$\fra g=\fra {so}(2,q)$: $\fra p^+=\Bobb C^q$ is regular;
$K^{\Bobb C} = \roman{SO}(2,\Bobb C) \times \roman{SO}(q,\Bobb C)$, 
the action on
$\Bobb C^q$ being given by
the ordinary 
$\roman{SO}(q,\Bobb C)$-action
on 
$\Bobb C^q$
together with the
multiplication by non-zero complex numbers,
when 
$\roman{SO}(2,\Bobb C)$ is identified with
$\Bobb C^*$.
\smallskip
Thus, $\fra p^+$ is regular if and only if
$(\fra g,z)$ is regular (cf. (3.3) above), that is,
if and only if its relative root system is of type
$(C_r)$ (where $r$ is the real rank).
Furthermore,
the property of being regular is equivalent to the corresponding
positive definite hermitian JTS to arise from a 
formally real Jordan algebra, cf. \cite\satakboo\ I.8 p.~31
(also referred to as {\it euclidean\/} Jordan algebra
in the literature \cite\farakora),
and the complexification of these 
formally real Jordan algebras gives the
simple complex Jordan algebra structure on
the constituent $\fra p^+$.
Plainly, in these cases,
whenever
$\fra p^+$ is regular, it is irreducible.
Recall that the {\it  rank\/}
of a non-zero element of a 
euclidean as well as of a positive definite hermitian
Jordan algebra,
referred to henceforth as its {\it Jordan rank\/}, is the number of
non-zero eigenvalues in its spectral decomposition,
with multiplicities counted
\cite\farakora\ (p.~77).

\proclaim{Theorem 7.1}
In the regular classical cases,
under the homeomorphism 
(induced by the projection from $\fra g$ to $\fra p$)
from the closure
$\overline {\Cal O}$ of the principal holomorphic nilpotent orbit
$\Cal O$
onto the complex vector space $\fra p^+$,
the $G$-orbit stratification of $\overline {\Cal O}$
passes to the stratification
of $\fra p^+$
by Jordan rank, 
with reference to the 
complex simple Jordan algebra structure on $\fra p^+$, 
and that stratification, in turn, coincides
with the $K^{\Bobb C}$-orbit stratification of
$\fra p^+$. 
In the standard cases,
when $\fra p^+$ is realized as a space
of complex matrices as indicated above,
the fundamental relative
invariant 
(Bernstein-Sato polynomial)
is the determinant.
For $\fra g=\fra {so}(2,q)$, in suitable 
complex coordinates $w_1,\dots, w_q$
on 
$\overline{\Cal O_2} \cong \fra p^+ \cong \Bobb C^q$
(e.g. those introduced in {\rm (6.2)}),
the fundamental relative invariant
(Bernstein-Sato polynomial)
$f$ 
is the quadratic form
$$
f(w_1,\dots, w_q) = w_1^2 + \dots + w_q^2 .
$$
Thus the stratified symplectic Poisson structure detects the stratification
by Jordan rank.
\endproclaim

\demo{Proof}
The statement involving the Jordan rank is just a reformulation
of the description of the strata in Theorems 3.5.5
and 3.6.3.
To see that the fundamental relative invariants
have the asserted form,
we note that,
by Theorem 5.4,
in each of the cases
$(\fra g,\fra p^+,r)=
(\fra{sp}(\ell,\Bobb R),\fra p^+=\roman S_{\Bobb C}^2[\Bobb C^{\ell}],\ell)$,
$(\fra g,\fra p^+,r)
=(\fra {su}(p,p),\roman M_{p,p}(\Bobb C),p)$, 
$(\fra g,\fra p^+,r)
=(\fra {so}^*(4\ell),\Lambda^2[\Bobb C^{2\ell}],\ell)$,
the complement
of the top stratum 
$\Cal O_r$
is the hypersurface
$\overline{\Cal O_{r-1}}$
which is given by a single determinantal equation.
The determinant is clearly a relative invariant.
For $\fra g=\fra {so}(2,q)$,
in the coordinates $w_1,\dots, w_q$
on 
$\overline{\Cal O_2} \cong \fra p^+ \cong \Bobb C^q$
introduced in {\rm (6.2)},
the 
fundamental relative invariant
$f$ on 
$\overline{\Cal O_2} \cong \fra p^+ \cong \Bobb C^q$
has the form
$$
f(w_1,\dots, w_q) = w_1^2 + \dots + w_q^2 .
$$
This is truly a relative
$K^{\Bobb C}$-invariant,
the requisite character 
for $K^{\Bobb C}$
being given by the projection
to 
$\roman{SO}(2,\Bobb C)\cong \Bobb C^*$,
followed by the squaring map from $\Bobb C^*$ to itself. 
Finally,
it is a standard fact that a stratified symplectic Poisson algebra
detects the strata. \qed
\enddemo

In the regular standard cases, the fundamental relative invariants
for the complexification of the compact constituent
in the Cartan decomposition of $\HG$
therefore arise
{\sl from the relations among the $K$-invariants
describing 
the complex analytic model\/}
$W_J\big /\big / K^{\Bobb C}$
of $\overline {\Cal O_{r-1}}$ since
these 
yield the defining equations
for
$\overline {\Cal O_r}$
as a complex affine variety.
Thus the fundamental relative invariants
are seen
as resulting from an application of the
\lq\lq Second Main Theorem\rq\rq\ of Invariant Theory.

\smallskip
\noindent
{\smc Remark 7.2.} 
See \cite{\farakora,\,\mccrione,\,\satakboo}
for notation and details.
For $\Bobb K = \Bobb R,\,\Bobb C,\, \Bobb H$,
consider the {\it real special Jordan algebra\/} 
of hermitian matrices
$\Cal H_n(\Bobb K)$ over $\Bobb K$
($\Cal H_n(\Bobb R) = S^2_{\Bobb R}[\Bobb R^n]$),
with Jordan product 
$\,\circ\,$ given by
$x \circ y = \frac 12 (xy + yx)$
($x,y \in \Cal H_n(\Bobb K)$).
The {\it Kantor-Koecher-Tits\/} construction 
of superstructure algebra for $\Cal H_n(\Bobb K)$ over $\Bobb K$ yields
the regular standard simple
Lie algebras of hermitian type, together with the appropriate
real symmetric Lie algebra decompositions (3.3.4($r$))
(where $r$ is the real rank), as follows:
$$
\alignat 2
\Bobb K &= \Bobb R \colon
\quad\Cal H_n(\Bobb R) \oplus \fra {gl}(n,\Bobb R) \oplus \Cal H_n(\Bobb R)
&&= \fra{sp}(n,\Bobb R)
\\
\Bobb K &= \Bobb C \colon
\quad
{\Cal H_n(\Bobb C)} 
\oplus \fra {gl}(n,\Bobb C) \oplus \Cal H_n(\Bobb C)
&&= \fra{u}(n,n)
\\
\Bobb K &= \Bobb H \colon
\quad
{\Cal H_n(\Bobb H)} 
\oplus \fra {gl}(n,\Bobb H) \oplus \Cal H_n(\Bobb H)
&&= \fra{so}^*(4n)
\endalignat
$$
Likewise, 
the 
superstructure algebra construction for 
the euclidean Jordan algebra $J(1,q-1)$
(arising from the Lorentz form of type $(1,q-1)$ on $\Bobb R^q$)
yields
the Lie algebra $\fra{so}(2,q)$, together with the corresponding
real symmetric Lie algebra decomposition (3.3.4($2$)),
which has the form
$$
J(1,q-1) \oplus \Bobb Rh^2 \oplus \fra {so}(1,q-1) 
\oplus J(1,q-1)= \fra{so}(2,q).
$$

\smallskip
We will now show that, in the non-regular classical cases,
the stratification can still
be described by Jordan rank:
The notion of rank of an element of
a positive definite hermitian JTS
(as well as of a real positive definite JTS)  
has been introduced in \cite\loobothr\ (11.6), see also \cite\neherone\
(where a general notion of rank, written $r$, is introduced without a name
on p. 2625);
accordingly we will continue to refer to the {\it Jordan rank\/}
of such an element.
For $\fra g = \fra{su}(p,q)$, when $p >q =r$
(so that $(\fra g,z)$ is non-regular), 
the projection from $\fra g$ to $\fra p$, followed by the 
isomorphism from $\fra p$ onto $\fra p^+= \roman M_{q,p}(\Bobb C)$,
identifies the closure $\overline{\Cal O_q}$  of the
principal holomorphic nilpotent orbit 
$\Cal O_q$ with 
the complex vector space underlying the
positive definite hermitian Jordan triple system 
$\roman M_{q,p}(\Bobb C)$
(the triple product  being given by
$\{x,y,z\} = x \overline y^t z + z \overline y^t x$).
In this fashion,
for $1 \leq s < r$,
$\Cal O_s$ amounts the space
of complex $(q \times p)$-matrices of rank $s$
and $\overline{\Cal O_s}$ with that
of complex $(q \times p)$-matrices of rank at most
$s$, and the ordinary rank of a matrix  
coincides with its Jordan rank.
Likewise,
in case $\fra g = \fra{so}^*(2n)$, when $n=2\ell+1$ 
(so that $(\fra g,z)$ is non-regular), 
the 
corresponding
map from $\fra g$ to 
$\fra p^+= \Lambda^2[\Bobb C^n]$
identifies the closure $\overline{\Cal O_{\ell}}$  of the
principal holomorphic nilpotent orbit 
$\Cal O_{\ell}$ with 
the complex vector space underlying the
positive definite hermitian Jordan triple system $\Lambda^2[\Bobb C^n]$.
For $1 \leq s \leq r=\ell$,
$\Cal O_s$ thus amounts to the space of
skew-symmetric
complex $(n \times n)$-matrices of rank 
$2s$ and
$\overline{\Cal O_s}$ to that of
skew-symmetric
complex $(n \times n)$-matrices of rank at most
$2s$, and the ordinary rank of a matrix  
is twice its Jordan rank.
Summarizing, we spell out the following.

\proclaim{Theorem 7.3}
In the  classical cases,
under the homeomorphism 
from the closure
$\overline {\Cal O}$ of the principal holomorphic nilpotent orbit
$\Cal O$
onto the complex vector space $\fra p^+$,
the $G$-orbit stratification of $\overline {\Cal O}$
passes to the stratification
of $\fra p^+$
by Jordan rank 
with reference to the  (positive definite hermitian) JTS-structure
on $\fra p^+$,
and that stratification, in turn, coincides
with the $K^{\Bobb C}$-orbit stratification of
$\fra p^+$. 
Thus the induced stratified symplectic Poisson structure 
on $\fra p^+$ 
turns the latter into a normal K\"ahler space
and detects the stratification by Jordan rank. \qed
\endproclaim

All geometric information about holomorphic nilpotent orbits
in a simple classical Lie algebra of hermitian type
is therefore encoded in the corresponding JTS.
Theorem 7.3 entails in particular that
the pre-homogeneous spaces
arising from positive definite hermitian JTS's
and classified in \cite\muruschi\  
carry much more structure than that made explicit in that reference.

\medskip\noindent
{\bf 8. The exceptional cases}
\smallskip\noindent
In this section we recollect some information relevant 
for the exceptional cases; 
see e.~g. \cite\farakora,\ \cite\jacobexc,
\ \cite\peteraci, and \cite\satakboo\ III.4 p.~117/118 for details
and notation. We then describe the geometry of  
the corresponding holomorphic nilpotent orbits.
Until the end of this section, $(\fra g,z)$ will be one of the two
simple exceptional Lie algebras of hermitian type;
as before, we 
write $\fra g = \fra k \oplus \fra p$ for the Cartan decomposition
and $\fra g_{\Bobb C} = \fra p^- \oplus \fra k_{\Bobb C} \oplus \fra p^+$
for the complexification thereof; further, we
denote the real rank by $r$,
the adjoint group of $\fra g$ by $G$,
and the maximal compact subgroup
determined by the choice of $z$
by $K$.
\medskip\noindent
{\smc 8.1. Octonions and Albert algebras.}
We shall use the neutral notation $\Bobb K$ for
$\Bobb R$ and/or $\Bobb C$ if need be.
Let $\Bobb O=\Bobb O_{\Bobb R}$ be 
the real {\it octonion algebra of Cayley numbers\/},
also called {\it Cayley division algebra\/},
and write
$\Bobb O_{\Bobb C}$
for its complexification;
this is a complex octonion algebra.
As usual, 
for $a \in  \Bobb O$ as well as for $a \in \Bobb O_{\Bobb C}$,
we write
$a \mapsto \overline a$ for the operation of conjugation,
and the {\it norm\/}
$n(a)\in \Bobb K$ and {\it trace\/}
$t(a)\in \Bobb K$ are defined by
$n(a)= a \overline a$
and 
$t(a)= a +\overline a$
where $\Bobb K$ is viewed as a subalgebra
of $\Bobb O_{\Bobb K}$ in the obvious fashion.
For $m=2$ and $m=3$,
let $\Cal H_m(\Bobb O)$
be the euclidean Jordan algebra
of positive definite hermitian $(m \times m)$-matrices
over $\Bobb O$.  
The Jordan algebra
$\Cal H_2(\Bobb O)$
has rank 2 and is isomorphic
to $J(1,9)$---in view of the classification,
every rank two hermitian Jordan algebra
arises from a quadratic form
in this fashion---and
$\Cal H_3(\Bobb O)$ is the
real reduced exceptional central simple rank 3 Jordan algebra
of dimension 27, a real {\it Albert algebra\/}.
The complexification of the latter
is the (reduced central simple) rank 3
{\it complex exceptional Jordan algebra\/}
$\Cal H_3(\Bobb O_{\Bobb C})$,
a complex {\it Albert algebra\/};
cf. e.~g. VIII.5 on p.~159 of \cite\farakora,\ p.~121 of \cite\jacobexc.
Here the property of being {\it reduced\/} 
amounts to the existence of three non-zero orthogonal idempotents,
and an Albert algebra
being {\it exceptional\/}
means that it cannot be embedded into any associative algebra
whatsoever.
Whenever we wish to 
treat the two Jordan algebras simultaneously,
we write $\Cal H_3(\Bobb O_{\Bobb K})$,
for $\Bobb K = \Bobb R$ and
$\Bobb K = \Bobb C$;
in \cite\jacobexc\ (p.~17),  
$\Cal H_3(\Bobb O_{\Bobb K})$ is denoted by
$\roman H(\Cal C_3,\gamma)$ with hermitian form $\gamma = \roman{diag}(1,1,1)$
where $\Cal C$ refers to the Cayley division algebra over $\Bobb R$
as well as to the complexification thereof.
Thus the Jordan algebra $\Cal H_3(\Bobb O_{\Bobb K})$
consists of 
$(3 \times 3)$-matrices
$$
A=\left[\matrix 
      \alpha_1 & a_3            & \overline{a_2}\\
\overline{a_3} & \alpha_2       & a_1\\
a_2            & \overline{a_1} & \alpha_3
\endmatrix \right]
$$
where $\alpha_1,\alpha_2,\alpha_3 \in \Bobb K$ and 
$a_1,a_2,a_3 \in \Bobb O_{\Bobb K}$,
and the Jordan product $\,\circ\,$
is given by
$x \circ y = \frac 12 (xy + yx)$
($x,y \in \Cal H_3(\Bobb O_{\Bobb K})$).
Such a matrix $A$ 
has {\it generic norm} $\nu(A)$ given by
$$
\nu(A) = \alpha_1\alpha_2\alpha_3 + t(a_3a_1a_2)
-\alpha_1 n(a_1)
-\alpha_2 n(a_2)
-\alpha_3 n(a_3).
$$
The generic norm $\nu$ amounts to Freudenthal's \lq\lq formally defined
determinant\rq\rq\ \cite\freudone.

\smallskip\noindent
{\smc 8.2.} $\fra g =\fra e_{7(-25)}$. 
This Lie algebra arises
from the
superstructure algebra construction
(cf. e.~g. \cite\mccrione)
for the positive definite hermitian Jordan algebra
$\Cal H_3(\Bobb O)$,
cf. \cite\jacobexc\ (No. 3 p.~121),\,\cite \mccrione.
The real rank $r$ of $\fra e_{7(-25)}$ equals 3, 
the maximal compact subgroup $K$ of the adjoint group 
is locally isomorphic to $E_{6(-78)} \times \roman{SO}(2,\Bobb R)$
($\roman E_{6(-78)}$ being a compact form),
$\dim \fra k = 79$, and $\dim \fra p = 54$.
\smallskip
The stratification of the closure of the
principal holomorphic nilpotent orbit
may be described in terms of 
the corresponding decompositions (3.3.4($k$)) 
($1\leq k \leq 3$);
in particular,
$\dim \fra z_{\fra g}(e^3)
= 52 + 27 = 79 = \dim \fra k$
whence
$\dim \Cal O_3$ equals $\dim \fra p$
($\Cal O_3=Ge^3$).
Indeed, direct inspection
confirms that the projection from the closure
$\overline {\Cal O_3}$ 
of the principal holomorphic nilpotent orbit
$\Cal O_3$
to $\fra p$ is a homeomorphism {\it onto\/} $\fra p$,
as claimed in Theorem 3.3.11.

\smallskip\noindent
{\smc 8.3.} $\fra g =\fra e_{6(-14)}$: 
With respect to the primitive idempotent $e_{3,3}$,
the {\it Peirce decomposition\/} 
(cf. e.~g. \cite\satakboo\ V.6 p.~242 ff.)
of
$\Cal H_3(\Bobb O)$, viewed as a JTS, 
has the form
$$
\Cal H_3(\Bobb O) = V_{0} \oplus V_{\frac 12}\oplus V_{1};
$$
here
$V_{0}$ is 
the rank 2 Jordan algebra which is
a copy of
$\Cal H_2(\Bobb O)$,
$V_{1}$ is 
the rank 1 Jordan algebra which is
the (real) span
of $e_{3,3}$,
and
$V_{\frac 12}$
is the JTS of 
$(3 \times 3)$-matrices
$A=\left[\matrix 
   0 & 0            & \overline{a_2}\\
   0 & 0       & a_1\\
a_2  & \overline{a_1} & 0
\endmatrix \right]
$
where 
$a_1,a_2 \in \Bobb O$.
The triple product on
$\Cal H_3(\Bobb O)$
induces a triple product on
$V_{\frac 12}$
turning the latter into a positive definite real JTS
of dimension 16. We identify the elements of
$V_{\frac 12}$
with the space $\roman M_{1,2}(\Bobb O)$
of $(1\times 2)$-matrices over $\Bobb O$
by means of the assignment to
$[a_1,a_2]^{\tau}$ of
$\left[\matrix 
   0 & 0            & \overline{a_2}\\
   0 & 0       & a_1\\
a_2  & \overline{a_1} & 0
\endmatrix \right].
$
The  superstructure algebra  
construction
for this positive definite 
real JTS
yields the exceptional Lie algebra 
$\fra e_{6(-14)}$.
The real rank $r$ of $\fra e_{6(-14)}$ equals 2, 
the maximal compact subgroup
$K$
of the adjoint
group $G$
is isomorphic to $\roman{Spin}_{\roman c}(10,\Bobb R)$,
and $\dim \fra k = 46$, $\dim \fra p = 32$.
To identify the representation
of
$\roman{Spin}_{\roman c}(10,\Bobb R)$
on 
$\fra p$ we note that,
$\fra p$ being endowed with the complex structure $J_z$,
the restriction
of this representation to 
$\roman{Spin}(10,\Bobb R) \subseteq\roman{Spin}_{\roman c}(10,\Bobb R)$
is the positive half-spin representation
$\Cal D_+^5$ of complex dimension 16;
here
$\Cal D_+^5$
is one of the two irreducible modules
for the complex Clifford algebra $C_9$.

\smallskip
As pointed out earlier,
the stratification of the closure of the principal holomorphic nilpotent orbit
may be described in terms of 
the corresponding decompositions (3.3.4($k$)) 
($1\leq k \leq 2$);
in particular,
$\dim \fra z_{\fra g}(e^2)= \dim \fra k$,
whence $\dim \Cal O_2 =\dim \fra p$ ($\Cal O_2=Ge^2$).
Again, direct inspection
confirms that
the projection from the closure
$\overline {\Cal O_2}$ 
of the principal holomorphic nilpotent orbit
$\Cal O_2$
to $\fra p$ is a homeomorphism {\it onto\/} $\fra p$,
as claimed in Theorem 3.3.11. 

\smallskip\noindent
{\smc 8.4. 
The statement of Theorem 7.3 for the exceptional cases.}
The generic norm
$\nu$ for
$\Cal H_3(\Bobb O_{\Bobb C})$
is a homogeneous complex cubic form.
The notion of
{\it Jordan rank\/}
has been recalled in Section 7 above.
For intelligibility, we note that
$\fra g_{\Bobb C} =\fra e_7(\Bobb C)$
is the superstructure algebra of  
$\Cal H_3(\Bobb O_{\Bobb C})$,
and that
the structure group $K^{\Bobb C}$
of
$\Cal H_3(\Bobb O_{\Bobb C})$
is locally isomorphic to  $\Bobb C^* \times E_6(\Bobb C)$.

\proclaim{Theorem 8.4.1}
For  $\fra g =\fra e_{7(-25)}$,
the pre-homogeneous space representation on
$\fra p^+$
is the resulting
$K^{\Bobb C}$-representation 
on $\Cal H_3(\Bobb O_{\Bobb C})$, and
the relative invariant 
(Bernstein-Sato polynomial)
is the generic norm $\nu$. Furthermore,
under the 
projection 
from $\fra g = \fra k \oplus \fra p$ to
$\fra p$, followed by the canonical complex linear isomorphism from
$\fra p$ to $\fra p^+$,
the closure
$\overline{\Cal O_3}$ 
of the top holomorphic nilpotent orbit
is identified with $\fra p^+=\Cal H_3(\Bobb O_{\Bobb C})$ 
in such a way that,
for $0 \leq s \leq 3$,
$\Cal O_s$ 
is mapped diffeomorphically onto
the smooth complex affine 
variety 
which consists of 
Jordan rank $s$ elements
whence,
for $1 \leq s < 3$,
as a complex analytic space, 
$\overline {\Cal O_s}$ 
amounts to the complex affine 
variety 
in 
$\fra p^+$
which consists of 
elements in $\fra p^+$
which have Jordan rank at most $s$.
Thus requiring that
the generic norm vanish yields an explicit equation
defining
$\overline {\Cal O_2}$ as an affine complex variety 
in $\fra p^+=\Cal H_3(\Bobb O_{\Bobb C})$ $(\cong \Bobb C^{27})$.
As  affine complex varieties,
$\overline {\Cal O_1}$  and
$\overline {\Cal O_2}$ are normal
(and so is $\overline {\Cal O_3}$, being an affine space).
\endproclaim

The top stratum
$\Cal O_3$ 
is thus identified
with the Jordan-invertible elements
in $\Cal H_3(\Bobb O_{\Bobb C})$.
In the literature, 
the generic norm hypersurface
$\overline{\Cal O_2} = \nu^{-1}(0) \subseteq \Cal H_3(\Bobb O_{\Bobb C})$ 
is sometimes referred to as {\it Freudenthal's cubic\/};
we write it as $\Cal F$.
We note that,
for $1 \leq s \leq 3$,
the group $E_6(\Bobb C)$ ($\subseteq K^{\Bobb C}$)
acts transitively on
(the image of)
each orbit $\Cal O_s$ 
(in $\fra p^+$),
and that,
in view of the infinitesimal structure given in
Table XIII in \cite\djokotwo\ (p.~516),
the point stabilizers 
are the following subgroups
of $E_6(\Bobb C)$:
\newline
\noindent
$F_4(\Bobb C)$ (for $\Cal O_3$),
$\roman{Spin}(9,\Bobb C) \cdot (\Bobb C^*)^{16}$ (for $\Cal O_2$),
$\roman{Spin}(10,\Bobb C) \cdot (\Bobb C^*)^{16}$ (for $\Cal O_1$).

\demo{Proof of Theorem {\rm 8.4.1}}
The first statement is a special case of Theorem 3.3.11. 
Since $\Cal F$ is an irreducible complex affine hypersurface which
is non-singular in codimension 1,
in view of standard algebraic geometry results,
$\overline{\Cal O_2}=\Cal F$ is normal.
The normality of
the lowest non-trivial stratum 
$\overline{\Cal O_1}$
will be justified in (10.7) below. \qed
\enddemo

\smallskip
Let $V$ be a complex vector space of even (complex) dimension $n$,
endowed with a non-degenerate quadratic form $q$, let
$C(V,q)$ be the corresponding complex Clifford algebra,
and let $\Cal S^{\Bobb C}$ be 
an irreducible complex spinor space, i.~e.
a fundamental Clifford module.
Recall that a spinor 
$\sigma \in \Cal S^{\Bobb C}$ 
is said to be {\it pure\/}
provided the 
kernel of the
$\Bobb C$-linear map
$$
j_{\sigma} \colon V @>>> \Cal S^{\Bobb C},
\quad
j_{\sigma}(v) = v \cdot \sigma
\quad \text{(Clifford multiplication)}
$$
has dimension $\frac n2$;
see e.~g. \cite\lawmiboo\ (IV.9~p.~336) for details.
For $n=10$,
$\Cal S^{\Bobb C} = \Cal D^5$ and, as a 
$\roman{Spin}(10,\Bobb R)$-representation,
$\Cal S^{\Bobb C}$
decomposes into the positive and negative
complex half-spin representations
$\Cal S^{\Bobb C}_+=\Cal D^5_+$  
and
$\Cal S^{\Bobb C}_-=\Cal D^5_-$.

\proclaim{Theorem 8.4.2}
For  $\fra g =\fra e_{6(-14)}$,
under the 
projection 
from $\fra g = \fra k \oplus \fra p$ to
$\fra p$, followed by the canonical complex linear isomorphism from
$\fra p$ to $\fra p^+$,
the orbit $\Cal O_1$ 
is 
identified with the space of positive pure spinors
in $\fra p^+\cong\Cal S^{\Bobb C}_+$.
As a complex analytic space, 
$\overline {\Cal O_1}$ 
amounts to the complex affine 
variety 
in 
$\fra p^+$
which consists of 
elements in $\fra p^+$
which have Jordan rank at most $1$.
As an  affine complex variety,
$\overline {\Cal O_1}$ is normal
(and so is $\overline {\Cal O_2}$, being an affine space).
\endproclaim

We note that
the group $\roman{Spin}(10, \Bobb C)$ 
($\subseteq K^{\Bobb C}$)
acts transitively on
the images of $\Cal O_1$ 
and $\Cal O_2$ 
in $\fra p^+$,
and that,
in view of Proposition 2 in \cite\junigusa\ (p.~1011), 
the point stabilizers 
for the orbits
$\Cal O_2$
and $\Cal O_1$
are the subgroups
$\roman{Spin}(7,\Bobb C) \cdot (\Bobb C^*)^8$
and $\roman{SL}(5,\Bobb C) \cdot (\Bobb C^*)^{10}$,
respectively,
of $\roman{Spin}(10, \Bobb C)$.

\demo{Proof}
The first statement is a special case of Theorem 3.3.11. 
Since the representation 
$\Cal S^{\Bobb C}_+$
has only two
non-zero $\roman{Spin}(10, \Bobb C)$-orbits,
the orbit
$\Cal O_1$ is necessarily that of pure spinors.
The normality of
the stratum 
$\overline{\Cal O_1}$
will be justified in (10.7) below. \qed
\enddemo

\medskip\noindent
{\bf 9. Contraction of semisimple  holomorphic orbits} 
\smallskip\noindent
The idea of contraction 
(as the operation inverse to a deformation)
goes back to \cite\inonwign;
a leisurely introduction 
into the notion of contraction we will use
may be found in Chap. V of \cite\guistebo.
Let $(\fra g,z)$ be a simple Lie algebra of Hermitian type
and, as before, write $G$ for its adjoint group and $K$ for the
maximal compact subgroup of $G$ determined by 
the choice
$z$ of $H$-element 
(i.~e. $\roman{Lie}(K) = \fra k = \roman{ker}(\roman{ad}_z)$).
Thus $K$ is the stabilizer $Z_G(z)$ of $z$,
and the Lie algebra $\fra k$ of $K$
equals the stabilizer algebra $\fra z_{\fra g}(z)$.
With reference to the 
decomposition $\fra g = \fra n^-_r \oplus \fra l_r \oplus\fra n^+_r$,
cf. (3.3.4($r$)) where, as before,  $r$ is the real rank of $\fra g$,
let
$\fra k_0 = \fra k \cap \fra l_r$, and consider the filtered Lie algebra
$\fra k_0 \subseteq \fra k$;
its associated graded Lie algebra
$E^0(\fra k) =\fra k_0 \lltimes \fra n^+_r$
is the semi-direct product of 
its compact constituent $\fra k_0$ with 
the nilpotent Lie algebra $\fra n^+_r$,
the $\fra k_0$-action
on $\fra n^+_r$
 being the obvious one
resulting from restriction of the
$\fra l_r$-action on $\fra n^+_r$ to $\fra k_0$.
The Lie algebra $E^0(\fra k)$
is the stabilizer algebra $\fra z_{\fra g}(e^r)$
of the element $e^r$ generating the principal holomorphic nilpotent orbit
$\Cal O_r$,
and the
stabilizer $\roman Z_{\fra g}(e^r)$
of $e^r$ is
the corresponding subgroup of the adjoint group $G$.
\smallskip
For $\varepsilon >0$, 
let $\Cal O^+(\varepsilon)=G(2\varepsilon z)$;
these orbits, being semisimple and holomorphic,
are ordinary K\"ahler manifolds,
the requisite complex structure being induced from the projection
to $\fra p \cong \fra p^+$,
cf. (3.7.6).
When $\varepsilon$ goes to zero, they tend to 
or, equivalently, {\it contract\/} to
the stratified K\"ahler space
$\overline{\Cal O_r}$ 
(the closure of the principal holomorphic nilpotent orbit $\Cal O_r$)
in such a way that the complex structure
remains fixed; indeed, the projection
to $\fra p^+$ is a diffeomorphism from
$\Cal O^+(\varepsilon)$ onto 
$\fra p^+$ for $\varepsilon >0$
and a 
homeomorphism 
from $\overline{\Cal O_r}$
onto 
$\fra p^+$ 
in the limit case $\varepsilon = 0$.
Writing
$(C^{\infty},\{\cdot,\cdot\})_{\varepsilon}$ ($\varepsilon >0$)
for the corresponding smooth symplectic Poisson structure
on
$\fra p^+$
and $(C^{\infty},\{\cdot,\cdot\})_0$ 
for the stratified symplectic Poisson structure
on
$\fra p^+$
resulting from the projection from
$\overline{\Cal O_r}$
onto 
$\fra p^+$,
we obtain a 1-parameter family 
$(C^{\infty},\{\cdot,\cdot\})_{\varepsilon}$ ($\varepsilon \geq 0$)
of
Poisson structures on
$\fra p^+$,
each of which is compatible with the complex structure
on 
$\fra p^+$
and thus a complex analytic stratified K\"ahler structure,
an ordinary smooth K\"ahler structure for $\varepsilon >0$.
{\sl This family is therefore a resolution of singularities\/}
for the
stratified K\"ahler structure
on $\overline{\Cal O_r}$, with deformation parameter $\varepsilon$.
See (3.2.2) for a very special case.
We have defined $\Cal O^+(\varepsilon)$ to be $G(2\varepsilon z)$
(rather than $G(\varepsilon z)$)
since, as explained in (5.6),
symplectic reduction  
then identifies $\Cal O^+(1)$ with the homogeneous space $G/K$
in a natural fashion;
see also (3.2.2) and (3.4) above. 
\smallskip
By Theorem 5.3.3, in each of the standard cases
($\fra g =\fra{sp}(\ell,\Bobb R),
\fra{su}(p,q),\fra{so}^*(2n)$),
the closure $\overline{\Cal O_r}$
of the principal holomorphic nilpotent orbit
arises via K\"ahler reduction
relative to the corresponding dual group
$\KKK$ (which equals, respectively, $\roman O(\ell,\Bobb R)$, 
$\roman U(p+q,\Bobb R)$,
$\roman {Sp}(n)$).
The semisimple holomorphic orbits
$\Cal O^+(\varepsilon)$
arise as 
well by reduction
relative to the dual group,
but with reference to the corresponding element
$\varepsilon J_V \in \fra {\kkk}$,
as explained in (5.6) above and,
in this fashion,
inherit K\"ahler structures which
coincide with the K\"ahler structures 
whose underlying complex analytic structure is
induced
from the projection  to the constituent $\fra p$
in the Cartan decomposition of $\fra g$.
When the deformation parameter $\varepsilon$ 
tends to zero,
$\varepsilon J_V \in \fra {\kkk}$ tends to the origin;
this shows once again how 
the K\"ahler manifolds $\Cal O^+(\varepsilon)$
contract to the stratified K\"ahler space
$\overline{\Cal O_r}$ in such a way that the complex structure
remains fixed.

\smallskip\noindent
{\smc Remark.}
For $\varepsilon >0$, the smooth K\"ahler structure
on the orbit $\Cal O^+(\varepsilon)$
(arising from the projection to
$\fra p$) is 
{\it not\/} the corresponding hermitian symmetric space
structure;
in particular, 
only $K$ acts by isometries,
the curvature tensor is not
parallel,
and the stereographic projection to the corresponding
symmetric domain 
(cf. 3.2.2)
is only a symplectic change of coordinates,
not a holomorphic one. 
Moreover, we conjecture that
the reasoning in (3.2.2) extends to the general case
to the effect that the  sectional curvature of the hermitian connection
on $\Cal O^+(\varepsilon)$ 
is negative
and that,
in the limit,
the hermitian connection
of each stratum of $\overline{\Cal O_r}$
is flat, i.~e. has zero curvature.

%\beginsection 10. Projectivization and exotic projective varieties
\medskip\noindent {\bf 10. Projectivization and exotic projective varieties}
\smallskip\noindent
Let $E$ be a unitary representation  of a compact Lie group $\KKK$,
$E$ being endowed with the standard symplectic form.
The $\KKK$-action on $E$ is well known to be hamiltonian
and to descend to a hamiltonian $\KKK$-action on
the complex projective space $\roman P(E)$
arising from $E$ by projectivization,
$\roman P(E)$ being endowed with the appropriate symplectic form.
One way to see this is to observe that the central copy 
of the circle group $S^1$ in
the unitary group $\roman U(E)$ of $E$
acts on $E$ 
in a hamiltonian fashion, with a momentum mapping
which in physics
arises as
the {\it harmonic oscillator\/} energy hamiltonian.
Symplectic reduction with respect to 
a constant positive value then yields $\roman P(E)$
endowed with the 
appropriate
symplectic form 
(a multiple of the negative of the imaginary part of
the Fubini-Study metric), and the $\KKK$-momentum mapping
on $E$ passes to a
$\KKK$-momentum mapping
on $\roman P(E)$.
\smallskip
We now apply this observation to the circumstances
of Theorem 5.3.3 and maintain the notation established there.
Thus, 
$\fra g$ is one of the classical Lie algebras
$\fra{sp}(\ell,\Bobb R)$, $\fra {su}(p,q)$ ($p \geq q$),
$\fra{so}^*(2n)$, and $r$ denotes its real rank.
\smallskip
The hamiltonian action of
$\KKK= \roman U(V^s,(\cdot,\cdot))$
on 
$W_J = \roman{Hom}_{\Bobb K}(V^s,V)$
descends to a hamiltonian action
on
$\roman P(W_J)$, and the latter carries a smooth (positive) K\"ahler
structure;
we remind the reader that
the requisite momentum mapping on
$W_J$
is that denoted by $\mu_{\KKK}$.
in Section 5 above. 
For $1 \leq s \leq r$,
let $Q_s$ be the resulting $\KKK$-reduced space;
in view of Proposition 4.2,
since
$\roman P(W_J)$ is a K\"ahler manifold,
$Q_s$ inherits a (positive) normal K\"ahler structure
whose
(stratified symplectic) Poisson structure 
is determined
by a choice of positive (energy) value $\eta$ 
for the $S^1$-momentum mapping.
The next result refers to the standard cases (cf. (3.5) above).

\proclaim{Theorem 10.1}
Each
$Q_s$ inherits a 
normal K\"ahler structure in such a way that
\newline\noindent
{\rm (i)} the injections
$Q_1\subseteq Q_2 \subseteq \dots\subseteq Q_r$
are morphisms of
(complex analytic) stratified K\"ahler
spaces; that 
\newline\noindent
{\rm (ii)} as a complex analytic space,
$Q_r$ is the complex projective space $\roman P(\fra p^+)$;
and that 
\newline\noindent
{\rm (iii)}, for $s < r$,
as a complex analytic space,
$Q_s$ is a projective determinantal variety in
$Q_r$.
\endproclaim

\demo{Proof} This is an immediate consequence of Theorem 5.3.3,
combined with Proposition 4.2. \qed 
\enddemo

Even though,
as a complex analytic space, $Q_r$
is ordinary complex projective space
(of complex dimension
$\dim(\fra p^+)-1$)
and hence non-singular,
as a {\it stratified symplectic space\/},
$Q_r$ has {\it singularities\/}, the {\it singular locus 
in the sense of stratified symplectic spaces\/}
being the subspace $Q_{r-1}$.
More precisely,
the reduced Poisson bracket
is symplectic on the top stratum
of $Q_r$
but 
has smaller rank
at every point of 
$Q_{r-1}$.
Hence the reduced Poisson structure 
on $Q_r$
does {\it not\/}
amount to the standard symplectic Poisson structure
on
complex projective space.
We refer to such a projective space
which is endowed with a normal K\"ahler structure
whose underlying stratified symplectic structure
is {\it not\/} standard
as an {\it exotic\/} projective space, and
we call the closure of a stratum
an
{\it exotic projective K\"ahler variety\/}.
The above construction yields
the exotic projective spaces
$Q_r=\roman P(\fra p^+) = \Bobb C \roman P^d$
according to the three standard cases
$(\fra g,r,d)=(\fra{sp}(\ell,\Bobb R),\ell,\frac 12 \ell(\ell+1) -1 )$,
$(\fra g,r,d)=(\fra {su}(p,q), q, pq-1)\, (p\geq q)$, 
$(\fra g,r,d)=(\fra {so}^*(2n),[\frac n2],\frac 12 n(n-1) -1 )$,
and in each case, the 
closures of the strata
constitute  a sequence
$Q_1\subseteq Q_2 \subseteq \dots\subseteq Q_r$
of exotic projective normal K\"ahler varieties
which, as projective varieties, are determinantal varieties
and thus quite explicitly given.

\proclaim{Theorem 10.2}
As a stratified symplectic subspace
of $Q_r$,
each $Q_s$ is a determinantal  subset of
$Q_r$ in the sense that it may be described
by finitely many  real homogeneous determinantal equations.
\endproclaim

\demo{Proof} This is a consequence of Theorem 5.4.1.
We leave the details to the reader. \qed 
\enddemo

\noindent
{\smc Remark 10.3.}
For the case where $\fra g=\fra {sp}(\ell,\Bobb R)$,
as a stratified symplectic space,
for $1 \leq s \leq \ell=r$,
$Q_s$ 
is the reduced phase space
of a system of $\ell$ harmonic oscillators in
$\Bobb R^s$
with 
constant energy and
total angular momentum zero,
and the injection of
$Q_{s-1}$ into  $Q_s$ 
amounts to 
such a system in
$\Bobb R^{s-1}$
being viewed
as a system in
$\Bobb R^s$
via the standard inclusion
of
$\Bobb R^{s-1}$
into 
$\Bobb R^s$.
Furthermore, for $s > \ell$, the reduced phase space
of a system of $\ell$ harmonic oscillators in
$\Bobb R^s$
with 
constant energy and
total angular momentum zero
is, as a stratified symplectic space,
diffeomorphic 
to the reduced phase space
$Q_\ell$
of a system of $\ell$ 
harmonic oscillators in
$\Bobb R^{\ell}$
with 
constant energy and
total angular momentum zero,
the diffeomorphism being induced by
the standard inclusion
of
$\Bobb R^{\ell}$
into 
$\Bobb R^s$.
This diffeomorphism is an isomorphism of 
(complex analytic)
stratified K\"ahler spaces.
\smallskip
Under the circumstances of Theorem 10.1,
the space $Q_s$ arises as well from the 
closure 
$\overline {\Cal O_s}$
of the corresponding holomorphic nilpotent orbit
$\Cal O_s$
by {\it stratified symplectic reduction\/}
in the sense explained in Section 4 above.
This observation carries over to {\it all\/}
Lie algebras of hermitian type, in the following fashion:
Let $(\fra g,z)$ be
a simple Lie algebra of hermitian type,
with Cartan decomposition $\fra g = \fra k \oplus \fra p$.
As before,
let $\HG$ be the adjoint group of $\fra g$ and $K$ the 
maximal compact subgroup thereof
with $\fra k = \roman{Lie}(K)$.
For $1 \leq s \leq r$,
consider $\overline {\Cal O_s}$,
endowed with the complex analytic stratified K\"ahler structure
given in Subsection 3.7.
The (co)adjoint action of
$\HG$ on $\fra {\hg} \cong\fra {\hg}^*$
induces a 
$\HG$-action on
$\overline {\Cal O_s}$ preserving the stratified symplectic structure,
and 
the injection 
$\overline {\Cal O_s} \subseteq \fra {\hg}^*$
is a stratified symplectic space momentum mapping
for this action.
The restriction  
of this
action to 
$K$ 
preserves the complex analytic stratified K\"ahler structure, 
and the composite
of the inclusion with
the projection
$\fra {\hg}^*\to \fra k^*$
is a stratified symplectic space momentum mapping
$\overline {\Cal O_s} @>>> \fra k^*$
for this action.
In particular,
the central 
copy of the
circle group $S^1$  in 
$K \subseteq G$
acts in this way 
on $\overline {\Cal O_s}$
in a hamiltonian fashion;
since
the 1-parameter subgroup of $G=\roman{Ad}(\fra g)$
generated by $z$ is this central copy of $S^1$
(it consists of the points 
$\roman{exp}(tz) \in G$
with $0 \leq t < 2\pi$ where
$\roman{exp} \colon \fra g \to G$
is the exponential map for the adjoint group $G$),
appropriately scaled, the momentum mapping
$\overline {\Cal O_s} @>>> \Bobb R \cong (\roman{Lie}(S^1))^*$
for this $S^1$-action
is 
simply
given by the composite of 
the inclusion of
$\overline {\Cal O_s}$ into 
$\fra {\hg}^*$,
followed by $z$, 
viewed as a linear map from
$\fra {\hg}^*$ to the reals.
For example,
for $\fra g = \fra{sl}(2,\Bobb R)$,
$\roman{exp}(tz) = \rho(\roman e^{tz})$ where
$\roman e^{tz} = \left[\matrix 
\cos{\frac t2} &-\sin{\frac t2} \\
\sin{\frac t2} & \cos{\frac t2} 
\endmatrix
\right] \in \roman{SL}(2,\Bobb R)
$
and where $\rho \colon \roman{SL}(2,\Bobb R) \to \roman{PSL}(2,\Bobb R)$
is the canonical projection.
Thus, in general, for each $1\leq s\leq r$,
stratified symplectic reduction with respect to 
the appropriate
positive 
(energy) value $\eta$
then yields
a complex analytic stratified K\"ahler space
$Q_s$, with its stratified symplectic
Poisson structure determined by the choice of $\eta$;
complex analytically,
the space
$\overline {\Cal O_s}$
is the (affine complex) cone on 
$Q_s$.
In the standard cases, the resulting spaces
$Q_s$
coincide with those described in
Theorem 10.1 (and denoted by $Q_s$ as well).
\smallskip
For the remaining classical case
$\fra g = \fra{so}(2,q)$ ($q \geq 3$),
this construction
leads to an exceedingly attractive
result:
Recall from Theorem 3.6.2 
that $\fra {so}(2,q)$ has two non-trivial
holomorphic nilpotent orbits
$\Cal O_1$ and $\Cal O_2$
whose closures constitute a sequence
${
0 \subseteq \overline{\Cal O_1} 
\subseteq \overline{\Cal O_2} \cong \Bobb C^q
}$
of normal K\"ahler spaces;
here
$\overline{\Cal O_1}$ 
is the complex quadric 
$w_1^2 + \dots + w_q^2 = 0$
in 
$\overline{\Cal O_2} = \Bobb C^q$.
Now
$\HG=\roman{SO}(2,q)^0$,
$K= \roman{SO}(2,\Bobb R)\times \roman{SO}(q,\Bobb R)$,
and 
$\roman{SO}(2,\Bobb R)$
is the requisite
central copy of the
circle group in $K$.
By construction, the $S^1$-reduced space
$Q_1$ sits inside 
the $S^1$-reduced space
$Q_2= \Bobb C \roman P^{q-1}$
as the projective quadric in
$\Bobb C \roman P^{q-1}$
given by the homogeneous equation
$w_1^2 + \dots + w_q^2 = 0$.
This quadric,
being a complex analytic stratified 
K\"ahler space with a single stratum, is necessarily
a smooth K\"ahler manifold.
Summing up, we obtain the following.

\proclaim{Theorem 10.4}
For $m \geq 2$,
$\Bobb C \roman P^m$
carries an exotic normal K\"ahler structure
whose singular locus (in the sense of stratified
K\"ahler spaces)
is a quadric $Q$ 
in $\Bobb C \roman P^m$
in such a way that
the exotic  K\"ahler structure
on $\Bobb C \roman P^m$
restricts to an ordinary
K\"ahler structure
on  $Q$. \qed
\endproclaim

The exceptional cases
lead 
again to 
{\it new\/} geometric insight:
Recall from Section 8
that $\fra g =\fra{e}_{7(-25)}$  has three non-trivial
holomorphic nilpotent orbits
$\Cal O_1$, $\Cal O_2$, and $\Cal O_3$
whose closures constitute a sequence
${
0 \subseteq \overline{\Cal O_1} 
\subseteq \overline{\Cal O_2} 
\subseteq \overline{\Cal O_3} 
}$
of complex analytic stratified
K\"ahler spaces;
furthermore, under the 
projection from $\fra g$ to $\fra p$,
followed by the complex linear isomorphism
from $\fra p$ onto $\fra p^+$,
as a stratified space,
the closure $\overline{\Cal O_3}$
of the principal holomorphic nilpotent orbit
is identified with $\Cal H_3(\Bobb O_{\Bobb C})$. Now
$\HG=\roman{E}_{7(-25)}$,
$K$ is locally isomorphic to 
$\roman{E}_{6(-78)}\times \roman{SO}(2,\Bobb R)$, 
and the constituent
$\roman{SO}(2,\Bobb R)$
is the requisite
central copy of the
circle group in $K$
(here \lq\lq requisite\rq\rq\ is intended to refer to the
fact that, in order for $K$ to be the stabilizer of
a point of a symmetric domain, it must have a 1-dimensional center).
By construction, the $S^1$-reduced space
$Q_2$ 
(here and below \lq\lq reduced\rq\rq\ is understood
in the sense of stratified symplectic reduction)
sits inside 
the $S^1$-reduced space
$Q_3= \roman P(\Cal H_3(\Bobb O^{\Bobb C})) \cong \Bobb C \roman P^{26}$
as the projective complex 
generic norm (cubic) hypersurface; it is the projectivization
of the affine Freudenthal cubic $\Cal F$ mentioned in Section 8 above.
Furthermore,
the $S^1$-reduced space
$Q_1$ sits inside 
$Q_2$
as its singular locus
in the sense of stratified K\"ahler spaces.
Thus the projective cubic $Q_2$
inherits a complex analytic stratified 
K\"ahler structure with two strata.
We note that the generic norm hypersurface 
has been shown by Jacobson to be rational
\cite\jacobstw.
Summing up, we obtain the following.

\proclaim{Theorem 10.5}
The complex projective space 
of dimension {\rm 26}
carries an exotic K\"ahler structure
whose singular locus (in the sense of stratified
K\"ahler spaces)
is 
the projective generic norm (cubic) hypersurface,
and the 
latter
inherits a complex analytic stratified 
K\"ahler structure with two strata,
having as singular stratum $Q_1$
(in the sense of stratified K\"ahler spaces)
the projectivization
$\roman P(\widetilde{\Cal O_1}) 
\subseteq \roman P(\fra p^+) \cong \Bobb C \roman P^{26}$
of the lowest non-zero stratum
$\widetilde{\Cal O_1} \subseteq \Cal H_3(\Bobb O_{\Bobb C})$
(i.~e. Jordan rank 1 stratum or, equivalently, lowest 
non-zero $K^{\Bobb C}$-stratum).
The latter
amounts to the compact hermitian symmetric space
$\roman{Ad}(\fra e_{6(-78)})
\big/(\roman{SO}(10,\Bobb R) \cdot \roman{SO}(2,\Bobb R))$
of complex dimension {\rm 16}, and 
the exotic  K\"ahler structure
on $\Bobb C \roman P^{26}$
restricts to an ordinary
K\"ahler structure
on  $Q_1$. \qed
\endproclaim

Likewise $\fra{e}_{6(-14)}$  has two non-trivial
holomorphic nilpotent orbits
$\Cal O_1$ and $\Cal O_2$
whose closures constitute a sequence
${
0 \subseteq \overline{\Cal O_1} 
\subseteq \overline{\Cal O_2} \cong \fra p^+
}$
of complex analytic stratified K\"ahler spaces,
and the pre-homogeneous space structure
on
$\fra p^+ \cong \Bobb C^{16}$
restricts to the positive half-spin representation
$\Cal S_+^{\Bobb C}$ of $\roman{Spin}(10,\Bobb C)$
(cf. Section 8 above for the notation).

\proclaim{Theorem 10.6}
The complex projective space 
of dimension {\rm 15}
carries an exotic K\"ahler structure
whose singular locus $Q$ (in the sense of stratified
K\"ahler spaces)
is 
the projectivization
$\roman P(\widetilde{\Cal O_1}) \subseteq \roman P(\fra p^+) 
\cong \Bobb C \roman P^{15}$
of the space $\widetilde{\Cal O_1}$ of pure spinors
in $\fra p^+\cong \Cal S_+^{\Bobb C}$,
and $Q$ amounts to the compact hermitian symmetric space
$\roman{SO}(10,\Bobb R)
\big/\roman{U}(5)$
of complex dimension {\rm 10}. 
The exotic  K\"ahler structure
on $\Bobb C \roman P^{15}$
restricts to an ordinary
K\"ahler structure
on $Q$. \qed
\endproclaim

Indeed, it is well known that
the projective space associated with the pure spinors
for the  (16-dimensional)
positive half-spin representation
$\Cal S_+^{\Bobb C}$ 
of $\roman{Spin}(10,\Bobb C)$
is the compact hermitian symmetric space
$\roman{SO}(10,\Bobb R)
\big/\roman{U}(5)$;
see e.~g.
\cite\lawmiboo\ (IV.9.7~p.~337).
For intelligibility, we recall that,
by Theorem 8.5.3,
under the homeomorphism from
$\overline {\Cal O_2}$ to
$\fra p^+$,
the orbit 
$\Cal O_1$ is identified with the space 
$\widetilde{\Cal O_1}$
of pure spinors
in 
$\fra p^+\cong \Cal S_+^{\Bobb C}$.

\smallskip
\noindent
{\smc 10.7. Projective normality.}
A projective variety $X \subseteq \Bobb C \roman P^d$ is said to be
{\it projectively normal\/}
(for the given embedding)
provided its homogeneous coordinate ring is integrally closed.
This is equivalent to requiring that
the cone on $X$ (in 
$\Bobb C^{d+1}$) be normal as an affine variety.
In the classical cases since,
for $1 \leq s \leq r$,
complex analytically, the space
$\overline {\Cal O_s}$ is a normal affine variety
(Theorems 3.5.5, 3.6.3, 5.3.3),
each $Q_s$ is projectively normal.
Furthermore, 
cf. Theorem 8.4.2,
the affine Freudenthal
cubic 
is normal whence the projective Freudenthal cubic is projectively normal.
We now settle the remaining cases, the lowest non-trivial stratum
in each of the two exceptional cases.
We are indebted to L. Manivel for having communicated to us the following.

\proclaim{Proposition 10.7.1}
The closed orbit $X$ in $\roman P(V)$
of an irreducible complex representation
$V$ of a complex simple Lie group 
is projectively normal.
\endproclaim

\demo{Proof}
Indeed, in view of the Borel-Weil theorem,
if $v$ is the highest weight of the irreducible representation
$V^* = \Gamma(X,\Cal O_X(1))$,
for $n \geq 1$,
the highest weight of the irreducible representation
$V^* = \Gamma(X,\Cal O_X(n))$ equals $nv$.
Consequently 
for $n \geq 1$
the natural map
$\Gamma(\roman P(V),\Cal O_{\roman P(V)}(n))
\to \Gamma(X,\Cal O_X(n))$
is surjective whence, cf.
Ex. II.5.14 on p.~126 of \cite\hartsboo, 
the projective variety $X$ is projectively normal. \qed
\enddemo

Given a simple Lie algebra of hermitian type $(\fra g,z)$,
we may apply this Proposition
to the irreducible representation of $K^{\Bobb C}$ on
$\fra p \cong \fra p^+$,
the closed orbit $X$ being the orbit denoted earlier by $Q_1$.
Thus we conclude that $Q_1$ is projectively normal
(which we already know in the classical cases).
Consequently,
for $\fra{e}_{6(-14)}$,  
the projectivization
$Q_1=\roman P(\widetilde{\Cal O_1}) \subseteq \roman P(\fra p^+) 
\cong \Bobb C \roman P^{15}$
of the space $\widetilde{\Cal O_1}$ of pure spinors
in $\Cal S_+^{\Bobb C}$ 
is projectively normal 
(for the embedding into
$\roman P(\fra p^+)$)
whence,
complex analytically,
$\overline{\Cal O_1}$ is a normal affine variety.
Likewise, for $\fra{e}_{7(-25)}$,  
the projectivization
$Q_1$ of
$\overline{\Cal O_1}$ is 
projectively normal
(for the embedding into
$Q_3= \roman P(\Cal H_3(\Bobb O^{\Bobb C})) \cong \Bobb C \roman P^{26}$)
whence $\overline{\Cal O_1}$
is normal as an affine complex variety.

\smallskip
\noindent
{\smc Remark 10.8.}
A Kodaira embedding 
of a smooth projective variety $N$
carrying a K\"ahler structure will in general not 
be compatible with the symplectic structures,
where projective  space is endowed with the Fubini-Studi metric.
Furthermore,
given an arbitrary smooth symplectic manifold $X$ and a smooth
symplectic submanifold $Y$,
the symplectic Poisson algebra on $X$ will {\it not\/}
restrict to the symplectic Poisson algebra on $Y$;
indeed this observation prompted {\smc Dirac\/} to introduce 
what is nowadays called the {\it Dirac bracket\/}.
The above results indicate that, perhaps,
the right question to ask
is this:
Given $N$ and an embedding thereof into
some complex projective space 
is there
an exotic K\"ahler structure 
on that projective space
which restricts to the
K\"ahler structure 
on $N$?
More generally,
we could ask this question for an arbitray
compact complex analytic stratified
K\"ahler space $N$,
and the above examples show that this situation indeed occurs.
In particular:
A  moduli space of the kind 
explained in Example 4.3 above
is known to be complex analytically a normal projective variety.
Thus we are led to
the following question:
Given such a
space $N$, 
{\sl is there
an embedding of $N$ into some
complex projective space
endowed with an exotic K\"ahler structure 
such that 
this
exotic K\"ahler structure 
restricts to the
normal K\"ahler structure 
on $N$\/} explained in Example 4.3? 

\medskip\noindent
{\bf 11. Comparison with other notions of K\"ahler space with 
singularities}
\smallskip\noindent
Let $X$ be a reduced complex analytic space
(where \lq\lq reduced\rq\rq\ is understood in the sense of complex analytic
spaces, {\it not\/} in the sense of symplectic geometry), 
let $N$ be its singular locus,
and consider a hermitian metric $g$ on $X\setminus N$.
In \cite\graueone\ (p.~346) {\smc Grauert}
refers to such a metric as
a {\it K\"ahler metric\/}
on $X$ provided every point $x$ of $X$ has a neighborhood
$U$
and a strongly plurisubharmonic function
$p$ defined on $U$
in such a way that, on $U \setminus N$,
the metric $g$ coincides with the Levi form of $p$,
that is to say,
in suitable coordinates
$z_1,\dots, z_m$,
$$
g_{j,\overline k} = \frac {\partial^2 p}{\partial z_j\partial \overline z_k}.
$$
As usual,
the function $p$ is then referred to as a {\it local K\"ahler
potential\/} for $g$.
\smallskip
On the other hand,
given a stratified K\"ahler space
of complex dimension $n$,
on its top stratum,
in local coordinates
$z_1,\dots, z_m$, 
the Poisson structure and
K\"ahler metric $g$ 
satisfy the identity
$$
\sum_{\ell} \{\overline z_k, z_{\ell} \} g_{\ell,\overline m} 
= 2 i \delta_{k,m},
\quad \text{for every} \quad 1 \leq k,m \leq n,
$$
that is, the functions
$
\frac 1{2 i}\{\overline z_k, z_{\ell} \}
$
yield the coefficients of 
the inverse matrix function
$g^{-1}$ 
(in these coordinates).
Since the Poisson structure is defined on the whole space,
even though the K\"ahler metric is {\it not\/},
we can somewhat loosely interpret the above formula by saying that
the \lq\lq inverse\rq\rq\ of the K\"ahler metric
{\it is\/} defined, in fact, even in such a way
that the restriction of this \lq\lq inverse\rq\rq\  to any lower stratum
is there invertible 
so that the inverse of the \lq\lq inverse\rq\rq, in turn,
defines a K\"ahler metric  on the stratum as well.
The Poisson structure then encapsulates the mutual positions
of these K\"ahler structures.
\smallskip
Grauert's notion has been refined
in \cite\varoucha, cf. also
\cite\heihuclo,
to that of
a {\it stratified K\"ahlerian space\/};
such a structure encapsulates the mutual positions of
K\"ahler structures on pieces of a decomposition,
but in a way different from 
our notion of complex analytic stratified K\"ahler space.
To reproduce 
stratified K\"ahlerian spaces briefly,
let $X$ be a reduced complex analytic space,
with a complex analytic stratification
by complex analytic manifolds
(in general genuinely finer than the ordinary complex analytic
stratification).
Given such a space,
a {\it K\"ahlerian\/} structure on $X$ is 
a collection of 
\lq\lq strictly plurisubharmonic functions\rq\rq\ 
$\{\phi_\alpha\}$
which are defined and continuous on the constituents $U_\alpha$ of an open  
covering $\{U_\alpha\}$ of $X$
and smooth on the intersection of each stratum with
$U_\alpha$
in such a way that the differences $u_{\alpha,\beta}$ are pluriharmonic
in the strict sense on the intersections
$U_{\alpha,\beta} =U_\alpha \cap U_\beta$, that is, there is 
$f_{\alpha,\beta} \in \Cal O(U_{\alpha,\beta})$
such that 
$u_{\alpha,\beta}= \roman{Re}(f_{\alpha,\beta})$.
There are various ways to make precise
the idea of
strictly plurisubharmonic function on $X$:
the strongest one is the notion of
\lq\lq strictly plurisubharmonic in the sense of perturbations\rq\rq\ 
(which we do not reproduce here) 
and a somewhat weaker on is that of
\lq\lq stratum-wise
strictly plurisubharmonic\rq\rq\ 
which means that
the restriction to 
(the appropriate open subset of)
any stratum
of $X$
is smooth in such a way
that its Levi-form
is positive definite 
(i.~e. a hermitian
structure);
see \cite\heihuclo\ for details.

\smallskip
To explain the
relationship between stratified K\"ahler spaces
and stratified K\"ahlerian spaces,
let
$(X, C^{\infty}(X), \{\cdot,\cdot\})$ be a 
stratified symplectic space.
Write $A=C^{\infty}(X)$, let
$(A,L)= (C^{\infty}(X), \roman{Der}(A))$,
with its Lie-Rinehart algebra structure,
and consider the 
differential graded algebra
$(\roman{Alt}_A(L,A),d)$ determined by it
\cite\poiscoho;
when $X$ is a smooth manifold
and
$C^{\infty}(X)$
the algebra of ordinary smooth functions,
$(\roman{Alt}_A(L,A),d)$
is the differential graded algebra of ordinary differential forms for $X$.
For a general
stratified symplectic space
$(X, C^{\infty}(X), \{\cdot,\cdot\})$,
the Poisson structure
induces a morphism
${
\pi^{\sharp}_{\{\cdot,\cdot\}}
\colon (A,\Omega^1(X)_{\{\cdot,\cdot\}})
@>>>
(A,L)
}$
of Lie-Rinehart algebras
(cf. (1.1) above)
and hence a morphism
$$
(\roman{Alt}_A(L,A),d)
@>>>
(\roman{Alt}_A(\Omega^1(X)_{\{\cdot,\cdot\}},A),d_{\{\cdot,\cdot\}})
$$
of differential graded algebras
\cite\poiscoho;
the Poisson 2-form $\pi_{\{\cdot,\cdot\}}$ is then an element
of $\roman{Alt}^2_A(\Omega^1(X)_{\{\cdot,\cdot\}},A)$,
cf. Section 1 above.
\smallskip 
On the other hand,
consider a stratified K\"ahlerien space $Y$,
with its 
ordinary
complex analytic stratification.
The latter determines a
smooth structure $C_{\roman{min}}^{\infty}(Y)$
which we refer to as the {\it minimal one\/};
this structure is, roughly speaking, locally determined
by the embedding of the corresponding complex analytic germ into
a complex domain \cite\varoucha;
the requisite smooth structure on
$Y$ is then obtained 
(locally)
by restricting the ordinary smooth structure
on the domain.
(Varouchas explains this in terms of sheaves.)
This smooth structure determines a
(complex) de Rham algebra $\Omega_Y = \Lambda \Omega^1_Y$,
$\Omega^1_Y$ being defined by an exact sequence of the kind
$$
I/I^2 @>>> 
C_{\roman{min}}^{\infty}(Y,\Bobb C)\otimes_{C^{\infty}(X,\Bobb C)}\Omega^1_X 
@>>>
\Omega^1_Y
@>>> 0
$$
where $Y \subseteq X$ is an embedding into a complex domain $X$
and where $I$ is the ideal (sheaf) of (germs of) functions
which vanish on $Y$.
The complex analytic structure of $Y$
entails a bigrading
$\Omega^{*,*}_Y$, just as for ordinary smooth complex manifolds,
and
the stratified K\"ahlerien structure determines a 2-form
$\omega \in \Omega^{1,1}_Y$ which,
when $Y$ is an ordinary smooth K\"ahler manifold, boils down to
its standard symplectic form.
\smallskip
Return to the circumstances of Proposition 4.2.
The reduced space $N^{\roman{red}}$ coming into play there
inherits,
on the one hand, 
a smooth structure
$C^{\infty}(N^{\roman{red}},\Bobb R)$
which underlies its
stratified symplectic Poisson algebra
and, furthermore,
a complex analytic stratified K\"ahler space structure.
On the other hand,
by (3.4) of \cite\heinloos,
$N^{\roman{red}}$
also inherits a stratified K\"ahlerien structure.
The corresponding minimal smooth structure 
$C_{\roman{min}}^{\infty}(N^{\roman{red}})$
(explained above)
will in general be a proper subalgebra
of  $C^{\infty}(N^{\roman{red}},\Bobb R)$,
and the space $\Omega^1_{N^{\roman{red}}}$
of 1-forms with reference to the complexified algebra
$C_{\roman{min}}^{\infty}(N^{\roman{red}},\Bobb C)$
will decompose into a direct sum of
$\Omega^{0,1}_{N^{\roman{red}}}$
and $\Omega^{1,0}_{N^{\roman{red}}}$.
Notice that, 
cf. Remark 3.7.3, 
a general complex analytic stratified K\"ahler structure
does {\it not\/}
give rise to such a decomposition; 
indeed, cf. (3.7.4), the sum of
$P$ and $\overline P$ does not exhaust the space of 1-forms.
To explore the mutual relationship
between the 
complex analytic stratified K\"ahler 
and stratified K\"ahlerien
structures (if any),
write 
$A=C^{\infty}(N^{\roman{red}},\Bobb C)$
and $L=\roman{Der}(A))$,
and consider the 
differential graded algebra
$(\roman{Alt}_A(L,A),d)$ determined by it,
that is, the Rinehart complex for the
Lie-Rinehart algebra
$(C^{\infty}(N^{\roman{red}},\Bobb C), \roman{Der}(A))$.
With reference to the
ordinary complex analytic stratification,
the stratified K\"ahlerien
structure 
determines a 2-form
$\omega \in \Omega^{1,1}_{N^{\roman{red}}}
\subseteq
\Omega^2_{N^{\roman{red}}}
$
whereas
the complex analytic stratified
K\"ahler structure
determines the
2-form
$\pi_{\{\cdot,\cdot\}} \in 
\roman{Alt}^2_A(\Omega^1(N^{\roman{red}})_{\{\cdot,\cdot\}},A)$.
In the case of an ordinary complex K\"ahler manifold,
we can identify
$\omega$ with $\pi_{\{\cdot,\cdot\}}$
but for 
a general
stratified K\"ahlerien structure
(i.~e. having underlying stratification finer than the
ordinary complex analytic one),
there is no obvious 
replacement for
the constituent
$\Omega^{1,1}_{N^{\roman{red}}}$
and hence no obvious
candidate for 
a
2-form of the kind $\omega$
whereas
given a complex analytic stratified
K\"ahler structure,
the 2-form
$\pi_{\{\cdot,\cdot\}}$
is always defined.
Furthermore,
when there is a
2-form of the kind $\omega$
with reference to the minimal structure
$C_{\roman{min}}^{\infty}(N^{\roman{red}})$,
its image in
$\roman{Alt}^2_A(\Omega^1(N^{\roman{red}})_{\{\cdot,\cdot\}},A)$
will not in general coincide with
the 2-form
$\pi_{\{\cdot,\cdot\}}$.
To justify this claim,
return to the circumstances of Theorem 3.5.4.
In view of 
Theorem 3.5.5,
complex analytically, $\overline {\Cal O_r}$ amounts to
$\fra p^+$ and, for $s<r$,
the complex analytic structure on
$\overline {\Cal O_s}$ arises from the embedding of
$\overline {\Cal O_s}$ into
$\fra p^+$.
Moreover,
the minimal smooth structure 
$C_{\roman{min}}^{\infty}(\overline {\Cal O_s})$
(i.~e. algebra of continuous functions)
on $\overline {\Cal O_s}$
is determined by the embedding 
of $\overline {\Cal O_s}$
into
$\fra p^+$, too;
this is the structure coming into play in \cite\varoucha.
On the other hand,
the stratified symplectic structure
of $\overline {\Cal O_s}$
results from its embedding  into
$\fra g\cong \fra g^*$;
a choice of basis
of $\fra g$
determines a system of coordinate functions
on $\overline {\Cal O_s}$,
and the last statement of
Theorem 5.4.1 says that
the stratified symplectic Poisson  algebra 
on $\overline {\Cal O_s}$
cannot be generated by
fewer functions.
Furthermore,
$C_{\roman{min}}^{\infty}(\overline {\Cal O_s})$
being viewed as a subalgebra of
$C^{\infty}(\overline {\Cal O_s},\Bobb C)$,
Poisson brackets of elements of
$C_{\roman{min}}^{\infty}(\overline {\Cal O_s})$
do not necessarily lie                          
in 
$C_{\roman{min}}^{\infty}(\overline {\Cal O_s})$.
Summing up, we see that
the stratified symplectic structure
of $\overline {\Cal O_s}$
results from its embedding  into
$\fra g$ 
whereas its complex analytic structure
comes from  the composite
of this embedding
with the 
projection from
$\fra g$ to $\fra p$.
Even though the projection 
from
$\fra g$ to $\fra p$
is {\it not\/}
injective,
the restriction thereof to
$\overline {\Cal O_s}$ is still injective and 
turns
$\overline {\Cal O_s}$
into a complex analytic space.
Thus, the real and complex structures of
$\overline {\Cal O_s}$
are obtained from the embeddings
thereof into $\fra g$
and
$\fra p$,
respectively,
and the relationship
between the two structures
comes from
the projection from $\fra g$ to $\fra p$.
Perhaps this kind of situation should be taken
as a local model for a complex analytic stratified K\"ahler space
in general.
\smallskip
This discussion raises the following questions:
\newline\noindent
1) Given a complex analytic stratified K\"ahler structure
on a stratified symplectic space $X$,
do there always exist 
local K\"ahler potentials which determine
a stratified K\"ahlerian structure?
\newline\noindent
2) Consider a complex analytic space $X$,
endowed with a general stratified K\"ahlerien structure
where \lq\lq general\rq\rq\ means that the underlying stratification
may be finer than
the ordinary complex analytic one,
and let
$C^{\infty}(X)$
be  
the smooth structure 
whose elements are continuous functions
on $X$ which, restricted to each stratum, are smooth.
Given two functions $f,h \in C^{\infty}(X)$,
taking their Poisson bracket
$\{f,h\}$
stratum-wise,
we obtain a function on $X$ which is necessarily smooth on every stratum.
The question is whether $\{f,h\}$ yields
a {\it continuous\/} function on $X$;
if so,
$C^{\infty}(X)$
would be closed under the Poisson bracket
and, in view of Theorem 2.5, the latter
would turn $X$
into a complex analytic stratified K\"ahler space.
\smallskip
The answers to both questions are \lq\lq yes\rq\rq\ for a space $X$ which
arises by symplectic reduction from a (positive) K\"ahler manifold
with respect to a compact group.
Thus to find a potential counterexample to (1), we would have to
look for a complex analytic stratified K\"ahler space
which does not arise in this way from an ordinary K\"ahler manifold;
for $q\geq 5$
{\sl the exotic complex projective space
described in \/} Theorem 10.4 {\sl above
having as its singular locus a quadric
seems indeed to be such an example,
and
we do 
not know
whether there  exist 
local K\"ahler potentials in this case which determine
a stratified K\"ahlerian structure\/}.
The same kind of question can be asked for the
two exceptional cases
$\fra e_{6(-14)}$ and $\fra e_{7(-25)}$.
Likewise,
for a potential counterexample to (2), we would need
a stratified K\"ahlerien 
space which does not arise from an ordinary K\"ahler manifold
by reduction.

\bigskip
%\vfill\eject
\centerline{\smc References}
\medskip
\widestnumber\key{999}

\ref \no  \armcusgo
\by J. M. Arms,  R. Cushman, and M. J. Gotay
\paper  A universal reduction procedure for Hamiltonian group actions
\paperinfo in: The geometry of Hamiltonian systems, T. Ratiu, ed.
\jour MSRI Publ. 
\vol 20
\pages 33--51
\yr 1991
\publ Springer
\publaddr Berlin $\cdot$ Heidelberg $\cdot$ New York $\cdot$ Tokyo
\endref

\ref \no \armamonc
\by J. M. Arms, J. E. Marsden, and V. Moncrief
\paper  Symmetry and bifurcation of moment mappings
\jour Comm. Math. Phys.
\vol 78
\yr 1981
\pages  455--478
\endref

\ref \no  \atibottw
\by M. F. Atiyah and R. Bott
\paper The Yang-Mills equations over Riemann surfaces
\jour Phil. Trans. R. Soc. London  A
\vol 308
\yr 1982
\pages  523--615
\endref

\ref \no \barbsepa
\by D. Barbasch and M. R. Sepanski
\paper Closure ordering and the Kostant-Sekiguchi correspondence
\jour Proc. Amer. Math. Soc.
\vol 126
\yr 1998
\pages 311--317
\endref

\ref \no \beardboo
\by A. F. Beardon
\book The Geometry of Discrete Groups
\bookinfo Graduate Texts in Mathematics
\publ Springer
\publaddr Berlin $\cdot$  Heidelberg $\cdot$ New York $\cdot$ Tokyo
\yr 1983
\endref

\ref \no \brunvett
\by W. Bruns and U. Vetter
\book Determinantal Rings
\bookinfo Lecture Notes in Mathematics, Vol. 1327
\yr 1988
\publ Springer 
\publaddr Berlin $\cdot$ Heidelberg $\cdot$ New York
\endref

\ref \no \burcush
\by N. Bourgoyne and R. Cushman
\paper Conjugacy classes in linear groups
\jour J. of Algebra
\vol 44
\yr 1977
\pages 339--362
\endref

\ref \no  \chevaone
\by C. Chevalley
\paper Invariants of finite groups generated by reflections
\jour Amer. J. of Math. 
\vol 77
\yr 1955
\pages 778--782
\endref

\ref \no \collmcgo
\by D.~H.~Collingwood and W. M. McGovern
\book Nilpotent orbits in semisimple Lie algebras
\bookinfo Mathematics Series
\publ Van Nostrand Reinhold
\publaddr New York
\yr 1993
\endref

\ref \no  \cushmone
\by R. Cushman
\paper The momentum mapping of the harmonic oscillator
\jour Symposia Math.
\vol 14
\yr 1975
\pages 323--342
\endref

\ref \no  \djokotwo
\by D. Djokovi\'c
\paper Classification of nilpotent elements in simple
exceptional real Lie algebras of inner type and descriptions
of their centralizers
\jour J. of Algebra
\vol 112
\yr 1988
\pages 503--524
\endref

\ref \no  \elkinone
\by G. B. Elkington
\paper Centralizers of unipotent elements in semisimple algebraic
groups
\jour J. of Algebra
\vol 23
\yr 1972
\pages 137--163
\endref

\ref \no \farakora
\by J. Faraut and A. Koranyi
\book Analysis On Symmetric Cones
\bookinfo Oxford Mathematical Monographs
 No. 114
\publ Oxford University Press
\publaddr Oxford U.K.
\yr 1994
\endref

\ref \no \freudone
\by H. Freudenthal
\paper Beziehungen der $E_7$ und $E_8$ zur Oktavenebene
\jour Indagationes Math.
\vol 16
\yr 1954
\pages 218--230
\endref

\ref \no \graueone
\by H. Grauert 
\paper \"Uber Modifikationen und exzeptionelle
analytische Mengen
\jour  Math. Ann.
\vol 146
\yr 1962
\pages 331--368
\endref

\ref \no \grauremm
\by H. Grauert and R. Remmert
\book Coherent Analytic Sheaves
\bookinfo Grundlehren der mathematischen Wissenschaften, vol. 265
\publ Springer
\publaddr Berlin $\cdot$ Heidelberg $\cdot$ New York $\cdot$ Tokyo
\yr 1984
\endref

\ref \no  \guistebo
\by V. Guillemin and S. Sternberg
\book Symplectic techniques in Physics
\publ Cambridge University Press
\publaddr London/New York
\yr 1984
\endref

\ref \no  \guistetw
\by V. W. Guillemin and S. Sternberg
\paper Geometric quantization and multiplicities of group representations
\jour Invent. Math.
\vol 67
\yr 1982
\pages 515--538
\endref

\ref \no \guhujewe
\by K. Guruprasad, J. Huebschmann, L. Jeffrey, and A. Weinstein
\paper Group systems, groupoids, and moduli spaces
of parabolic bundles
\jour Duke Math. J.
\vol 89
\yr 1997
\pages 377--412
\finalinfo {\tt dg-ga/9510006}
\endref

\ref \no \hartsboo
\by  R. Hartshorne
\book Algebraic Geometry
\bookinfo Graduate texts in Mathematics
 No. 52
\publ Springer
\publaddr Berlin-G\"ottingen-Heidelberg
\yr 1977
\endref

\ref \no \heihuclo
\by P. Heinzner, A. Huckleberry, and F. Loose
\paper K\"ahlerien extensions of the symplectic reduction 
\jour J. f\"ur die reine und angewandte Mathematik
\vol 455
\yr 1994
\pages 123--140
\endref

\ref \no  \heinloos
\by P. Heinzner and F. Loose
\paper Reduction of complex Hamiltonian $G$-spaces
\jour  Geom. and Func. Analysis
\vol 4
\yr 1994
\pages  288--297
\endref

\ref \no \helgaboo
\by S. Helgason
\book Differential Geometry, Lie Groups, and Symmetric Spaces
\bookinfo Pure and Applied Mathematics, 
A series of monographs and textbooks
\publ Academic Press
\publaddr New York $\cdot$ London $\cdot$ Toronto 
$\cdot$ Sydney $\cdot$ San Francisco
\yr 1978     
\endref

\ref \no \hilneors
\by J. Hilgert, K.-H. Neeb, and B. Oersted
\paper The geometry of nilpotent coadjoint orbits
of convex type in hermitian Lie algebras
\jour J. of Lie Theory
\vol 4
\yr 1994
\pages 185--235
\endref

\ref \no \hocheago
\by M. Hochster and J. A. Eagon
\paper Cohen-Macaulay rings, invariant theory, and the generic perfection
of  determinantal loci
\jour Amer. J. of Math.
\vol 93
\yr 1971
\pages 1020--1058
\endref

\ref \no  \howeone
\by R. Howe
\paper Remarks on classical invariant theory
\jour  Trans. Amer. Math. Soc.
\vol 313
\yr 1989
\pages  539--570
\endref

\ref \no \howetwo
\by R. Howe
\paper Transcending classical invariant theory
\jour  J. Amer. Math. Soc.
\vol 2
\yr 1989
\pages  535--552
\endref

\ref \no \poiscoho
\by J. Huebschmann
\paper Poisson cohomology and quantization
\jour 
J. f\"ur die reine und angewandte Mathematik
\vol  408 
\yr 1990
\pages 57--113
\endref

\ref \no  \souriau
\by J. Huebschmann
\paper On the quantization of Poisson algebras
\paperinfo Symplectic Geometry and Mathematical Physics,
Actes du colloque en l'honneur de Jean-Marie Souriau,
P. Donato, C. Duval, J. Elhadad, G.M. Tuynman, eds.;
Progress in Mathematics, Vol. 99
\publ Birkh\"auser
\publaddr Boston $\cdot$ Basel $\cdot$ Berlin
\yr 1991
\pages 204--233
\endref

\ref \no \singula
\by J. Huebschmann
\paper The singularities of Yang-Mills connections
for bundles on a surface. I. The local model
\jour Math. Z. 
\vol 220
\yr 1995
\pages 595--605
\finalinfo {\tt dg-ga/9411006}
\endref

\ref \no \singulat
\by J. Huebschmann
\paper The singularities of Yang-Mills connections
for bundles on a surface. II. The stratification
\jour Math. Z. 
\vol 221
\yr 1996
\pages 83--92
\finalinfo {\tt dg-ga/9411007}
\endref

\ref \no \locpois
\by J. Huebschmann
\paper Poisson geometry of flat connections 
for {\rm SU(2)}-bundles on surfaces
\jour Math. Z.
\vol 221
\yr 1996
\pages 243--259
\finalinfo {\tt hep-th/9312113}
\endref

\ref \no \poisson
\by J. Huebschmann
\paper 
Poisson
structures on certain
moduli spaces 
for bundles on a surface
\jour Annales de l'Institut Fourier
\vol 45
\yr 1995
\pages 65--91
\finalinfo {\tt dg-ga/9411009}
\endref

\ref \no \modus
\by J. Huebschmann
\paper Symplectic and Poisson structures of certain moduli spaces
\jour Duke Math. J.
\vol 80
\yr 1995
\pages 737--756
\finalinfo {\tt hep-th/9312112}
\endref

\ref \no \modustwo
\by J. Huebschmann
\paper Symplectic and Poisson structures of certain moduli spaces. II.
Projective
representations of cocompact planar discrete groups
\jour Duke Math. J.
\vol 80
\yr 1995
\pages 757--770
\finalinfo {\tt dg-ga/9412003}
\endref

\ref \no \poismod
\by J. Huebschmann
\paper On the variation of the Poisson structures of certain moduli spaces
\jour Math. Ann.
\vol 319
\yr 2001
\pages 267--310
\finalinfo {\tt dg-ga/9710033}
\endref

\ref \no \smooth
\by J. Huebschmann
\paper 
Smooth structures on moduli spaces of central Yang-Mills connections 
for bundles on a surface
\jour J. of Pure and Applied Algebra
\vol 126
\yr 1998
\pages 183--221
\finalinfo {\tt dg-ga/9411008}
\endref

\ref \no \srni
\by J. Huebschmann
\paper
Poisson geometry of certain
moduli spaces
\paperinfo
Lectures delivered at the \lq\lq 14th Winter School\rq\rq, Srni,
Czeque Republic,
January 1994
\jour Rendiconti del Circolo Matematico di Palermo, Serie II
\vol 39
\yr 1996
\pages 15--35
\endref

\ref \no \claustha
\by J. Huebschmann
\paper On the Poisson geometry of certain moduli spaces
\paperinfo in: Proceedings of an international workshop on
\lq\lq Lie theory and its applications in physics\rq\rq,
Clausthal, 1995
H. D. Doebner, V. K. Dobrev, J. Hilgert, eds.
\publ World Scientific
\publaddr Singapore $\cdot$
New Jersey $\cdot$
London $\cdot$
Hong Kong 
\pages 89--101
\yr 1996
\endref

\ref \no  \oberwork
\by J. Huebschmann
\paper Singularities and Poisson geometry of certain representation spaces
\paperinfo in: Quantization of Singular Symplectic Quotients,
N. P. Landsman, M. Pflaum, M. Schlichenmaier, eds.,
Workshop, Oberwolfach,
August 1999,
Progress in Mathematics, Vol. 198
\publ Birkh\"auser
\publaddr Boston $\cdot$ Basel $\cdot$ Berlin
\yr 2001
\pages 119--135
\finalinfo{\tt math.dg/0012184}
\endref

\ref \no \qr
\by J. Huebschmann
\paper K\"ahler reduction and quantization
\paperinfo {\tt math.SG/0207166}
%\jour \vol \yr \pages 
\endref

\ref \no \severi
\by J. Huebschmann
\paper Severi varieties and holomorphic nilpotent orbits
\linebreak
\paperinfo {\tt math.DG/0206143}
%\jour \vol \yr \pages 
\endref

\ref \no \junigusa
\by J.-I. Igusa
\paper A classification of spinors up to dimension twelve
\jour Amer. J. of Math.
\vol 92
\yr 1970
\pages 997--1028
\endref

\ref \no \inonwign
\by E. Inonu and E.P. Wigner
\paper On the contraction of groups and their representations
\jour Proc. Nat. Acad. Sci. USA
\vol 39
\yr 1953
\pages 510--524
\endref

\ref \no \jacobexc
\by  N. Jacobson
\book Exceptional Lie algebras
\bookinfo Lecture Notes in Pure and Applied Mathematics
\publ Marcel Dekker
\publaddr New York
\yr 1971
\endref

\ref \no \jacobstw
\by N. Jacobson
\paper Some projective varieties defined by Jordan algebras
\jour J. of Algebra
\vol 97
\yr 1985
\pages 556--598
\endref

\ref \no \kemneone
\by G. Kempf and L. Ness
\paper The length of vectors in representation spaces
\jour Lecture Notes in Mathematics
\vol 732
\yr 1978
\pages 233--244
\paperinfo Algebraic geometry, Copenhagen, 1978
\publ Springer 
\publaddr Berlin $\cdot$ Heidelberg $\cdot$ New York
\endref

\ref \no \kirwaboo
\by F. Kirwan
\book Cohomology of quotients in symplectic and algebraic geometry
\publ Princeton University Press
\publaddr Princeton, New Jersey
\yr 1984
\endref

\ref \no \kobaswan
\by P. Z. Kobak and A. Swann
\paper Classical nilpotent orbits as hyperk\"ahler quotients
\jour Int. J. Math.
\vol 7
\yr 1996
\pages 193--210
\endref

\ref \no \kodnispe
\by M. Kodaira, L. Nirenberg and D. C. Spencer
\paper On the existence of deformations of complex analytic structures
\jour Ann. of Math.
\vol 68
\yr 1958
\pages 450--457
\endref

\ref \no \kostasix
\by B. Kostant
\paper The principal three-dimensional subgroup
and the Betti numbers of a complex simple Lie group
\jour Amer. J. of Math.
\vol 81
\yr 1959
\pages 973--1032
\endref

\ref \no \kostafou
\by B. Kostant
\paper Lie group representations on polynomial rings
\jour Amer. J. of Math.
\vol 85
\yr 1963
\pages 327--404
\endref

\ref \no \kraprotw
\by H. Kraft and C. Procesi
\paper On the geometry of
conjugacy classes in classical groups
\jour Comment. Math. Helv.
\vol 57
\yr 1982
\pages 529--602
\endref

\ref \no \kralyvin
\by I. S. Krasil'shchik, V. V. Lychagin, and A. M. Vinogradov
\book Geometry of Jet Spaces and Nonlinear Partial Differential Equations
\bookinfo Advanced Studies in Contemporary Mathematics, vol. 1
\publ Gordon and Breach Science Publishers
\publaddr New York, London, Paris, Montreux, Tokyo
\yr 1986
\endref

\ref \no \kronhtwo
\by P. B. Kronheimer
\paper Instantons and the geometry of the nilpotent variety
\jour J. Diff. Geom.
\vol 32
\yr 1990
\pages 473--490
\endref

\ref \no \kumnarra
\by S. Kumar, M. S. Narasimhan, and A. Ramanathan
\paper Infinite Grassmannian and moduli space of $G$-bundles
\jour Math. Ann.
\vol 300
\yr 1994
\pages 41-75
\endref

\ref \no \lawmiboo
\by H. B. Lawson and M. L. Michelsohn
\book Spin Geometry
\bookinfo Princeton Mathematical Series
\vol 34
\publ Princeton University Press
\publaddr Princeton, N. J.
\yr 1989
\endref

\ref \no \lermonsj
\by E. Lerman, R. Montgomery, and R. Sjamaar
\paper Examples of singular reduction
\paperinfo Symplectic Geometry,
Warwick, 1990,  D. A. Salamon, editor, 
London Math. Soc. Lecture Note 
Series, vol. 192
\yr 1993
\pages  127--155
\publ Cambridge University Press
\publaddr Cambridge, UK
\endref

\ref \no \loobotwo
\by O. Loos
\book Jordan pairs
\bookinfo Lecture Notes in Mathematics, Vol. 460
\yr 1975
\publ Springer 
\publaddr Berlin $\cdot$ Heidelberg $\cdot$ New York
\endref

\ref \no \loobothr
\by O. Loos
\book Bounded symmetric domains and Jordan pairs
\bookinfo Math. Lectures
\yr 1977
\publ University of California
\publaddr Irvine, California
\endref

\ref \no \marsrati
\by J. Marsden and T. Ratiu
\paper Reduction of Poisson manifolds 
\jour Lett. in Math. Phys.
\vol 11
\yr 1985
\pages 161--170
\endref

\ref \no \marswein
\by J. Marsden and A. Weinstein
\paper Reduction of symplectic manifolds with symmetries
\jour Rep. on Math. Phys.
\vol 5
\yr 1974
\pages 121--130
\endref

\ref \no \mccrione
\by  K. McCrimmon
\paper Jordan algebras and their applications
\jour Bull. Amer. Math. Soc.
\vol 84
\yr 1978
\pages 612--627
\endref

\ref \no \mehtsesh
\by V. Mehta and C. Seshadri
\paper Moduli of vector bundles on curves with parabolic structure
\jour Math. Ann.
\vol 248
\yr 1980
\pages 205--239
\endref

\ref \no \mizuntwo
\by K. Mizuno
\paper The conjugate classes of unipotent elements
of the Chevalley groups of type $E_7$ and $E_8$
\jour Tokyo J. Math.
\vol 3
\yr 1980
\pages 391--461
\endref

\ref \no \muruschi
\by I. Muller, H. Rubenthaler, and G. Schiffmann
\paper Structures des espaces pr\'ehomo-\linebreak
g\`enes
associ\'es \`a certaines alg\`ebres de Lie gradu\'ees
\jour Math. Ann.
\vol 274
\yr 1986
\pages 95--123
\endref

\ref \no \narasesh
\by M. S. Narasimhan and C. S. Seshadri
\paper Stable and unitary vector bundles on a compact Riemann surface
\jour Ann. of Math.
\vol 82
\yr 1965
\pages  540--567
\endref

\ref \no \naramntw
\by M. S. Narasimhan and S. Ramanan
\paper Moduli of vector bundles on a compact Riemann surface
\jour Ann. of Math.
\vol 89
\yr 1969
\pages  19--51
\endref

\ref \no \naramnth
\by M. S. Narasimhan and S. Ramanan
\paper 2$\theta$-linear systems on abelian varieties
\jour Bombay colloquium
\yr 1985
\pages  515--427
\endref

\ref \no \neherone
\by E. Neher
\paper On the classification of Lie and Jordan triple systems
\jour Comm. Algebra
\vol 13
\yr 1985
\pages 2615--2667
\endref

\ref \no \neherboo
\by E. Neher
\book Jordan triple systems by the grid approach
\bookinfo Lecture Notes in Mathematics, Vol. 1280
\yr 1987
\publ Springer 
\publaddr Berlin $\cdot$ Heidelberg $\cdot$ New York
\endref

\ref \no \nessone
\by L. Ness
\paper A stratification of the null cone via the moment map
\jour Amer. J. of Math.
\vol 106
\yr 1984
\pages 1231--1329
\endref

\ref \no \ohtaone
\by T. Ohta
\paper The closures of nilpotent orbits
in the classical symmetric pairs
and their singularities
\jour T\hataccent ohoku Math. J.
\vol 43
\yr 1991
\pages 161--211
\endref

\ref \no \ortratwo
\by J.-P. Ortega and T. Ratiu
\paper Singular reduction of Poisson manifolds
\jour Letters in Math. Physics
\vol 46
\yr 1998
\pages 359--372
\endref

\ref \no  \peteraci
\by H. P. Petersson and M. L. Racine
\paper Albert Algebras
\paperinfo in: Jordan Algebras,
Proceedings of the Conference
held in Oberwolfach, Germany, August 9 -- 15, 1992,
W. Kaup, K. McCrimmon, H. P. Petersson, editors
\publ Walter de Gruyter
\publaddr Berlin $\cdot$ New York
\yr 1994
\pages 197--207
\endref

\ref \no \rinehart
\by G. Rinehart
\paper Differential forms for general commutative algebras
\jour  Trans. Amer. Math. Soc.
\vol 108
\yr 1963
\pages 195--222
\endref

\ref \no \satakboo
\by I. Satake
\book Algebraic structures of symmetric domains
\bookinfo Publications of the Math. Soc. of Japan, vol. 14
\publ Princeton University Press
\publaddr Princeton, NJ
\yr 1980
\endref

\ref \no \satokimu
\by M. Sato and T. Kimura
\paper A classification of irreducible prehomogeneous vector spaces
and their relative invariants
\jour Nagoya Math. J.
\vol 65
\yr 1977
\pages 1--155
\endref

\ref \no \gwschwar
\by G. W. Schwarz
\paper Smooth functions invariant under the action of
a compact Lie group
\jour Topology 
\vol 14
\yr 1975
\pages 63--68
\endref

\ref \no \gwschwat
\by G. W. Schwarz
\paper The topology of algebraic quotients
\paperinfo In: Topological methods in algebraic
transformation groups,
Progress in Mathematics, Vol. 80
\yr 1989
\pages 135--152
\publ Birkh\"auser
\publaddr Boston $\cdot$ Basel $\cdot$ Berlin
\endref

\ref \no \sekiguch
\by J. Sekiguchi
\paper Remarks on real nilpotent orbits of a symmetric pair
\jour J. Math. Soc. Japan
\vol 39
\yr 1987
\pages 127--138
\endref

\ref \no \sjamatwo
\by R. Sjamaar
\paper Holomorphic slices, symplectic reduction, and multiplicities of 
representations
\jour Ann. of Math.
\vol 141
\yr 1995
\pages 87--129
\endref

\ref \no \sjamlerm
\by R. Sjamaar and E. Lerman
\paper Stratified symplectic spaces and reduction
\jour Ann. of Math.
\vol 134
\yr 1991
\pages 375--422
\endref

\ref \no \slodoboo
\by P. Slodowy
\book Simple singularities and simple algebraic groups
\bookinfo Lecture Notes in Mathematics, Vol. 815
\yr 1980
\publ Springer 
\publaddr Berlin $\cdot$ Heidelberg $\cdot$ New York
\endref

\ref \no \spristei
\by T. Springer and R. Steinberg
\paper Conjugacy classes
\paperinfo in: Seminar on Algebraic Groups and Related Finite Groups;
Lecture Notes in Mathematics
No. 131
\publ Springer Verlag
\publaddr Berlin $\cdot$ Heidelberg $\cdot$
New York
\yr 1970
\pages 167--266
\endref

\ref \no \varoucha
\by J. Varouchas
\paper 
K\"ahler spaces and proper open morphisms
\jour Math. Ann.
\vol 283
\yr 1989
\pages 13--52
\endref

\ref \no \vergnsix
\by M. Vergne
\paper 
Instantons et correspondance de Kostant-Sekiguchi
\jour C. R. Acad. Sci. Paris S\'erie I
\vol 320
\yr 1995
\pages 901--906
\endref

\ref \no \vinberg
\by E. B. Vinberg
\paper 
On the classification of the nilpotent elements of graded
Lie algebras
\jour Sov. Math. Dokl.
\vol 16
\yr 1975
\pages 1517--1520
\endref

\ref \no \weinstwo
\by A. Weinstein
\paper The local structure of Poisson manifolds
\jour J. of Diff. Geom.
\vol 18
\yr 1983
\pages 523--557
\endref

\ref \no \weinsone
\by A. Weinstein
\paper Poisson structures
\jour Ast\'erisque,
\vol hors-serie 
\yr 1985
\pages 421--434
\paperinfo in: E. Cartan et les Math\'ematiciens \linebreak
d'aujourd'hui, 
Lyon, 25--29 Juin, 1984
\endref

\ref \no \weylbook
\by H. Weyl
\book The classical groups
\publ Princeton University  Press
\publaddr Princeton, New Jersey
\yr 1946
\endref

\ref \no \whitnboo
\by H. Whitney
\book Complex analytic varieties
\publ Addison-Wesley Pub. Comp.
\publaddr Reading, Ma, Menlo Park, Ca, London, Don Mills, Ontario
\yr 1972
\endref
\ref \no \woodhous
\by N. M. J. Woodhouse
\book Geometric quantization
\bookinfo Second edition
\publ Clarendon Press
\publaddr Oxford
\yr 1991
\endref

\enddocument